\documentclass[11pt,a4paper]{article}
\usepackage[margin=1in]{geometry}
\usepackage{amsmath,amssymb,amsthm,mathtools}
\usepackage{mathrsfs}
\usepackage{bm}
\usepackage{enumitem}
\usepackage{fancyhdr}
\usepackage[hidelinks]{hyperref}
\usepackage[capitalise,nameinlink]{cleveref}

\usepackage{hyperref}
\hypersetup{
  colorlinks=true,
  linkcolor=blue,
  citecolor=blue,
  urlcolor=blue
}

\numberwithin{equation}{section}
\theoremstyle{plain}
\newtheorem{theorem}{Theorem}[section]
\newtheorem{lemma}{Lemma}[section]
\newtheorem{proposition}{Proposition}[section]
\newtheorem{corollary}{Corollary}[section]
\theoremstyle{definition}
\newtheorem{definition}{Definition}[section]
\theoremstyle{remark}
\newtheorem{remark}{Remark}[section]

\newcommand{\R}{\mathbb{R}}
\newcommand{\T}{\mathbb{T}}
\newcommand{\Tl}{\mathbb{T}_l^2}
\newcommand{\N}{\mathbb{N}}
\newcommand{\Z}{\mathbb{Z}}
\newcommand{\bx}{\mathbf{x}}

\newcommand{\bU}{\mathbf{U}}
\newcommand{\bX}{\mathbf{X}}
\newcommand{\rs}{\rho_{\mathrm{s}}}
\newcommand{\ps}{p_{\mathrm{s}}}
\newcommand{\rin}{\rho_{\mathrm{in}}}
\newcommand{\vin}{\varrho_{\mathrm{in}}}
\newcommand{\uin}{v_{\mathrm{in}}}
\newcommand{\qin}{q_{\mathrm{in}}}
\newcommand{\rinf}{\rho_\infty}
\newcommand{\vinf}{\varrho_\infty}
\newcommand{\dd}{\,\mathrm{d}}
\newcommand{\dx}{\,\mathrm{d}\mathbf{x}}
\newcommand{\dxh}{\,\mathrm{d}x_h}
\newcommand{\dt}{\,\mathrm{d}t}
\newcommand{\dtau}{\,\mathrm{d}\tau}
\newcommand{\ddt}{\frac{\mathrm{d}}{\mathrm{d}t}}
\newcommand{\AB}{A_B}
\newcommand{\ABh}{A_B^{\frac12}}
\newcommand{\ABmh}{A_B^{-\frac12}}
\newcommand{\ABmi}{A_B^{-1}}
\newcommand{\Es}{\mathcal{E}_\sharp}
\newcommand{\Eb}{\mathcal{E}_\flat}
\newcommand{\Ef}{E_\flat}
\newcommand{\Ff}{F_\flat}
\newcommand{\Esh}{E_\sharp}
\newcommand{\Fsh}{F_\sharp}
\newcommand{\Fn}{F_\natural}
\newcommand{\Ev}{E_\varrho}
\newcommand{\Fv}{F_\varrho}
\newcommand{\Eu}{E_u}
\newcommand{\Fu}{F_u}
\newcommand{\Eut}{E_{u_t}}
\newcommand{\Fut}{F_{u_t}}
\newcommand{\Eutt}{E_{u_{tt}}}
\newcommand{\Futt}{F_{u_{tt}}}
\newcommand{\Fuz}{F_{u,z}}
\newcommand{\cEv}{\mathcal{E}_\varrho}
\newcommand{\cEu}{\mathcal{E}_u}
\newcommand{\cEut}{\mathcal{E}_{u_t}}
\newcommand{\cEutt}{\mathcal{E}_{u_{tt}}}
\newcommand{\cT}{\mathcal{T}}
\newcommand{\cN}{\mathcal{N}}
\newcommand{\cG}{\mathcal{G}}
\newcommand{\cR}{\mathcal{R}}
\newcommand{\cJ}{\mathcal{J}}
\newcommand{\cK}{\mathcal{K}}
\newcommand{\cC}{\mathcal{C}}
\newcommand{\cM}{\mathcal{M}}
\newcommand{\cA}{\mathcal{A}}
\newcommand{\cB}{\mathcal{B}}
\newcommand{\cD}{\mathcal{D}}

\newcommand{\cc}[1]{\mathsf{c}_{#1}}
\newcommand{\norm}[1]{\left\|#1\right\|}
\newcommand{\nrm}[1]{\|#1\|}
\newcommand{\ip}[1]{\left\langle #1\right\rangle}
\newcommand{\ol}[1]{\overline{#1}}
\newcommand{\wt}[1]{\widetilde{#1}}
\newcommand{\eqdef}{\stackrel{\mathrm{def}}{=}}
\newcommand{\Lap}{\Delta}
\newcommand{\Laph}{\Delta_h}
\newcommand{\nh}{\nabla_h}
\newcommand{\nhp}{\nabla_h^{\perp}}
\newcommand{\pa}{\partial}
\newcommand{\px}{\partial_x}
\newcommand{\py}{\partial_y}
\newcommand{\pz}{\partial_z}
\newcommand{\pt}{\partial_t}
\newcommand{\sgm}{\sigma}

\DeclareMathOperator{\esssup}{ess\,sup}

\DeclareMathOperator{\divv}{div}
\DeclareMathOperator{\curl}{curl}
\newcommand{\divh}{\divv_h}
\newcommand{\curlh}{\curl_h}

\title{\textbf{Asymptotic stability and Rayleigh--Taylor instability of stratified states for the three-dimensional rotating primitive equations with gravity}}
\author{Yanlong Fan$^{a}$\thanks{fanyanlong@swjtu.edu.cn, (Y. Fan)} \qquad 
Fengjie Luo$^{b}$\thanks{Luofengjie664@163.com, (F. Luo)}  \qquad 
Yiqiu Mao$^{c}$\thanks{yqmao@gzhu.edu.cn, (Y. Mao)}  \qquad 
Quan Wang$^{b}$\thanks{Corresponding author: xihujunzi@scu.edu.cn, (Q. Wang)}\\[4pt]
{\small\itshape $^a$ School of Mathematics, Southwest Jiaotong University, Chengdu, Sichuan, 611756, China}\\
{\small\itshape $^b$ College of Mathematics, Sichuan University, Chengdu, Sichuan, 610065, China}\\
{\small\itshape $^c$ School of Mathematics and Information Science, Guangzhou University,
Guangzhou, 510000,  China}
}
\date{}

\begin{document}
\maketitle
\thispagestyle{empty}

\begin{abstract}
We study the three-dimensional rotating primitive equations with gravity in a horizontally periodic channel. The horizontal velocity is viscous, while the density is transported without diffusion. For every value of the Coriolis parameter, stably stratified equilibria are asymptotically stable. Small compatible $H^4$ perturbations generate global classical solutions with uniform $H^4$ bounds; the oscillatory part of the density and the velocity decay like $(1+t)^{-1}$, and the density converges at the same rate to a stratified profile, its horizontal mean converging faster. For small data, the limit profile is the decreasing vertical rearrangement of the initial density, that is, the unique minimizer of the potential energy among all rearrangements of the initial density; it depends neither on the rotation, nor on the viscosity, nor on the initial velocity. The kinetic energy plus the available potential energy is dissipated exactly by viscosity and decays like $(1+t)^{-2}$, and every fluid particle settles at the height where the limit density equals its own density. The rates for the density and for the energy are optimal for every rotation rate. The linearized problem has real, slowly decaying normal modes, which are internal gravity waves of high vertical wavenumber overdamped by viscosity, and small solutions built on them show that the density rates $(1+t)^{-(1-r/4)}$ in $H^r$, $0\le r<4$, the pressure rate $(1+t)^{-1}$ and the energy rate $(1+t)^{-2}$ cannot be improved uniformly for small data. The main difficulties are the loss of a horizontal derivative in the vertical velocity, the failure of the embedding $H^1(\mathbb T^2)\subset L^\infty$, and the coupling of divergent and vortical motion by rotation; at top order, the rotational coupling is controlled by a positive power of the initial low-order energy. When the stratification is unstable somewhere, every horizontal mode is linearly unstable with real growth rates, and the equilibrium is nonlinearly unstable. Rotation strictly lowers each growth rate, with an explicit bound in terms of a Taylor number built on the horizontal wavelength, but it leaves the supremum of the growth rates unchanged.
\end{abstract}



\newpage
\tableofcontents
\newpage

\section{Introduction}

Large-scale motions of the ocean and the atmosphere are described, to a good approximation, by the primitive equations, in which the vertical momentum equation is replaced by hydrostatic balance \cite{Pedlosky1987,Gill1982,Majda2003,Vallis2017}. Hydrostatic balance is the formal limit of small aspect ratio, and this limit has been justified rigorously \cite{AzeradGuillen2001,LiTiti2019}. The mathematical theory of the primitive equations goes back to Lions, Temam and Wang \cite{LionsTemamWang1992}; we refer to \cite{TemamZiane2004,LiTiti2018} for surveys. At these scales two effects dominate: the rotation of the Earth, which enters through the Coriolis force, and gravity acting on a density that varies with height. When heavier fluid lies beneath lighter fluid, the fluid is stably stratified and supports internal gravity waves \cite{Gill1982,Sutherland2010}; when heavier fluid lies above lighter fluid somewhere, it is subject to the Rayleigh--Taylor instability. In this paper we address the basic long-time questions for such equilibria. Does a small perturbation of a stably stratified state return to a stratified state, at what rate, and to which profile? And what changes when the stratification is unstable somewhere?

We consider the three-dimensional rotating primitive equations in a horizontally periodic channel. The horizontal velocity is dissipated by the full Laplacian, while the density is transported without any diffusion. This regime models stratified fluids in which heat or salt diffuses much more slowly than momentum, and it is delicate from the mathematical point of view: the density equation provides neither smoothing nor decay, so the stabilizing effect of gravity has to be extracted from the coupling between density and velocity. The same difficulty lies at the heart of the stability theory for the Boussinesq equations without thermal diffusion \cite{DoeringWuZhaoZheng2018,TaoWuZhaoZheng2020,DongSun2022,JangKim2023}, the incompressible porous medium equation \cite{Elgindi2017,CastroCordobaLear2019} and the Stokes--transport system \cite{DalibardGuillodLeblond2025}. Rotation adds a further feature. The Coriolis force does no work, but it couples the divergent part of the horizontal velocity, which moves the density, to the vortical part, which gravity does not see. All our results hold for every value of the Coriolis parameter, including the non-rotating case.

\subsection{Formulation of the problem}

Let $l>0$, $\Tl\eqdef(\R/2\pi l\Z)^2$ and
\begin{equation*}
\Omega\eqdef\Tl\times(0,1),\qquad \bx=(x_h,z),\quad x_h=(x,y)\in\Tl,\quad z\in(0,1).
\end{equation*}
We write $v=(v_1,v_2)$ for the horizontal velocity, $w$ for the vertical velocity and $\bU=(v,w)$; $\nh=(\px,\py)$, $\divh$ and $\Laph$ are the horizontal gradient, divergence and Laplacian, $\nabla=(\nh,\pz)$ and $\Lap=\Laph+\pz^2$. For a horizontal vector field $v$ we set
\begin{equation*}
v^\perp\eqdef(-v_2,v_1),\qquad \curlh v\eqdef\px v_2-\py v_1 .
\end{equation*}
The system reads
\begin{equation*}
\begin{cases}
\rho_t+\bU\cdot\nabla\rho=0,\\
\rho\,(v_t+\bU\cdot\nabla v)+f\rho\, v^\perp+\nh p-\mu\Lap v=0,\\
\divh v+\pz w=0,\\
\pz p=-\rho g,
\end{cases}
\end{equation*}
where $\rho$ is the density, $p$ the pressure, $g>0$ the gravitational constant, $\mu>0$ the viscosity and $f\in\R$ the Coriolis parameter (the $f$-plane approximation, with the traditional approximation, which keeps only the horizontal component of the Coriolis force). The density equation carries no diffusion.

We consider perturbations of a hydrostatic equilibrium $(\rho,v,w,p)=(\rs(z),0,0,\ps(z))$ with $\ps'=-g\rs$. The equilibrium is \emph{stably stratified} if
\begin{equation}\label{eq:stable}
\rs\in C^6([0,1]),\qquad 0<\rs(z)<\infty,\qquad \sup_{z\in[0,1]}\rs'(z)<0 .
\end{equation}
Under \eqref{eq:stable} we write $\sgm\eqdef-1/\rs'$, so that $\sgm\in C^5([0,1])$. Six derivatives of $\rs$ are used only once, through $\sgm^{(5)}$ in the top-order transport estimate (\Cref{rem:C6}); local well-posedness needs only $\rs\in C^5$. With $\rho=\rs+\varrho$ and $p=\ps+q$, the perturbation satisfies
\begin{equation}\label{PE}\tag{PE}
\begin{cases}
\varrho_t+\bU\cdot\nabla\varrho=-w\rs',\\
\rho\,(v_t+\bU\cdot\nabla v)+f\rho\, v^\perp+\nh q-\mu\Lap v=0,\\
\divh v+\pz w=0,\\
\pz q=-g\varrho,
\end{cases}
\end{equation}
with initial data
\begin{equation}\label{eq:ID}
\varrho(0,\bx)=\vin(\bx),\qquad v(0,\bx)=\uin(\bx),
\end{equation}
and boundary conditions
\begin{equation}\label{eq:BC}
v|_{\pa\Omega}=0,\qquad \varrho|_{\pa\Omega}=\pz\varrho|_{\pa\Omega}=0,
\end{equation}
where $\pa\Omega=\Tl\times\{0,1\}$. The velocity satisfies the no-slip condition and the density is clamped at the lids. The vertical velocity is diagnostic:
\begin{equation}\label{eq:w}
w(t,x_h,z)=-\int_0^z\divh v(t,x_h,z')\dd z' .
\end{equation}
Thus $w|_{z=0}=0$, and the no-penetration condition $w|_{z=1}=0$ is the lid constraint
\begin{equation}\label{eq:lid}
\divh\int_0^1v(t,x_h,z)\dd z\equiv0 .
\end{equation}
At the initial time it reads
\begin{equation}\label{eq:lid0}
\divh\int_0^1\uin(x_h,z)\dd z\equiv0,
\end{equation}
and the first-order compatibility condition for the fixed boundary value of $w$ is
\begin{equation}\label{eq:wt0}
w_t(0)|_{\pa\Omega}=0 .
\end{equation}
By \Cref{lem:propagation}, \eqref{eq:lid0} propagates, so $w|_{\pa\Omega}=0$ at all later times, and
\begin{equation}\label{eq:wtt}
w_t|_{\pa\Omega}=w_{tt}|_{\pa\Omega}=0
\end{equation}
at positive times at which these traces exist. The surface pressure that enforces \eqref{eq:lid} is determined by the solution (\Cref{lem:propagation}); in particular \eqref{eq:wt0} fixes the initial pressure $\qin$ in terms of $(\vin,\uin)$, and the genuine compatibility condition on the data is \eqref{eq:compat} below. Throughout the paper, a solution of \eqref{PE}--\eqref{eq:BC} always means a solution of \eqref{PE} with \eqref{eq:ID}--\eqref{eq:lid}: $w$ is given by \eqref{eq:w} and the lid constraint \eqref{eq:lid} holds. The pressure $q$ is then determined by $(\varrho,v)$ up to an additive function of $t$; we fix it by requiring the surface pressure $P=q+g\int_0^z\varrho\dd s$ to have zero horizontal mean.

Three structural facts are used throughout. The Coriolis force does no work, $v^\perp\cdot v=0$, and its horizontal derivatives satisfy $\pa_jv^\perp\cdot\pa_jv=0$ pointwise. It exchanges divergence and vorticity, $\divh v^\perp=-\curlh v$ and $\curlh v^\perp=\divh v$. And the total horizontal force
\begin{equation}\label{eq:Gdef}
G\eqdef\rho\big(\dot v+fv^\perp\big),\qquad \dot v\eqdef v_t+\bU\cdot\nabla v,
\end{equation}
enters the hydrostatic structure in exactly the way the inertial term enters it without rotation; see \eqref{eq:wrep} and \eqref{eq:decomp2}.

\subsection{Main results}

Our first result is global existence for small data, with uniform bounds in $H^4$ and algebraic decay, for every $f\in\R$.

\begin{theorem}\label{thm:main}
Let $\rs$ satisfy \eqref{eq:stable} and let $f\in\R$. Suppose that the initial data satisfy
\begin{equation*}
\vin\in H^4,\quad \rs+\vin>0\ \text{in }\Omega,\quad \vin=\pz\vin=0\ \text{on }\pa\Omega,\quad\text{and}\quad \uin\in H^4\cap H^1_0,
\end{equation*}
together with the constraint \eqref{eq:lid0} and the compatibility condition \eqref{eq:compat}. There exist $\varepsilon_0>0$ and $C>0$, depending only on $\rs$, $\mu$, $g$, $f$ and $l$, such that if
\begin{equation*}
\nrm{\vin}_{H^4}+\nrm{\uin}_{H^4}\le\varepsilon_0,
\end{equation*}
then \eqref{PE} with \eqref{eq:ID}--\eqref{eq:lid} has a unique global classical solution $(\varrho,v,q)$ on $[0,\infty)\times\Omega$ (with $q$ normalized as above), satisfying
\begin{equation}\label{eq:unif}
\nrm{\varrho(t)}_{H^4}+\nrm{v(t)}_{H^4}\le C\big(\nrm{\vin}_{H^4}+\nrm{\uin}_{H^4}\big),\qquad t\ge0,
\end{equation}
the velocity regularity
\begin{equation}\label{eq:ureg}
\sup_{t\ge0}\nrm{v(t)}_{H^4}^2+\int_0^\infty\nrm{v(t)}_{H^5}^2\dt\le C\big(\nrm{\vin}_{H^4}+\nrm{\uin}_{H^4}\big)^2,
\end{equation}
the clamped conditions $\varrho(t)|_{\pa\Omega}=\pz\varrho(t)|_{\pa\Omega}=0$ for all $t\ge0$, and the decay rates
\begin{equation}\label{eq:decay}
\nrm{\wt\varrho(t)}_{L^2}+\nrm{\nh\varrho(t)}_{L^2}+\nrm{v(t)}_{H^2}+\nrm{v_t(t)}_{L^2}+\nrm{v_{tt}(t)}_{L^2}\le C\big(\nrm{\vin}_{H^4}+\nrm{\uin}_{H^4}\big)(1+t)^{-1},
\end{equation}
for all $t\ge0$. Here $\wt\varrho=\varrho-\bar\varrho$ is the oscillatory part of the density, $\bar\varrho(t,z)=|\Tl|^{-1}\int_{\Tl}\varrho(t,x_h,z)\dxh$.
\end{theorem}

Our second result identifies the terminal state. The horizontal mean $\bar\varrho$ of the density perturbation does not decay: by \eqref{eq:meaneq}, it changes only through the vertical flux $\ol{w\wt\varrho}$. Since the velocity is divergence-free and tangent to the lids, transport preserves the volume of every superlevel set of the density. For small data the limit is strictly stratified, and therefore it is the decreasing vertical rearrangement of the initial density. Neither the initial velocity, nor the viscosity, nor the rotation enters.

\begin{theorem}\label{thm:limit}
Let the hypotheses of \Cref{thm:main} hold, let $(\varrho,v,q)$ be the global solution, and write $\mathsf N\eqdef\nrm{\vin}_{H^4}+\nrm{\uin}_{H^4}$. There is $\vinf=\vinf(z)\in H^4(0,1)$ with $\vinf(0)=\vinf(1)=0$ such that
\begin{equation*}
\nrm{\varrho(t)-\vinf}_{L^2(\Omega)}\le C\mathsf N(1+t)^{-1},\qquad t\ge0,
\end{equation*}
and the assertions (i)--(iv) below hold. Moreover there is $\cc{*}=\cc{*}(\rs,\mu,g,f,l)$, independent of $\varepsilon_0$, with $\sup_{t\ge0}\nrm{\pz\varrho(t)}_{L^\infty}+\nrm{\pz\vinf}_{L^\infty(0,1)}\le\cc{*}\varepsilon_0$, and if, after decreasing $\varepsilon_0$ if necessary (the constants in \eqref{eq:unif} do not depend on $\varepsilon_0$),
\begin{equation}\label{eq:slope}
\varepsilon_0\le\frac1{2\cc{*}}\min_{z\in[0,1]}\big(-\rs'(z)\big),
\end{equation}
then (v) and (vi) hold as well.
\begin{enumerate}[label=(\roman*)]
\item The limit is
\begin{equation}\label{eq:limdef}
\vinf(z)\eqdef\ol\varrho_{\mathrm{in}}(z)-\int_0^\infty\pz\big(\ol{w\wt\varrho}\big)(\tau,z)\dtau,
\end{equation}
where the time integral converges as a Bochner integral in $L^2(0,1)$, and
\begin{equation}\label{eq:meanconv}
\bar\varrho(t)\to\vinf\quad\text{in }L^2(0,1)\quad\text{as }t\to\infty;
\end{equation}
moreover $\nrm{\vinf}_{H^4(0,1)}\le C\mathsf N$ and $\pz\vinf(0)=\pz\vinf(1)=0$.
\item The mean converges at a faster rate, quadratic in the data:
\begin{equation}\label{eq:rate1eps}
\nrm{\bar\varrho(t)-\vinf}_{L^2(0,1)}\le\int_t^\infty\nrm{\pt\bar\varrho(\tau)}_{L^2(0,1)}\dtau\le C\mathsf N^2(1+t)^{-\frac54},
\end{equation}
and for $s\in[0,4)$,
\begin{equation}\label{eq:rateHs}
\nrm{\bar\varrho(t)-\vinf}_{H^s(0,1)}\le C_s\mathsf N(1+t)^{-\frac54(1-\frac s4)} .
\end{equation}
\item With $\rinf\eqdef\rs+\vinf$,
\begin{equation}\label{eq:ratefull}
\nrm{\rho(t)-\rinf}_{L^2(\Omega)}=\nrm{\varrho(t)-\vinf}_{L^2(\Omega)}\le C\mathsf N(1+t)^{-1},
\end{equation}
and, for $s\in[0,4)$,
\begin{equation}\label{eq:rateLinf}
\begin{aligned}
&\nrm{\rho(t)-\rinf}_{H^s}\le C_s\mathsf N(1+t)^{-(1-\frac s4)},\\
&\nrm{\rho(t)-\rinf}_{L^\infty}\le C\mathsf N(1+t)^{-\frac56},\qquad \nrm{\nabla(\rho(t)-\rinf)}_{L^\infty}\le C\mathsf N(1+t)^{-\frac12}.
\end{aligned}
\end{equation}
With $q_\infty(z)\eqdef-g\int_0^z\vinf(s)\dd s$, the hydrostatic pressure perturbation of the limit, and the surface pressure normalized to zero horizontal mean,
\begin{equation}\label{eq:ratepress}
\nrm{q(t)-q_\infty}_{H^1(\Omega)}\le C\mathsf N(1+t)^{-1}.
\end{equation}
Together with \eqref{eq:decay}, the state $(\rho,v,w,p)$ converges at the rate $(1+t)^{-1}$ to the hydrostatic state $(\rinf,0,0,\ps+q_\infty)$.
\item $\rinf(0)=\rs(0)$ and $\rinf(1)=\rs(1)$.
\item $\pz\rinf<0$ on $[0,1]$, and $\rinf$ is the decreasing vertical rearrangement of $\rin\eqdef\rs+\vin$:
\begin{equation}\label{eq:rearr}
\rinf(z)=\rin^*(z)\eqdef\int_0^\infty\mathbf 1_{\{|\{\rin>\lambda\}|\ge|\Tl|\,z\}}\dd\lambda,\qquad0<z<1,
\end{equation}
where $|\{\rin>\lambda\}|$ is the three-dimensional Lebesgue measure and $|\Tl|=4\pi^2l^2$. Moreover, the density stays strictly stratified, $-\pz\rho(t)\ge\frac12\min_{[0,1]}(-\rs')$ on $\ol\Omega$ for all $t\ge0$, its rearrangement is $\rho(t)^*=\rinf$ at every time, and the mean converges at the faster rate
\begin{equation}\label{eq:rate138}
\nrm{\bar\varrho(t)-\vinf}_{L^2(0,1)}\le C\mathsf N^2(1+t)^{-\frac{13}8},\qquad t\ge0 .
\end{equation}
\item With $\sgm=-1/\rs'$,
\begin{equation}\label{eq:quadexp}
\Big\|\vinf-\ol\varrho_{\mathrm{in}}-\tfrac12\pz\big(\sgm\,\ol{\wt\varrho_{\mathrm{in}}^{\,2}}\big)\Big\|_{C([0,1])}\lesssim\nrm{\vin}_{H^4}^3 .
\end{equation}
\end{enumerate}
\end{theorem}

The relaxation to $\rinf$ also has a variational and a Lagrangian description, which is our third result. Gravity acts through the potential energy $g\int_\Omega z\rho\dx$, and transport preserves the distribution of the density, so every attainable density field is a rearrangement of $\rin$; among all rearrangements the potential energy is minimized exactly by $\rin^*=\rinf$. The excess $g\int_\Omega z(\rho-\rinf)\dx$ is the available potential energy of Lorenz \cite{Lorenz1955,WintersLombardRileyDAsaro1995}. The Coriolis force does no work, so the kinetic energy plus the available potential energy is dissipated exactly by viscosity, for every $f$, and each fluid particle comes to rest at the height where the limit density equals its own density.

\begin{theorem}\label{thm:energy}
Let the hypotheses of \Cref{thm:limit} hold, including \eqref{eq:slope}. Let $\mathsf N=\nrm{\vin}_{H^4}+\nrm{\uin}_{H^4}$, $m\eqdef\min_{[0,1]}(-\rs')$, $M_1\eqdef\max_{[0,1]}(-\rs')$, and
\begin{equation}\label{eq:E-H}
\mathcal H(t)\eqdef\frac12\nrm{\sqrt{\rho(t)}\,v(t)}_{L^2}^2+g\int_\Omega z\big(\rho(t)-\rinf\big)\dx .
\end{equation}
\begin{enumerate}[label=(\roman*)]
\item \emph{Minimal potential energy.} For every measurable $\phi$ on $\Omega$ that is equimeasurable with $\rin$, that is, $|\{\phi>\lambda\}|=|\{\rin>\lambda\}|$ for all $\lambda\in\R$,
\begin{equation}\label{eq:E-min}
\frac g{3M_1}\nrm{\phi-\rinf}_{L^2}^2\le g\int_\Omega z\big(\phi-\rinf\big)\dx\le\frac gm\nrm{\phi-\rinf}_{L^2}^2 .
\end{equation}
In particular, $\rinf$ is the unique minimizer of the potential energy $\phi\mapsto g\int_\Omega z\phi\dx$ among all rearrangements of $\rin$.
\item \emph{Energy identity.} For $0\le t\le T$,
\begin{equation}\label{eq:E-identity}
\mathcal H(T)+\mu\int_t^T\nrm{\nabla v(\tau)}_{L^2}^2\dtau=\mathcal H(t),
\end{equation}
and there are $c,C>0$ depending only on $\rs$ and $g$ with
\begin{equation}\label{eq:E-coerc}
c\big(\nrm{v(t)}_{L^2}^2+\nrm{\rho(t)-\rinf}_{L^2}^2\big)\le\mathcal H(t)\le C\big(\nrm{v(t)}_{L^2}^2+\nrm{\rho(t)-\rinf}_{L^2}^2\big).
\end{equation}
Consequently
\begin{equation}\label{eq:E-tail}
\mathcal H(t)=\mu\int_t^\infty\nrm{\nabla v(\tau)}_{L^2}^2\dtau\le C\mathsf N^2(1+t)^{-2},\qquad t\ge0,
\end{equation}
and the total viscous dissipation is determined by the data, independently of $f$ and $\mu$:
\begin{equation}\label{eq:E-total}
\mu\int_0^\infty\nrm{\nabla v(\tau)}_{L^2}^2\dtau=\frac12\nrm{\sqrt{\rin}\,\uin}_{L^2}^2+g\int_\Omega z\big(\rin-\rin^*\big)\dx .
\end{equation}
\item \emph{Lagrangian picture.} Let $\bX(t,\mathbf a)$ be the flow of $\bU=(v,w)$,
\begin{equation}\label{eq:E-flow}
\pt\bX(t,\mathbf a)=\bU\big(t,\bX(t,\mathbf a)\big),\qquad \bX(0,\mathbf a)=\mathbf a,
\end{equation}
on the lift $\R^2\times[0,1]$ of $\Omega$, and write $\bX=(X^h,X^z)$. There is a volume-preserving $C^1$ diffeomorphism $\bX_\infty$, bi-Lipschitz, with $\bX_\infty(x_h+2\pi lk,z)=\bX_\infty(x_h,z)+(2\pi lk,0)$ for $k\in\Z^2$, fixing every point of $\pa\Omega$, such that
\begin{equation}\label{eq:E-flowconv}
\nrm{\bX(t)-\bX_\infty}_{C^1}\le C\mathsf N(1+t)^{-\frac13},\qquad t\ge0 .
\end{equation}
More precisely, $\nrm{X^h(t)-X^h_\infty}_{C^1}\le C\mathsf N(1+t)^{-1/3}$ and $\nrm{X^z(t)-X^z_\infty}_{C^1}\le C\mathsf N(1+t)^{-1/2}$.
The final height of a particle is determined by its density:
\begin{equation}\label{eq:E-height}
X^z_\infty(\mathbf a)=\rinf^{-1}\big(\rin(\mathbf a)\big).
\end{equation}
Moreover, for $1\le p\le\infty$,
\begin{equation}\label{eq:E-disp}
\frac m2\nrm{X^z(t)-X^z_\infty}_{L^p}\le\nrm{\rho(t)-\rinf}_{L^p}\le\frac32M_1\nrm{X^z(t)-X^z_\infty}_{L^p},
\end{equation}
so that $\nrm{X^z(t)-X^z_\infty}_{L^2}\le C\mathsf N(1+t)^{-1}$ and $\nrm{X^z(t)-X^z_\infty}_{L^\infty}\le C\mathsf N(1+t)^{-5/6}$.
\item \emph{The terminal map.} If two sets of data satisfy these hypotheses, then
\begin{equation}\label{eq:E-contr}
\nrm{\rinf^{(1)}-\rinf^{(2)}}_{L^p(\Omega)}\le\nrm{\rin^{(1)}-\rin^{(2)}}_{L^p(\Omega)},\qquad1\le p\le\infty .
\end{equation}
In particular, $\rinf=\rs$ if and only if $\rin$ is equimeasurable with $\rs$.
\end{enumerate}
\end{theorem}

Neither the limit $\rinf$, nor the final heights \eqref{eq:E-height}, nor the total dissipation \eqref{eq:E-total} depends on the rotation; rotation only changes how the energy is dissipated in time. In item (iii), the bound for the vertical component of the flow is better than that for the horizontal one: it follows from the transported density, whereas the horizontal bound comes from the time integral of $\nrm{v}_{W^{1,\infty}}$. Item (iii) gives no formula for the final horizontal position $X^h_\infty$, and we do not claim that the exponents $\frac13$ and $\frac12$ are sharp.

Our fourth result shows that the rate $(1+t)^{-1}$ for the density cannot be improved, for any $f$. It rests on a one-dimensional eigenvalue problem. For $\nu>0$ let $\theta_1\ge\theta_2\ge\dots>0$ be the eigenvalues of the clamped problem, understood in weak form,
\begin{equation}\label{eq:theta-intro}
g\nu^2\sgm^{-1}\phi=\theta\,\mu\big(\phi''''-\nu^2\phi''\big),\qquad \phi\in H^2_0(0,1).
\end{equation}
This problem balances buoyancy against viscosity, and $\theta_n$ is of order $n^{-4}$ (\Cref{prop:O-modes}(a)). We also write
\begin{equation}\label{eq:vartheta}
\vartheta_\nu\eqdef\frac{\max\rs}{\mu(\pi^2+\nu^2)} .
\end{equation}

\begin{theorem}\label{thm:opt}
Let $\rs$ satisfy \eqref{eq:stable}, let $f\in\R$, fix $k\in(l^{-1}\Z)^2\setminus\{0\}$, and let $\nu=|k|$.
\begin{enumerate}[label=(\roman*)]
\item \emph{Slow normal modes.} For every sufficiently large $n$, the linearization of \eqref{PE} at $(\rs,0)$ has a real solution $(\varrho_n,v_n)$ with $\varrho_n=e^{-s_n^*t}\cos(k\cdot x_h)\,\hat\varrho_n(z)$ and $\hat\varrho_n\not\equiv0$. Its decay rate $s_n^*$ is real and satisfies
\begin{equation*}
\frac{\theta_n}{1+\frac43f^2\vartheta_\nu^2}\le s_n^*\le\frac43\theta_n,\qquad \lim_{n\to\infty}\frac{s_n^*}{\theta_n}=1 .
\end{equation*}
Its data satisfy $c(s_n^*)^{-1}\nrm{\varrho_n(0)}_{L^2}\le\nrm{\varrho_n(0)}_{H^4}+\nrm{v_n(0)}_{H^4}\le C(s_n^*)^{-1}\nrm{\varrho_n(0)}_{L^2}$, where $c>0$ and $C$ do not depend on $n$. Rotation slows these modes: $|f|\mapsto s_n^*$ is strictly decreasing.
\item \emph{Optimality.} There are $c>0$, data $(\vin^{(n)},\uin^{(n)})$ satisfying the hypotheses of \Cref{thm:energy}, with $\mathsf N_n>0$ and $\mathsf N_n\to0$, and times $t_n\to\infty$, such that the corresponding solutions satisfy
\begin{equation*}
\begin{aligned}
&\nrm{\rho^{(n)}(t_n)-\rinf^{(n)}}_{L^2}\ge\nrm{\wt\varrho^{(n)}(t_n)}_{L^2}\ge c(1+t_n)^{-1}\mathsf N_n,\qquad \nrm{\nh\varrho^{(n)}(t_n)}_{L^2}\ge c(1+t_n)^{-1}\mathsf N_n,\\
&\nrm{q^{(n)}(t_n)-q^{(n)}_\infty}_{H^1}\ge c(1+t_n)^{-1}\mathsf N_n,\qquad \mathcal H^{(n)}(t_n)\ge c(1+t_n)^{-2}\mathsf N_n^2,\\
&\nrm{X^{z,(n)}(t_n)-X^{z,(n)}_\infty}_{L^2}\ge c(1+t_n)^{-1}\mathsf N_n .
\end{aligned}
\end{equation*}
\item \emph{Sobolev norms.} For every $0\le r<4$ there is $c_r>0$ such that the data and times in (ii) can be chosen, depending on $r$, so that, in addition,
\begin{equation*}
\nrm{\wt\varrho^{(n)}(t_n)}_{H^r}\ge c_r(1+t_n)^{-(1-\frac r4)}\mathsf N_n,\qquad \nrm{\rho^{(n)}(t_n)-\rinf^{(n)}}_{H^r}\ge c_r(1+t_n)^{-(1-\frac r4)}\mathsf N_n .
\end{equation*}
In $H^4$ there is no uniform decay: there is $c>0$ such that, for every $t\ge0$, there are data of arbitrarily small size $\mathsf N$ with $\nrm{\wt\varrho(t)}_{H^4}\ge c\mathsf N$.
\end{enumerate}
\end{theorem}

\Cref{thm:opt} collects \Cref{prop:O-modes,thm:O-lin,thm:O-nonlin,thm:O-Hr} and \Cref{cor:O-H4}. By \Cref{thm:O-nonlin}(i) and \Cref{thm:O-Hr}(ii), the lower bounds in (ii) and (iii) also hold at each fixed time $t$, for data of arbitrarily small size. Consequently, none of the following rates can be improved by a factor $(1+t)^{-\delta}$, $\delta>0$, uniformly for small data: the rate $(1+t)^{-(1-r/4)}$ for $\wt\varrho$ and for $\rho-\rinf$ in $H^r$, $0\le r<4$ (by interpolation of \eqref{eq:decay} and \eqref{eq:unif}, see \Cref{thm:O-Hr}, and by \eqref{eq:rateLinf}); the rate $(1+t)^{-1}$ for $\nh\varrho$ in \eqref{eq:decay}, for the pressure in \eqref{eq:ratepress} and for the vertical displacement in \eqref{eq:E-disp}; and the rate $(1+t)^{-2}$ for $\mathcal H$ in \eqref{eq:E-tail}. The faster rates for the mean in \Cref{thm:limit}, the pointwise rates in \eqref{eq:rateLinf} and the flow-map rates in \Cref{thm:energy}(iii) are upper bounds, and we do not claim that they are sharp.

The slow modes have a simple physical meaning. They are internal gravity waves of high vertical wavenumber $n$ that are overdamped by viscosity: buoyancy restores them only on the time scale $\theta_n^{-1}\sim n^4$, while their $H^4$ norm exceeds their $L^2$ norm by a factor of order $n^4$. Rotation adds the Ekman-type drag $f^2R_\nu$ of \Cref{sec:RT-lin} and slows them down further, but by a factor that tends to $1$ as $n\to\infty$. These modes carry little velocity and have zero horizontal mean, so they do not directly decide whether the velocity rates in \eqref{eq:decay}, or the faster rates for the mean in \Cref{thm:limit}, are sharp.

Finally, we turn to unstable stratification. For a horizontal wave vector $k\in(l^{-1}\Z)^2\setminus\{0\}$ the linearized problem decouples (\Cref{sec:RT}); its growth rates depend on $k$ only through $\nu=|k|$, and we write $\lambda(\nu)=\lambda(\nu;f)$ for the largest one. Let
\begin{equation}\label{eq:Lambdainf}
\Lambda_\infty\eqdef\sup_{\phi\in H^1_0(0,1)\setminus\{0\}}\frac{g\int_0^1\rs'|\phi|^2\dd z}{\mu\int_0^1|\phi'|^2\dd z}.
\end{equation}

\begin{theorem}\label{thm:RT}
Let $\rs\in C^\infty([0,1])$ with $\min\rs>0$ and $\max\rs'>0$, and let $f\in\R$.
\begin{enumerate}[label=(\roman*)]
\item Every horizontal mode $k\ne0$ is linearly unstable: there is a growing normal mode $e^{\lambda(|k|)t}$ with $\lambda(|k|)>0$, and every normal mode growing in time has a real growth rate $\le\lambda(|k|)$.
\item $\nu\mapsto\lambda(\nu;f)$ is strictly increasing, $\lambda(\nu;f)<\Lambda_\infty$, and $\lambda(\nu;f)\to\Lambda_\infty$ as $\nu\to\infty$. The supremum $\Lambda_\infty$ does not depend on $f$ and is not attained.
\item Rotation lowers every growth rate: $|f|\mapsto\lambda(\nu;f)$ is strictly decreasing, $\lambda(\nu;f)\to0$ as $|f|\to\infty$ for each fixed $\nu$, and, with $\vartheta_\nu$ as in \eqref{eq:vartheta},
\begin{equation*}
\frac{\lambda(\nu;0)}{1+f^2\vartheta_\nu^2}<\lambda(\nu;f)<\lambda(\nu;0)\quad(f\ne0),
\end{equation*}
\item $(\rs,0)$ is nonlinearly unstable: there are $\varepsilon_1>0$ and smooth data of arbitrarily small $H^m$-norm, for every $m$, whose solutions leave the $\varepsilon_1$-ball of $H^4$. More precisely, fix $k_*\ne0$ with $\lambda_*\eqdef\lambda(|k_*|;f)>\Lambda_\infty/2$ (possible by (ii)). There are $\bar\varepsilon_0,\varepsilon_1,\delta_0>0$ such that for every $\delta\in(0,\delta_0)$ there are smooth compatible data of size $O(\delta)$ in every $H^m$ for which either the maximal existence time is at most $T^\delta\eqdef\lambda_*^{-1}\log(\bar\varepsilon_0/\delta)$, or $\nrm{\varrho(T^\delta)}_{L^2}\ge\varepsilon_1$.
\end{enumerate}
\end{theorem}

\Cref{thm:RT} collects \Cref{thm:RTlinear,thm:RTnonlinear}. Together with \Cref{thm:main}, it gives a dichotomy for every rotation rate, apart from the borderline case $\max\rs'=0$: stable stratification leads to algebraic convergence to a stratified state, and unstable stratification leads to instability. Rotation slows the instability mode by mode, and it suppresses each fixed mode as $|f|\to\infty$, but it cannot suppress the instability as a whole. Indeed, $f^2\vartheta_\nu^2$ is a Taylor number built on the horizontal wavelength, so the short horizontal waves, which grow at nearly the maximal rate $\Lambda_\infty$, do not feel the Coriolis force.

\emph{The two-dimensional problem.} For $f=0$ and data that do not depend on $y$ and have $v_2=0$, the solution keeps these properties, by uniqueness, and \eqref{PE} reduces to the non-rotating primitive equations in the vertical slab $\T_l\times(0,1)$, where $\T_l\eqdef\R/2\pi l\Z$. All the results above apply to this problem. In two dimensions the pressure also converges in $L^\infty$, because $H^1(\T_l)\subset L^\infty(\T_l)$; see \Cref{rem:2D}, which also discusses smooth data.

\subsection{Ideas of the proofs}\label{sec:ideas}

\emph{The hydrostatic loss of derivatives and the inherited dissipation.} The vertical velocity $w=-\int_0^z\divh v$ costs one horizontal derivative of $v$, so the transport terms $w\pz v$ and $w\pz\varrho$ do not fit into the usual $H^m$ energy calculus. We bypass this loss by an elliptic representation of $w$. Applying $\pz\divh$ to the momentum equation and using $\pz q=-g\varrho$, we obtain
\begin{equation*}
\mu\AB w=g\Laph\varrho-\pz\divh G,\qquad \AB=\Lap\pz^2,
\end{equation*}
with clamped boundary conditions, where $G=\rho(\dot v+fv^\perp)$ is the force in \eqref{eq:Gdef}. Testing the density equation with $-\sgm\Laph\varrho$ and inserting this representation produces the dissipation $\frac g\mu\nrm{\ABmh\Laph\varrho}_{L^2}^2$. In this way the density, which does not diffuse, inherits a dissipation from the gravitational restoring force; this mechanism was introduced by Dong, Li and Sun \cite{DongLiSun2026}. The inherited dissipation does not act on the horizontal mean, which therefore converges to a new profile instead of decaying.

\emph{Rotation.} The Coriolis force contributes the linear term $f\rs v^\perp$ to $G$. It drops out of the basic energy identity, since $v^\perp\cdot v=0$, but it enters the representation of $w$ through $(f/\mu)\ABmi\pz(\rs\curlh v)$: the vorticity feeds the vertical velocity. At low order we compare this term with the velocity dissipation $\mu\nrm{\nabla\nh v}^2$ and lower the weight of the density estimate accordingly (the constant $\kappa$ in \eqref{eq:lowenergy}). At top order the Coriolis term in the estimate with multiplier $\AB(\sgm\varrho)$ is bounded by $\Esh\Ff$, the product of the top-order energy and the decaying low-order dissipation. This term can be closed because $\int_0^\infty\Ff\dt$ is controlled by a positive power of the initial low-order energy (\Cref{lem:intFflat}). Finally, the horizontal mean of the velocity obeys the rotating balance $-\mu\pz^2\bar v+f\rs\bar v^\perp=-\ol{\rho\dot v}-f\ol{\varrho v^\perp}$, whose operator is coercive because the Coriolis term is skew-symmetric (\Cref{lem:Lf}).

\emph{Three dimensions.} The horizontal variable lives on the two-dimensional torus, where $H^1(\T^2)\not\subset L^\infty(\T^2)$. We compensate for this by three devices. We put the $L^\infty$ norm on a factor that has spare regularity in $H^4$; we put a decaying factor in $L^p_h$ for a large $p<\infty$, using $H^1(\T^2)\subset L^p(\T^2)$, at the price of a small fraction of a horizontal derivative on the other factor; and we bound the $L^\infty$ norm of decaying quantities through $H^{1+s}(\T^2)\subset L^\infty(\T^2)$ with $s>0$ small. With these devices the regularity level stays at $H^4$, and $\nrm{\nabla\bU}_{L^\infty}$ remains integrable in time: it is controlled by $\Ff^{7/12}\Fsh^{5/12}$.

\emph{Faster rates without $H^1(\T^2)\subset L^\infty$.} Two refinements circumvent the failure of $H^1(\T^2)\subset L^\infty(\T^2)$ in the long-time asymptotics. For the mean, the difference between the horizontal mean of the density and its rearrangement is quadratic in the oscillation, and we measure this quadratic error in the mixed norm $L^2_hL^\infty_z$. This norm is controlled by a vertical interpolation for the clamped oscillation alone, which gives the rate $(1+t)^{-13/8}$. For pointwise rates, an anisotropic interpolation inequality (\Cref{lem:endpoint}) exploits the decay of the horizontal gradient of the density, and not only of its $L^2$ norm. The same inequality, applied to the velocity gradient, improves the bound for the horizontal component of the flow map; the vertical component obeys a better bound, because it follows from the transported density.

\emph{Hierarchy of energies.} The estimates are organized in two tiers. A low-order energy $\Ef$, built on $\nh\varrho$, $\nh v$, $v_t$ and $w$, satisfies a closed dissipative inequality and decays like $(1+t)^{-1}$. A high-order functional $\Es$ controls $\varrho$ and $v$ in $H^4$, $v_t$ in $H^2$ and $v_{tt}$ in $L^2$. The top-order density estimate uses the three multipliers $\AB(\sgm\varrho)$, $\Laph^2\varrho$ and $\Laph\pz\nabla(\sgm\varrho)$, together with a perfect square that absorbs the linear inertial coupling.

\emph{Optimality.} In each horizontal Fourier mode the linearized problem leads to a quadratic eigenvalue problem, the pencil \eqref{eq:pencil}, in which rotation appears through a nonnegative Ekman-type operator $f^2R_\nu(\lambda)$. For stable stratification and real negative rates $\lambda=-s$ with $0<s\le\mu\pi^2/(4\max\rs)$, the pencil is a self-adjoint family that decreases strictly in $s$; the Ekman term preserves this monotonicity because $sf^2R_\nu(-s)$ is nonnegative and nondecreasing in $s$. This is the overdamped situation, in which min--max principles are available \cite{Duffin1955,Markus1988}. A comparison with the eigenvalues $\theta_n$ of \eqref{eq:theta-intro} produces exact real normal modes $e^{-s_n^*t}$ with $s_n^*/\theta_n\to1$. Elliptic regularity shows that their $H^4$ norm is of order $(s_n^*)^{-1}$ times their $L^2$ norm, which gives the lower bound $(1+t)^{-1}$ at $t=1/s_n^*$. For the nonlinear problem we take $\varepsilon$ times such a mode, corrected at second order in $\varepsilon$ to make the data compatible; the resulting solution stays within $O(\varepsilon^2(1+t))$ of the linear mode in $L^2$.

\emph{Instability.} When $\rs'>0$ somewhere, every nonzero horizontal mode has a real growing mode, obtained from a variational problem for the same pencil. The Ekman term lowers every growth rate, but it becomes negligible for short horizontal waves; hence the supremum of the growth rates is the same as without rotation, and it is still not attained. Since there is no fastest growing mode, we prove nonlinear instability with high-order approximate solutions in the spirit of Grenier \cite{Grenier2000}.

\subsection{Related work}\label{sec:lit}

\emph{Primitive equations.} Global well-posedness of strong solutions of the three-dimensional viscous primitive equations was proved by Cao and Titi \cite{CaoTiti2007} and, independently, by Kobelkov \cite{Kobelkov2006}; see also \cite{KukavicaZiane2007,HieberKashiwabara2016}. Cao, Li and Titi studied models with partial dissipation \cite{CaoLiTiti2014a,CaoLiTiti2014b,CaoLiTiti2016}; in particular, global well-posedness holds with only horizontal viscosity and diffusion \cite{CaoLiTiti2016}. Without viscosity the primitive equations may blow up in finite time \cite{CaoIbrahimNakanishiTiti2015}, and the hydrostatic Euler equations are ill-posed in Sobolev spaces \cite{Renardy2009,HanKwanNguyen2016}. This ill-posedness reflects the loss of a horizontal derivative in the vertical velocity, which also governs our analysis. The nonhomogeneous two-dimensional system is studied in \cite{WangJiuXu2021}, the stochastic system in \cite{DebusscheGlattHoltzTemamZiane2012}, and the stability of the slow manifold in \cite{TemamWirosoetisno2010}. In these works either the density (or temperature) diffuses, or gravity does not act on it; in \eqref{PE} the density does not diffuse, and it drives the flow through gravity.

\emph{Stratified fluids without diffusion.} Gravity stabilizes a stably stratified fluid even when the density does not diffuse, but the mechanism is indirect, since the density obeys a pure transport equation. For the Boussinesq equations without thermal diffusion, stability near hydrostatic equilibrium and the long-time behaviour were studied in \cite{DoeringWuZhaoZheng2018,TaoWuZhaoZheng2020,DongSun2022,JangKim2023,KangLeeNguyen2024}; see also \cite{HouLi2005,ElgindiWidmayer2015,CastroCordobaLear2019b} for related results, and \cite{YangLin2018,MasmoudiSaidHouariZhao2022,BedrossianBianchiniCotiZelatiDolce2023} for stratified shear flows. Dong, Li and Sun \cite{DongLiSun2026} proved asymptotic stability for the three-dimensional inhomogeneous Navier--Stokes equations driven by gravity, and their inherited dissipation is a key tool here. The models closest to \eqref{PE} are the incompressible porous medium equation and the Stokes--transport system, in which the velocity is slaved to the density. Asymptotic stability of stratified states for the porous medium equation is due to Elgindi \cite{Elgindi2017} and to Castro, C\'ordoba and Lear \cite{CastroCordobaLear2019}; see also \cite{KiselevYao2023}. For the Stokes--transport system in a channel, Dalibard, Guillod and Leblond \cite{DalibardGuillodLeblond2025} proved convergence at the rate $(1+t)^{-1}$ and identified the limit as the decreasing vertical rearrangement of the initial density; see also \cite{Leblond2022}. In \eqref{PE}, by contrast, the velocity has inertia, its horizontal mean is not slaved to the density, the vertical velocity costs a horizontal derivative, and rotation couples the divergent and the vortical parts of the flow.

\emph{Rotation.} At high rotation rates the Coriolis force is a source of dispersion and of improved regularity \cite{BabinMahalovNicolaenko1999,CheminDesjardinsGallagherGrenier2006,Majda2003}. In our problem it acts instead as a bounded skew-symmetric perturbation. Its effect on the long-time behaviour is quantitative: it changes rates and smallness thresholds, but not the limit state.

\emph{Rayleigh--Taylor instability.} Nonlinear Rayleigh--Taylor instability of inhomogeneous fluids was proved by Hwang and Guo \cite{HwangGuo2003} for the Euler equations and by F.~Jiang and S.~Jiang \cite{JiangJiang2014} for viscous flows (see also \cite{JiangJiangNi2013}), using the bootstrap argument of Guo and Strauss \cite{GuoStrauss1995}. That argument relies on a fastest growing mode. In the hydrostatic setting the maximal growth rate is not attained, and we use instead the high-order approximate solutions of Grenier \cite{Grenier2000}.

\emph{Energy and rearrangements.} The available potential energy of Lorenz \cite{Lorenz1955} is the potential energy in excess of that of the rearranged, stably stratified state; Winters, Lombard, Riley and D'Asaro \cite{WintersLombardRileyDAsaro1995} used it to separate reversible motion from irreversible mixing. Rearrangements are classical in analysis \cite{Kawohl1985,LiebLoss2001}. We use the $L^p$ contraction property of rearrangements due to Crowe, Zweibel and Rosenbloom \cite{CroweZweibelRosenbloom1986}.

\subsection{Organization of the paper}

\Cref{sec:prelim} collects the three-dimensional analytic toolbox: anisotropic embeddings, the operator $\AB$, the hydrostatic Stokes estimate and the rotating mean-flow operator. \Cref{sec:low} proves the low-order energy inequality, \Cref{sec:stokes} the Stokes estimates, and \Cref{sec:decay} the decay estimates. \Cref{sec:high} contains the higher-order estimates and the proof of \Cref{thm:main}. \Cref{sec:limit} proves the convergence to the limit state (\Cref{thm:limit}), and \cref{sec:energy} the variational and Lagrangian description of the relaxation (\Cref{thm:energy}). \Cref{sec:RT} treats unstable stratification (\Cref{thm:RT}); it also sets up the per-mode linear analysis that \Cref{sec:opt} uses to prove the optimality of the decay rate (\Cref{thm:opt}). \Cref{sec:opt} ends with a remark on the two-dimensional problem (\Cref{rem:2D}). \Cref{app:local} proves local well-posedness.

\subsection{Notation}

$\dot f\eqdef(\pt+\bU\cdot\nabla)f$ is the material derivative. For a function $f$ on $\Omega$, $\bar f(z)\eqdef|\Tl|^{-1}\int_{\Tl}f(x_h,z)\dxh$ is the horizontal mean and $\wt f\eqdef f-\bar f$ the oscillatory part. For an integer $m\ge1$, $\nh^mf$ denotes the tensor of all horizontal derivatives $\pa_{i_1}\cdots\pa_{i_m}f$, $i_j\in\{1,2\}$, and $|\nh^mf|^2=\sum|\pa_{i_1}\cdots\pa_{i_m}f|^2$; by Parseval in $x_h$, $\nrm{\nh^mf}_{L^2}=\nrm{|\nh|^mf}_{L^2}=\nrm{(-\Laph)^{m/2}f}_{L^2}$, and for real $s\ge0$ we write $|\nh|^s$ for the Fourier multiplier $|k|^s$. Mixed norms are written $L^p_hL^r_z$ (inner norm in $z$) and $L^r_zL^p_h$ (inner norm in $x_h$). We write
\begin{equation}\label{eq:Lambda}
\Lambda^{-1}\eqdef\nabla\pz\ABmi ,
\end{equation}
so that $\nrm{\Lambda^{-1}g}_{L^2}=\nrm{\ABmh g}_{L^2}$ (\Cref{lem:ABid}).

$\varrho_{\mathrm{in}}$ is the initial density \emph{perturbation} and $\rin\eqdef\rs+\vin$ the initial total density. By transport, $\inf_\Omega\rho(t)=\inf_\Omega\rin$, and on the data class used in the proof of \Cref{thm:main} one has $\rin\ge\frac12\min\rs$. The symbol $\lesssim$ means $\le C$ for a constant $C>0$ that may change from line to line and depends only on $\mu$, $g$, $f$, $l$, $\rs$, $\nrm{\rs+\vin}_{L^\infty}$ and a positive lower bound for $\rs+\vin$; on the data class $\nrm{\vin}_{H^4}+\nrm{\uin}_{H^4}\le1$, $\rin\ge\frac12\min\rs$, the last two quantities are bounded in terms of $\rs$. The numbered constants $\cc0,\dots,\cc7$ are fixed numbers of the same kind. Constants are uniform for $|f|$ in bounded sets.

\section{Preliminaries}\label{sec:prelim}

\subsection{Boundary conditions and elementary inequalities}

The system \eqref{PE} contains no evolution equation for the surface pressure, and it does not by itself enforce the boundary conditions on $w$ and on the density. The next lemma shows that the surface pressure is the Lagrange multiplier of the lid constraint. Once the pressure is chosen accordingly, the lid constraint and the clamped conditions on the density, imposed at $t=0$, persist for all times.

\begin{lemma}\label{lem:propagation}
Let $\varrho$ and $v$ be classical fields on $\Omega\times[0,T]$ with $v|_{\pa\Omega}=0$ and $\inf_\Omega\rho>0$, where $\rho\eqdef\rs+\varrho$. Define $w$ by \eqref{eq:w}, and suppose that \eqref{PE}$_1$ and \eqref{PE}$_2$ hold with $q=P(t,x_h)-g\int_0^z\varrho\dd s$ for a surface pressure $P$ of zero horizontal mean.
\begin{enumerate}[label=(\roman*)]
\item One has $\pt\divh\int_0^1v\dd z=\cN-\divh(a\nh P)$ with $\cN$ and $a$ given in \eqref{eq:Phieq} below; hence $\pt\divh\int_0^1v\dd z\equiv0$ on $(0,T)$ if and only if $P$ is the unique zero-mean solution of \eqref{eq:Peq}. If $P$ solves \eqref{eq:Peq} and \eqref{eq:lid0} holds, then $w|_{\pa\Omega}=0$ on $[0,T]$, and $w_t|_{\pa\Omega}=w_{tt}|_{\pa\Omega}=0$ whenever these traces exist.
\item If $P$ solves \eqref{eq:Peq}, \eqref{eq:lid0} holds and $\vin|_{\pa\Omega}=\pz\vin|_{\pa\Omega}=0$, then $\varrho|_{\pa\Omega}=\pz\varrho|_{\pa\Omega}=0$ on $[0,T]$.
\end{enumerate}
\end{lemma}

\begin{proof}
(i) By \eqref{eq:w}, $w(t,x_h,0)=0$. Set $\Phi(t,x_h)\eqdef\divh\int_0^1v\dd z$, so that $w(t,x_h,1)=-\Phi$. Dividing \eqref{PE}$_2$ by $\rho$ and inserting $q=P-g\int_0^z\varrho$,
\begin{equation*}
v_t+v\cdot\nh v+w\pz v+fv^\perp+\frac1\rho\nh P-\frac g\rho\int_0^z\nh\varrho\dd s-\frac\mu\rho\Lap v=0 .
\end{equation*}
Since $\divh v+\pz w=0$, $w|_{z=0}=0$ and $v|_{z=1}=0$, we have $\int_0^1(v\cdot\nh v+w\pz v)\dd z=\divh\int_0^1v\otimes v\dd z+[wv]_{z=0}^{z=1}=\divh\int_0^1v\otimes v\dd z$. With $a\eqdef\int_0^1\rho^{-1}\dd z$, integration in $z$ and application of $\divh$ give
\begin{equation}\label{eq:Phieq}
\begin{aligned}
&\Phi_t+\divh(a\nh P)=\cN,\\
&\cN\eqdef-\divh\divh\int_0^1v\otimes v\dd z-f\divh\int_0^1v^\perp\dd z\\&\qquad+g\,\divh\int_0^1\frac1\rho\int_0^z\nh\varrho\dd s\dd z+\mu\,\divh\int_0^1\frac{\Lap v}\rho\dd z .
\end{aligned}
\end{equation}
Each summand of $\cN$ is a horizontal divergence of a periodic field, so $\int_{\Tl}\cN\dxh=0$. Since $(\sup\rho)^{-1}\le a\le(\inf\rho)^{-1}$, the form $\int_{\Tl}a\nh\phi\cdot\nh\psi\dxh$ is coercive on zero-mean functions in $H^1(\Tl)$, and the Lax--Milgram lemma gives a unique zero-mean $P$ with
\begin{equation}\label{eq:Peq}
\divh(a\nh P)=\cN .
\end{equation}
This is the Lagrange multiplier of the lid constraint. Indeed, \eqref{eq:Phieq} shows that $\Phi_t=\cN-\divh(a\nh P)$; hence the lid constraint is preserved if and only if $P$ solves \eqref{eq:Peq}. At $t=0$, \eqref{eq:wt0} at $z=1$ is equivalent to \eqref{eq:Peq} at $t=0$; it determines the initial surface pressure, hence $\qin$, from $(\vin,\uin)$ and is not a restriction on the data. With \eqref{eq:Peq}, $\Phi_t\equiv0$, so $\Phi(t)=\Phi(0)=0$ by \eqref{eq:lid0} and $w(t,x_h,1)=0$. Differentiating the time-independent identity $w|_{\pa\Omega}=0$ gives the statements on $w_t,w_{tt}$.

(ii) On $\pa\Omega$ we have $v=0$, $w=0$ and $\pz w=-\divh v=0$. Hence \eqref{PE}$_1$ reduces to $\varrho_t=0$ on $\pa\Omega$, so $\varrho|_{\pa\Omega}=0$ for all $t$, and therefore $\nh\varrho|_{\pa\Omega}=0$. Differentiating \eqref{PE}$_1$ in $z$ and evaluating on $\pa\Omega$,
\begin{equation*}
\pz\varrho_t=-\pz v\cdot\nh\varrho-v\cdot\nh\pz\varrho-\pz w\,\pz\varrho-w\pz^2\varrho-\pz w\,\rs'-w\rs''=0 .
\end{equation*}
With $\pz\vin|_{\pa\Omega}=0$ this gives $\pz\varrho|_{\pa\Omega}=0$.
\end{proof}

We next record the Poincar\'e inequalities used throughout. The oscillatory part of a function is controlled by its horizontal gradient, and a function that vanishes on a lid is controlled by its vertical derivative. This is how $L^2$ norms are bounded by the dissipative norms in the energy estimates.

\begin{lemma}\label{lem:poincare}
Let $f\in H^1(\Omega)$.
\begin{enumerate}[label=(\arabic*)]
\item If $\bar f\equiv0$, then $\nrm{f}_{L^2}\le l\nrm{\nh f}_{L^2}$.
\item If $f$ vanishes on $\pa\Omega$ (or only on one of the lids $z=0$, $z=1$), then $\nrm{f}_{L^2}\le\pi^{-1}\nrm{\pz f}_{L^2}$ (respectively $\le\frac2\pi\nrm{\pz f}_{L^2}$).
\end{enumerate}
We write $C_P\eqdef\max\{l,1\}$; thus $\nrm{f}_{L^2}\le C_P\nrm{\nabla f}_{L^2}$ in both cases.
\end{lemma}

\begin{proof}
(1) follows from Parseval, since the nonzero wave vectors $k\in(l^{-1}\Z)^2$ satisfy $|k|\ge l^{-1}$. (2) is the one-dimensional Poincar\'e inequality on $(0,1)$ with Dirichlet (respectively mixed) conditions, applied for a.e.\ $x_h$.
\end{proof}

\begin{lemma}\label{lem:GN}
\begin{enumerate}[label=(\arabic*)]
\item Let $0\le s_1\le s_2$, $\theta\in[0,1]$ and $s=(1-\theta)s_1+\theta s_2$. Then $\nrm{f}_{H^s(\Omega)}\lesssim\nrm{f}_{H^{s_1}}^{1-\theta}\nrm{f}_{H^{s_2}}^\theta$.
\item For $f,g\in H^m(\Omega)$, $m\ge2$: $\nrm{fg}_{H^m}\lesssim\nrm{f}_{L^\infty}\nrm{g}_{H^m}+\nrm{f}_{H^m}\nrm{g}_{L^\infty}$.
\item (Three-dimensional embeddings.) $\nrm{f}_{L^6}\lesssim\nrm{f}_{H^1}$, $\nrm{f}_{L^4}^2\lesssim\nrm{f}_{L^2}^{\frac12}\nrm{f}_{H^1}^{\frac32}$, $\nrm{f}_{L^\infty}\lesssim\nrm{f}_{H^2}$; in particular $H^4(\Omega)\hookrightarrow W^{2,\infty}(\Omega)$ and $H^3(\Omega)\hookrightarrow W^{1,\infty}(\Omega)$.
\end{enumerate}
\end{lemma}

These are standard on the Lipschitz domain $\Omega$ (by extension to $\Tl\times\R$). The horizontal variable, however, lives on the two-dimensional torus, where $H^1$ is not embedded in $L^\infty$, and the isotropic embeddings of \Cref{lem:GN} are often too crude. The next lemma treats the horizontal and the vertical variables separately. It bounds traces on horizontal planes, gives $L^p$ bounds on $\Tl$ for $p<\infty$ and $L^\infty$ bounds at the cost of a fraction of a derivative, and uses the vanishing of a function on a lid to gain an $L^\infty$ bound in $z$.

\begin{lemma}\label{lem:aniso}
\begin{enumerate}[label=(\alph*)]
\item For $f\in H^1(\Omega)$, $\sup_{z}\nrm{f(\cdot,z)}_{L^2(\Tl)}^2\le\nrm{f}_{L^2}^2+2\nrm{f}_{L^2}\nrm{\pz f}_{L^2}$. If $f$ vanishes on one lid, $\sup_z\nrm{f(\cdot,z)}_{L^2(\Tl)}^2\le2\nrm{f}_{L^2}\nrm{\pz f}_{L^2}$.
\item For $g\in H^1(\Tl)$ and $2\le p<\infty$, $\nrm{g}_{L^p(\Tl)}\le C_p\nrm{g}_{L^2}^{\frac2p}\nrm{g}_{H^1}^{1-\frac2p}$. If $\int_{\Tl}g=0$, then $\nrm{g}_{L^p(\Tl)}\le C_p\nrm{g}_{L^2}^{\frac2p}\nrm{\nh g}_{L^2}^{1-\frac2p}\le C_p'\nrm{\nh g}_{L^2}$.
\item For $0<s\le1$ and $g\in H^{1+s}(\Tl)$, $\nrm{g}_{L^\infty(\Tl)}\le C_s\nrm{g}_{H^{1+s}}$. If $\int_{\Tl}g=0$, then $\nrm{g}_{L^\infty(\Tl)}\le C_s\nrm{\nh g}_{L^2}^{1-s}\nrm{\nh^2g}_{L^2}^s$.
\item If $f\in H^1(\Omega)$ vanishes on one lid, then $|f(x_h,z)|\le\nrm{\pz f(x_h,\cdot)}_{L^2(0,1)}$, and for $2\le p\le\infty$,
\begin{equation*}
\nrm{f}_{L^p_hL^\infty_z}\le\nrm{\pz f}_{L^p_hL^2_z}\le\nrm{\pz f}_{L^2_zL^p_h}.
\end{equation*}
\item For $f\in L^2_zH^{1+s}_h$, $0<s\le1$: $\nrm{f}_{L^\infty_hL^2_z}\le\nrm{f}_{L^2_zL^\infty_h}\le C_s\nrm{f}_{L^2_zH^{1+s}_h}\le C_s\nrm{f}_{H^{1+s}(\Omega)}$.
\item For $2\le q<\infty$ and $f\in L^2_zH^1_h$: $\nrm{f}_{L^2_zL^q_h}\le C_q\nrm{f}_{L^2_zH^{1-2/q}_h}\le C_q\nrm{f}_{H^{1-2/q}(\Omega)}$.
\end{enumerate}
\end{lemma}

\begin{proof}
(a) $\nrm{f(\cdot,z)}^2_{L^2_h}=\nrm{f(\cdot,z_0)}_{L^2_h}^2+2\int_{z_0}^z\int_{\Tl}f\pz f$; average in $z_0\in(0,1)$, or take $z_0$ on the lid where $f$ vanishes. (b) is the Gagliardo--Nirenberg (Ladyzhenskaya) inequality on the two-dimensional torus, combined with \Cref{lem:poincare} in the zero-mean case. (c) is the embedding $H^{1+s}(\T^2)\subset L^\infty(\T^2)$; in the zero-mean case $\nrm{g}_{H^{1+s}}\simeq\nrm{|\nh|^{1+s}g}_{L^2}\le\nrm{\nh g}^{1-s}\nrm{\nh^2g}^s$ by Parseval and H\"older. (d) follows from $f(x_h,z)=\int_{z_0}^z\pz f$ and the Cauchy--Schwarz inequality on $(0,1)$, then Minkowski's inequality ($p\ge2$). (e) is Minkowski's inequality followed by (c) for a.e.\ $z$. (f) follows from the embedding $H^{1-2/q}(\T^2)\subset L^q(\T^2)$ for a.e.\ $z$.
\end{proof}

The vertical velocity is not an independent unknown. It is recovered from the horizontal velocity through \eqref{eq:w}, at the cost of one horizontal derivative. The next lemma turns this formula into bounds on $w$ in the mixed norms of \Cref{lem:aniso}. These bounds are the source of the loss of derivatives described in \Cref{sec:ideas}.

\begin{lemma}\label{lem:w}
Let $v\in H^1(\Omega)$ and $w=-\int_0^z\divh v\dd z'$. Then, for $2\le p\le\infty$, $j\ge0$ and $0<s\le1$,
\begin{equation*}
\begin{aligned}
&\nrm{\nh^jw}_{L^p_hL^\infty_z}\le\nrm{\nh^j\divh v}_{L^2_zL^p_h},\qquad \nrm{\nh^jw}_{L^2}\le\nrm{\nh^jw}_{L^2_hL^\infty_z}\le\nrm{\nh^j\divh v}_{L^2},\\
&\nrm{w}_{L^\infty_zL^p_h}\le C_p\nrm{\nh\divh v}_{L^2}\quad(p<\infty),\qquad \nrm{\nh^jw}_{L^\infty}\le C_s\nrm{\nh^j\divh v}_{L^2_zH^{1+s}_h},
\end{aligned}
\end{equation*}
and $\nrm{w}_{H^k}\lesssim\nrm{v}_{H^{k+1}}$ for $k\ge0$.
\end{lemma}

\begin{proof}
$|\nh^jw(x_h,z)|\le\int_0^1|\nh^j\divh v(x_h,z')|\dd z'\le\nrm{\nh^j\divh v(x_h,\cdot)}_{L^2_z}$; then use Minkowski's inequality, and \Cref{lem:aniso}(b) and (c) (the function $\divh v(\cdot,z)$ has zero horizontal mean). The last bound follows from $\pz w=-\divh v$ and the first bounds.
\end{proof}

\begin{lemma}\label{lem:Linf}
\begin{enumerate}[label=(\roman*)]
\item If $f\in H^2(\Omega)$ vanishes on $\pa\Omega$, then $\nrm{f}_{H^2}\lesssim\nrm{\nabla^2f}_{L^2}$, and hence $\nrm{f}_{L^\infty}\lesssim\nrm{\nabla^2f}_{L^2}$.
\item If $\bar f\equiv0$ and $f$ vanishes on one lid, then for $p<\infty$ and $0<s\le1$,
\begin{equation*}
\nrm{f}_{L^\infty_zL^p_h}\le C_p\nrm{\nh\pz f}_{L^2},\qquad \nrm{f}_{L^\infty}\le C_s\nrm{\nh\pz f}_{L^2}^{1-s}\nrm{\nh^2\pz f}_{L^2}^s .
\end{equation*}
\end{enumerate}
\end{lemma}

\begin{proof}
(i) By \Cref{lem:poincare}, $\nrm{f}\le C_P\nrm{\nabla f}$, and $\nrm{\nabla f}^2=-\int_\Omega f\Lap f\le C_P\nrm{\nabla f}\nrm{\Lap f}$; together with $\nrm{\Lap f}\le\sqrt3\nrm{\nabla^2f}$ this gives the claim, and \Cref{lem:GN}(3) gives the $L^\infty$ bound. (ii) By \Cref{lem:aniso}(d), $\nrm{f}_{L^\infty_zL^p_h}\le\nrm{\pz f}_{L^2_zL^p_h}$, and $\pz f(\cdot,z)$ has zero mean; apply \Cref{lem:aniso}(b). Similarly $\nrm{f}_{L^\infty}\le\nrm{\pz f}_{L^2_zL^\infty_h}$, and \Cref{lem:aniso}(c) for a.e.\ $z$ followed by H\"older's inequality in $z$ gives the second bound.
\end{proof}

For pointwise decay rates we need an interpolation inequality that uses the horizontal gradient. In the long-time analysis, a function and its horizontal gradient decay at the same rate, while only its $H^4$ norm stays bounded; the vertical derivatives need not decay at this rate. The isotropic inequality $\nrm{h}_{L^\infty}\lesssim\nrm{h}_{L^2}^{5/8}\nrm{h}_{H^4}^{3/8}$ does not exploit the decay of $\nh h$. The next lemma does, and it gives larger exponents of the decaying norms. No boundary condition or mean-value condition is needed.

\begin{lemma}\label{lem:endpoint}
There is $C=C(l)$ such that every $h\in H^4(\Omega)$ satisfies
\begin{equation}\label{eq:endpoint}
\nrm{h}_{L^\infty(\Omega)}\le C\big(\nrm{h}_{L^2}+\nrm{\nh h}_{L^2}\big)^{\frac56}\nrm{h}_{H^4}^{\frac16},\qquad \nrm{\nabla h}_{L^\infty(\Omega)}\le C\big(\nrm{h}_{L^2}+\nrm{\nh h}_{L^2}\big)^{\frac12}\nrm{h}_{H^4}^{\frac12}.
\end{equation}
The same bounds hold componentwise for vector fields and tensors.
\end{lemma}

\begin{proof}
\emph{Extension.} Choose $a_1,\dots,a_4$ with $\sum_{j=1}^4a_j(-j)^m=1$ for $m=0,1,2,3$, and extend $h$ to $\Tl\times(-\frac14,\frac54)$ by $Eh(x_h,z)\eqdef\sum_ja_jh(x_h,-jz)$ for $z<0$ and $Eh(x_h,z)\eqdef\sum_ja_jh(x_h,1-j(z-1))$ for $z>1$. The derivatives of order at most three match at $z=0,1$, so $E$ is bounded from $H^m(\Omega)$ to $H^m(\Tl\times(-\frac14,\frac54))$ for $0\le m\le4$, and it commutes with horizontal derivatives. Let $\chi\in C_c^\infty(-\frac14,\frac54)$ with $\chi=1$ on $[0,1]$. We regard $\chi Eh$ as a $2\pi l$-periodic function on $\R^2\times\R$, and let $\zeta\in C_c^\infty(\R^2)$ equal $1$ on the period cell $[0,2\pi l]^2$. Set $H\eqdef\zeta\chi Eh$. Then $H=h$ on the cell, and, since the support of $\zeta$ meets a fixed finite number of cells,
\begin{equation*}
\nrm{\nh H}_{L^2(\R^3)}\le C\big(\nrm{h}_{L^2(\Omega)}+\nrm{\nh h}_{L^2(\Omega)}\big),\qquad \nrm{H}_{H^4(\R^3)}\le C\nrm{h}_{H^4(\Omega)};
\end{equation*}
the $L^2$ term comes from $\nh\zeta$.

\emph{Fourier estimate.} Let $A\eqdef\nrm{\nh H}_{L^2(\R^3)}$ and $B\eqdef\nrm{|D|^4H}_{L^2(\R^3)}\le C\nrm{H}_{H^4}$, where $|D|^4$ has the symbol $|\xi|^4$. If $A=0$, then $H$ does not depend on $x_h$ and has compact support, so $H=0$; if $B=0$, then $H=0$. Otherwise, for $R>0$ and $j=0,1$, Fourier inversion and the Cauchy--Schwarz inequality with the weight $R^6|\xi_h|^2+|\xi|^8$ give
\begin{equation*}
\nrm{D^jH}_{L^\infty}\le C\int_{\R^3}|\xi|^j|\hat H(\xi)|\dd\xi\le C\big(R^6A^2+B^2\big)^{\frac12}\Big(\int_{\R^3}\frac{|\xi|^{2j}\dd\xi}{R^6|\xi_h|^2+|\xi|^8}\Big)^{\frac12}.
\end{equation*}
The substitution $\xi=R\eta$ turns the last integral into $R^{2j-5}J_j$, with $J_j\eqdef\int_{\R^3}|\eta|^{2j}(|\eta_h|^2+|\eta|^8)^{-1}\dd\eta$. In polar coordinates $\eta=r\omega$, with $u=\omega_3$ and $|\omega_h|^2=1-u^2\ge1-|u|$,
\begin{equation*}
J_j=\int_0^\infty r^{2j}\int_{\mathbb S^2}\frac{\dd\omega}{|\omega_h|^2+r^6}\dd r,\qquad \int_{\mathbb S^2}\frac{\dd\omega}{|\omega_h|^2+r^6}\le4\pi\int_0^1\frac{\dd u}{1-u+r^6}\le C\min\{1+|\log r|,\ r^{-6}\}.
\end{equation*}
Hence $J_j\le C\int_0^1r^{2j}(1+|\log r|)\dd r+C\int_1^\infty r^{2j-6}\dd r<\infty$ for $j=0,1$. Choosing $R=(B/A)^{1/3}$, so that $R^6A^2=B^2$, we obtain $\nrm{H}_{L^\infty}\le CA^{5/6}B^{1/6}$ and $\nrm{\nabla H}_{L^\infty}\le CA^{1/2}B^{1/2}$. Restriction to the period cell proves \eqref{eq:endpoint}.
\end{proof}

\subsection{The anisotropic operator \texorpdfstring{$\AB=\Lap\pz^2$}{A_B}}

Throughout, functions on $\Omega$ are expanded as $f=\sum_{k\in(l^{-1}\Z)^2}e^{ik\cdot x_h}f_k(z)$, and we write $\nu=|k|$.

\begin{lemma}\label{lem:aniso4}
Let $\varphi\in H^m(\Omega)$, $m\ge0$. The problem $\Lap\pz^2h=\varphi$ in $\Omega$, $h=\pz h=0$ on $\pa\Omega$, has a unique solution $h\in H^{m+2}(\Omega)$, and $\nrm{h}_{H^{m+2}}\lesssim\nrm{\varphi}_{H^m}$.
\end{lemma}

\begin{proof}
Mode by mode, $(\pz^2-\nu^2)\pz^2h_k=\varphi_k$ on $(0,1)$ with $h_k=\pz h_k=0$ at $z=0,1$. For $k=0$ this is the clamped problem $\pz^4h_0=\varphi_0$, with $\nrm{h_0}_{H^{m+4}(0,1)}\lesssim\nrm{\varphi_0}_{H^m(0,1)}$. For $k\ne0$, testing with $\bar h_k$ and using $\nrm{h_k}\lesssim\nrm{\pz h_k}$ gives $\nrm{\pz^2h_k}^2+\nu^2\nrm{\pz h_k}^2\le\nrm{\varphi_k}\nrm{h_k}\lesssim\nrm{\varphi_k}\nrm{\pz h_k}$, hence
\begin{equation*}
\nu^2\nrm{h_k}+\nu^2\nrm{\pz h_k}+\nu\nrm{\pz^2h_k}\lesssim\nrm{\varphi_k},
\end{equation*}
where we used $\nu\ge l^{-1}$. For higher $z$-derivatives we use $\pz^jh_k=\pz^{j-4}\varphi_k+\nu^2\pz^{j-2}h_k$ ($j\ge4$) and $\nrm{\pz^3h_k}^2\lesssim\nrm{\pz^2h_k}\nrm{\pz^4h_k}+\nrm{\pz^2h_k}^2$. By induction, with $r\eqdef\max(a+j-2,0)$,
\begin{equation*}
\nu^a\nrm{\pz^jh_k}\lesssim\sum_{i\le r}\nu^{r-i}\nrm{\pz^i\varphi_k}\qquad\text{for }a+j\le m+2 .
\end{equation*} Summation over $k$ (Parseval) gives the claim; uniqueness follows from the energy identity.
\end{proof}

\begin{remark}\label{rem:onlytwo}
The gain is two derivatives, not four: the symbol $(|\xi_h|^2+\xi_z^2)\xi_z^2$ degenerates at $\xi_z=0$, and horizontal regularity is recovered only through the clamped conditions. For $\varphi=e^{ik\cdot x_h}$ one finds $\nrm{\pz^4h_k}_{L^2(0,1)}\sim\frac12\nu^{1/2}$ as $\nu\to\infty$ (boundary layers of width $\nu^{-1}$), so $\Lap\pz^2$ restricted to $H^4$ with clamped conditions is not closed. This is why $\AB$ is defined through its quadratic form.
\end{remark}

\begin{definition}\label{def:AB}
Let
\begin{equation*}
\begin{gathered}V\eqdef\big\{f\in L^2(\Omega):\pz f,\ \nh\pz f,\ \pz^2f\in L^2(\Omega),\ f|_{\pa\Omega}=\pz f|_{\pa\Omega}=0\big\},\\ \mathfrak a(f,g)\eqdef\int_\Omega\nabla\pz f\cdot\nabla\pz\bar g\dx .\end{gathered}
\end{equation*}
By the vertical Poincar\'e inequality, $\mathfrak a(f,f)\ge\pi^4\nrm{f}^2$ on $V$, and $V$ is complete for $\mathfrak a(f,f)^{1/2}$. Let $\AB$ be the positive self-adjoint operator associated with $\mathfrak a$ (Friedrichs extension): $f\in\cD(\AB)$ iff $f\in V$ and $g\mapsto\mathfrak a(f,g)$ is $L^2$-bounded; then $\AB f=\Lap\pz^2f$ in $\mathcal D'$. By \Cref{lem:aniso4}, $\cD(\AB)\subset H^2$, $\ABmi H^m\subset H^{m+2}$, and $\mathcal C\eqdef\{f\in H^4:f|_{\pa\Omega}=\pz f|_{\pa\Omega}=0\}$ is a core. Fractional powers are defined by functional calculus; equivalently $\AB^{-\alpha}=\Gamma(\alpha)^{-1}\int_0^\infty t^{\alpha-1}e^{-t\AB}\dt$ for $\alpha>0$, and $\AB^\alpha=(\AB^{-\alpha})^{-1}$.
\end{definition}

\begin{lemma}\label{lem:ABprop}
$\AB$ is self-adjoint and positive definite, $\ip{\AB f,f}=\nrm{\nabla\pz f}_{L^2}^2$ for $f\in\cD(\AB)$, $\cD(\ABh)=V$, $\ABmh$ is bounded on $L^2$, and every $\AB^{\pm\alpha}$ is self-adjoint.
\end{lemma}

\begin{lemma}\label{lem:ABmean}
\begin{enumerate}[label=(\roman*)]
\item (Fibre decomposition.) $\AB=\bigoplus_kA_k$ with $A_k=(\pz^2-|k|^2)\pz^2$ on $(0,1)$ with clamped conditions; every function of $\AB$ acts mode by mode. In particular, $\ABmi$ and $\AB^{-1/2}$ preserve functions of zero horizontal mean.
\item For $j=1,2$ and $f$ with $\pa_jf\in L^2$: $\pa_j\ABmi f=\ABmi\pa_jf$ and $\pa_j\ABmh f=\ABmh\pa_jf$; if $f,\pa_jf\in\cD(\ABh)$, then $\pa_j\ABh f=\ABh\pa_jf$. The same holds for $\Laph$ and every horizontal Fourier multiplier.
\end{enumerate}
\end{lemma}

\begin{proof}
The space $V$ and the form $\mathfrak a$ split orthogonally over $k$: $\mathfrak a(f,g)=|\Tl|\sum_k\int_0^1(\pz^2f_k\,\ol{\pz^2g_k}+|k|^2\pz f_k\,\ol{\pz g_k})\dd z$. Hence $\AB=\bigoplus_kA_k$, which gives (i), and horizontal derivatives act on the $k$-th fibre as multiplication by $ik_j$, which gives (ii).
\end{proof}

\begin{lemma}\label{lem:ABid}
\begin{enumerate}[label=(\roman*)]
\item For $g\in\cD(\ABh)$, $\nrm{\ABh g}_{L^2}=\nrm{\nabla\pz g}_{L^2}$; for $g\in L^2$, $\nrm{\nabla\pz\ABmh g}_{L^2}=\nrm{g}_{L^2}$ and $\nrm{\Lambda^{-1}g}_{L^2}=\nrm{\ABmh g}_{L^2}$.
\item For a horizontal vector field $Y\in L^2(\Omega;\R^2)$ and $g\in L^2$,
\begin{equation}\label{eq:ABdual}
\nrm{\ABmh\pz\divh Y}_{L^2}\le\nrm{Y}_{L^2},\quad \nrm{\ABmh\pz\pa_jg}_{L^2}\le\nrm{g}_{L^2},\quad \nrm{\ABmh\pz^2g}_{L^2}\le\nrm{g}_{L^2},
\end{equation}
where the left-hand sides are defined by duality, e.g.\ $\ip{\ABmh\pz\divh Y,\phi}\eqdef\ip{Y,\nh\pz\ABmh\phi}$.
\item For $g\in L^2$, $\nrm{\pz\ABmi g}_{L^2}\le\pi^{-1}\nrm{\ABmh g}_{L^2}$ and $\nrm{\ABmi g}_{L^2}\le\pi^{-2}\nrm{\ABmh g}_{L^2}$.
\end{enumerate}
\end{lemma}

\begin{proof}
(i) For $g\in\cD(\AB)$, two integrations by parts (with $g=\pz g=0$ on $\pa\Omega$) give $\nrm{\ABh g}^2=\ip{\AB g,g}=\nrm{\nabla\pz g}^2$, and density extends this to $\cD(\ABh)$. With $g$ replaced by $\ABmh g$ we get the second identity, and the third is the second applied to $\ABmh g$. (ii) For $\phi\in L^2$ let $\psi=\ABmh\phi$. Then $|\ip{Y,\nh\pz\psi}|\le\nrm{Y}\nrm{\nabla\pz\psi}=\nrm{Y}\nrm{\phi}$ by (i); the other two bounds are proved in the same way. (iii) $\pz\ABmi g$ vanishes on $\pa\Omega$, so $\nrm{\pz\ABmi g}\le\pi^{-1}\nrm{\pz^2\ABmi g}\le\pi^{-1}\nrm{\nabla\pz\ABmi g}=\pi^{-1}\nrm{\ABmh g}$ by (i); one more Poincar\'e inequality gives the second bound.
\end{proof}

\begin{lemma}\label{lem:meanosc}
For $f\in H^s(\Omega)$, $s\ge0$: $\int_{\Tl}\wt f\dxh=0$, $\ip{\bar f,\wt f}=0$, $\nh\bar f=0$, $\nh f=\nh\wt f$, and $\nrm{\wt f}_{H^s}+\nrm{\bar f}_{H^s}\le2\nrm{f}_{H^s}$. Moreover $\nrm{\bar f}_{L^2(\Omega)}\le\nrm{f}_{L^2(\Omega)}$ and $\nrm{\bar f}_{L^2(\Omega)}=|\Tl|^{1/2}\nrm{\bar f}_{L^2(0,1)}$.
\end{lemma}

The dissipation that the density inherits from gravity controls only the weak norm $\nrm{\ABmh\Laph\varrho}$. To turn it into decay of $\nrm{\nh\varrho}$, and to compare it with the dissipation of higher horizontal derivatives, we interpolate between this weak norm and stronger norms that remain bounded. The next lemma provides the required inequalities for clamped functions.

\begin{lemma}\label{lem:interp}
Let $f\in H^4(\Omega)$ with $f=\pz f=0$ on $\pa\Omega$. Then
\begin{align}
&\nrm{\nh f}_{L^2}^2\le\nrm{\nabla\pz f}_{L^2}\nrm{\ABmh\Laph f}_{L^2},\label{eq:2.9}\\
&\nrm{\nh f}_{L^2}^2\lesssim\nrm{\pz f}_{H^3}^{\frac23}\nrm{\ABmh\Laph f}_{L^2}^{\frac43},\label{eq:2.10}\\
&\nrm{\nabla\nh f}_{L^2}^2\lesssim\nrm{\ABmh\Laph f}_{L^2}\nrm{f}_{H^4}.\label{eq:2.11}
\end{align}
Moreover, for $f\in L^2$ and $2\le a\le5$,
\begin{equation}\label{eq:interpA}
\nrm{\ABmh|\nh|^af}_{L^2}\le\nrm{\ABmh\Laph f}_{L^2}^{\frac{5-a}3}\nrm{\ABmh|\nh|^5f}_{L^2}^{\frac{a-2}3}.
\end{equation}
\end{lemma}

\begin{proof}
\eqref{eq:2.9}: $\nrm{\nh f}^2=-\ip{f,\Laph f}=-\ip{\ABh f,\ABmh\Laph f}$ by \Cref{lem:ABmean}, and $\nrm{\ABh f}=\nrm{\nabla\pz f}$.

\eqref{eq:2.10}: it suffices to treat $\wt f$, since $\nh f=\nh\wt f$. For zero-mean $f$, $\nrm{f}\le l\nrm{\nh f}$, and the clamped conditions give $\nrm{\nabla\pz f}^2=\ip{\AB f,f}\le\nrm{f}\nrm{\Lap\pz^2f}\le\sqrt2\nrm{f}\nrm{\pz\nabla^3f}$. Hence $\nrm{\nabla\pz f}\lesssim\nrm{\nh f}^{1/2}\nrm{\pz f}_{H^3}^{1/2}$, and inserting this in \eqref{eq:2.9} gives \eqref{eq:2.10}. (This is not a plain Gagliardo--Nirenberg inequality: the intermediate bound $\nrm{\nabla\pz f}\lesssim\nrm{\nh f}^{1/2}\nrm{\pz f}_{H^3}^{1/2}$ fails for the clamped, horizontally constant function $f=z^2(1-z)^2$, which is why we first pass to $\wt f$.)

\eqref{eq:2.11}: $\nrm{\nabla\nh f}^2=\nrm{\nh^2f}^2+\nrm{\nh\pz f}^2$. First, $\nrm{\nh^2f}^2=\nrm{\Laph f}^2=\ip{\ABh\Laph f,\ABmh\Laph f}\le\nrm{f}_{H^4}\nrm{\ABmh\Laph f}$. Second, for each $k\ne0$, with $\nu=|k|$, the clamped conditions give $\nrm{\pz f_k}^2\le\nrm{f_k}\nrm{\pz^2f_k}$ and $\nrm{\pz^2f_k}^2\le\nrm{f_k}\nrm{\pz^4f_k}$, so
\begin{equation*}
X_k\eqdef\nrm{\pz^2f_k}^2+\nu^2\nrm{\pz f_k}^2\le\nrm{f_k}\,Y_k,\qquad Y_k\eqdef\nrm{\pz^4f_k}+\nu^2\nrm{\pz^2f_k},
\end{equation*}
while self-adjointness gives $\nrm{f_k}^2\le X_k^{1/2}\ip{A_k^{-1}f_k,f_k}^{1/2}$. Combining, $X_k^{3/4}\le\ip{A_k^{-1}f_k,f_k}^{1/4}Y_k$ and
\begin{equation*}
\nu^2\nrm{\pz f_k}^2\le\nu^2\nrm{f_k}X_k^{1/2}\le\nu^2X_k^{3/4}\ip{A_k^{-1}f_k,f_k}^{1/4}\le\nu^2\ip{A_k^{-1}f_k,f_k}^{1/2}Y_k .
\end{equation*}
By Parseval and the Cauchy--Schwarz inequality in $k$, $\nrm{\nh\pz f}^2\le\nrm{\ABmh\Laph f}\big(|\Tl|\sum_kY_k^2\big)^{1/2}\lesssim\nrm{\ABmh\Laph f}\nrm{f}_{H^4}$.

\eqref{eq:interpA}: with $a_k\eqdef\ip{A_k^{-1}f_k,f_k}\ge0$, $\nrm{\ABmh|\nh|^af}^2=|\Tl|\sum_k|k|^{2a}a_k$; write $|k|^{2a}=(|k|^4)^{\frac{5-a}3}(|k|^{10})^{\frac{a-2}3}$ and apply H\"older's inequality in $k$.
\end{proof}

In the higher-order estimates, the operator $\ABmh$ acts on products, and in some terms its output must be bounded in $L^\infty$. In three dimensions $\ABmh$ does not map $L^2$ into $L^\infty$ (\Cref{rem:Linfty}). The next lemma provides substitutes. It gives an $L^\infty$ bound in $z$ combined with an $L^p$ bound in $x_h$ for $p<\infty$, a duality bound for $\ABmh$ of a product, and a full $L^\infty$ bound that costs a fraction $s$ of a horizontal derivative.

\begin{lemma}\label{lem:product}
Let $2\le p<\infty$ and $q=\frac{2p}{p-2}$ ($q=\infty$ if $p=2$).
\begin{enumerate}[label=(\roman*)]
\item For $\psi\in L^2$ and $h=\ABmh\psi$,
\begin{equation*}
\nrm{h}_{L^\infty_zL^p_h}+\nrm{\pz h}_{L^2_zL^p_h}+\nrm{\nabla h}_{L^2}+\nrm{\nabla\pz h}_{L^2}\le C_p\nrm{\psi}_{L^2}.
\end{equation*}
\item For $\phi\in L^2$ and $\chi\in L^2_zL^q_h$, $\ABmh(\phi\chi)$, defined by duality as in \Cref{lem:ABid}(ii), satisfies $\nrm{\ABmh(\phi\chi)}_{L^2}\le C_p\nrm{\phi}_{L^2}\nrm{\chi}_{L^2_zL^q_h}$.
\item For $g\in L^2$: $\nrm{\ABmi g}_{L^\infty_zL^p_h}\le C_p\nrm{\ABmh g}_{L^2}$. If moreover $\bar g\equiv0$ and $0<s\le1$, then $\nrm{\ABmi g}_{L^\infty}\le C_s\nrm{\ABmh g}_{L^2}^{1-s}\nrm{\ABmh\nh g}_{L^2}^s$.
\end{enumerate}
\end{lemma}

\begin{proof}
(i) $h\in V$, so $h$ and $\pz h$ vanish on $\pa\Omega$, and $\nrm{\nabla\pz h}=\nrm{\psi}$ by \Cref{lem:ABid}. By \Cref{lem:aniso}(d) and (b), $\nrm{h}_{L^\infty_zL^p_h}\le\nrm{\pz h}_{L^2_zL^p_h}\le C_p(\nrm{\pz h}+\nrm{\nh\pz h})\le C_p\nrm{\psi}$, using $\nrm{\pz h}\le\pi^{-1}\nrm{\pz^2h}$. Finally $\nrm{\nabla h}\le\pi^{-1}\nrm{\nabla\pz h}$, since $\nabla h$ vanishes on $\pa\Omega$.

(ii) By self-adjointness,\begin{equation*}\nrm{\ABmh(\phi\chi)}=\sup_{\nrm\psi=1}\Big|\int_\Omega\phi\chi\,\ABmh\psi\Big|\le\nrm{\phi}\nrm{\chi}_{L^2_zL^q_h}\sup_{\nrm\psi=1}\nrm{\ABmh\psi}_{L^\infty_zL^p_h},\end{equation*} and (i) applies.

(iii) Apply (i) to $\psi=\ABmh g$. If $\bar g=0$, then $\ABmi g$ has zero horizontal mean and vanishes on $\pa\Omega$; \Cref{lem:Linf}(ii) gives $\nrm{\ABmi g}_{L^\infty}\le C_s\nrm{\nh\pz\ABmi g}^{1-s}\nrm{\nh^2\pz\ABmi g}^s$, and $\nrm{\nh\pz\ABmi g}\le\nrm{\nabla\pz\ABmi g}=\nrm{\ABmh g}$, $\nrm{\nh^2\pz\ABmi g}=\nrm{\nh\pz\ABmi\nh g}\le\nrm{\ABmh\nh g}$ by \Cref{lem:ABmean,lem:ABid}.
\end{proof}

\begin{remark}\label{rem:Linfty}
In three dimensions the bound $\nrm{\ABmh\psi}_{L^\infty}\lesssim\nrm{\psi}_{L^2}$ fails, since $H^1(\T^2)\not\subset L^\infty(\T^2)$, and \Cref{lem:product} is the replacement: $L^\infty_h$ is weakened to $L^p_h$ with $p<\infty$, and a full $L^\infty$ bound costs the fraction $s$ of a horizontal derivative.
\end{remark}

\subsection{Elliptic tools}

The horizontal velocity is recovered from the momentum equation by solving a Stokes problem. In it, the pressure depends only on $x_h$, and the incompressibility condition is replaced by the lid constraint. The next lemma provides elliptic regularity for this hydrostatic Stokes problem; see \cite{WangJiuXu2021,Ziane1995}. We use it to turn bounds on the time derivative of the velocity and on the density into $H^m$ bounds for $v$.

\begin{lemma}\label{lem:stokes}
Let $\zeta\in H^m(\Omega;\R^2)$, $m\ge0$. The problem
\begin{equation*}
-\mu\Lap v+\nh P=\zeta\ \text{in }\Omega,\qquad \divh\int_0^1v\dd z=0,\qquad v|_{\pa\Omega}=0,
\end{equation*}
with $P=P(x_h)$, has a solution $(v,P)$, unique up to an additive constant in $P$, with
\begin{equation*}
\nrm{v}_{H^{m+2}}+\nrm{\nh P}_{H^m(\Tl)}\le C\nrm{\zeta}_{H^m}.
\end{equation*}
\end{lemma}

This is the case $n=3$, constant depth, of the result of Ziane \cite{Ziane1995}. The constraint is the lid condition \eqref{eq:lid}.

\begin{theorem}[{\cite[Corollary 2.21]{GazzolaGrunauSweers2010}}]\label{thm:biharm}
For $f\in L^2(\Omega)$ the clamped problem $\Lap^2u=f$ in $\Omega$, $u=\pz u=0$ on $\pa\Omega$, has a unique solution $u\in H^4(\Omega)$, and $\nrm{u}_{H^4}\le C\nrm{f}_{L^2}$.
\end{theorem}

The horizontal mean of the velocity obeys a one-dimensional problem in $z$, the second line of \eqref{eq:decomp2} below, in which rotation adds the term $f\rs\phi^\perp$. The next lemma shows that this rotating mean-flow operator is coercive and invertible, because the Coriolis term is skew-symmetric. Hence the mean velocity is controlled for every value of $f$.

\begin{lemma}\label{lem:Lf}
For $\phi\in H^2(0,1;\R^2)$ with $\phi(0)=\phi(1)=0$ set $L_f\phi\eqdef-\mu\phi''+f\rs\phi^\perp$. Then $\ip{L_f\phi,\phi}_{L^2(0,1)}=\mu\nrm{\phi'}^2$, $L_f$ is an isomorphism onto $L^2(0,1;\R^2)$, and $\nrm{\phi}_{H^{m+2}(0,1)}\le C\nrm{L_f\phi}_{H^m(0,1)}$ for $0\le m\le3$, with $C$ depending on $\mu$, $f$ and $\nrm{\rs}_{C^3}$.
\end{lemma}

\begin{proof}
The identity holds because $\phi^\perp\cdot\phi=0$. Hence $\mu\pi^2\nrm{\phi}^2\le\mu\nrm{\phi'}^2\le\nrm{L_f\phi}\nrm{\phi}$, which gives injectivity, $\nrm{\phi}_{H^1}\lesssim\nrm{L_f\phi}$, and surjectivity by the Lax--Milgram lemma applied to the coercive form $\mu\ip{\phi',\psi'}+f\ip{\rs\phi^\perp,\psi}$ on $H^1_0$. Then $\mu\phi''=f\rs\phi^\perp-L_f\phi$ gives the $H^{m+2}$ bounds by induction.
\end{proof}

\subsection{Decomposition of the equations}

Averaging \eqref{PE}$_1$ over $\Tl$ and using $\bar w=-\int_0^z\ol{\divh v}=0$, $\divh v+\pz w=0$ and periodicity,
\begin{equation*}
\ol{\bU\cdot\nabla\varrho}=\ol{\divh(v\varrho)}+\pz\ol{w\varrho}=\pz\ol{w\wt\varrho},
\end{equation*}
so that
\begin{equation}\label{eq:meaneq}
\pt\bar\varrho=-\pz\big(\ol{w\wt\varrho}\big),
\end{equation}
and, subtracting,
\begin{equation*}
\pt\wt\varrho+\bU\cdot\nabla\wt\varrho+w\big(\pz\bar\varrho+\rs'\big)=\pz\big(\ol{w\wt\varrho}\big).
\end{equation*}
For the momentum equation we use the hydrostatic split $q=P(t,x_h)-g\int_0^z\varrho\dd s$ (\Cref{lem:propagation}) and the force $G$ of \eqref{eq:Gdef}. Since $\ol{\nh P}=0$ and $\ol{\rs v^\perp}=\rs\bar v^\perp$,
\begin{equation}\label{eq:decomp2}
\begin{cases}
-\mu\Laph v-\mu\wt{\pz^2v}+\nh P=-\wt G+g\int_0^z\nh\varrho(x_h,s)\dd s,\\
\mu\pz^2\bar v=\ol G,\quad\text{equivalently}\quad L_f\bar v=-\ol{\rho\dot v}-f\,\ol{\varrho v^\perp}.
\end{cases}
\end{equation}

The mean velocity $\bar v$ is not reached by the horizontal Poincar\'e inequality, and the horizontal dissipation does not act on it. The next lemma controls it through the second line of \eqref{eq:decomp2} and \Cref{lem:Lf}: by the material acceleration and, when $f\ne0$, by the oscillatory part of the velocity.

\begin{lemma}\label{lem:vbar}
For a classical solution of \eqref{PE}--\eqref{eq:BC},
\begin{equation*}
\nrm{\bar v}_{H^2(\Omega)}\le C\big(1+\nrm{\varrho}_{L^\infty}\big)^2\big(\nrm{\rho\dot v}_{L^2}+|f|\nrm{\nh v}_{L^2}\big),
\end{equation*}
and, for $m=1,2,3$, $\nrm{\bar v}_{H^{m+2}(\Omega)}\le C\nrm{\ol G}_{H^m(\Omega)}+C\nrm{\bar v}_{H^2(\Omega)}$.
\end{lemma}

\begin{proof}
Test the second line of \eqref{eq:decomp2} with $\bar v$ in $L^2(\Omega)$. Since $\bar v$ depends only on $z$, $\ip{\ol{\varrho v^\perp},\bar v}=\ip{\varrho v^\perp,\bar v}=\ip{\varrho\wt v^\perp,\bar v}$, because $\bar v^\perp\cdot\bar v=0$. Thus
\begin{equation*}
\mu\nrm{\pz\bar v}^2\le\big(\nrm{\rho\dot v}+|f|\nrm{\varrho}_{L^\infty}\nrm{\wt v}\big)\nrm{\bar v}\le\big(\nrm{\rho\dot v}+|f|l\nrm{\varrho}_{L^\infty}\nrm{\nh v}\big)\pi^{-1}\nrm{\pz\bar v}.
\end{equation*}
Hence $\nrm{\bar v}_{H^1}\lesssim\nrm{\rho\dot v}+|f|\nrm{\varrho}_{L^\infty}\nrm{\nh v}$. Then $\mu\pz^2\bar v=\ol G=\ol{\rho\dot v}+f\rs\bar v^\perp+f\ol{\varrho(\bar v+\wt v)^\perp}$ gives $\nrm{\pz^2\bar v}\lesssim\nrm{\rho\dot v}+|f|(1+\nrm{\varrho}_{L^\infty})(\nrm{\bar v}+\nrm{\wt v})$, which yields the first claim. The second follows by differentiating $\mu\pz^2\bar v=\ol G$.
\end{proof}

\section{Lower-order energy estimates}\label{sec:low}

\subsection{Energy functionals}

The higher-order energy is
\begin{equation}\label{eq:energies}
\Es^2(t)=\cEv^2(t)+\cEu^2(t)+\cEut^2(t)+\cEutt^2(t),
\end{equation}
where each component consists of an energy and a dissipation part:
\begin{equation*}
\begin{aligned}
&\cEv^2(t)\eqdef\sup_{\tau\in(0,t)}\Ev^2+\int_0^t\Fv^2\dtau,\quad \Ev\eqdef\nrm{\varrho}_{H^4},\\&\qquad \Fv\eqdef\nrm{\nabla\nh\pz\varrho}_{L^2}+\norm{\ABmh|\nh|^5\varrho}_{L^2}+\nrm{\nh^3\varrho}_{L^2},\\
&\cEu^2(t)\eqdef\sup_{\tau\in(0,t)}\Eu^2+\int_0^t\Fu^2\dtau,\quad \Eu\eqdef\nrm{v}_{H^4},\quad \Fu\eqdef\nrm{v}_{H^5},\\
&\cEut^2(t)\eqdef\sup_{\tau\in(0,t)}\Eut^2+\int_0^t\Fut^2\dtau,\quad \Eut\eqdef\nrm{v_t}_{H^2},\ \ \Fut\eqdef\nrm{v_t}_{H^3},\\
&\cEutt^2(t)\eqdef\sup_{\tau\in(0,t)}\Eutt^2+\int_0^t\Futt^2\dtau,\quad \Eutt\eqdef\nrm{v_{tt}}_{L^2},\ \ \Futt\eqdef\nrm{v_{tt}}_{H^1}.
\end{aligned}
\end{equation*}
We also write
\begin{equation}\label{eq:Esharpdef}
\Esh\eqdef\Ev+\Eu+\Eut+\Eutt,\qquad \Fsh\eqdef\Fv+\Fu+\Fut+\Futt,\qquad \Es^2(t)\simeq\sup_{\tau\in(0,t)}\Esh^2(\tau)+\int_0^t\Fsh^2\dtau .
\end{equation}
For the top-order velocity dissipation we use
\begin{equation*}
\Fuz\eqdef\nrm{\pz v}_{H^4},\qquad D_4\eqdef\nrm{\nabla\nh^4v}_{L^2},\qquad \cT_4^2\eqdef\nrm{\sqrt\rho\,\nh^4v}_{L^2}^2+g\nrm{\sqrt\sgm\,\nh^4\varrho}_{L^2}^2,
\end{equation*}
with $\sgm=-1/\rs'$. By positivity and boundedness of $\rho$ and $\sgm$, $\cT_4\lesssim\Eu+\Ev$. Since $v|_{\pa\Omega}=0$, the vertical Poincar\'e inequality applied to purely horizontal derivatives gives
\begin{equation}\label{eq:3.3}
\Eu^2\simeq\nrm{\pz v}_{H^3}^2+\nrm{\nh^4v}_{L^2}^2,\qquad \Fu^2\simeq\Fuz^2+\nrm{\nh^5v}_{L^2}^2\simeq\Fuz^2+D_4^2,
\end{equation}
and also $\nrm{v}_{H^m}\lesssim\nrm{\pz v}_{H^{m-1}}+\nrm{\nh^mv}_{L^2}\lesssim\nrm{\pz v}_{H^m}$ for $m\ge1$ (pure horizontal derivatives are not controlled by $\nrm{\pz v}_{H^{m-1}}$ alone, but $\nrm{\nh^jv}\le\pi^{-1}\nrm{\pz\nh^jv}$).

For the lower-order estimates let
\begin{equation*}
M\eqdef\pi^{-1}\nrm{\rs}_{L^\infty},\qquad \kappa\eqdef\min\Big\{\frac{\mu^2}{4M^2},\ \frac{g\pi^2\mu^2}{4f^2M^2}\Big\}\quad(\kappa\eqdef\tfrac{\mu^2}{4M^2}\ \text{if}\ f=0),
\end{equation*}
and define the low-order energy and dissipation
\begin{equation}\label{eq:lowenergy}
\begin{aligned}
\Ef^2(t)&\eqdef(\kappa+g)\nrm{\sqrt\sgm\nh\varrho}_{L^2}^2+\nrm{\sqrt{\rs}\nh v}_{L^2}^2+\frac1g\nrm{\sqrt\rho\,v_t}_{L^2}^2+\nrm{\sgm^{-1/2}w}_{L^2}^2,\\
\Ff^2(t)&\eqdef\mu\Big(\nrm{\nabla\nh v}_{L^2}^2+\frac1g\nrm{\nabla v_t}_{L^2}^2\Big)+\frac{g\kappa}\mu\norm{\ABmh\Laph\varrho}_{L^2}^2,
\end{aligned}
\end{equation}
and
\begin{equation}\label{eq:3.5}
\Eb^2(t)\eqdef\sup_{\tau\in(0,t)}\Ef^2+\int_0^t\Ff^2\dtau .
\end{equation}
The second entry in the definition of $\kappa$ is the only place where the rotation affects the structure of the low-order energy; it balances the vortical part of the vertical velocity against the velocity dissipation (see \eqref{eq:J23}). We finally collect degree-one quantities whose time integral is bounded (\Cref{prop:Fnat}):
\begin{equation}\label{eq:Fnat}
\Fn\eqdef\nrm{\nabla\bU}_{L^\infty}+\nrm{\nabla^2v}_{L^2}+\Ff^{\frac7{12}}\Fsh^{\frac5{12}} .
\end{equation}

\begin{remark}\label{rem:simplify}
The following elementary bounds are used throughout \cref{sec:low}--\cref{sec:high}, sometimes without reference. By definition,
\begin{equation}\label{eq:3.7}
\Ev,\Eu,\Eut,\Eutt\le\Esh,\qquad \Fv,\Fu,\Fut,\Futt\le\Fsh,\qquad \Fuz\le\Fu .
\end{equation}
Since $v$, $v_t$ vanish on $\pa\Omega$,
\begin{equation}\label{eq:3.8}
\nrm{\nabla\nh v}_{L^2}\le\nrm{v}_{H^2}\le\Eu\le\Fu,\qquad \nrm{\nabla v_t}_{L^2}\le\Eut\le\Fut .
\end{equation}
By \Cref{lem:ABprop} and \Cref{lem:poincare}, using that $\Laph\varrho$ has zero horizontal mean and $\nrm{\Laph\varrho}=\nrm{\nh^2\varrho}$,
\begin{equation}\label{eq:3.10}
\nrm{\ABmh\Laph\varrho}_{L^2}\le\pi^{-2}\nrm{\nh^2\varrho}_{L^2}\le\pi^{-2}\min\big\{\Ev,\ l\,\nrm{\nh^3\varrho}_{L^2}\big\},
\end{equation}
hence
\begin{equation}\label{eq:3.11}
\Ff\lesssim\Esh,\qquad \Ff\lesssim\Fsh .
\end{equation}
Since $\rs\in C^6$ and $H^4(\Omega)\hookrightarrow W^{2,\infty}(\Omega)$,
\begin{equation}\label{eq:3.12}
\nrm{\nabla\rho}_{L^\infty}+\nrm{\nabla^2\rho}_{L^\infty}\lesssim1+\Ev,\qquad \nrm{\varrho}_{W^{2,\infty}}\lesssim\Ev .
\end{equation}
By \Cref{lem:Linf}, \Cref{lem:w} and \Cref{lem:aniso},
\begin{equation}\label{eq:3.13}
\begin{aligned}
&\nrm{v}_{L^\infty}\lesssim\nrm{\nabla^2v}_{L^2},\qquad \nrm{w}_{L^\infty}+\nrm{\nabla v}_{L^\infty}\lesssim\nrm{v}_{H^3},\qquad \nrm{\nabla\bU}_{L^\infty}\lesssim\nrm{v}_{H^4}=\Eu,\\
&\nrm{w}_{L^2}+\nrm{w}_{L^2_hL^\infty_z}\le\nrm{\divh v}_{L^2},\qquad \nrm{w}_{H^1}+\nrm{w}_{L^\infty_zL^p_h}\le C_p\nrm{\nabla\nh v}_{L^2}\lesssim_p\Ff\quad(p<\infty).
\end{aligned}
\end{equation}
Finally, by \Cref{lem:vbar}, \Cref{lem:poincare} and the bound $\nrm{\bU\cdot\nabla v}_{L^2}\le\nrm{v}_{L^\infty}\nrm{\nh v}+\nrm{w}_{L^2_hL^\infty_z}\nrm{\pz v}_{L^\infty_hL^2_z}\lesssim\nrm{v}_{H^3}\nrm{\nh v}_{L^2}$ (\Cref{lem:aniso}(e)),
\begin{equation}\label{eq:3.13b}
\nrm{\bU\cdot\nabla v}_{L^2}\lesssim\Eu\Ff,\qquad \nrm{\bar v}_{H^2}+\nrm{v}_{H^1}\lesssim(1+\Esh)^3\Ff .
\end{equation}
Indeed $\nrm{\rho\dot v}\lesssim(1+\Ev)(\nrm{v_t}+\Eu\Ff)$, $\nrm{\wt v}_{H^1}\le C_P\nrm{\nabla\nh v}$, and $\nrm{v_t}\le\pi^{-1}\nrm{\nabla v_t}$.
\end{remark}

\subsection{The basic energy identity}

\begin{proposition}\label{prop:basic}
For a classical solution of \eqref{PE}--\eqref{eq:BC} on $(0,T)$, $\nrm{\rho(t)}_{L^\infty}=\nrm{\rin}_{L^\infty}$ and
\begin{equation}\label{eq:3.15}
\frac12\ddt\Big(g\nrm{\sqrt\sgm\varrho}_{L^2}^2+\nrm{\sqrt\rho\,v}_{L^2}^2\Big)+\mu\nrm{\nabla v}_{L^2}^2=\frac g2\int_\Omega w\,\pz\sgm\,\varrho^2\dx .
\end{equation}
Consequently, there is $\cc0>0$ depending only on $\rs$ and $g$ such that
\begin{equation}\label{eq:3.16}
g\nrm{\sqrt\sgm\varrho(t)}_{L^2}^2+\nrm{\sqrt\rho\,v(t)}_{L^2}^2\le\Big(g\nrm{\sqrt\sgm\vin}_{L^2}^2+\nrm{\sqrt{\rin}\,\uin}_{L^2}^2\Big)\exp\Big(\cc0\int_0^t\nrm{w}_{L^\infty}\dtau\Big).
\end{equation}
\end{proposition}

\begin{proof}
Testing \eqref{PE}$_1$ with $\sgm\varrho$ and using $\divv\bU=0$, $\bU\cdot n=0$ gives $\frac12\ddt\nrm{\sqrt\sgm\varrho}^2=\int_\Omega\varrho w+\frac12\int_\Omega w\pz\sgm\varrho^2$. Testing \eqref{PE}$_2$ with $v$: the Coriolis term vanishes since $v^\perp\cdot v=0$; $\int\rho\dot v\cdot v=\frac12\ddt\nrm{\sqrt\rho v}^2$ because $\rho_t+\bU\cdot\nabla\rho=0$; and $\int\nh q\cdot v=-\int q\divh v=\int q\pz w=-\int\pz q\,w=g\int\varrho w$, the boundary terms vanishing by \Cref{lem:propagation}. Adding $g$ times the first identity to the second, the coupling $\int\varrho w$ cancels. The bound \eqref{eq:3.16} follows from $|\frac g2\int w\pz\sgm\varrho^2|\le\frac{\cc0}2\nrm{w}_{L^\infty}\cdot g\nrm{\sqrt\sgm\varrho}^2$ and Gronwall's inequality.
\end{proof}

\begin{remark}\label{rem:c0}
Since $w$ vanishes on the lids and $\pz w=-\divh v$, $\nrm{w}_{L^\infty}\le\nrm{\pz w}_{L^\infty}=\nrm{\divh v}_{L^\infty}\le\nrm{\nabla\bU}_{L^\infty}\le\Fn$; the time integral of $\nrm{w}_{L^\infty}$ in \eqref{eq:3.16} is therefore bounded by $\int_0^t\Fn\dtau$, which is finite under the hypotheses of \Cref{prop:Fnat}.
\end{remark}

\subsection{The first-order energy identities}

Since the horizontal mean $\bar\varrho$ need not decay, we estimate $\nh\varrho$.

\emph{Momentum--density coupling.} Testing \eqref{PE}$_1$ with $-g\sgm\Laph\varrho$ gives, since $\sgm$ depends only on $z$,
\begin{equation*}
\frac g2\ddt\nrm{\sqrt\sgm\nh\varrho}_{L^2}^2=g\int_\Omega\sgm\,\bU\cdot\nabla\varrho\,\Laph\varrho\dx-g\int_\Omega w\Laph\varrho\dx .
\end{equation*}
Testing \eqref{PE}$_2$ with $-\Laph v$: since $\rs$ depends only on $z$, $\int\rs v_t\cdot(-\Laph v)=\frac12\ddt\nrm{\sqrt\rs\nh v}^2$; the linear Coriolis term vanishes, $\int f\rs v^\perp\cdot(-\Laph v)=f\int\rs\pa_jv^\perp\cdot\pa_jv=0$; and, by \Cref{lem:propagation},
\begin{equation*}
\int_\Omega\nh q\cdot(-\Laph v)\dx=\int_\Omega q\,\Laph\divh v\dx=-\int_\Omega q\,\Laph\pz w\dx=\int_\Omega\pz q\,\Laph w\dx=-g\int_\Omega w\Laph\varrho\dx .
\end{equation*}
Adding the two identities, the coupling $g\int w\Laph\varrho$ cancels:
\begin{equation}\label{eq:3.19}
\frac12\ddt\Big(\nrm{\sqrt\rs\nh v}_{L^2}^2+g\nrm{\sqrt\sgm\nh\varrho}_{L^2}^2\Big)+\mu\nrm{\nabla\nh v}_{L^2}^2=J_{u,1}+J_{u,2}+J_{u,3}+gJ_{\varrho,1},
\end{equation}
with
\begin{equation*}
\begin{gathered}J_{u,1}\eqdef\int_\Omega\varrho\,v_t\cdot\Laph v\dx,\qquad J_{u,2}\eqdef\int_\Omega\rho(\bU\cdot\nabla v)\cdot\Laph v\dx,\\ J_{u,3}\eqdef-f\int_\Omega\pa_j\varrho\,v^\perp\cdot\pa_jv\dx,\qquad J_{\varrho,1}\eqdef\int_\Omega\sgm\,\bU\cdot\nabla\varrho\,\Laph\varrho\dx .\end{gathered}
\end{equation*}
Here $J_{u,3}$ is the nonlinear part of the Coriolis force: $\int f\rho v^\perp\cdot(-\Laph v)=f\int\pa_j\varrho\, v^\perp\cdot\pa_jv$, because $\rho\,\pa_jv^\perp\cdot\pa_jv=0$ and $\pa_j\rs=0$ for horizontal $j$.

\emph{Inherited dissipation.} Apply $\pz\divh$ to \eqref{PE}$_2$, written as $G+\nh q-\mu\Lap v=0$ with $G$ from \eqref{eq:Gdef}, and use $\divh v=-\pz w$ and $\pz q=-g\varrho$:
\begin{equation*}
\pz\divh G-g\Laph\varrho+\mu\Lap\pz^2w=0 .
\end{equation*}
Since $w=\pz w=0$ on $\pa\Omega$, this gives the representation
\begin{equation}\label{eq:wrep}
\begin{aligned}
w&=\frac1\mu\ABmi\big[g\Laph\varrho-\pz\divh G\big]=W^{\mathrm L_1}+W^{\mathrm L_2}+W^{\mathrm L_3}+W^{\mathrm{NL}},\\
W^{\mathrm L_1}&\eqdef\frac g\mu\ABmi\Laph\varrho,\qquad W^{\mathrm L_2}\eqdef-\frac1\mu\ABmi\pz\divh(\rs v_t),\qquad W^{\mathrm L_3}\eqdef\frac f\mu\ABmi\pz(\rs\,\curlh v),\\
W^{\mathrm{NL}}&\eqdef-\frac1\mu\ABmi\pz\divh\big(\varrho\dot v+\rs\bU\cdot\nabla v+f\varrho v^\perp\big),
\end{aligned}
\end{equation}
where we used $\divh(\rs v^\perp)=-\rs\curlh v$. Differentiating \eqref{PE}$_2$ in $z$ also gives
\begin{equation}\label{eq:3.22}
g\nh\varrho=-\mu\Lap\pz v+\pz G .
\end{equation}
Testing \eqref{PE}$_1$ with $-\sgm\Laph\varrho$, inserting \eqref{eq:wrep} into $-\int w\Laph\varrho$ and using $\ip{\ABmi\Laph\varrho,\Laph\varrho}=\nrm{\ABmh\Laph\varrho}^2$, we obtain
\begin{equation}\label{eq:3.23}
\frac12\ddt\nrm{\sqrt\sgm\nh\varrho}_{L^2}^2+\frac g\mu\nrm{\ABmh\Laph\varrho}_{L^2}^2=J_{\varrho,1}+J_{\varrho,21}+J_{\varrho,22}+J_{\varrho,23},
\end{equation}
where, after moving $\ABmi$ to the other factor and integrating by parts once in $z$ (using $\ABmi\Laph\varrho=0$ on $\pa\Omega$) and once horizontally,
\begin{equation*}
\begin{aligned}
J_{\varrho,21}&\eqdef\frac1\mu\int_\Omega\rs v_t\cdot\nh\pz\ABmi\Laph\varrho\dx,\qquad J_{\varrho,23}\eqdef\frac f\mu\int_\Omega\rs\curlh v\;\pz\ABmi\Laph\varrho\dx,\\
J_{\varrho,22}&\eqdef\frac1\mu\int_\Omega\big(\varrho\dot v+\rs\bU\cdot\nabla v+f\varrho v^\perp\big)\cdot\nh\pz\ABmi\Laph\varrho\dx .
\end{aligned}
\end{equation*}

\emph{The time derivative.} Differentiating \eqref{PE}$_2$ in $t$,
\begin{equation*}
\rho v_{tt}+\rho_tv_t+\pt(\rho\bU\cdot\nabla v)+f\rho_tv^\perp+f\rho v_t^\perp+\nh q_t-\mu\Lap v_t=0 .
\end{equation*}
Test with $v_t$. The term $f\rho v_t^\perp\cdot v_t$ vanishes, and $\int\nh q_t\cdot v_t=g\int\varrho_tw_t=g\int(\sgm^{-1}w-\bU\cdot\nabla\varrho)w_t$ by \Cref{lem:propagation}. Hence
\begin{equation}\label{eq:3.24}
\frac12\ddt\Big(\nrm{\sqrt\rho\,v_t}_{L^2}^2+g\nrm{\sgm^{-1/2}w}_{L^2}^2\Big)+\mu\nrm{\nabla v_t}_{L^2}^2=J_{t,1}+J_{t,2}+J_{t,3}+J_{t,4},
\end{equation}
with
\begin{equation*}
\begin{gathered}J_{t,1}\eqdef-\frac12\int_\Omega\varrho_t|v_t|^2\dx,\qquad J_{t,2}\eqdef-\int_\Omega\pt(\rho\bU\cdot\nabla v)\cdot v_t\dx,\\ J_{t,3}\eqdef g\int_\Omega\bU\cdot\nabla\varrho\,w_t\dx,\qquad J_{t,4}\eqdef-f\int_\Omega\varrho_t\,v^\perp\cdot v_t\dx .\end{gathered}
\end{equation*}

\emph{The combination.} Multiplying \eqref{eq:3.23} by $\kappa$, adding \eqref{eq:3.19} and $\frac1g$\eqref{eq:3.24}, we obtain
\begin{equation}\label{eq:3.25}
\frac12\ddt\Ef^2+\Ff^2=(\kappa+g)J_{\varrho,1}+\kappa\big(J_{\varrho,21}+J_{\varrho,22}+J_{\varrho,23}\big)+J_{u,1}+J_{u,2}+J_{u,3}+\frac1g\big(J_{t,1}+J_{t,2}+J_{t,3}+J_{t,4}\big).
\end{equation}
By \Cref{lem:poincare} and the bounds on $\rho$, $\sgm$,
\begin{equation}\label{eq:3.27}
\Ef^2\simeq\nrm{\nh\varrho}_{L^2}^2+\nrm{\nh v}_{L^2}^2+\nrm{v_t}_{L^2}^2+\nrm{w}_{L^2}^2,\qquad \Ff^2\simeq\nrm{\nabla\nh v}_{L^2}^2+\nrm{\nabla v_t}_{L^2}^2+\nrm{\ABmh\Laph\varrho}_{L^2}^2 .
\end{equation}

\subsection{The linear terms}

The two linear terms on the right of \eqref{eq:3.25} are $\kappa J_{\varrho,21}$ (inertia) and $\kappa J_{\varrho,23}$ (rotation). By \Cref{lem:ABid}, $\nrm{\nh\pz\ABmi\Laph\varrho}\le\nrm{\ABmh\Laph\varrho}$ and $\nrm{\pz\ABmi\Laph\varrho}\le\pi^{-1}\nrm{\ABmh\Laph\varrho}$. With $\nrm{v_t}\le\pi^{-1}\nrm{\nabla v_t}$ and $\nrm{\curlh v}\le\nrm{\nh v}\le\pi^{-1}\nrm{\nabla\nh v}$,
\begin{align}
\kappa|J_{\varrho,21}|&\le\frac{\kappa M}\mu\nrm{\nabla v_t}\nrm{\ABmh\Laph\varrho}\le\frac{\kappa g}{4\mu}\nrm{\ABmh\Laph\varrho}^2+\frac{\kappa M^2}{g\mu}\nrm{\nabla v_t}^2,\label{eq:J21}\\
\kappa|J_{\varrho,23}|&\le\frac{\kappa|f|M}{\pi\mu}\nrm{\nabla\nh v}\nrm{\ABmh\Laph\varrho}\le\frac{\kappa g}{4\mu}\nrm{\ABmh\Laph\varrho}^2+\frac{\kappa f^2M^2}{g\pi^2\mu}\nrm{\nabla\nh v}^2.\label{eq:J23}
\end{align}
By the choice of $\kappa$, $\frac{\kappa M^2}{g\mu}\le\frac\mu{4g}$ and $\frac{\kappa f^2M^2}{g\pi^2\mu}\le\frac\mu4$. Hence
\begin{equation}\label{eq:linear}
\kappa\big(|J_{\varrho,21}|+|J_{\varrho,23}|\big)\le\frac12\Ff^2 .
\end{equation}
These are the only terms in \eqref{eq:3.25} that are not small.

\subsection{The nonlinear terms}

Every remaining term is bounded by a polynomial in $\Esh$ times $\Ff^2$. We use \eqref{eq:3.13}, \eqref{eq:3.13b} and the interpolation bounds of \Cref{lem:interp},
\begin{equation}\label{eq:B1}
\nrm{\nh\varrho}_{L^2}^2\le\Ev\nrm{\ABmh\Laph\varrho}\lesssim\Esh\Ff,\qquad \nrm{\nabla\nh\varrho}_{L^2}^2\lesssim\Esh\Ff,\qquad \nrm{\Lambda^{-1}\nh\varrho}_{L^2}\le l\nrm{\ABmh\Laph\varrho}\lesssim\Ff,
\end{equation}
where the last one follows from \Cref{lem:poincare} applied to the zero-mean functions $\Lambda^{-1}\pa_i\varrho$ and from $\sum_{i,j}\nrm{\ABmh\pa_i\pa_j\varrho}^2=\nrm{\ABmh\Laph\varrho}^2$.

\emph{The term $J_{\varrho,1}$.} Write $J_{\varrho,1}=J_{\varrho,11}+J_{\varrho,12}$ with the horizontal and vertical transport. Integrating by parts horizontally,
\begin{equation*}
J_{\varrho,11}=\int_\Omega\sgm(v\cdot\nh\varrho)\Laph\varrho=-\int_\Omega\sgm\,\pa_jv_i\,\pa_i\varrho\,\pa_j\varrho+\frac12\int_\Omega\sgm\divh v\,|\nh\varrho|^2 .
\end{equation*}
By \Cref{lem:GN}(3) and \eqref{eq:B1}, $\nrm{\nh\varrho}_{L^4}^2\lesssim\nrm{\nh\varrho}^{1/2}(\nrm{\nh\varrho}+\nrm{\nabla\nh\varrho})^{3/2}\lesssim\Esh\Ff$, so
\begin{equation*}
|J_{\varrho,11}|\lesssim\nrm{\nh v}_{L^2}\nrm{\nh\varrho}_{L^4}^2\lesssim\Esh\Ff^2 .
\end{equation*}
For $J_{\varrho,12}$ we use that $\sgm w\pz\varrho$ and its $z$-derivative vanish on $\pa\Omega$ (because $w=\pz w=\pz\varrho=0$ there), so that
\begin{equation*}
J_{\varrho,12}=\int_\Omega(\sgm w\pz\varrho)\,\AB\ABmi\Laph\varrho=\int_\Omega\nabla\pz(\sgm w\pz\varrho)\cdot\Lambda^{-1}\Laph\varrho .
\end{equation*}
Expanding $\nabla\pz(\sgm w\pz\varrho)=\nabla\pz w\,\sgm\pz\varrho+\pz w\,\nabla(\sgm\pz\varrho)+\nabla w\,\pz(\sgm\pz\varrho)+w\,\nabla\pz(\sgm\pz\varrho)$, we estimate the four products in $L^2$:
\begin{itemize}
\item $\nrm{\nabla\pz w}\nrm{\sgm\pz\varrho}_{L^\infty}\lesssim\nrm{\nabla\divh v}\Ev\lesssim\Esh\Ff$;
\item $\nrm{\pz w}_{L^4}\nrm{\nabla(\sgm\pz\varrho)}_{L^4}\lesssim\nrm{\divh v}_{H^1}\Ev\lesssim\Esh\Ff$;
\item $\nrm{\nabla w}_{L^2_hL^\infty_z}\nrm{\pz(\sgm\pz\varrho)}_{L^\infty_hL^2_z}\lesssim(\nrm{\nh\divh v}+\nrm{\pz\divh v})\nrm{\varrho}_{H^{7/2}}\lesssim\Esh\Ff$, by \Cref{lem:w}, \Cref{lem:aniso}(d) for $\divh v$ (which vanishes on $\pa\Omega$) and \Cref{lem:aniso}(e);
\item $\nrm{w}_{L^\infty_zL^4_h}\nrm{\nabla\pz(\sgm\pz\varrho)}_{L^2_zL^4_h}\lesssim\Ff\nrm{\varrho}_{H^{7/2}}\lesssim\Esh\Ff$, by \eqref{eq:3.13} and \Cref{lem:aniso}(f).
\end{itemize}
Hence $|J_{\varrho,12}|\lesssim\Esh\Ff\nrm{\Lambda^{-1}\Laph\varrho}\lesssim\Esh\Ff^2$, and $|J_{\varrho,1}|\lesssim\Esh\Ff^2$.

\emph{The term $J_{\varrho,22}$.} By \eqref{eq:3.13b} and \Cref{lem:ABid},
\begin{equation*}
|J_{\varrho,22}|\lesssim\big(\nrm{\varrho}_{L^\infty}\nrm{v_t}+(1+\Ev)\nrm{\bU\cdot\nabla v}+|f|\nrm{\varrho}_{L^\infty}\nrm{v}\big)\nrm{\ABmh\Laph\varrho}\lesssim\Esh(1+\Esh)^3\Ff^2 .
\end{equation*}

\emph{The terms $J_{u,i}$.} Using $\nrm{\Laph v}\le\sqrt2\nrm{\nabla\nh v}$, \eqref{eq:3.13b} and $H^4\hookrightarrow W^{1,\infty}$,
\begin{equation*}
\begin{aligned}
|J_{u,1}|&\le\nrm{\varrho}_{L^\infty}\nrm{v_t}\nrm{\Laph v}\lesssim\Esh\Ff^2,\qquad |J_{u,2}|\le\nrm{\rho}_{L^\infty}\nrm{\bU\cdot\nabla v}\nrm{\Laph v}\lesssim(1+\Esh)\Esh\Ff^2,\\
|J_{u,3}|&\le|f|\nrm{\nh\varrho}_{L^\infty}\nrm{v}\nrm{\nh v}\lesssim\Esh(1+\Esh)^3\Ff^2 .
\end{aligned}
\end{equation*}

\emph{The terms $J_{t,1}$, $J_{t,2}$ and $J_{t,4}$.} From \eqref{PE}$_1$, $\varrho_t=\sgm^{-1}w-\bU\cdot\nabla\varrho$, so $\nrm{\varrho_t}_{H^1}\lesssim\nrm{w}_{H^1}+\nrm{\bU}_{H^1}\nrm{\nabla\varrho}_{W^{1,\infty}}\lesssim\Esh(1+\Esh)$, and by \Cref{lem:GN}(3), $\nrm{v_t}_{L^4}^2+\nrm{v_t}_{L^6}^2\lesssim\nrm{\nabla v_t}^2$. Hence
\begin{equation*}
|J_{t,1}|\lesssim\nrm{\varrho_t}\nrm{v_t}_{L^4}^2\lesssim\Esh(1+\Esh)\Ff^2 .
\end{equation*}
Next, $J_{t,2}=-\int\varrho_t(\bU\cdot\nabla v)\cdot v_t-\int\rho(\bU_t\cdot\nabla v)\cdot v_t-\int\rho(\bU\cdot\nabla v_t)\cdot v_t$. For the first, $\nrm{\varrho_t}_{L^3}\nrm{\bU\cdot\nabla v}\nrm{v_t}_{L^6}\lesssim\Esh^2(1+\Esh)\Ff^2$. For the second, $\nrm{v_t\cdot\nh v}\le\nrm{v_t}_{L^6}\nrm{\nh v}_{L^3}\lesssim\Ff\Eu$ and $\nrm{w_t\pz v}\le\nrm{w_t}_{L^2_hL^\infty_z}\nrm{\pz v}_{L^\infty_hL^2_z}\lesssim\nrm{\divh v_t}\Eu$, so this term is $\lesssim(1+\Esh)\Esh\Ff^2$. For the third, $\divv\bU=0$ and $\bU\cdot n=0$ give $-\int\rho(\bU\cdot\nabla v_t)\cdot v_t=\frac12\int(\bU\cdot\nabla\rho)|v_t|^2=\frac12\int(\bU\cdot\nabla\varrho+\rs'w)|v_t|^2$, bounded by $(\nrm{\bU\cdot\nabla\varrho}+C\nrm{w})\nrm{v_t}_{L^4}^2\lesssim\Esh(1+\Esh)\Ff^2$. Finally,
\begin{equation*}
|J_{t,4}|\le|f|\Big(\nrm{\bU\cdot\nabla\varrho}_{L^\infty}\nrm{v}\nrm{v_t}+C\nrm{w}_{L^3}\nrm{v}\nrm{v_t}_{L^6}\Big)\lesssim\Esh(1+\Esh)^4\Ff^2,
\end{equation*}
where we used $\nrm{w}_{L^3}\lesssim\nrm{w}_{H^1}\lesssim\Ff$, $\nrm{v}\le\Eu$ and, for the first product, $\nrm{v}\lesssim(1+\Esh)^3\Ff$.

\emph{The term $J_{t,3}$.} Write $J_{t,3}=g(J_{t,31}+J_{t,32})$ with $J_{t,31}=\int(v\cdot\nh\varrho)w_t$ and $J_{t,32}=\int w\pz\varrho\,w_t$. Since $v_iw_t$ and $\pz(v_iw_t)$ vanish on $\pa\Omega$,
\begin{equation*}
J_{t,31}=\sum_i\int_\Omega v_iw_t\,\AB\ABmi\pa_i\varrho=\sum_i\int_\Omega\nabla\pz(v_iw_t)\cdot\Lambda^{-1}\pa_i\varrho,
\end{equation*}
and $\nabla\pz(v_iw_t)=\nabla\pz v_i\,w_t+\pz v_i\,\nabla w_t+\nabla v_i\,\pz w_t+v_i\nabla\pz w_t$. We treat the four parts in turn; recall $\pz w_t=-\divh v_t$ and $\nrm{\Lambda^{-1}\nh\varrho}\lesssim\Ff$.
\begin{itemize}
\item $\nrm{\nabla\pz v}_{L^\infty_hL^2_z}\nrm{w_t}_{L^2_hL^\infty_z}\nrm{\Lambda^{-1}\nh\varrho}\lesssim\Eu\nrm{\divh v_t}\Ff\lesssim\Esh\Ff^2$.
\item For the horizontal part of $\pz v_i\nabla w_t$ we integrate by parts horizontally: with $(\Lambda^{-1}\pa_i\varrho)_j=\pa_j\pz\ABmi\pa_i\varrho$,
\begin{equation*}
\int\pz v_i\,\pa_jw_t\,(\Lambda^{-1}\pa_i\varrho)_j=-\int\pa_j\pz v_i\,w_t\,(\Lambda^{-1}\pa_i\varrho)_j-\int\pz v_i\,w_t\,\Laph\pz\ABmi\pa_i\varrho,
\end{equation*}
and $\nrm{\Laph\pz\ABmi\pa_i\varrho}\le\sum_j\nrm{\ABmh\pa_j\pa_i\varrho}\lesssim\Ff$ by \Cref{lem:ABid}; both integrals are $\lesssim\nrm{v}_{H^{7/2}}\nrm{w_t}_{L^2_hL^\infty_z}\Ff\lesssim\Esh\Ff^2$. The vertical part is $-\int\pz v_i\divh v_t(\Lambda^{-1}\pa_i\varrho)_z\lesssim\nrm{\pz v}_{L^\infty}\nrm{\nabla v_t}\Ff\lesssim\Esh\Ff^2$.
\item $\nrm{\nabla v}_{L^\infty}\nrm{\divh v_t}\nrm{\Lambda^{-1}\nh\varrho}\lesssim\Esh\Ff^2$.
\item For $v_i\nabla\pz w_t=-v_i\nabla\divh v_t$ we integrate by parts horizontally:
\begin{equation*}
\begin{aligned}-\int v_i\,\nabla\pa_k(v_t)_k\cdot\Lambda^{-1}\pa_i\varrho&=\int\pa_kv_i\,\nabla(v_t)_k\cdot\Lambda^{-1}\pa_i\varrho+\int v_i\,\nabla(v_t)_k\cdot\Lambda^{-1}\pa_k\pa_i\varrho\\&\lesssim\nrm{v}_{W^{1,\infty}}\nrm{\nabla v_t}\Ff\lesssim\Esh\Ff^2 .\end{aligned}
\end{equation*}
\end{itemize}
For $J_{t,32}$, integrating by parts in $z$ and using $\pz w=-\divh v$,
\begin{equation*}
J_{t,32}=\int_\Omega\divh v\,\varrho\,w_t+\int_\Omega w\,\varrho\divh v_t\lesssim\nrm{\varrho}_{L^\infty}\big(\nrm{\divh v}\nrm{w_t}+\nrm{w}\nrm{\divh v_t}\big)\lesssim\Esh\Ff^2 .
\end{equation*}

Collecting \eqref{eq:3.25}, \eqref{eq:linear} and the bounds of this subsection, we have proved:

\begin{proposition}\label{prop:low}
There is a constant $\cc1\ge1$, depending only on $\mu$, $g$, $f$, $l$, $\rs$, $\nrm{\rs+\vin}_{L^\infty}$ and a positive lower bound for $\rs+\vin$, such that every classical solution of \eqref{PE} with \eqref{eq:ID}--\eqref{eq:BC} satisfies
\begin{equation}\label{eq:3.33}
\frac12\ddt\Ef^2(t)+\frac14\Ff^2(t)\le\cc1\Esh(t)\big(1+\Esh(t)\big)^4\Ff^2(t).
\end{equation}
\end{proposition}

(In fact \eqref{eq:3.25} and \eqref{eq:linear} give the left-hand side with $\frac12\Ff^2$; we keep the weaker form for convenience.)

\section{Stokes estimates}\label{sec:stokes}

With the surface pressure $P(t,x_h)\eqdef q+g\int_0^z\varrho\dd s$, the momentum equation becomes the hydrostatic Stokes system
\begin{equation}\label{eq:4.1}
-\mu\Lap v+\nh P=\zeta\eqdef-G+g\int_0^z\nh\varrho(x_h,s)\dd s,\qquad \divh\int_0^1v\dd z=0,\qquad v|_{\pa\Omega}=0,
\end{equation}
with $G=\rho(\dot v+fv^\perp)$. Compared with the non-rotating case, $\zeta$ contains the linear term $-f\rs v^\perp$. It is of lower order, but it is not small, and it is controlled through the decomposition $v=\bar v+\wt v$, \Cref{lem:poincare} and \Cref{lem:vbar}.

\begin{proposition}\label{prop:stokes}
Let $(\varrho,v,q)$ be a classical solution of \eqref{PE}--\eqref{eq:BC}. The constants below depend on $\mu,g,f,l,\rs$, $\nrm{\rho}_{L^\infty}$ and a positive lower bound for $\rho=\rs+\varrho$ (by transport, $\inf\rho=\inf(\rs+\vin)$).
\begin{enumerate}[label=(\roman*)]
\item One has
\begin{equation}\label{eq:4.2}
\nrm{v}_{H^2}+\nrm{\nh P}_{L^2}\le C(1+\Ev)^4\big(\nrm{v_t}_{L^2}+\nrm{v}_{H^3}\nrm{\nh v}_{L^2}+\nrm{\nh v}_{L^2}+\nrm{\nh\varrho}_{L^2}\big),
\end{equation}
and consequently
\begin{equation}\label{eq:4.3}
\nrm{v}_{H^2}+\nrm{\nh P}_{L^2}\le C(1+\Esh)^4\Ef,
\end{equation}
\begin{equation}\label{eq:4.4}
\nrm{v}_{H^2}\le C\Big(\Ff^{\frac23}\Fuz^{\frac13}+(1+\Esh)^4\Ff\Big).
\end{equation}
\item One has
\begin{equation}\label{eq:4.5}
\nrm{v}_{H^3}\le C\Big(\Ff^{\frac13}\Fuz^{\frac23}+(1+\Esh)^4(1+\Fn)\Ff+\nrm{\nabla^2v}_{L^2}^2\Big).
\end{equation}
\end{enumerate}
\end{proposition}

\begin{proof}
(i) By \Cref{lem:stokes} with $m=0$, $\nrm{v}_{H^2}+\nrm{\nh P}\le C\nrm{\zeta}$, and
\begin{equation*}
\nrm{\zeta}\le\nrm{\rho}_{L^\infty}\big(\nrm{v_t}+\nrm{\bU\cdot\nabla v}+|f|\nrm{v}\big)+g\nrm{\nh\varrho}.
\end{equation*}
By \eqref{eq:3.13b}, $\nrm{\bU\cdot\nabla v}\lesssim\nrm{v}_{H^3}\nrm{\nh v}$. By \Cref{lem:poincare} and \Cref{lem:vbar}, $\nrm{v}\le l\nrm{\nh v}+\nrm{\bar v}\lesssim(1+\Ev)^2(\nrm{v_t}+\nrm{v}_{H^3}\nrm{\nh v}+\nrm{\nh v})$. This proves \eqref{eq:4.2}; \eqref{eq:4.3} follows from \eqref{eq:3.27} and $\nrm{v}_{H^3}\le\Esh$.

For \eqref{eq:4.4} we split $\nabla^2v$. First, $\nrm{\nh^2v}+\nrm{\nh\pz v}\le\nrm{\nabla\nh v}\lesssim\Ff$. Second, by the horizontal Poincar\'e inequality and interpolation between $L^2$ and $H^3$,
\begin{equation}\label{eq:4.6}
\nrm{\wt{\pz^2v}}^2\le l^2\nrm{\nh\pz^2v}^2\lesssim\nrm{\nh\pz v}^{\frac43}\nrm{\nh\pz v}_{H^3}^{\frac23}\lesssim\Ff^{\frac43}\Fuz^{\frac23}.
\end{equation}
Third, by \Cref{lem:vbar} and \eqref{eq:3.13b},
\begin{equation}\label{eq:4.7}
\nrm{\pz^2\bar v}\le\nrm{\bar v}_{H^2}\lesssim(1+\Esh)^3\Ff .
\end{equation}
Since $v$ vanishes on $\pa\Omega$, $\nrm{v}_{H^2}\lesssim\nrm{\nabla^2v}$ (\Cref{lem:Linf}), and \eqref{eq:4.4} follows.

(ii) We split $\nabla^3v$ in the same way. The derivatives containing at least one horizontal derivative satisfy, by interpolation and \eqref{eq:3.3},
\begin{equation}\label{eq:4.10}
\nrm{\nh\nabla^2v}^2\lesssim\nrm{\nh\nabla v}\nrm{\nh\nabla^3v}\lesssim\Ff\Fuz,
\end{equation}
since $\nrm{\nh^4v}\le\pi^{-1}\nrm{\nh^4\pz v}$ and every other component of $\nh\nabla^3v$ contains $\pz v$. As in \eqref{eq:4.6},
\begin{equation}\label{eq:4.8}
\nrm{\wt{\pz^3v}}^2\le l^2\nrm{\nh\pz^3v}^2\lesssim\nrm{\nh\pz v}^{\frac23}\nrm{\nh\pz v}_{H^3}^{\frac43}\lesssim\Ff^{\frac23}\Fuz^{\frac43}.
\end{equation}
For the mean, differentiating \eqref{eq:decomp2}$_2$ gives $\mu\pz^3\bar v=\pz\ol G$, and
\begin{equation*}
\nrm{\pz G}\le\nrm{\nabla\rho}_{L^\infty}\nrm{\dot v+fv^\perp}+\nrm{\rho}_{L^\infty}\big(\nrm{\pz v_t}+\nrm{\pz(\bU\cdot\nabla v)}+|f|\nrm{\pz v}\big).
\end{equation*}
Here $\nrm{\dot v}+\nrm{v}_{H^1}\lesssim(1+\Esh)^3\Ff$ by \eqref{eq:3.13b}, $\nrm{\pz v_t}\lesssim\Ff$, and
\begin{equation*}
\begin{aligned}\nrm{\pz(\bU\cdot\nabla v)}&\le\nrm{\nabla\bU}_{L^\infty}\nrm{\nabla v}+\nrm{v}_{L^\infty}\nrm{\nabla\nh v}+\nrm{w}_{L^\infty_zL^4_h}\nrm{\pz\nabla v}_{L^2_zL^4_h}\\&\lesssim\Fn(1+\Esh)^3\Ff+\nrm{\nabla^2v}\Ff+\Ff\nrm{v}_{H^{5/2}},\end{aligned}
\end{equation*}
using \eqref{eq:3.13} and \Cref{lem:aniso}(f). With $\nrm{v}_{H^{5/2}}\le\nrm{v}_{H^3}\le\Esh$ we obtain $\nrm{\pz^3\bar v}\lesssim(1+\Esh)^4(1+\Fn)\Ff$. Taking square roots in \eqref{eq:4.10}, \eqref{eq:4.8} and using $\Ff^{1/2}\Fuz^{1/2}\lesssim\Ff^{1/3}\Fuz^{2/3}+\Ff$ bounds $\nrm{\nabla^3v}$. Finally $\nrm{v}_{H^3}\lesssim\nrm{\nabla^3v}+\nrm{v}_{H^2}$, and the lower-order part is bounded by \eqref{eq:4.4} and $\Ff^{2/3}\Fuz^{1/3}\le\Ff^{1/3}\Fuz^{2/3}+\Ff$; this gives \eqref{eq:4.5} (the last term $\nrm{\nabla^2v}^2$ in \eqref{eq:4.5} is not needed, and is kept only for later convenience).
\end{proof}

\begin{corollary}\label{cor:stokes}
Assume \eqref{eq:small}. Then
\begin{equation}\label{eq:4.11}
\nrm{v}_{H^2}\lesssim\Ff^{\frac23}\Fuz^{\frac13}+\Ff,\qquad \nrm{v}_{H^3}\lesssim\Ff^{\frac13}\Fuz^{\frac23}+\Ff+\Ff^{\frac7{12}}\Fsh^{\frac5{12}},
\end{equation}
hence
\begin{equation}\label{eq:4.12}
\nrm{v}_{H^2}\nrm{v}_{H^3}\lesssim\Fsh^2,
\end{equation}
and
\begin{equation}\label{eq:4.13}
\Eu\lesssim\Eut+\Ev+\Ef+\Esh^2 .
\end{equation}
\end{corollary}

\begin{proof}
Under \eqref{eq:small}, $\Esh\lesssim1$, so \eqref{eq:4.4} gives the first bound. Then $\nrm{\nabla^2v}^2\lesssim\Ff^{4/3}\Fuz^{2/3}+\Ff^2\lesssim\Ff^{1/3}\Fuz^{2/3}+\Ff$, because $\Ff\lesssim\Esh\lesssim1$ by \eqref{eq:3.11} and $\Ff^{4/3}\Fuz^{2/3}=\Ff\cdot\Ff^{1/3}\Fuz^{2/3}$. By \eqref{eq:3.13}, $\Fn\lesssim\Esh+\Ff^{7/12}\Fsh^{5/12}\lesssim1+\Ff^{7/12}\Fsh^{5/12}$, and $\Ff\cdot\Ff^{7/12}\Fsh^{5/12}\lesssim\Ff^{7/12}\Fsh^{5/12}$; inserting these into \eqref{eq:4.5} gives the second bound. \eqref{eq:4.12} follows from $\Ff,\Fuz\le C\Fsh$.

For \eqref{eq:4.13}, \Cref{lem:stokes} with $m=2$ gives $\Eu\lesssim\nrm{G}_{H^2}+\nrm{\int_0^z\nh\varrho}_{H^2}$. By \Cref{lem:GN}(2), $\nrm{G}_{H^2}\lesssim(1+\Ev)(\nrm{v_t}_{H^2}+\nrm{\bU\cdot\nabla v}_{H^2}+|f|\nrm{v}_{H^2})$, with $\nrm{\bU\cdot\nabla v}_{H^2}\lesssim\nrm{v}_{H^3}\nrm{v}_{H^4}$ (using $\nrm{\pa^\beta w\,\pa^\gamma\pz v}\le\nrm{\pa^\beta w}_{L^2_hL^\infty_z}\nrm{\pa^\gamma\pz v}_{L^\infty_hL^2_z}$ for $|\beta|=2$) and $\nrm{v}_{H^2}\lesssim\Ef$ by \eqref{eq:4.3}. Finally $\nrm{\int_0^z\nh\varrho}_{H^2}\lesssim\Ev$.
\end{proof}

\begin{remark}\label{rem:twoest}
The two bounds \eqref{eq:4.3} and \eqref{eq:4.4} serve different purposes: \eqref{eq:4.3} gives pointwise decay of $\nrm{v}_{H^2}$ through $\Ef$, and \eqref{eq:4.4} controls its time integral through $\Ff$.
\end{remark}

Differentiating \eqref{eq:4.1} in $z$ eliminates the surface pressure and gives \eqref{eq:3.22}:
\begin{equation}\label{eq:4.14}
-\mu\Lap\pz v=-\pz G+g\nh\varrho .
\end{equation}
Since $\pz v$ satisfies no boundary condition, \Cref{lem:stokes} does not apply to \eqref{eq:4.14}; the following recovery lemma supplies the anisotropic estimate.

\begin{lemma}\label{lem:recovery}
Let $V=\pz v$, where $v$ is a classical solution of \eqref{PE}--\eqref{eq:BC} satisfying \eqref{eq:lid0}, \eqref{eq:wt0} and \eqref{eq:compat}. Then
\begin{equation}\label{eq:4.15}
\Fuz=\nrm{V}_{H^4}\lesssim\nrm{\nh^4V}_{L^2}+\nrm{\pz G}_{H^2}+\nrm{\nh\varrho}_{H^2}.
\end{equation}
\end{lemma}

\begin{proof}
Set $\cG\eqdef-\pz G+g\nh\varrho$, so that $\pz^2V=-\Laph V-\mu^{-1}\cG$. For the oscillatory part, the horizontal Poincar\'e inequality gives $\nrm{\nh^j\wt V}\lesssim\nrm{\nh^4V}$ for $j\le4$. The one-dimensional inequality $\nrm{\phi'}^2\lesssim\nrm{\phi}\nrm{\phi''}+\nrm{\phi}^2$ on $(0,1)$, applied to the Fourier coefficients, gives
\begin{equation*}
\nrm{\nh^3\pz\wt V}^2\lesssim\nrm{\nh^4V}\nrm{\nh^2\pz^2\wt V}+\nrm{\nh^4V}^2,\qquad \nrm{\nh^2\pz^2\wt V}\lesssim\nrm{\nh^4V}+\nrm{\nh^2\wt\cG},
\end{equation*}
hence $\nrm{\nh^3\pz\wt V}\lesssim\nrm{\nh^4V}+\nrm{\wt\cG}_{H^2}$, and $\nrm{\nh^j\pz\wt V}\lesssim\nrm{\nh^3\pz\wt V}$ for $j\le3$. Every derivative of order at most four with at least two $z$-derivatives reduces, through $\pz^2V=-\Laph V-\mu^{-1}\cG$, to a derivative already estimated or to a derivative of $\wt\cG$ of order at most two. Thus $\nrm{\wt V}_{H^4}\lesssim\nrm{\nh^4V}+\nrm{\wt\cG}_{H^2}$.

For the mean, \eqref{eq:decomp2}$_2$ gives $\mu\pz\bar V=\ol G$ and $\int_0^1\bar V\dd z=0$. The no-slip condition and its time compatibility give $\dot v|_{\pa\Omega}=0$ and $v|_{\pa\Omega}=0$, hence $G|_{\pa\Omega}=0$ and $\ol G|_{z=0,1}=0$. By Poincar\'e inequalities in $z$, $\nrm{\bar V}_{H^4(0,1)}\lesssim\nrm{\ol G}_{H^3(0,1)}\lesssim\nrm{\pz\ol G}_{H^2(0,1)}=\nrm{\bar\cG}_{H^2(0,1)}$, since $\ol{\nh\varrho}=0$.
\end{proof}

\begin{remark}\label{rem:fH3}
In \eqref{eq:4.15}, $\nrm{\pz G}_{H^2}\le\nrm{\rho\dot v}_{H^3}+|f|\nrm{\rho v}_{H^3}$. The rotational term is linear and not small; it is controlled by interpolation. Since $\nrm{v}_{H^4}\lesssim\nrm{\pz v}_{H^4}=\Fuz$ by \eqref{eq:3.3}, for every $\epsilon>0$,
\begin{equation}\label{eq:fH3}
|f|\nrm{\rho v}_{H^3}\lesssim(1+\Ev)\nrm{v}_{H^2}^{\frac12}\nrm{v}_{H^4}^{\frac12}\le\epsilon\Fuz+C_\epsilon(1+\Ev)^2\nrm{v}_{H^2},
\end{equation}
and $\nrm{v}_{H^2}$ is controlled by \eqref{eq:4.3} or \eqref{eq:4.4}.
\end{remark}

\begin{proposition}\label{prop:stokes_t}
Let $(\rho,v,q)$ be a strong solution of \eqref{eq:4.1}. For every integer $k\ge0$, if $\zeta_t\in H^k$, then
\begin{equation}\label{eq:4.17}
\nrm{v_t}_{H^{k+2}}+\nrm{\nh P_t}_{H^k}\le C\Big(\nrm{G_t}_{H^k}+\norm{\int_0^z\nh\varrho_t\dd s}_{H^k}\Big).
\end{equation}
\end{proposition}

\begin{proof}
Differentiate \eqref{eq:4.1} in $t$ and apply \Cref{lem:stokes}.
\end{proof}

\begin{remark}\label{rem:gradq}
\Cref{lem:stokes} with $m=2$ and $\nrm{\zeta}_{H^2}\lesssim(1+\nrm{\varrho}_{H^4})(\nrm{v_t}_{H^2}+\nrm{v}_{H^3}\nrm{v}_{H^4}+\nrm{v}_{H^2}+\nrm{\varrho}_{H^3})$ give $\nrm{\nabla q}_{H^2}\lesssim\nrm{\zeta}_{H^2}+\nrm{\varrho}_{H^3}$; this is the regularity $\nabla q\in C([0,T];H^2)$ in \eqref{eq:localclass}.
\end{remark}

\section{Decay estimates}\label{sec:decay}

By \eqref{eq:2.10}, \Cref{lem:poincare} and \Cref{lem:w},
\begin{equation*}
\nrm{\nh\varrho}^2\lesssim\nrm{\pz\varrho}_{H^3}^{\frac23}\nrm{\ABmh\Laph\varrho}^{\frac43},\qquad \nrm{\nh v}+\nrm{w}\lesssim\nrm{\nabla\nh v},\qquad \nrm{v_t}\lesssim\nrm{\nabla v_t}.
\end{equation*}
Hence \eqref{eq:3.27} gives
\begin{equation}\label{eq:5.1}
\Ef\lesssim\Es^{\frac13}\Ff^{\frac23}.
\end{equation}

\begin{proposition}\label{prop:decay1}
Suppose that
\begin{equation}\label{eq:small}
\Es(t)\le\frac1{256\cc1}\qquad\text{for all }t\in[0,T],
\end{equation}
with $\cc1$ from \Cref{prop:low}. Then, for $t\in[0,T]$,
\begin{equation}\label{eq:5.3}
\nrm{\wt\varrho}_{L^2}+\nrm{\nh\varrho}_{L^2}+\nrm{\nh v}_{L^2}+\nrm{v_t}_{L^2}+\nrm{w}_{L^2}\lesssim\min\big\{\Ef(0),\Es(t)(1+t)^{-1}\big\}.
\end{equation}
\end{proposition}

\begin{remark}\label{rem:c1}
(i) \eqref{eq:3.33} remains true if $\cc1$ is enlarged. We fix $\cc1\ge1$ once and for all so large that the finitely many absorption conditions of the form ``$C\Es\le\frac12$'' (with $C$ depending only on $\mu,g,f,l,\rs$, $\nrm{\rs+\vin}_{L^\infty}$ and a positive lower bound for $\rs+\vin$, hence only on $\mu,g,f,l,\rs$ on the data class used in the proof of \Cref{thm:main}) in \cref{sec:decay} and \cref{sec:high} follow from \eqref{eq:small}. (ii) The statements of \cref{sec:decay} and \cref{sec:high} use \eqref{eq:small} only on $[0,t]$, and $t\mapsto\Es(t)$ is nondecreasing; they hold on every interval $[0,T]$ on which the solution exists and \eqref{eq:small} holds. This is how they are used in the bootstrap.
\end{remark}

\begin{proof}[Proof of \Cref{prop:decay1}]
Since $\Esh\le2\Es$ by \eqref{eq:Esharpdef}, \eqref{eq:small} and $\cc1\ge1$ give $\Esh\le(128\cc1)^{-1}\le\frac1{128}$, so $\cc1\Esh(1+\Esh)^4\le2\cc1\Esh\le\frac1{64}\le\frac18$, and \eqref{eq:3.33} gives $\frac12\ddt\Ef^2+\frac18\Ff^2\le0$, in particular
\begin{equation}\label{eq:5.2}
\ddt\Ef^2+\frac18\Ff^2\le0\qquad\text{on }[0,T].
\end{equation}
In particular $\Ef(\tau)\le\Ef(0)$. By \eqref{eq:5.1}, $\Ff^2\ge c\,\Es(t)^{-1}\Ef^3$ on $[0,t]$, so $y=\Ef^2$ satisfies $y'+c\Es^{-1}(t)y^{3/2}\le0$ and
\begin{equation*}
\Ef(\tau)\le\Big(\Ef(0)^{-1}+\tfrac c2\Es^{-1}(t)\tau\Big)^{-1}\lesssim\min\big\{\Ef(0),\Es(t)(1+\tau)^{-1}\big\},\qquad0\le\tau\le t,
\end{equation*}
where we used $\Ef(0)\lesssim\Esh(0)\le2\Es(t)$. By \eqref{eq:3.27} the left-hand side of \eqref{eq:5.3}, without $\wt\varrho$, is bounded by $C\Ef$, and $\nrm{\wt\varrho}\le l\nrm{\nh\varrho}$. This is the mechanism of \cite[Remark 3.6]{DongLiSun2026}.
\end{proof}

\begin{proposition}\label{prop:decay2}
Assume \eqref{eq:small}. Then
\begin{equation}\label{eq:5.4}
\Ef^2(t)+\frac18\int_0^t\Ff^2\dtau\le\Ef^2(0),\qquad \nrm{v(t)}_{H^2}\lesssim\Es(t)(1+t)^{-1},
\end{equation}
and for $0<\delta<2$,
\begin{align}
\int_0^t(1+\tau)^{2-\delta}\Ff^2(\tau)\dtau&\lesssim\delta^{-1}\Es^2(t),\label{eq:5.5}\\
\int_0^t\Ff(\tau)\dtau&\lesssim\Es(t).\label{eq:5.6}
\end{align}
The dissipation has the tail bound
\begin{equation}\label{eq:5.tail}
\int_\tau^t\Ff^2(s)\dd s\le8\Ef^2(\tau)\lesssim\Es^2(t)(1+\tau)^{-2},\qquad0\le\tau\le t .
\end{equation}
Moreover, for $0<\delta'<\frac43$,
\begin{equation*}
\int_0^t(1+\tau)^{\frac43-\delta'}\nrm{v(\tau)}_{H^2}^2\dtau\lesssim\delta'^{-1}\Es^2(t).
\end{equation*}
\end{proposition}

\begin{proof}
\eqref{eq:5.4} follows from \eqref{eq:5.2}, \eqref{eq:4.3} and \eqref{eq:5.3}, and integrating \eqref{eq:5.2} over $[\tau,t]$ and using the bound $\Ef(\tau)\lesssim\Es(t)(1+\tau)^{-1}$ from the proof of \Cref{prop:decay1} gives \eqref{eq:5.tail}. Multiplying \eqref{eq:5.2} by $(1+\tau)^{2-\delta}$ and integrating, with \eqref{eq:5.3},
\begin{equation*}
(1+t)^{2-\delta}\Ef^2(t)+\frac18\int_0^t(1+\tau)^{2-\delta}\Ff^2\dtau\le\Ef^2(0)+(2-\delta)\int_0^t(1+\tau)^{1-\delta}\Ef^2\dtau\lesssim\Ef^2(0)+\delta^{-1}\Es^2(t),
\end{equation*}
which gives \eqref{eq:5.5}, since $\Ef(0)\lesssim\Es(t)$. Then \eqref{eq:5.6} follows from the Cauchy--Schwarz inequality with $\delta=\frac12$. Finally, by \eqref{eq:4.4} and \eqref{eq:small}, $\nrm{v}_{H^2}^2\lesssim\Ff^2+\Ff^{4/3}\Fuz^{2/3}$; with $k=\frac43-\delta'$ and Young's inequality $(1+\tau)^k\Ff^{4/3}\Fuz^{2/3}\le\frac23(1+\tau)^{3k/2}\Ff^2+\frac13\Fuz^2$, where $\frac32k=2-\frac32\delta'$, we conclude by \eqref{eq:5.5} and $\int_0^t\Fuz^2\le\Es^2$.
\end{proof}

The rotational term in the top-order density estimate is linear in $\Esh$ and is multiplied by $\Ff$. It is controlled by the following refinement of \eqref{eq:5.6}, which measures $\int\Ff$ by a positive power of the initial low-order energy.

\begin{lemma}\label{lem:intFflat}
Assume \eqref{eq:small}. Then
\begin{equation}\label{eq:intFflat}
\int_0^t\Ff(\tau)\dtau\le\cc4\,\Ef(0)^{\frac14}\Es(t)^{\frac34},
\end{equation}
with $\cc4$ depending only on $\mu,g,f,l,\rs$, $\nrm{\rs+\vin}_{L^\infty}$ and a positive lower bound for $\rs+\vin$.
\end{lemma}

\begin{proof}
Write $a=\Ef(0)$ and $b=\Es(t)$; then $a\le Cb$. If $a=0$, then $\Ef\equiv0$ and $\Ff\equiv0$ on $[0,t]$ by \eqref{eq:5.2}, and there is nothing to prove; let $a>0$. By \eqref{eq:5.3}, $\Ef(\tau)\le C\min\{a,b(1+\tau)^{-1}\}$. Multiplying \eqref{eq:5.2} by $(1+\tau)^{3/2}$ and integrating,
\begin{equation*}
\frac18\int_0^t(1+\tau)^{\frac32}\Ff^2\dtau\le a^2+\frac32\int_0^\infty(1+\tau)^{\frac12}C^2\min\{a^2,b^2(1+\tau)^{-2}\}\dtau .
\end{equation*}
Let $\tau_0\ge0$ be defined by $1+\tau_0=\max\{1,b/a\}$. The integral over $[0,\tau_0]$ is at most $\frac23a^2(1+\tau_0)^{3/2}\le\frac23a^{1/2}b^{3/2}$, and the integral over $[\tau_0,\infty)$ is at most $2b^2(1+\tau_0)^{-1/2}\le2a^{1/2}b^{3/2}$ if $b\ge a$; if $b<a$ the whole integral is at most $2b^2\le2a^{1/2}b^{3/2}$. Hence $\int_0^t(1+\tau)^{3/2}\Ff^2\lesssim a^{1/2}b^{3/2}$, and the Cauchy--Schwarz inequality gives
\begin{equation*}
\int_0^t\Ff\dtau\le\Big(\int_0^t(1+\tau)^{\frac32}\Ff^2\dtau\Big)^{\frac12}\Big(\int_0^\infty(1+\tau)^{-\frac32}\dtau\Big)^{\frac12}\lesssim a^{\frac14}b^{\frac34}.\qedhere
\end{equation*}
\end{proof}

\begin{proposition}\label{prop:Fnat}
Assume \eqref{eq:small}. Then $\int_0^t\Fn(\tau)\dtau\lesssim\Es(t)$.
\end{proposition}

\begin{proof}
\emph{Step 1: $\nrm{\nabla\wt v}_{L^\infty}$.} For a.e.\ $z$, the function $\nabla\wt v(\cdot,z)$ has zero horizontal mean, so by \Cref{lem:aniso}(c) with $s=\frac12$,
\begin{equation*}
\nrm{\nabla\wt v}_{L^\infty}\lesssim\sup_z\nrm{\nh\nabla v(\cdot,z)}_{L^2_h}^{\frac12}\sup_z\nrm{\nh^2\nabla v(\cdot,z)}_{L^2_h}^{\frac12}.
\end{equation*}
By \Cref{lem:aniso}(a), $\sup_z\nrm{\phi(\cdot,z)}_{L^2_h}\lesssim\nrm{\phi}^{1/2}\nrm{\phi}_{H^1}^{1/2}$. Applied to $\phi=\nh\nabla v$ and $\phi=\nh^2\nabla v$,
\begin{equation*}
\nrm{\nabla\wt v}_{L^\infty}\lesssim\nrm{\nh v}_{H^1}^{\frac14}\nrm{\nh v}_{H^2}^{\frac12}\nrm{\nh v}_{H^3}^{\frac14}.
\end{equation*}
Interpolating $\nrm{\nh v}_{H^m}\lesssim\nrm{\nh v}_{H^1}^{(4-m)/3}\nrm{\nh v}_{H^4}^{(m-1)/3}$ for $m=2,3$, and using $\nrm{\nh v}_{H^1}\lesssim\nrm{\nabla\nh v}\lesssim\Ff$ (\Cref{lem:poincare}) and $\nrm{\nh v}_{H^4}\le\Fu$,
\begin{equation}\label{eq:gradvtilde}
\nrm{\nabla\wt v}_{L^\infty}\lesssim\Ff^{\frac23}\Fu^{\frac13}.
\end{equation}

\emph{Step 2: $\nrm{\nh w}_{L^\infty}$ and $\nrm{\pz w}_{L^\infty}$.} By \Cref{lem:w} with $s=\frac14$, and interpolation,
\begin{equation}\label{eq:gradhw}
\nrm{\nh w}_{L^\infty}\lesssim\nrm{\nh\divh v}_{L^2_zH^{5/4}_h}\lesssim\nrm{\nh v}_{H^{9/4}}\lesssim\nrm{\nh v}_{H^1}^{\frac7{12}}\nrm{\nh v}_{H^4}^{\frac5{12}}\lesssim\Ff^{\frac7{12}}\Fu^{\frac5{12}},
\end{equation}
and $\nrm{\pz w}_{L^\infty}=\nrm{\divh v}_{L^\infty}\le\nrm{\nabla v}_{L^\infty}$.

\emph{Step 3: $\nrm{\pz\bar v}_{L^\infty}$.} By \Cref{lem:vbar} and \eqref{eq:3.13b}, $\nrm{\pz\bar v}_{L^\infty}\lesssim\nrm{\bar v}_{H^2}\lesssim\Ff$.

\emph{Step 4: time integrals.} For $a\in(\frac12,1]$ and $0<\delta<2-\frac1a$, H\"older's inequality and \eqref{eq:5.5} give
\begin{equation}\label{eq:Fa}
\int_0^t\Ff^a\Fsh^{1-a}\dtau\le\Big(\int_0^t(1+\tau)^{2-\delta}\Ff^2\Big)^{\frac a2}\Big(\int_0^t\Fsh^2\Big)^{\frac{1-a}2}\Big(\int_0^\infty(1+\tau)^{-(2-\delta)a}\dtau\Big)^{\frac12}\lesssim\Es(t),
\end{equation}
since $(2-\delta)a>1$. We apply \eqref{eq:Fa} with $a=1$ (Step 3), $a=\frac23$ (Step 1) and $a=\frac7{12}$ (Step 2 and the last term of $\Fn$). Since $\nabla v=\nabla\wt v+\pz\bar v\otimes e_z$, Steps 1--3 bound $\int\nrm{\nabla\bU}_{L^\infty}$. Finally, $\nrm{\nabla^2v}\le\nrm{v}_{H^2}$ and \Cref{prop:decay2} give, with $k=\frac43-\delta'\in(1,\frac43)$,
\begin{equation*}
\int_0^t\nrm{\nabla^2v}\dtau\le\Big(\int_0^\infty(1+\tau)^{-k}\Big)^{\frac12}\Big(\int_0^t(1+\tau)^k\nrm{v}_{H^2}^2\Big)^{\frac12}\lesssim\Es(t).\qedhere
\end{equation*}
\end{proof}

\section{Higher-order energy estimates}\label{sec:high}

In this section we prove closed a priori bounds for the four components of the higher-order energy $\Es$ in \eqref{eq:energies}. The estimates are derived in the order
\begin{equation*}
\cEutt\ \longrightarrow\ \cEut\ \longrightarrow\ \cT_4\ \longrightarrow\ \cEu\ \longrightarrow\ \cEv .
\end{equation*}
The bound for $\cEutt$ (\Cref{prop:utt}) follows from the energy identity of the twice time-differentiated momentum equation, and the decay of $v_{tt}$ (\Cref{prop:uttdecay}) from the same identity with a time weight. The temporal Stokes estimate \eqref{eq:4.17} converts it into a bound for $\cEut$ (\Cref{prop:ut}). The tangential energy $\cT_4$ supplies the dissipation $D_4$ of the fourth horizontal derivatives (\Cref{prop:T4}), which, together with the recovery \Cref{lem:recovery}, controls $\cEu$ (\Cref{prop:u}). The density estimate (\Cref{prop:rho}), based on the multipliers $\AB(\sgm\varrho)$, $\Laph^2\varrho$ and $\Laph\pz\nabla(\sgm\varrho)$, closes the hierarchy. In \Cref{prop:total} the four estimates are combined into a single inequality for $\Es$, from which \Cref{thm:main} follows by a bootstrap argument (\cref{sec:proofmain}).

The Coriolis force enters in two ways. At leading order it drops out of every energy identity, because $v\mapsto v^\perp$ commutes with time derivatives and horizontal derivatives, and $\pa^\alpha v^\perp\cdot\pa^\alpha v=0$ pointwise. The remaining rotational terms are either cubic, when a derivative falls on $\rho$, or linear and of lower order. The linear ones are routed through $D_4$, $\nrm{v}_{H^2}$ or $\Ff$ (see \Cref{rem:fH3}), never through an $O(1)$ multiple of $\Fuz$. The only linear rotational term that is not of lower order occurs in the top-order density estimate; it is bounded by $\Esh\Ff$ and is controlled by \Cref{lem:intFflat}.

\emph{Standing conventions.} Throughout this section we assume \eqref{eq:small} on $[0,t]$ (see \Cref{rem:c1}), and $\tau$ denotes a time in $[0,t]$. We write $\Es(0)\eqdef\lim_{t\downarrow0}\Es(t)\simeq\Esh(0)$. By \eqref{eq:Esharpdef} and the Cauchy--Schwarz inequality, $\Esh(\tau)\le2\Es(t)$, so under \eqref{eq:small} $\Esh(\tau)\le2\Es(t)\le1$; hence $\nrm{\varrho}_{W^{2,\infty}}+\nrm{\nabla\rho}_{L^\infty}\lesssim1$ by \eqref{eq:3.12}, $\rho\ge\inf\rin>0$ by transport (so the constants below may also depend on a positive lower bound for $\rs+\vin$), and \eqref{eq:3.11}, \eqref{eq:4.3} give
\begin{equation}\label{eq:A6-std}
\Ff\lesssim\Esh,\qquad \nrm{v}_{H^2}\lesssim\Ef\lesssim\Esh .
\end{equation}
We use freely the bounds of \Cref{rem:simplify}, in particular
\begin{gather*}
\nrm{v}_{L^\infty}+\nrm{\nabla v}_{L^6}+\nrm{w}_{H^1}+\nrm{w}_{L^\infty_zL^4_h}\lesssim\nrm{v}_{H^2},\\
\nrm{\bU}_{L^\infty}+\nrm{\nabla v}_{L^\infty}+\nrm{w}_{H^2}\lesssim\nrm{v}_{H^3},\qquad \nrm{\nabla\bU}_{L^\infty}+\nrm{w}_{W^{1,\infty}}\lesssim\Eu,
\end{gather*}
and the bounds of \Cref{lem:w} for $w_t=-\int_0^z\divh v_t\dd z'$ and $w_{tt}=-\int_0^z\divh v_{tt}\dd z'$; for instance $\nrm{w_t}_{L^2}+\nrm{w_t}_{L^2_hL^\infty_z}\le\nrm{\divh v_t}_{L^2}\lesssim\Ff$ and $\nrm{w_{tt}}_{L^2}+\nrm{w_{tt}}_{L^2_hL^\infty_z}\le\nrm{\divh v_{tt}}_{L^2}$. By \Cref{lem:poincare}, $\nrm{v_t}_{H^1}\lesssim\nrm{\nabla v_t}_{L^2}\lesssim\Ff$.

\begin{remark}\label{rem:justify}
\emph{Justification of the energy identities.} In the class \eqref{eq:localclass} several intermediate expressions below are not in $L^2$ (for instance $\AB(\sgm\bU\cdot\nabla\varrho)$ contains fifth derivatives of $\varrho$, and testing \eqref{eq:6.2} with $v_{tt}$ involves $v_{ttt}$), although the final identities contain only controlled quantities. They are justified as follows.

(a) For the density, each of $f=\pa^\alpha\varrho$ ($|\alpha|\le4$), $\AB(\sgm\varrho)$, $\Laph^2\varrho$ and $\Laph\pz\nabla(\sgm\varrho)$ solves a linear transport equation $\pt f+\bU\cdot\nabla f=S_f$, in which the commutators have been moved into $S_f\in L^1(0,T;L^2)$; this is the content of the commutator expansions in the proofs of \Cref{prop:T4,prop:rho}. In three dimensions these source terms contain $\nh^4w$, which lies in $L^2(0,T;L^2_hL^\infty_z)$ by \Cref{lem:w}, since $v\in L^2(0,T;H^5)$. The transporting field satisfies $\divv\bU=0$, $\bU\cdot n=0$ on $\pa\Omega$, and $\bU\in L^\infty(0,T;W^{1,\infty})$: by \eqref{eq:3.13}, $\nrm{\bU}_{L^\infty}\lesssim\nrm{v}_{H^3}$ and $\nrm{\nabla\bU}_{L^\infty}\lesssim\nrm{v}_{H^4}$, where the horizontal derivatives of $w$ require $\nrm{\nh\divh v}_{L^2_zH^{1+s}_h}\lesssim\nrm{v}_{H^{3+s}}$ (\Cref{lem:w}). By the renormalization theory of DiPerna and Lions \cite{DiPernaLions1989}, $f\in C([0,T];L^2)$ and $\ddt\int_\Omega|f|^2\dx=2\int_\Omega fS_f\dx$; with the weight $\sgm(z)$ one obtains in addition the term $\int_\Omega w\,\sgm'|f|^2\dx$.

(b) For the momentum equation and its time derivatives we use the Lions--Magenes lemma \cite[Lemma III.1.2]{Temam2001} in the lid-constrained space
\begin{equation*}
V_\sigma\eqdef\Big\{\varphi\in H^1_0(\Omega;\R^2):\ \divh\int_0^1\varphi\dd z=0\Big\},
\end{equation*}
with dual $V_\sigma'$. For $\varphi\in V_\sigma$ and $P=P(x_h)$, $\int_\Omega\nh P\cdot\varphi\dx=-\int_{\Tl}P\,\divh\int_0^1\varphi\dd z\dxh=0$, so the surface pressure drops out. By \eqref{eq:localclass} and \eqref{eq:lid}, $v_{tt}\in L^2(0,T;V_\sigma)$, and \eqref{eq:6.2}, read in $V_\sigma'$, gives $\pt(\rho v_{tt})\in L^2(0,T;V_\sigma')$; each term on its right-hand side is bounded in the proof of \Cref{prop:utt}. Since $\rho\in W^{1,\infty}((0,T)\times\Omega)$ by \eqref{eq:3.12} and \eqref{eq:A6-rhotH2}, a mollification in time shows that $t\mapsto\int_\Omega\rho|v_{tt}|^2\dx$ is absolutely continuous, with
\begin{equation*}
\frac12\ddt\int_\Omega\rho|v_{tt}|^2\dx=\ip{\pt(\rho v_{tt}),v_{tt}}-\frac12\int_\Omega\varrho_t|v_{tt}|^2\dx .
\end{equation*}
The identities for $v$ and $v_t$ are obtained in the same way, using $v\in L^2(0,T;H^5)$ and $v_t\in L^2(0,T;H^3)$. The same applies to the horizontal derivatives $\pa_h^\alpha v$, $|\alpha|=4$, used in \Cref{prop:T4}. Horizontal derivatives preserve the no-slip condition, and $\pa_h^\alpha\divh\int_0^1v\dd z\equiv0$, so $\pa_h^\alpha v\in L^2(0,T;V_\sigma)$ because $v\in L^2(0,T;H^5)$. Not every term of the differentiated momentum equation lies in $L^2(0,T;L^2)$: the hydrostatic part $-g\pa_h^\alpha\nh\int_0^z\varrho\dd s$ of the pressure gradient contains five horizontal derivatives of $\varrho$. Writing $\pa_h^\alpha=\pa_i\pa_h^{\alpha-e_i}$ for some $i$ with $\alpha_i\ge1$, it is the horizontal derivative $\pa_i$ of $-g\pa_h^{\alpha-e_i}\nh\int_0^z\varrho\dd s\in C([0,T];L^2)$, hence it lies in $C([0,T];H^{-1})\subset L^2(0,T;V_\sigma')$, which is all the Lions--Magenes lemma needs. Likewise $\rho\,\pa_h^\alpha v_t=\pa_i(\rho\,\pa_h^{\alpha-e_i}v_t)-\pa_i\varrho\,\pa_h^{\alpha-e_i}v_t\in L^2(0,T;H^{-1})$, since $v_t\in L^2(0,T;H^3)$. The surface pressure drops out by a mode-by-mode density argument: the functions in $V_\sigma$ with finitely many horizontal Fourier modes are dense in $V_\sigma$ (horizontal Fourier truncation preserves $H^1_0$ and the lid constraint), and for such $\varphi$, $\ip{\pa_h^\alpha\nh P,\varphi}=\ip{P,\pa_h^\alpha\divh\int_0^1\varphi\dd z}_{\Tl}=0$. Finally, the identity actually used in \Cref{prop:T4},
\begin{equation*}
\ip{\pa_h^\alpha\nh q,\pa_h^\alpha v}=g\int_\Omega\pa_h^\alpha\varrho\,\pa_h^\alpha w\dx ,
\end{equation*}
obtained as in the proof of \Cref{prop:utt} (horizontal integration by parts and Fubini's theorem in $z$, first for finitely many modes), has both factors on the right in $L^2(0,T;L^2)$: $\pa_h^\alpha\varrho\in C([0,T];L^2)$, and $\nrm{\pa_h^\alpha w}_{L^2}\le\nrm{\pa_h^\alpha\divh v}_{L^2}\le\nrm{v}_{H^5}$ by \Cref{lem:w}.

(c) Expressions such as $\ABmh|\nh|^5\varrho$ or $\ABmi\Laph^3\varrho$ are defined mode by mode via \Cref{lem:ABmean}. On the mode $k\ne0$, $\ip{A_k\phi,\phi}\ge|k|^2\nrm{\pz\phi}^2\ge\pi^2|k|^2\nrm{\phi}^2$, so $\nrm{A_k^{-1/2}}\le(\pi|k|)^{-1}$ and $\nrm{A_k^{-1}}\le(\pi|k|)^{-2}$. Hence $\nrm{\ABmh|\nh|^5\varrho}_{L^2}+\nrm{\ABmi\Laph^3\varrho}_{L^2}\le2\pi^{-1}\nrm{\nh^4\varrho}_{L^2}$, which is finite in the class \eqref{eq:localclass}.
\end{remark}

\subsection{Estimate of \texorpdfstring{$\cEutt$}{E utt}}\label{sec:A6-utt}

By \eqref{PE}$_1$ and $\sgm=-1/\rs'$,
\begin{equation}\label{eq:A6-rhot}
\varrho_t=\sgm^{-1}w-\bU\cdot\nabla\varrho,\qquad \varrho_{tt}=\sgm^{-1}w_t-\bU_t\cdot\nabla\varrho-\bU\cdot\nabla\varrho_t,
\end{equation}
and $\rho_t=\varrho_t$, $\rho_{tt}=\varrho_{tt}$. We first collect bounds for these time derivatives and for two convective products.

\begin{lemma}\label{lem:A6-rhot}
Assume \eqref{eq:small}. Then, on $[0,t]$,
\begin{align}
\nrm{\varrho_t}_{H^1}&\lesssim\nrm{v}_{H^2},\label{eq:6.4}\\
\nrm{\varrho_t}_{L^\infty}+\nrm{\varrho_t}_{H^2}&\lesssim\nrm{v}_{H^3},\label{eq:A6-rhotH2}\\
\nrm{\varrho_{tt}}_{L^2}&\lesssim\Ff+\nrm{v}_{H^2}\nrm{v}_{H^3},\label{eq:6.5}\\
\nrm{\bU\cdot\nabla v}_{L^3}&\lesssim\nrm{v}_{H^2}^2,\label{eq:A6-Uv}\\
\nrm{\bU_t\cdot\nabla v}_{L^2}&\lesssim\Ff\nrm{v}_{H^3}.\label{eq:A6-Utv}
\end{align}
\end{lemma}

\begin{proof}
Recall $\nrm{\nabla\varrho}_{L^\infty}+\nrm{\nabla^2\varrho}_{L^\infty}\lesssim\Ev\lesssim1$. By \eqref{eq:A6-rhot}, $\nrm{\varrho_t}\le\nrm{\sgm^{-1}}_{L^\infty}\nrm{w}+\nrm{\bU}\nrm{\nabla\varrho}_{L^\infty}\lesssim\nrm{v}_{H^1}$, and
\begin{equation*}
\nabla\varrho_t=\nabla(\sgm^{-1}w)-(\nabla\bU)^{\top}\nabla\varrho-(\bU\cdot\nabla)\nabla\varrho
\end{equation*}
gives $\nrm{\nabla\varrho_t}\lesssim\nrm{w}_{H^1}+\nrm{\bU}_{H^1}\big(\nrm{\nabla\varrho}_{L^\infty}+\nrm{\nabla^2\varrho}_{L^\infty}\big)\lesssim\nrm{v}_{H^2}$. This is \eqref{eq:6.4}. For \eqref{eq:A6-rhotH2}, \Cref{lem:w} gives $\nrm{w}_{H^2}\lesssim\nrm{v}_{H^3}$, and \Cref{lem:GN}(2) with $m=2$ gives
\begin{equation*}
\nrm{\bU\cdot\nabla\varrho}_{H^2}\lesssim\nrm{\bU}_{L^\infty}\nrm{\nabla\varrho}_{H^2}+\nrm{\bU}_{H^2}\nrm{\nabla\varrho}_{L^\infty}\lesssim\nrm{v}_{H^3}\Ev ;
\end{equation*}
the $L^\infty$ bound follows from $H^2\subset L^\infty$. For \eqref{eq:6.5}, by \eqref{eq:A6-rhot}, \Cref{lem:w}, \Cref{lem:poincare} and \eqref{eq:6.4},
\begin{equation*}
\nrm{\varrho_{tt}}\lesssim\nrm{w_t}+\big(\nrm{v_t}+\nrm{w_t}\big)\nrm{\nabla\varrho}_{L^\infty}+\nrm{\bU}_{L^\infty}\nrm{\nabla\varrho_t}\lesssim\nrm{\nabla v_t}+\nrm{v}_{H^3}\nrm{v}_{H^2},
\end{equation*}
and $\nrm{\nabla v_t}\lesssim\Ff$. Next, by $H^1\subset L^6$ and $\nrm{w}_{H^1}\lesssim\nrm{v}_{H^2}$,
\begin{equation*}
\nrm{\bU\cdot\nabla v}_{L^3}\le\nrm{v}_{L^\infty}\nrm{\nh v}_{L^3}+\nrm{w}_{L^6}\nrm{\pz v}_{L^6}\lesssim\nrm{v}_{H^2}^2 .
\end{equation*}
Finally, by \Cref{lem:w} and \Cref{lem:aniso}(e) with $s=\frac14$,
\begin{equation*}
\begin{aligned}
\nrm{\bU_t\cdot\nabla v}&\le\nrm{v_t}_{L^6}\nrm{\nh v}_{L^3}+\nrm{w_t}_{L^2_hL^\infty_z}\nrm{\pz v}_{L^\infty_hL^2_z}\\
&\lesssim\nrm{\nabla v_t}\nrm{v}_{H^2}+\nrm{\divh v_t}\nrm{v}_{H^{9/4}}\lesssim\Ff\nrm{v}_{H^3}.
\end{aligned}\qedhere
\end{equation*}
\end{proof}

\begin{proposition}\label{prop:utt}
Assume \eqref{eq:small}. Then
\begin{equation}\label{eq:6.1}
\cEutt^2(t)=\sup_{\tau\in(0,t)}\nrm{v_{tt}(\tau)}_{L^2}^2+\int_0^t\nrm{v_{tt}(\tau)}_{H^1}^2\dtau\lesssim\Es(0)^2+\Es^4(t).
\end{equation}
\end{proposition}

\begin{proof}
\emph{The energy identity.} Differentiating \eqref{PE}$_2$ twice in time and using $\rho_t=\varrho_t$, $\rho_{tt}=\varrho_{tt}$, we obtain
\begin{equation}\label{eq:6.2}
\rho v_{ttt}+\nh q_{tt}-\mu\Lap v_{tt}=-2\varrho_tv_{tt}-\varrho_{tt}v_t-\pt^2\big(\rho\,\bU\cdot\nabla v\big)-f\big(\varrho_{tt}v^\perp+2\varrho_tv_t^\perp+\rho v_{tt}^\perp\big).
\end{equation}
We test \eqref{eq:6.2} with $v_{tt}$ (\Cref{rem:justify}(b)). First, $\int\rho v_{ttt}\cdot v_{tt}=\frac12\ddt\nrm{\sqrt\rho v_{tt}}^2-\frac12\int\varrho_t|v_{tt}|^2$, and $\rho v_{tt}^\perp\cdot v_{tt}=0$ pointwise. Next, write $q_{tt}=P_{tt}-g\int_0^z\varrho_{tt}\dd s$ (\Cref{lem:propagation}). The surface pressure drops out because $\divh\int_0^1v_{tt}\dd z=0$. For the remaining part we integrate by parts horizontally and use Fubini's theorem in $z$ together with $\int_z^1\divh v_{tt}\dd z'=\divh\int_0^1v_{tt}\dd z'+w_{tt}(z)=w_{tt}(z)$:
\begin{equation*}
\int_\Omega\nh q_{tt}\cdot v_{tt}\dx=g\int_\Omega\Big(\int_0^z\varrho_{tt}\dd s\Big)\divh v_{tt}\dx=g\int_\Omega\varrho_{tt}\,w_{tt}\dx .
\end{equation*}
By \eqref{eq:A6-rhot}, $g\int\varrho_{tt}w_{tt}=\frac g2\ddt\nrm{\sgm^{-1/2}w_t}^2-g\int(\bU_t\cdot\nabla\varrho)w_{tt}-g\int(\bU\cdot\nabla\varrho_t)w_{tt}$. Hence
\begin{equation}\label{eq:6.3}
\frac12\ddt\Big(\nrm{\sqrt\rho\,v_{tt}}_{L^2}^2+g\nrm{\sgm^{-1/2}w_t}_{L^2}^2\Big)+\mu\nrm{\nabla v_{tt}}_{L^2}^2=\sum_{i=1}^6R_i,
\end{equation}
where
\begin{equation*}
\begin{aligned}
&R_1=-\frac32\int_\Omega\varrho_t|v_{tt}|^2\dx,\qquad R_2=g\int_\Omega(\bU_t\cdot\nabla\varrho)\,w_{tt}\dx,\qquad R_3=g\int_\Omega(\bU\cdot\nabla\varrho_t)\,w_{tt}\dx,\\
&R_4=-\int_\Omega\varrho_{tt}\,v_t\cdot v_{tt}\dx,\qquad R_5=-\int_\Omega\pt^2(\rho\,\bU\cdot\nabla v)\cdot v_{tt}\dx,\\
&R_6=-f\int_\Omega\big(\varrho_{tt}v^\perp+2\varrho_tv_t^\perp\big)\cdot v_{tt}\dx .
\end{aligned}
\end{equation*}
The leading Coriolis term has dropped out; $R_6$ collects the rotational terms in which a time derivative falls on the density.

\emph{Two embeddings.} For a.e.\ $\tau$, $v_{tt}(\tau)\in H^1$ vanishes on $\pa\Omega$ by the no-slip condition, so $\nrm{v_{tt}}\le C_P\nrm{\nabla v_{tt}}$ (\Cref{lem:poincare}), and by $H^1\subset L^6$ and interpolation,
\begin{equation}\label{eq:A6-L3}
\nrm{v_{tt}}_{L^6}\lesssim\nrm{\nabla v_{tt}},\qquad \nrm{v_{tt}}_{L^3}^2\le\nrm{v_{tt}}\nrm{v_{tt}}_{L^6}\lesssim\nrm{v_{tt}}\nrm{\nabla v_{tt}} .
\end{equation}
Moreover, for $h\in H^1(\Omega)$, Minkowski's inequality, \Cref{lem:aniso}(b) for a.e.\ $z$ and the Cauchy--Schwarz inequality in $z$ give
\begin{equation}\label{eq:A6-L4h}
\nrm{h}_{L^4_hL^2_z}^2\le\nrm{h}_{L^2_zL^4_h}^2\le C\int_0^1\nrm{h(\cdot,z)}_{L^2_h}\nrm{h(\cdot,z)}_{H^1_h}\dd z\le C\nrm{h}_{L^2}\nrm{h}_{H^1}.
\end{equation}

In the rest of the proof we write $a\eqdef\nrm{v}_{H^2}$. By \eqref{eq:A6-std}, $a+\Ff+\nrm{v}_{H^3}\lesssim\Esh\lesssim1$. Let $\epsilon\in(0,1)$; $C_\epsilon$ denotes a constant depending on $\epsilon$ and on the parameters listed in the notation.

\emph{The term $R_1$.} By \eqref{eq:6.4} and \eqref{eq:A6-L3},
\begin{equation*}
|R_1|\le\frac32\nrm{\varrho_t}_{L^3}\nrm{v_{tt}}_{L^3}^2\lesssim\nrm{\varrho_t}_{H^1}\nrm{v_{tt}}\nrm{\nabla v_{tt}}\lesssim a\nrm{v_{tt}}\nrm{\nabla v_{tt}}\le\epsilon\nrm{\nabla v_{tt}}^2+C_\epsilon a^2\nrm{v_{tt}}^2 .
\end{equation*}

\emph{The terms $R_2$ and $R_3$.} Since $\nrm{w_{tt}}\le\nrm{\divh v_{tt}}\le\sqrt2\nrm{\nabla v_{tt}}$ and $\nrm{\bU_t}\le\nrm{v_t}+\nrm{w_t}\lesssim\Ff$,
\begin{equation*}
|R_2|\le g\nrm{\bU_t}\nrm{\nabla\varrho}_{L^\infty}\nrm{w_{tt}}\lesssim\Ff\Ev\nrm{\nabla v_{tt}}\le\epsilon\nrm{\nabla v_{tt}}^2+C_\epsilon\Ff^2 .
\end{equation*}
For $R_3$ we use $\nrm{\bU}_{L^\infty}\lesssim\nrm{v}_{H^3}$ and \eqref{eq:6.4}:
\begin{equation*}
|R_3|\le g\nrm{\bU}_{L^\infty}\nrm{\nabla\varrho_t}\nrm{w_{tt}}\lesssim\nrm{v}_{H^3}\,a\,\nrm{\nabla v_{tt}}\le\epsilon\nrm{\nabla v_{tt}}^2+C_\epsilon a^2\nrm{v}_{H^3}^2 .
\end{equation*}

\emph{The term $R_4$.} By \eqref{eq:6.5}, $\nrm{v_t}_{L^3}\lesssim\nrm{v_t}_{H^1}\lesssim\Ff$ and \eqref{eq:A6-L3},
\begin{equation*}
|R_4|\le\nrm{\varrho_{tt}}\nrm{v_t}_{L^3}\nrm{v_{tt}}_{L^6}\lesssim\big(\Ff+a\nrm{v}_{H^3}\big)\Ff\nrm{\nabla v_{tt}}\le\epsilon\nrm{\nabla v_{tt}}^2+C_\epsilon\Ff^2 ,
\end{equation*}
where we used $\Ff+a\nrm{v}_{H^3}\lesssim1$.

\emph{The term $R_5$.} Expanding $\pt^2(\rho\,\bU\cdot\nabla v)$,
\begin{equation*}
\begin{aligned}
R_5&=-\int_\Omega\varrho_{tt}(\bU\cdot\nabla v)\cdot v_{tt}-2\int_\Omega\varrho_t(\bU_t\cdot\nabla v)\cdot v_{tt}-2\int_\Omega\varrho_t(\bU\cdot\nabla v_t)\cdot v_{tt}\\
&\quad-\int_\Omega\rho(\bU_{tt}\cdot\nabla v)\cdot v_{tt}-2\int_\Omega\rho(\bU_t\cdot\nabla v_t)\cdot v_{tt}-\int_\Omega\rho(\bU\cdot\nabla v_{tt})\cdot v_{tt}\eqdef\sum_{i=1}^6R_{5i}.
\end{aligned}
\end{equation*}
By \eqref{eq:6.5}, \eqref{eq:A6-Uv} and \eqref{eq:A6-L3},
\begin{equation*}
|R_{51}|\le\nrm{\varrho_{tt}}\nrm{\bU\cdot\nabla v}_{L^3}\nrm{v_{tt}}_{L^6}\lesssim\big(\Ff+a\nrm{v}_{H^3}\big)a^2\nrm{\nabla v_{tt}}\le\epsilon\nrm{\nabla v_{tt}}^2+C_\epsilon\big(\Ff^2+a^2\nrm{v}_{H^3}^2\big),
\end{equation*}
since $a\lesssim1$. By \eqref{eq:6.4}, \eqref{eq:A6-Utv}, $\nrm{\bU}_{L^\infty}\lesssim\nrm{v}_{H^3}$ and \eqref{eq:A6-L3},
\begin{equation*}
|R_{52}|+|R_{53}|\le2\nrm{\varrho_t}_{L^3}\big(\nrm{\bU_t\cdot\nabla v}+\nrm{\bU}_{L^\infty}\nrm{\nabla v_t}\big)\nrm{v_{tt}}_{L^6}\lesssim a\,\nrm{v}_{H^3}\Ff\nrm{\nabla v_{tt}}\le\epsilon\nrm{\nabla v_{tt}}^2+C_\epsilon\Ff^2 .
\end{equation*}
For $R_{54}$ we split $\bU_{tt}\cdot\nabla v=v_{tt}\cdot\nh v+w_{tt}\pz v$. The first part is bounded by \eqref{eq:A6-L3}:
\begin{equation*}
\Big|\int_\Omega\rho(v_{tt}\cdot\nh v)\cdot v_{tt}\Big|\lesssim\nrm{\nh v}_{L^3}\nrm{v_{tt}}_{L^3}^2\lesssim a\nrm{v_{tt}}\nrm{\nabla v_{tt}} .
\end{equation*}
In the second part the factor $w_{tt}$ is controlled only in $L^2_hL^\infty_z$. We use H\"older's inequality in $z$ and then in $x_h$, with $\nrm{w_{tt}}_{L^2_hL^\infty_z}\le\nrm{\divh v_{tt}}$, and \eqref{eq:A6-L4h} for $\pz v$ and $v_{tt}$:
\begin{equation*}
\begin{aligned}
\Big|\int_\Omega\rho\,w_{tt}\,\pz v\cdot v_{tt}\Big|&\lesssim\nrm{w_{tt}}_{L^2_hL^\infty_z}\nrm{\pz v}_{L^4_hL^2_z}\nrm{v_{tt}}_{L^4_hL^2_z}\lesssim\nrm{\nabla v_{tt}}\,\nrm{v}_{H^2}\,\nrm{v_{tt}}^{\frac12}\nrm{v_{tt}}_{H^1}^{\frac12}\\
&\lesssim a\nrm{v_{tt}}^{\frac12}\nrm{\nabla v_{tt}}^{\frac32}\le\epsilon\nrm{\nabla v_{tt}}^2+C_\epsilon a^4\nrm{v_{tt}}^2 .
\end{aligned}
\end{equation*}
Hence $|R_{54}|\le2\epsilon\nrm{\nabla v_{tt}}^2+C_\epsilon(a^2+a^4)\nrm{v_{tt}}^2$. For $R_{55}$ we split $\bU_t\cdot\nabla v_t=v_t\cdot\nh v_t+w_t\pz v_t$. By \eqref{eq:A6-L3},
\begin{equation*}
\Big|\int_\Omega\rho(v_t\cdot\nh v_t)\cdot v_{tt}\Big|\lesssim\nrm{v_t}_{L^3}\nrm{\nabla v_t}\nrm{v_{tt}}_{L^6}\lesssim\Ff^2\nrm{\nabla v_{tt}},
\end{equation*}
and, as for $R_{54}$, using $\nrm{w_t}_{L^2_hL^\infty_z}\le\nrm{\divh v_t}\lesssim\Ff$, \eqref{eq:A6-L4h}, $\nrm{\pz v_t}\lesssim\Ff$, $\nrm{\pz v_t}_{H^1}\le\Eut$ and $\nrm{v_{tt}}\le C_P\nrm{\nabla v_{tt}}$,
\begin{equation*}
\Big|\int_\Omega\rho\,w_t\,\pz v_t\cdot v_{tt}\Big|\lesssim\nrm{w_t}_{L^2_hL^\infty_z}\nrm{\pz v_t}_{L^4_hL^2_z}\nrm{v_{tt}}_{L^4_hL^2_z}\lesssim\Ff\cdot\Ff^{\frac12}\Eut^{\frac12}\cdot\nrm{\nabla v_{tt}} .
\end{equation*}
Since $\Ff+\Eut\lesssim1$, $|R_{55}|\le\epsilon\nrm{\nabla v_{tt}}^2+C_\epsilon\Ff^2$. Finally, since $\divv\bU=0$, $\bU\cdot n=0$ on $\pa\Omega$ and $\bU\cdot\nabla\rho=-\rho_t=-\varrho_t$,
\begin{equation*}
R_{56}=\frac12\int_\Omega(\bU\cdot\nabla\rho)|v_{tt}|^2\dx=-\frac12\int_\Omega\varrho_t|v_{tt}|^2\dx,
\end{equation*}
which is bounded like $R_1$: $|R_{56}|\le\epsilon\nrm{\nabla v_{tt}}^2+C_\epsilon a^2\nrm{v_{tt}}^2$.

\emph{The rotational term $R_6$.} By \eqref{eq:6.5}, \eqref{eq:6.4}, \eqref{eq:A6-L3} and $\nrm{v}_{L^3}\lesssim\nrm{v}_{H^1}\lesssim\Ff$ (see \eqref{eq:3.13b}),
\begin{equation*}
\begin{aligned}
|R_6|&\le|f|\big(\nrm{\varrho_{tt}}\nrm{v}_{L^3}+2\nrm{\varrho_t}_{L^3}\nrm{v_t}\big)\nrm{v_{tt}}_{L^6}\lesssim\Big(\big(\Ff+a\nrm{v}_{H^3}\big)\Ff+a\Ff\Big)\nrm{\nabla v_{tt}}\\
&\lesssim\Ff\nrm{\nabla v_{tt}}\le\epsilon\nrm{\nabla v_{tt}}^2+C_\epsilon\Ff^2 .
\end{aligned}
\end{equation*}

\emph{Conclusion.} We insert these bounds into \eqref{eq:6.3}, choose $\epsilon$ so small that the total coefficient of $\nrm{\nabla v_{tt}}^2$ on the right is at most $\frac\mu2$, and use $a^4\lesssim a^2$. With a constant $C_{\mathrm{tt}}\ge1$ depending only on $\mu,g,f,l,\rs$, $\nrm{\rs+\vin}_{L^\infty}$ and a positive lower bound for $\rs+\vin$,
\begin{equation}\label{eq:6.6}
\frac12\ddt\Big(\nrm{\sqrt\rho\,v_{tt}}^2+g\nrm{\sgm^{-1/2}w_t}^2\Big)+\frac\mu2\nrm{\nabla v_{tt}}^2\le C_{\mathrm{tt}}\Big(\nrm{v}_{H^2}^2\nrm{v_{tt}}^2+\Ff^2+\nrm{v}_{H^2}^2\nrm{v}_{H^3}^2\Big).
\end{equation}
The first term on the right is absorbed. By the Poincar\'e inequality and \eqref{eq:A6-std}, $C_{\mathrm{tt}}\nrm{v}_{H^2}^2\nrm{v_{tt}}^2\le C\,C_P^2\Esh^2\nrm{\nabla v_{tt}}^2\le4C\,C_P^2\Es^2(t)\nrm{\nabla v_{tt}}^2\le\frac\mu4\nrm{\nabla v_{tt}}^2$, the last inequality being one of the absorption conditions fixed in \Cref{rem:c1}(i). Therefore
\begin{equation}\label{eq:6.8}
\frac12\ddt\Big(\nrm{\sqrt\rho\,v_{tt}}^2+g\nrm{\sgm^{-1/2}w_t}^2\Big)+\frac\mu4\nrm{\nabla v_{tt}}^2\le C_{\mathrm{tt}}\Big(\Ff^2+\nrm{v}_{H^2}^2\nrm{v}_{H^3}^2\Big).
\end{equation}
We integrate \eqref{eq:6.8} over $(0,\tau)$, $\tau\le t$. At the initial time, $\nrm{\sqrt\rho\,v_{tt}(0)}^2+g\nrm{\sgm^{-1/2}w_t(0)}^2\lesssim\Eutt^2(0)+\nrm{\nabla v_t(0)}^2\lesssim\Es(0)^2$. By \eqref{eq:5.4} and $\Ef(0)\lesssim\Esh(0)$, $\int_0^t\Ff^2\dtau\le8\Ef^2(0)\lesssim\Es(0)^2$. By \eqref{eq:A6-std} and $\nrm{v}_{H^3}\le\Fu$,
\begin{equation*}
\int_0^t\nrm{v}_{H^2}^2\nrm{v}_{H^3}^2\dtau\le\sup_{(0,t)}\nrm{v}_{H^2}^2\int_0^t\Fu^2\dtau\lesssim\Es^4(t).
\end{equation*}
Hence
\begin{equation}\label{eq:A6-uttint}
\sup_{\tau\in(0,t)}\Big(\nrm{\sqrt\rho\,v_{tt}}^2+g\nrm{\sgm^{-1/2}w_t}^2\Big)+\frac\mu2\int_0^t\nrm{\nabla v_{tt}}^2\dtau\lesssim\Es(0)^2+\Es^4(t).
\end{equation}
Since $\rho\ge\inf\rin>0$ and $\nrm{v_{tt}}_{H^1}\lesssim\nrm{\nabla v_{tt}}$, this gives \eqref{eq:6.1}; the implicit constant depends on $\mu,g,f,l,\rs$, $\nrm{\rs+\vin}_{L^\infty}$ and a positive lower bound for $\rs+\vin$.
\end{proof}

\begin{remark}\label{rem:A6-3D}
Two points in the proof of \Cref{prop:utt} are specific to three dimensions. (i) The cubic terms $\int|\varrho_t||v_{tt}|^2$ and $\int|\nh v||v_{tt}|^2$ are estimated with coefficients in $L^3$ and \eqref{eq:A6-L3}. The planar bound $\nrm{v_{tt}}_{L^4}^2\lesssim\nrm{v_{tt}}\nrm{\nabla v_{tt}}$ fails in three dimensions, where only $\nrm{v_{tt}}_{L^4}^2\lesssim\nrm{v_{tt}}^{1/2}\nrm{\nabla v_{tt}}^{3/2}$ holds, and Young's inequality then puts a coefficient of fourth order in front of $\nrm{v_{tt}}^2$. This occurs only in $R_{54}$, through the anisotropic product $w_{tt}\pz v\cdot v_{tt}$ (where $w_{tt}$ is controlled only in $L^2_hL^\infty_z$ and $H^1(\T^2)\not\subset L^\infty(\T^2)$ excludes $\nrm{v_{tt}}_{L^\infty_hL^2_z}$); the resulting term $\nrm{v}_{H^2}^4\nrm{v_{tt}}^2$ is absorbed by the Poincar\'e inequality. (ii) The vertical velocity is not bounded in $L^\infty$ by $\nrm{v}_{H^2}$, because $\nrm{w}_{L^\infty}$ requires $H^{1+s}_h$-control of $\divh v$. Therefore $R_3$ produces $\nrm{v}_{H^2}^2\nrm{v}_{H^3}^2$ instead of a fourth power of $\nrm{v}_{H^2}$. This term is integrable in time since $\nrm{v}_{H^3}\le\Fu$, and it is compatible with the time weight of \Cref{prop:uttdecay} by the decay of $\nrm{v}_{H^2}$. The rotational terms $R_6$ are cubic and bounded by $\Ff\nrm{\nabla v_{tt}}$.
\end{remark}

\begin{proposition}\label{prop:uttdecay}
Assume \eqref{eq:small} and let $0<\delta<2$. Then
\begin{equation*}
\nrm{v_{tt}(t)}_{L^2}\lesssim\delta^{-\frac12}\Es(t)(1+t)^{-1+\frac\delta2},
\end{equation*}
and
\begin{equation}\label{eq:6.7}
\int_0^t(1+\tau)^{2-\delta}\nrm{\nabla v_{tt}(\tau)}_{L^2}^2\dtau\lesssim\delta^{-1}\Es^2(t).
\end{equation}
Finally, the pointwise rate holds without loss: for $0\le\tau\le t$,
\begin{equation}\label{eq:6.7b}
\nrm{v_{tt}(\tau)}_{L^2}+\nrm{w_t(\tau)}_{L^2}\lesssim\Es(t)(1+\tau)^{-1},\qquad \int_\tau^t\nrm{\nabla v_{tt}(s)}_{L^2}^2\dd s\lesssim\Es^2(t)(1+\tau)^{-2}.
\end{equation}
\end{proposition}

\begin{proof}
Let $Y\eqdef\nrm{\sqrt\rho\,v_{tt}}^2+g\nrm{\sgm^{-1/2}w_t}^2$ and $S\eqdef\Ff^2+\nrm{v}_{H^2}^2\nrm{v}_{H^3}^2$, so that \eqref{eq:6.8} reads $\frac12Y'+\frac\mu4\nrm{\nabla v_{tt}}^2\le C_{\mathrm{tt}}S$. By \eqref{eq:5.4} and the monotonicity of $\Es$, $\nrm{v(\tau)}_{H^2}\lesssim\Es(t)(1+\tau)^{-1}$ for $0\le\tau\le t$. Hence
\begin{equation}\label{eq:6.9}
(1+\tau)^{2-\delta}\nrm{v(\tau)}_{H^2}^2\nrm{v(\tau)}_{H^3}^2\lesssim\Es^2(t)(1+\tau)^{-\delta}\nrm{v(\tau)}_{H^3}^2\le\Es^2(t)\Fu^2(\tau),
\end{equation}
and, together with \eqref{eq:5.5},
\begin{equation*}
\int_0^t(1+\tau)^{2-\delta}S\dtau\lesssim\delta^{-1}\Es^2(t)+\Es^4(t)\lesssim\delta^{-1}\Es^2(t).
\end{equation*}
Let $k=2-\delta\in(0,2)$. Multiplying \eqref{eq:6.8} by $(1+\tau)^k$ and integrating over $(0,t)$,
\begin{equation}\label{eq:A6-weighted}
\frac12(1+t)^kY(t)+\frac\mu4\int_0^t(1+\tau)^k\nrm{\nabla v_{tt}}^2\dtau\le\frac12Y(0)+C\delta^{-1}\Es^2(t)+\frac k2\int_0^t(1+\tau)^{k-1}Y\dtau .
\end{equation}
We bound the last integral. For the $w_t$-part, $\nrm{w_t}\le\nrm{\divh v_t}\lesssim\Ff$ gives
\begin{equation*}
\frac k2\int_0^t(1+\tau)^{k-1}g\nrm{\sgm^{-1/2}w_t}^2\dtau\lesssim\int_0^t(1+\tau)^{1-\delta}\Ff^2\dtau\lesssim\Es^2(t);
\end{equation*}
for $\delta<1$ this is \eqref{eq:5.5} with $\delta$ replaced by $1+\delta\in(1,2)$, and for $\delta\ge1$ the weight is at most $1$ and $\int_0^t\Ff^2\le8\Ef^2(0)\lesssim\Es^2(t)$. For the $v_{tt}$-part let $T_1\eqdef4kC_P^2\mu^{-1}\nrm{\rs+\vin}_{L^\infty}$. If $1+\tau\ge T_1$, the Poincar\'e inequality gives
\begin{equation*}
\begin{aligned}
\frac k2(1+\tau)^{k-1}\nrm{\sqrt\rho\,v_{tt}}^2&\le\frac k2\nrm{\rs+\vin}_{L^\infty}C_P^2(1+\tau)^{k-1}\nrm{\nabla v_{tt}}^2\\
&=\frac\mu8T_1(1+\tau)^{k-1}\nrm{\nabla v_{tt}}^2\le\frac\mu8(1+\tau)^k\nrm{\nabla v_{tt}}^2,
\end{aligned}
\end{equation*}
which is absorbed by the left-hand side of \eqref{eq:A6-weighted}. On the remaining set $\{\tau:1+\tau<T_1\}$, \eqref{eq:A6-uttint} gives, uniformly in $k\in(0,2)$,
\begin{equation*}
\frac k2\int_{\{1+\tau<T_1\}}(1+\tau)^{k-1}\nrm{\sqrt\rho\,v_{tt}}^2\dtau\le T_1\max\{1,T_1\}\sup_{(0,t)}Y\lesssim\Es(0)^2+\Es^4(t)\lesssim\Es^2(t).
\end{equation*}
Finally $Y(0)\lesssim\Es(0)^2$. Inserting these bounds into \eqref{eq:A6-weighted},
\begin{equation*}
(1+t)^{2-\delta}\Big(\nrm{\sqrt\rho\,v_{tt}(t)}^2+g\nrm{\sgm^{-1/2}w_t(t)}^2\Big)+\frac\mu4\int_0^t(1+\tau)^{2-\delta}\nrm{\nabla v_{tt}}^2\dtau\le\frac C\delta\Es^2(t),
\end{equation*}
with $C$ depending only on $\mu,g,f,l,\rs$, $\nrm{\rs+\vin}_{L^\infty}$ and a positive lower bound for $\rs+\vin$. This proves \eqref{eq:6.7}, and $\rho\ge\inf\rin>0$ gives $\nrm{v_{tt}(t)}\le C\delta^{-1/2}\Es(t)(1+t)^{-1+\delta/2}$.

\emph{The rate without loss.} By the Poincar\'e inequality for $v_{tt}$, whose trace vanishes for a.e.\ time, and $\nrm{w_t}\le\nrm{\divh v_t}\lesssim\Ff$, one has $Y\le C_Y(\nrm{\nabla v_{tt}}^2+\Ff^2)$. Moreover $\nrm{v}_{H^2}^2\nrm{v}_{H^3}^2\lesssim\Es^2(t)(1+\tau)^{-2}\Esh^2\lesssim\Es^2(t)(1+\tau)^{-2}$ on $[0,t]$, by \eqref{eq:5.4} and \eqref{eq:small}. Hence, with $c\eqdef\mu/(4C_Y)$, \eqref{eq:6.8} gives
\begin{equation*}
Y'+cY+\frac\mu4\nrm{\nabla v_{tt}}^2\le C\Ff^2+C\Es^2(t)(1+\tau)^{-2}\qquad\text{on }[0,t],
\end{equation*}
and Duhamel's formula gives, for $0\le\tau\le t$,
\begin{equation*}
Y(\tau)\le e^{-c\tau}Y(0)+C\int_0^\tau e^{-c(\tau-s)}\Ff^2(s)\dd s+C\Es^2(t)\int_0^\tau e^{-c(\tau-s)}(1+s)^{-2}\dd s .
\end{equation*}
Splitting the second integral at $s=\tau/2$, and using $\int_0^t\Ff^2\le8\Ef^2(0)\lesssim\Es^2(t)$ on $[0,\tau/2]$ and \eqref{eq:5.tail} on $[\tau/2,\tau]$, it is $\lesssim\Es^2(t)(1+\tau)^{-2}$; the third integral is $\lesssim(1+\tau)^{-2}$ in the same way, and $Y(0)\lesssim\Es(0)^2$. Hence $Y(\tau)\lesssim\Es^2(t)(1+\tau)^{-2}$, which gives the first bound in \eqref{eq:6.7b} because $\rho$ and $\sgm^{-1}$ are bounded below. Integrating \eqref{eq:6.8} over $[\tau,t]$ and using $Y\ge0$, this bound, \eqref{eq:5.tail} and $\int_\tau^t\nrm{v}_{H^2}^2\nrm{v}_{H^3}^2\lesssim\Es^2(t)(1+\tau)^{-2}\int_0^t\Fu^2\lesssim\Es^2(t)(1+\tau)^{-2}$ gives the second bound.
\end{proof}

\subsection{Estimate of \texorpdfstring{$\cEut$}{E ut}}\label{sec:A6-ut}

\begin{proposition}\label{prop:ut}
Assume \eqref{eq:small}. Then
\begin{equation}\label{eq:6.10}
\cEut^2(t)=\sup_{\tau\in(0,t)}\Eut^2+\int_0^t\Fut^2\dtau\lesssim\Es(0)^2+\Eb^2+\Eb\Es+\Es^3 .
\end{equation}
\end{proposition}

\begin{proof}
By \eqref{eq:4.3}, \eqref{eq:4.11} and \eqref{eq:small},
\begin{equation}\label{eq:6.11}
\nrm{v}_{H^2}\lesssim\Ef,\qquad \nrm{v}_{H^2}\lesssim\Ff^{\frac23}\Fuz^{\frac13}+\Ff .
\end{equation}
We use the temporal Stokes estimate \eqref{eq:4.17}, in which $\zeta_t=-G_t+g\int_0^z\nh\varrho_t\dd s$ and, by \eqref{eq:Gdef},
\begin{equation}\label{eq:A6-Gt}
G_t=\varrho_t\dot v+\rho\,\dot v_t+f\big(\varrho_tv^\perp+\rho\,v_t^\perp\big),\qquad \dot v_t\eqdef\pt\dot v=v_{tt}+\bU_t\cdot\nabla v+\bU\cdot\nabla v_t .
\end{equation}
The Coriolis part $f(\rho v^\perp)_t=f(\varrho_tv^\perp+\rho v_t^\perp)$ is of lower order: besides the bounded factor $\varrho_t$, it involves only $v$, $v_t$ and their first derivatives, which are controlled by $\Ef$ and $\Ff$.

\textbf{Step 1. Estimate of $\nrm{v_t}_{H^2}$.} By \eqref{eq:4.17} with $k=0$,
\begin{equation}\label{eq:6.12}
\nrm{v_t}_{H^2}\lesssim\nrm{G_t}_{L^2}+\norm{\int_0^z\nh\varrho_t\dd s}_{L^2}.
\end{equation}
By the Cauchy--Schwarz inequality in $z$, \eqref{eq:6.4} and \eqref{eq:6.11}, $\nrm{\int_0^z\nh\varrho_t\dd s}\le\nrm{\nh\varrho_t}\lesssim\nrm{v}_{H^2}\lesssim\Ef$. For $G_t$ we use \eqref{eq:A6-rhotH2}, $\nrm{\dot v}\le\nrm{v_t}+\nrm{\bU\cdot\nabla v}\lesssim\Ef+\nrm{v}_{H^2}^2\lesssim\Ef$ (by \eqref{eq:A6-Uv} and \eqref{eq:6.11}), and
\begin{equation*}
\begin{aligned}
\nrm{\bU_t\cdot\nabla v}&\le\nrm{v_t}_{L^6}\nrm{\nh v}_{L^3}+\nrm{w_t}_{L^\infty_zL^4_h}\nrm{\pz v}_{L^2_zL^4_h}\\
&\lesssim\big(\Ff+\nrm{\nh\divh v_t}\big)\nrm{v}_{H^2}\lesssim\Ef\big(\Ff+\nrm{v_t}_{H^2}\big),\\
\nrm{\bU\cdot\nabla v_t}&\le\nrm{v}_{L^\infty}\nrm{\nh v_t}+\nrm{w}_{L^6}\nrm{\pz v_t}_{L^3}\\
&\lesssim\nrm{v}_{H^2}\big(\Ff+\nrm{v_t}_{H^{3/2}}\big)\lesssim\Ef\big(\Ff+\nrm{v_t}_{H^2}\big),
\end{aligned}
\end{equation*}
where we used \Cref{lem:w} for $w_t$ (namely $\nrm{w_t}_{L^\infty_zL^4_h}\le C\nrm{\nh\divh v_t}$), \Cref{lem:aniso}(f) with $q=4$, and $H^{1/2}\subset L^3$. Hence, by \eqref{eq:A6-Gt} and $\nrm{\varrho_t}_{L^\infty}\lesssim\nrm{v}_{H^3}\lesssim1$,
\begin{equation*}
\begin{aligned}
\nrm{G_t}&\le\nrm{\varrho_t}_{L^\infty}\nrm{\dot v}+\nrm{\rho}_{L^\infty}\big(\nrm{v_{tt}}+\nrm{\bU_t\cdot\nabla v}+\nrm{\bU\cdot\nabla v_t}\big)+|f|\big(\nrm{\varrho_t}_{L^\infty}\nrm{v}+\nrm{\rho}_{L^\infty}\nrm{v_t}\big)\\
&\lesssim\Ef+\Eutt+\Ef\nrm{v_t}_{H^2},
\end{aligned}
\end{equation*}
where we also used $\Ff\lesssim1$, $\nrm{v}\le\nrm{v}_{H^2}\lesssim\Ef$ and $\nrm{v_t}\lesssim\Ef$. Inserting this into \eqref{eq:6.12} gives $\nrm{v_t}_{H^2}\le C_*(\Ef+\Eutt)+C_*\Ef\nrm{v_t}_{H^2}$. Since $C_*\Ef\le C_*C\Esh\le2C_*C\Es(t)\le\frac12$ by \eqref{eq:small} and \Cref{rem:c1}(i), we conclude
\begin{equation}\label{eq:6.13}
\nrm{v_t}_{H^2}\lesssim\Ef+\Eutt .
\end{equation}
(The products $w_t\pz v$ and $w\pz v_t$ are the reason for this absorption: in three dimensions $w_t$ cannot be put in $L^\infty_h$ at the level of $\nrm{\nabla v_t}$, so $\nrm{v_t}_{H^2}$ itself appears, multiplied by the small factor $\Ef$.) Squaring \eqref{eq:6.13} and using \eqref{eq:6.1} and $\Ef^2\le\Eb^2$,
\begin{equation}\label{eq:6.14}
\sup_{\tau\in(0,t)}\nrm{v_t}_{H^2}^2\lesssim\Es(0)^2+\Eb^2+\Es^4\lesssim\Es(0)^2+\Eb^2+\Es^3 .
\end{equation}

\textbf{Step 2. Estimate of $\nrm{v_t}_{H^3}$.} By \eqref{eq:4.17} with $k=1$ and the horizontal Poincar\'e inequality (the function $\nh\varrho_t$ has zero horizontal mean),
\begin{equation}\label{eq:6.15}
\nrm{v_t}_{H^3}\lesssim\nrm{G_t}_{H^1}+\norm{\int_0^z\nh\varrho_t\dd s}_{H^1}\lesssim\nrm{G_t}_{H^1}+\nrm{\nh^2\varrho_t}_{L^2}.
\end{equation}
For the density term we use $\varrho_t=\sgm^{-1}w-\bU\cdot\nabla\varrho$ and expand $\nh^2(v\cdot\nh\varrho+w\pz\varrho)$ by the Leibniz rule:
\begin{equation}\label{eq:6.16}
\begin{aligned}
\nrm{\nh^2\varrho_t}&\lesssim\nrm{\nh^2w}\big(1+\nrm{\pz\varrho}_{L^\infty}\big)+\nrm{\nh^2v}\nrm{\nh\varrho}_{L^\infty}+\nrm{\nh v}_{L^4}\nrm{\nh^2\varrho}_{L^4}+\nrm{v}_{L^\infty}\nrm{\nh^3\varrho}\\
&\quad+\nrm{\nh w}\nrm{\nh\pz\varrho}_{L^\infty}+\nrm{w}_{L^\infty_zL^4_h}\nrm{\nh^2\pz\varrho}_{L^2_zL^4_h}\\
&\lesssim\nrm{\nh^3v}+\nrm{v}_{H^2}\nrm{\varrho}_{H^4}\lesssim\Ff^{\frac12}\Fuz^{\frac12}+\Esh\Fsh .
\end{aligned}
\end{equation}
Here we used $\nrm{\nh^2w}\le\nrm{\nh^2\divh v}\le\sqrt2\nrm{\nh^3v}$ and $\nrm{\nh w}+\nrm{w}_{L^\infty_zL^4_h}\lesssim\nrm{\nh\divh v}$ (\Cref{lem:w}), $\nrm{\nh^2\pz\varrho}_{L^2_zL^4_h}\lesssim\nrm{\varrho}_{H^{7/2}}$ (\Cref{lem:aniso}(f)), $\nrm{v}_{H^2}\le\Fu\le\Fsh$, and, by Parseval in $x_h$, the Cauchy--Schwarz inequality and the vertical Poincar\'e inequality,
\begin{equation*}
\nrm{\nh^3v}^2\le\nrm{\nh^2v}\nrm{\nh^4v}\le\pi^{-1}\nrm{\nabla\nh v}\nrm{\pz\nh^4v}\lesssim\Ff\Fuz .
\end{equation*}
Combining \eqref{eq:6.15} and \eqref{eq:6.16},
\begin{equation}\label{eq:6.17}
\nrm{v_t}_{H^3}\lesssim\nrm{G_t}_{H^1}+\Ff^{\frac12}\Fuz^{\frac12}+\Esh\Fsh .
\end{equation}
For a scalar $\phi$ and a vector field $X$, $\nrm{\phi X}_{H^1}\le\nrm{\phi}_{L^\infty}\nrm{X}_{H^1}+\nrm{\nabla\phi}_{L^3}\nrm{X}_{L^6}\lesssim\nrm{\phi}_{H^2}\nrm{X}_{H^1}$. Applying this with $\phi=\varrho_t$ and \eqref{eq:A6-rhotH2}, and using $\nrm{v}_{H^1}+\nrm{v_t}_{H^1}\lesssim\Ff$ (by \eqref{eq:3.13b} and \Cref{lem:poincare}), we obtain from \eqref{eq:A6-Gt}
\begin{equation}\label{eq:6.18}
\begin{aligned}
\nrm{G_t}_{H^1}&\le\nrm{\varrho_t\dot v}_{H^1}+\nrm{\rho\,\dot v_t}_{H^1}+|f|\big(\nrm{\varrho_tv^\perp}_{H^1}+\nrm{\rho\,v_t^\perp}_{H^1}\big)\\
&\lesssim\nrm{\varrho_t}_{H^2}\big(\nrm{\dot v}_{H^1}+\nrm{v}_{H^1}\big)+\big(1+\nrm{\nabla\varrho}_{L^\infty}\big)\big(\nrm{\dot v_t}_{H^1}+\nrm{v_t}_{H^1}\big)\\
&\lesssim\nrm{\dot v}_{H^1}+\nrm{\dot v_t}_{H^1}+\Ff .
\end{aligned}
\end{equation}
In particular, the rotational part of $G_t$ contributes only $\Ff$. For $\dot v=v_t+v\cdot\nh v+w\pz v$ we have, with $\nabla w=(\nh w,-\divh v)$,
\begin{equation}\label{eq:6.19}
\begin{aligned}
\nrm{\dot v}_{H^1}&\lesssim\nrm{v_t}_{H^1}+\nrm{v}_{L^\infty}\nrm{\nh v}_{H^1}+\nrm{\nabla v}_{L^4}\nrm{\nh v}_{L^4}+\nrm{w}_{L^\infty_zL^4_h}\big(\nrm{\pz v}_{L^2_zL^4_h}+\nrm{\nabla\pz v}_{L^2_zL^4_h}\big)\\
&\quad+\nrm{\nh w}_{L^2_hL^\infty_z}\nrm{\pz v}_{L^\infty_hL^2_z}+\nrm{\divh v}_{L^4}\nrm{\pz v}_{L^4}\\
&\lesssim\Ff\big(1+\nrm{v}_{H^{5/2}}+\nrm{v}_{H^{9/4}}\big)\lesssim\Ff .
\end{aligned}
\end{equation}
Here $\nrm{\nh v}_{H^1}+\nrm{\nh v}_{L^4}+\nrm{\divh v}_{L^4}+\nrm{w}_{L^\infty_zL^4_h}+\nrm{\nh w}_{L^2_hL^\infty_z}\lesssim\nrm{\nabla\nh v}\lesssim\Ff$ by \Cref{lem:poincare} and \Cref{lem:w}; the factors $\pz v$ and $\nabla\pz v$ were placed in $L^\infty_hL^2_z$ or $L^2_zL^4_h$ by \Cref{lem:aniso}(e) (with $s=\frac14$) and (f), and $\nrm{v}_{L^\infty}+\nrm{\nabla v}_{L^4}\lesssim\nrm{v}_{H^2}\lesssim1$. For $\dot v_t=v_{tt}+v_t\cdot\nh v+w_t\pz v+v\cdot\nh v_t+w\pz v_t$ we have, with $\nabla w_t=(\nh w_t,-\divh v_t)$,
\begin{equation}\label{eq:6.20}
\begin{aligned}
\nrm{\dot v_t}_{H^1}&\lesssim\nrm{v_{tt}}_{H^1}+\nrm{v_t}_{H^1}\nrm{\nh v}_{L^\infty}+\nrm{v_t}_{L^6}\nrm{\nabla\nh v}_{L^3}+\nrm{\divh v_t}\nrm{\pz v}_{L^\infty}\\
&\quad+\big(\nrm{w_t}_{L^2_hL^\infty_z}+\nrm{\nh w_t}_{L^2_hL^\infty_z}\big)\big(\nrm{\pz v}_{L^\infty_hL^2_z}+\nrm{\nabla\pz v}_{L^\infty_hL^2_z}\big)\\
&\quad+\big(\nrm{v}_{W^{1,\infty}}+\nrm{w}_{W^{1,\infty}}\big)\nrm{v_t}_{H^2}\\
&\lesssim\Futt+\Ff\nrm{v}_{H^3}+\nrm{v_t}_{H^2}\nrm{v}_{H^{13/4}}+\Eu\nrm{v_t}_{H^2}\lesssim\Futt+\Esh\Fsh .
\end{aligned}
\end{equation}
Here we used $\nrm{w_t}_{L^2_hL^\infty_z}+\nrm{\nh w_t}_{L^2_hL^\infty_z}\le\nrm{\divh v_t}+\nrm{\nh\divh v_t}\lesssim\nrm{v_t}_{H^2}$ (\Cref{lem:w}), \Cref{lem:aniso}(e) with $s=\frac14$ for the factors $\pz v$, $\nabla\pz v$ in $L^\infty_hL^2_z$, $\nrm{v}_{W^{1,\infty}}\lesssim\nrm{v}_{H^3}$, $\nrm{w}_{W^{1,\infty}}\lesssim\Eu$, and $\nrm{v_t}_{H^2}\le\Fut\le\Fsh$, $\Ff\lesssim\Fsh$, $\nrm{v}_{H^{13/4}}\le\Eu\le\Esh$. Inserting \eqref{eq:6.19} and \eqref{eq:6.20} into \eqref{eq:6.18} gives
\begin{equation*}
\nrm{G_t}_{H^1}\lesssim\Ff+\Futt+\Esh\Fsh ,
\end{equation*}
and with \eqref{eq:6.17},
\begin{equation}\label{eq:6.21}
\nrm{v_t}_{H^3}^2\lesssim\Ff\Fuz+\Ff^2+\Futt^2+\Esh^2\Fsh^2 .
\end{equation}
We integrate \eqref{eq:6.21} over $(0,t)$. By \eqref{eq:6.1}, $\int_0^t\Futt^2\dtau\le\cEutt^2(t)\lesssim\Es(0)^2+\Es^4$. By the Cauchy--Schwarz inequality, $\int_0^t\Ff\Fuz\dtau\le\Eb\Es$, and $\int_0^t\Ff^2\dtau\le\Eb^2$. Since $\Esh\le2\Es\lesssim1$ and $\int_0^t\Fsh^2\dtau\le4\Es^2$, $\int_0^t\Esh^2\Fsh^2\dtau\le16\Es^4\lesssim\Es^3$. Hence
\begin{equation*}
\int_0^t\nrm{v_t}_{H^3}^2\dtau\lesssim\Es(0)^2+\Eb^2+\Eb\Es+\Es^3 .
\end{equation*}
Adding \eqref{eq:6.14} proves \eqref{eq:6.10}.
\end{proof}

\subsection{Fourth-order tangential momentum-density energy}\label{sec:6tang}

Horizontal derivatives are written with multi-indices. For $\alpha=(\alpha_1,\alpha_2)\in\N_0^2$ let $\pa^\alpha\eqdef\px^{\alpha_1}\py^{\alpha_2}$ and $|\alpha|\eqdef\alpha_1+\alpha_2$; for $\beta\le\alpha$ (componentwise) let $\binom\alpha\beta\eqdef\binom{\alpha_1}{\beta_1}\binom{\alpha_2}{\beta_2}$, and let $e_1=(1,0)$, $e_2=(0,1)$. Let
\begin{equation*}
c_\alpha\eqdef\frac{|\alpha|!}{\alpha_1!\,\alpha_2!}
\end{equation*}
be the number of ordered tuples $(i_1,\dots,i_{|\alpha|})\in\{1,2\}^{|\alpha|}$ with $\pa_{i_1}\cdots\pa_{i_{|\alpha|}}=\pa^\alpha$. Then $\sum_{|\alpha|=4}c_\alpha|\pa^\alpha\phi|^2=|\nh^4\phi|^2$ pointwise, and hence
\begin{equation}\label{eq:T4-tensor}
\cT_4^2=\sum_{|\alpha|=4}c_\alpha\Big(\nrm{\sqrt\rho\,\pa^\alpha v}^2+g\nrm{\sqrt\sgm\,\pa^\alpha\varrho}^2\Big),\qquad D_4^2=\sum_{|\alpha|=4}c_\alpha\nrm{\nabla\pa^\alpha v}^2 .
\end{equation}
Since $\rs$ and $\sgm$ depend only on $z$, $\pa^\beta\rho=\pa^\beta\varrho$ for $\beta\ne0$ and $\pa^\alpha(\sgm\phi)=\sgm\pa^\alpha\phi$. In this subsection and the next, unlabeled norms are $L^2(\Omega)$ norms.

\begin{proposition}\label{prop:T4}
Assume \eqref{eq:small}. Then
\begin{equation}\label{eq:6.22}
\sup_{\tau\in[0,t]}\cT_4^2(\tau)+\int_0^tD_4^2\dtau\lesssim\Es(0)^2+\Es^3 .
\end{equation}
Moreover,
\begin{equation}\label{eq:6.23}
\int_0^t\nrm{\rho\dot v}_{H^3}^2\dtau\lesssim\Es(0)^2+\Eb^2+\Eb\Es+\Es^3 .
\end{equation}
\end{proposition}

\begin{proof}
\textbf{Step 1. Setup.} By \eqref{eq:Esharpdef} and \eqref{eq:small}, $\Esh\le2\Es\le1$, so that $1+\Esh\lesssim1$; in particular $\nrm{\rho}_{W^{2,\infty}}+\nrm{\rho}_{H^3}\lesssim1$ by \eqref{eq:3.12} and $\rs\in C^5$.

\emph{The convective term.} We claim
\begin{equation}\label{eq:T4-conv}
\nrm{\bU\cdot\nabla v}_{H^3}\lesssim\nrm{v}_{H^3}\nrm{v}_{H^4}.
\end{equation}
For the horizontal part, \Cref{lem:GN}(2) with $m=3$ and $H^2\hookrightarrow L^\infty$ give $\nrm{v\cdot\nh v}_{H^3}\lesssim\nrm{v}_{L^\infty}\nrm{\nh v}_{H^3}+\nrm{v}_{H^3}\nrm{\nh v}_{L^\infty}\lesssim\nrm{v}_{H^3}\nrm{v}_{H^4}$. For the vertical part we use, in this paragraph only, full multi-indices $\delta,\gamma\in\N_0^3$, $\pa^\delta=\px^{\delta_1}\py^{\delta_2}\pz^{\delta_3}$, and $e_z\eqdef(0,0,1)$. We first note that, for every $\delta\in\N_0^3$ with $|\delta|\le3$,
\begin{equation}\label{eq:T4-wL2Linf}
\nrm{\pa^\delta w}_{L^2_hL^\infty_z}\lesssim\nrm{v}_{H^{|\delta|+1}} .
\end{equation}
Indeed, if $\pa^\delta$ is purely horizontal this is \Cref{lem:w}; otherwise $\pa^\delta w=-\pa^{\delta-e_z}\divh v$, and the embedding $H^1(0,1)\subset L^\infty(0,1)$ gives $\nrm{\pa^{\delta-e_z}\divh v}_{L^2_hL^\infty_z}\lesssim\nrm{\pa^{\delta-e_z}\divh v}_{H^1}\lesssim\nrm{v}_{H^{|\delta|+1}}$. By the Leibniz rule, $\pa^\gamma(w\pz v)$, $|\gamma|\le3$, is a sum of products $\pa^{\gamma_1}w\,\pa^{\gamma_2}\pz v$ with $\gamma_1+\gamma_2=\gamma$. If $|\gamma_1|=0$ we use $\nrm{w}_{L^\infty}\lesssim\nrm{v}_{H^3}$, and if $|\gamma_1|=1$ we use $\nrm{\nabla w}_{L^\infty}\lesssim\nrm{v}_{H^4}$ (both from \eqref{eq:3.13}); this gives the bound $\nrm{v}_{H^3}\nrm{v}_{H^4}$. If $|\gamma_1|\in\{2,3\}$, then $|\gamma_2|\le1$, and \eqref{eq:T4-wL2Linf} and \Cref{lem:aniso}(e) with $s=\frac12$ give
\begin{equation*}
\nrm{\pa^{\gamma_1}w\,\pa^{\gamma_2}\pz v}\le\nrm{\pa^{\gamma_1}w}_{L^2_hL^\infty_z}\nrm{\pa^{\gamma_2}\pz v}_{L^\infty_hL^2_z}\lesssim\nrm{v}_{H^{|\gamma_1|+1}}\nrm{v}_{H^{|\gamma_2|+\frac52}}\lesssim\nrm{v}_{H^3}\nrm{v}_{H^4},
\end{equation*}
since $(|\gamma_1|+1,|\gamma_2|+\frac52)\in\{(3,\frac72),(4,\frac52)\}$. This proves \eqref{eq:T4-conv}.

Since $v$ vanishes on $\pa\Omega$, $\nrm{\nh^jv}\le\pi^{-1}\nrm{\pz\nh^jv}$ for every $j$, while every other derivative of $v$ of order at most $m$ is a derivative of $\pz v$ of order at most $m-1$; hence
\begin{equation*}
\nrm{v}_{H^m}\lesssim\nrm{\pz v}_{H^m}\le\Fuz,\qquad m\le4 .
\end{equation*}
By \eqref{eq:T4-conv} (with $\nrm{v}_{H^4}\le\Eu$ and $\nrm{v}_{H^3}\lesssim\Fuz$) and \Cref{lem:GN}(2), (3),
\begin{equation}\label{eq:6.24}
\begin{aligned}
\nrm{\dot v}_{H^3}&\le\nrm{v_t}_{H^3}+\nrm{\bU\cdot\nabla v}_{H^3}\lesssim\Fut+\Eu\Fuz\le\Fut+\Esh\Fuz,\\
\nrm{\rho\dot v}_{H^3}&\lesssim\nrm{\rho}_{L^\infty}\nrm{\dot v}_{H^3}+\nrm{\rho}_{H^3}\nrm{\dot v}_{L^\infty}\lesssim\nrm{\dot v}_{H^3}\lesssim\Fut+\Esh\Fuz .
\end{aligned}
\end{equation}

\emph{Recovery of $\Fuz$.} By \eqref{eq:4.15}, $\Fuz\lesssim\nrm{\nh^4\pz v}+\nrm{\pz G}_{H^2}+\nrm{\nh\varrho}_{H^2}$. Here $\nrm{\nh^4\pz v}\le D_4$, and every derivative in $\nrm{\nh\varrho}_{H^2}$ has at least one horizontal derivative and at most two vertical ones, so the horizontal Poincar\'e inequality (\Cref{lem:poincare}) gives
\begin{equation*}
\nrm{\nh\varrho}_{H^2}\lesssim\nrm{\nh^3\varrho}+\nrm{\nh^2\pz\varrho}+\nrm{\nh\pz^2\varrho}\lesssim\nrm{\nh^3\varrho}+\nrm{\nabla\nh\pz\varrho}\le\Fv .
\end{equation*}
By \Cref{rem:fH3}, $\nrm{\pz G}_{H^2}\le\nrm{\rho\dot v}_{H^3}+|f|\nrm{\rho v}_{H^3}$, and by \eqref{eq:fH3} the linear Coriolis contribution satisfies, for every $\epsilon>0$ (after renaming $\epsilon$), $|f|\nrm{\rho v}_{H^3}\le\epsilon\Fuz+C_\epsilon(1+\Ev)^2\nrm{v}_{H^2}\le\epsilon\Fuz+C'_\epsilon\nrm{v}_{H^2}$, since $(1+\Ev)^2\lesssim1$ under \eqref{eq:small}. Choosing $\epsilon$ small and absorbing $\epsilon\Fuz$,
\begin{equation}\label{eq:6.25}
\Fuz\lesssim D_4+\nrm{\rho\dot v}_{H^3}+\Fv+\nrm{v}_{H^2}.
\end{equation}
Inserting \eqref{eq:6.24} into \eqref{eq:6.25} gives $\Fuz\le C(D_4+\Fv+\Fut+\nrm{v}_{H^2})+C\Esh\Fuz$, and $C\Esh\le\frac12$ by \eqref{eq:small} and \Cref{rem:c1}. Together with $\Fu\lesssim D_4+\Fuz$ (by \eqref{eq:3.3}) and \eqref{eq:6.24} this yields
\begin{equation}\label{eq:6.26}
\begin{aligned}
\Fuz+\Fu&\lesssim D_4+\Fv+\Fut+\nrm{v}_{H^2},\\
\nrm{\rho\dot v}_{H^3}&\lesssim\Fut+\Esh\big(D_4+\Fv+\nrm{v}_{H^2}\big).
\end{aligned}
\end{equation}
The term $\nrm{v}_{H^2}$ is the contribution of the linear Coriolis force; it is of lower order and is controlled by \eqref{eq:4.3} or \eqref{eq:4.4}. Squaring the second line of \eqref{eq:6.26},
\begin{equation}\label{eq:6.27}
\nrm{\rho\dot v}_{H^3}^2\lesssim\Fut^2+\Esh^2\big(D_4^2+\Fv^2+\nrm{v}_{H^2}^2\big).
\end{equation}

\emph{Horizontal derivatives and the Coriolis force.} Derivatives of $v$ that contain at least one horizontal derivative have zero horizontal mean, and $v$ together with all its horizontal derivatives vanishes on $\pa\Omega$. The horizontal and vertical Poincar\'e inequalities (\Cref{lem:poincare}) therefore give
\begin{equation}\label{eq:T4-P}
\nrm{\pa^\gamma\pz^mv}\lesssim D_4\qquad\text{for }\gamma\in\N_0^2,\ |\gamma|\ge1,\ m\in\{0,1\},\ |\gamma|+m\le5,
\end{equation}
and also $\nrm{\nh^4v}\le\pi^{-1}\nrm{\pz\nh^4v}\le D_4$. We use this to control the Coriolis part of $G$:
\begin{equation}\label{eq:T4-GH3}
\nrm{\nh^3(\rho v)}+\nrm{\Laph\pz(\rho v)}+\nrm{\nh^2\pz(\rho v)}\lesssim D_4+\Ev\nrm{v}_{H^2}\lesssim D_4+\Esh\Fuz .
\end{equation}
Indeed, $\nrm{\Laph\pz(\rho v)}=\nrm{\nh^2\pz(\rho v)}$ by Parseval in $x_h$. For the part with $\rs$, $\nh^3(\rs v)=\rs\nh^3v$ and $\nh^2\pz(\rs v)=\rs'\nh^2v+\rs\nh^2\pz v$ are bounded by $CD_4$ by \eqref{eq:T4-P}. For the part with $\varrho$ we expand by the Leibniz rule and use $\nrm{\varrho}_{W^{2,\infty}}\lesssim\Ev$, $\nrm{v}_{L^\infty}\lesssim\nrm{v}_{H^2}$ (\eqref{eq:3.12}, \eqref{eq:3.13}) and \eqref{eq:T4-P}: in the products in which at most two derivatives fall on $\varrho$, the factor of $v$ carries at least one derivative and at most one $\pz$, so it is either $\pz v$ or one of the derivatives in \eqref{eq:T4-P}; these products are bounded by $\nrm{\varrho}_{W^{2,\infty}}(D_4+\nrm{\pz v})\lesssim\Ev(D_4+\nrm{v}_{H^2})$. The products with three derivatives on $\varrho$ are bounded by $\nrm{\varrho}_{H^3}\nrm{v}_{L^\infty}\lesssim\Ev\nrm{v}_{H^2}$. Finally $\nrm{v}_{H^2}\lesssim\nrm{\pz v}_{H^2}\le\Fuz$. As a consequence of \eqref{eq:T4-GH3} and \eqref{eq:6.26}, since $G=\rho\dot v+f\rho v^\perp$ and $|\pa^\gamma(\rho v^\perp)|=|\pa^\gamma(\rho v)|$,
\begin{equation}\label{eq:T4-G3}
\nrm{\nh^3G}+\nrm{\nh^2\pz G}\lesssim\nrm{\rho\dot v}_{H^3}+D_4+\Ev\nrm{v}_{H^2}\lesssim\Fut+D_4+\Esh\big(\Fv+\nrm{v}_{H^2}\big).
\end{equation}
Thus the linear Coriolis force, measured with at least two horizontal derivatives and at most one vertical derivative, costs only $D_4$ and a quadratic term; it is never estimated by an $O(1)$ multiple of $\Fuz$.

\textbf{Step 2. Hydrostatic cancellation.} Fix $\alpha\in\N_0^2$ with $|\alpha|=4$ and define
\begin{equation*}
\cC_m^\alpha\eqdef\pa^\alpha(\rho\dot v)-\rho(\pt+\bU\cdot\nabla)\pa^\alpha v,\qquad \cC_f^\alpha\eqdef\pa^\alpha(\rho v^\perp)-\rho\,\pa^\alpha v^\perp,\qquad \cC_\varrho^\alpha\eqdef\pa^\alpha(\bU\cdot\nabla\varrho)-\bU\cdot\nabla\pa^\alpha\varrho .
\end{equation*}
Apply $\pa^\alpha$ to \eqref{PE}$_2$ and test with $\pa^\alpha v$, which vanishes on $\pa\Omega$. Since $\rho_t+\bU\cdot\nabla\rho=0$, $\divv\bU=0$ and $w|_{\pa\Omega}=0$, we have $\int_\Omega\rho(\pt+\bU\cdot\nabla)\pa^\alpha v\cdot\pa^\alpha v\dx=\frac12\ddt\int_\Omega\rho|\pa^\alpha v|^2\dx$. The leading Coriolis term vanishes pointwise, $f\rho\,\pa^\alpha v^\perp\cdot\pa^\alpha v=0$. Hence
\begin{equation}\label{eq:6.28}
\frac12\ddt\int_\Omega\rho|\pa^\alpha v|^2\dx+\mu\nrm{\nabla\pa^\alpha v}^2+\int_\Omega\nh\pa^\alpha q\cdot\pa^\alpha v\dx=-\int_\Omega\cC_m^\alpha\cdot\pa^\alpha v\dx-f\int_\Omega\cC_f^\alpha\cdot\pa^\alpha v\dx .
\end{equation}
By $\divh v=-\pz w$ and $\pz q=-g\varrho$, the pressure integral is
\begin{equation*}
\begin{aligned}
\int_\Omega\nh\pa^\alpha q\cdot\pa^\alpha v\dx&=-\int_\Omega\pa^\alpha q\,\divh\pa^\alpha v\dx=\int_\Omega\pa^\alpha q\,\pz\pa^\alpha w\dx\\
&=-\int_\Omega\pz\pa^\alpha q\,\pa^\alpha w\dx=g\int_\Omega\pa^\alpha\varrho\,\pa^\alpha w\dx,
\end{aligned}
\end{equation*}
where the boundary term $\int_{\Tl}[\pa^\alpha q\,\pa^\alpha w]_{z=0}^{z=1}\dxh$ vanishes because $w=0$ on $\pa\Omega$ (\Cref{lem:propagation}), hence $\pa^\alpha w=0$ on $\pa\Omega$. Next, \eqref{PE}$_1$ reads $\sgm(\pt+\bU\cdot\nabla)\varrho=w$. Applying $\pa^\alpha$ (which commutes with $\sgm$) and testing with $g\pa^\alpha\varrho$, we obtain, since $\int_\Omega\sgm\bU\cdot\nabla\frac{|\pa^\alpha\varrho|^2}2\dx=-\frac12\int_\Omega\sgm'w|\pa^\alpha\varrho|^2\dx$,
\begin{equation}\label{eq:6.29}
\frac g2\ddt\nrm{\sqrt\sgm\,\pa^\alpha\varrho}^2+g\int_\Omega\sgm\,\cC_\varrho^\alpha\,\pa^\alpha\varrho\dx=g\int_\Omega\pa^\alpha w\,\pa^\alpha\varrho\dx+\frac g2\int_\Omega\sgm'w|\pa^\alpha\varrho|^2\dx .
\end{equation}
Adding \eqref{eq:6.28} and \eqref{eq:6.29}, the two integrals $g\int_\Omega\pa^\alpha\varrho\,\pa^\alpha w\dx$ cancel. Multiplying by $c_\alpha$ and summing over $|\alpha|=4$, \eqref{eq:T4-tensor} gives
\begin{equation}\label{eq:6.30}
\begin{aligned}
\frac12\ddt\cT_4^2+\mu D_4^2={}&-\sum_{|\alpha|=4}c_\alpha\int_\Omega\cC_m^\alpha\cdot\pa^\alpha v\dx-f\sum_{|\alpha|=4}c_\alpha\int_\Omega\cC_f^\alpha\cdot\pa^\alpha v\dx\\
&+\frac g2\int_\Omega\sgm'w|\nh^4\varrho|^2\dx-g\sum_{|\alpha|=4}c_\alpha\int_\Omega\sgm\,\cC_\varrho^\alpha\,\pa^\alpha\varrho\dx .
\end{aligned}
\end{equation}
These identities are justified as in \Cref{rem:justify}: in \eqref{eq:6.28} the pressure enters only through the integral $g\int_\Omega\pa^\alpha\varrho\,\pa^\alpha w\dx$, whose factors belong to $L^2$, and in \eqref{eq:6.29} the commutator $\cC_\varrho^\alpha$ belongs to $L^2$. Since $c_\alpha\le c_{(2,2)}=6$ (and $\sum_{|\alpha|=4}c_\alpha=16$), it suffices to estimate each commutator for fixed $\alpha$.

\textbf{Step 3. Momentum and Coriolis commutators.} The commutators have the exact expansions
\begin{equation*}
\begin{aligned}
\cC_m^\alpha&=\rho\sum_{0<\beta\le\alpha}\binom\alpha\beta\big(\pa^\beta v\cdot\nh\pa^{\alpha-\beta}v+\pa^\beta w\,\pz\pa^{\alpha-\beta}v\big)+\sum_{0<\beta\le\alpha}\binom\alpha\beta\pa^\beta\varrho\,\pa^{\alpha-\beta}\dot v,\\
\cC_f^\alpha&=\sum_{0<\beta\le\alpha}\binom\alpha\beta\pa^\beta\varrho\,\pa^{\alpha-\beta}v^\perp .
\end{aligned}
\end{equation*}
We collect the three-dimensional bounds of horizontal derivatives. By \eqref{eq:T4-P}, $\nrm{\nh^kv}_{H^1}\lesssim\nrm{\nh^kv}+\nrm{\nh^{k+1}v}+\nrm{\pz\nh^kv}\lesssim D_4$ for $1\le k\le4$. By \Cref{lem:w}, $\nrm{\nh^jw}\le\nrm{\nh^j\divh v}\le\sqrt2\nrm{\nh^{j+1}v}$, and $\pz\nh^jw=-\nh^j\divh v$; hence $\nrm{\nh^mw}_{H^1}\lesssim D_4$ for $1\le m\le3$ and $\nrm{\nh^4w}\le\sqrt2\nrm{\nh^5v}\le\sqrt2D_4$. By \Cref{lem:aniso}(d) with $p=2$, $\nrm{\nh^4v}_{L^2_hL^\infty_z}\le\nrm{\pz\nh^4v}\le D_4$. For $j=1,2$, \Cref{lem:Linf}(ii) with $s=\frac12$ gives $\nrm{\nh^jv}_{L^\infty}\lesssim\nrm{\nh^{j+1}\pz v}^{1/2}\nrm{\nh^{j+2}\pz v}^{1/2}\lesssim D_4$, and \Cref{lem:w} with $s=1$ gives $\nrm{\nh^jw}_{L^\infty}\lesssim\nrm{\nh^j\divh v}_{L^2_zH^2_h}\lesssim\nrm{\nh^{j+3}v}\lesssim D_4$ (zero horizontal mean). With $H^1(\Omega)\hookrightarrow L^4(\Omega)$ (\Cref{lem:GN}(3)) we obtain
\begin{equation}\label{eq:6.31}
\begin{aligned}
&\nrm{\nh^kv}_{L^4}+\nrm{\nh^kv}_{H^1}+\nrm{\nh^mw}_{L^4}+\nrm{\nh^mw}_{H^1}\lesssim D_4,\qquad 1\le k\le4,\ 1\le m\le3,\\
&\nrm{\nh^4w}+\nrm{\nh^4v}_{L^2_hL^\infty_z}+\nrm{\nh^jv}_{L^\infty}+\nrm{\nh^jw}_{L^\infty}\lesssim D_4,\qquad j=1,2.
\end{aligned}
\end{equation}
Since $H^1(\T^2)\not\subset L^\infty(\T^2)$, $\nrm{\nh^3v}_{L^\infty}$ is not controlled by $D_4$ (only $\nrm{\nh^3v}_{L^\infty_zL^p_h}\le C_pD_4$, $p<\infty$, holds by \Cref{lem:Linf}(ii)); third horizontal derivatives are therefore always placed in $L^4$ below. For the transport part of $\cC_m^\alpha$, let $j=|\beta|$. Using $\nrm{\rho}_{L^\infty}\lesssim1$, $\nrm{\nabla\bU}_{L^\infty}+\nrm{\nabla^2v}_{L^\infty}\lesssim\Eu$ (\eqref{eq:3.13} and $H^2\hookrightarrow L^\infty$), $\nrm{\pz\nh v}_{L^4}\lesssim\nrm{v}_{H^3}\le\Eu$, $\nrm{\nh^4v}\le\Eu$ and \eqref{eq:6.31},
\begin{equation}\label{eq:6.32}
\begin{aligned}
&\left|\int_\Omega\rho\,\pa^\beta v\cdot\nh\pa^{\alpha-\beta}v\cdot\pa^\alpha v\dx\right|+\left|\int_\Omega\rho\,\pa^\beta w\,\pz\pa^{\alpha-\beta}v\cdot\pa^\alpha v\dx\right|\\
&\quad\lesssim\begin{cases}
\nrm{\nh v}_{L^\infty}\nrm{\nh^4v}^2+\nrm{\nh w}_{L^\infty}\nrm{\pz\nh^3v}\nrm{\nh^4v},&j=1,\\
\nrm{\nh^2v}_{L^\infty}\nrm{\nh^3v}\nrm{\nh^4v}+\nrm{\nh^2w}_{L^\infty}\nrm{\pz\nh^2v}\nrm{\nh^4v},&j=2,\\
\nrm{\nh^3v}_{L^4}\nrm{\nh^2v}_{L^4}\nrm{\nh^4v}+\nrm{\nh^3w}_{L^4}\nrm{\pz\nh v}_{L^4}\nrm{\nh^4v},&j=3,\\
\nrm{\nh^4v}\nrm{\nh v}_{L^\infty}\nrm{\nh^4v}+\nrm{\nh^4w}\nrm{\pz v}_{L^\infty}\nrm{\nh^4v},&j=4,
\end{cases}\\
&\quad\lesssim\Eu D_4^2 .
\end{aligned}
\end{equation}
For the density factors, \eqref{eq:3.12} and $H^1\hookrightarrow L^4$ give
\begin{equation}\label{eq:6.33}
\nrm{\varrho}_{W^{2,\infty}}+\nrm{\nh^3\varrho}_{L^4}+\nrm{\nh^4\varrho}+\nrm{\pz\nh^3\varrho}\lesssim\Ev .
\end{equation}
With $j=|\beta|$, \eqref{eq:6.31}, \eqref{eq:6.33} and $H^2\hookrightarrow L^\infty\cap W^{1,4}$ give
\begin{equation}\label{eq:6.34}
\left|\int_\Omega\pa^\beta\varrho\,\pa^{\alpha-\beta}\dot v\cdot\pa^\alpha v\dx\right|\lesssim\begin{cases}
\nrm{\nh^j\varrho}_{L^\infty}\nrm{\nh^{4-j}\dot v}\nrm{\nh^4v},&j=1,2,\\
\nrm{\nh^3\varrho}_{L^4}\nrm{\nh\dot v}_{L^4}\nrm{\nh^4v},&j=3,\\
\nrm{\nh^4\varrho}\nrm{\dot v}_{L^\infty}\nrm{\nh^4v},&j=4,
\end{cases}\quad\lesssim\Ev\nrm{\dot v}_{H^3}D_4 .
\end{equation}
By \eqref{eq:6.32} and \eqref{eq:6.34},
\begin{equation}\label{eq:6.35}
\left|\int_\Omega\cC_m^\alpha\cdot\pa^\alpha v\dx\right|\lesssim\Eu D_4^2+\Ev\nrm{\dot v}_{H^3}D_4 .
\end{equation}
Since $\nrm{\dot v}_{H^3}\lesssim\Fut+\Esh\Fuz$ by \eqref{eq:6.24}, $\Esh\lesssim1$, and $D_4+\Fuz\le2\Fu$, we obtain
\begin{equation}\label{eq:6.36}
\left|\sum_{|\alpha|=4}c_\alpha\int_\Omega\cC_m^\alpha\cdot\pa^\alpha v\dx\right|\lesssim\Esh\Fsh D_4 .
\end{equation}
The Coriolis commutator $\cC_f^\alpha$ is cubic. With $j=|\beta|$, the same placements as in \eqref{eq:6.34}, \eqref{eq:6.31} and $\nrm{v}_{L^\infty}\lesssim\nrm{v}_{H^2}$ give
\begin{equation*}
\left|\int_\Omega\pa^\beta\varrho\,\pa^{\alpha-\beta}v^\perp\cdot\pa^\alpha v\dx\right|\lesssim\begin{cases}
\nrm{\nh^j\varrho}_{L^\infty}\nrm{\nh^{4-j}v}\nrm{\nh^4v},&j=1,2,\\
\nrm{\nh^3\varrho}_{L^4}\nrm{\nh v}_{L^4}\nrm{\nh^4v},&j=3,\\
\nrm{\nh^4\varrho}\nrm{v}_{L^\infty}\nrm{\nh^4v},&j=4,
\end{cases}\quad\lesssim\Ev\big(D_4+\nrm{v}_{H^2}\big)D_4 .
\end{equation*}
Only the summand $j=4$ sees the horizontal mean $\bar v$, through $\nrm{v}_{L^\infty}$. Since $D_4+\nrm{v}_{H^2}\le2\Fu$,
\begin{equation}\label{eq:T4-Cf}
\left|f\sum_{|\alpha|=4}c_\alpha\int_\Omega\cC_f^\alpha\cdot\pa^\alpha v\dx\right|\lesssim\Ev\big(D_4+\nrm{v}_{H^2}\big)D_4\lesssim\Esh\Fsh D_4 .
\end{equation}

\textbf{Step 4. Density commutator.} The density commutator has the exact expansion
\begin{equation}\label{eq:6.37}
\cC_\varrho^\alpha=\sum_{0<\beta\le\alpha}\binom\alpha\beta\big(\pa^\beta v\cdot\nh\pa^{\alpha-\beta}\varrho+\pa^\beta w\,\pz\pa^{\alpha-\beta}\varrho\big).
\end{equation}
The terms containing $v$ and those containing $w$ are called the horizontal and the vertical summands; we write $j=|\beta|$. Besides \eqref{eq:6.31} and \eqref{eq:6.33} we use the dissipative bounds
\begin{equation}\label{eq:T4-Fv}
\nrm{\nh\varrho}_{H^1}+\nrm{\nh^2\varrho}_{H^1}+\nrm{\nh\pz\varrho}_{H^1}+\nrm{\nh^3\varrho}+\nrm{\nh^2\pz\varrho}\lesssim\Fv,
\end{equation}
which follow from the definition of $\Fv$ and the horizontal Poincar\'e inequality: every derivative occurring on the left has at least one horizontal derivative and is one of $\nh^k\varrho$ ($k\le3$), $\nh^k\pz\varrho$ ($k\le2$) or $\nh\pz^2\varrho$. By $H^1\hookrightarrow L^4$, the $H^1$ norms in \eqref{eq:T4-Fv} also bound the corresponding $L^4$ norms.

First, the horizontal summand with $j=1$ contains four horizontal derivatives of $\varrho$, and $\nrm{\nh v}_{L^\infty}\le\nrm{\nabla\bU}_{L^\infty}\le\Fn$; hence
\begin{equation}\label{eq:6.38}
\left|\int_\Omega\sgm\,\pa^\beta v\cdot\nh\pa^{\alpha-\beta}\varrho\,\pa^\alpha\varrho\dx\right|\lesssim\nrm{\nh v}_{L^\infty}\nrm{\nh^4\varrho}^2\lesssim\Esh^2\Fn,\qquad|\beta|=1 .
\end{equation}

Second, consider the vertical summand with $j=4$, that is, $\beta=\alpha$. Here
\begin{equation}\label{eq:6.39}
\nrm{\pa^\alpha w}\le\nrm{\nh^4\divh v}\le\sqrt2\nrm{\nh^5v}\le\sqrt2D_4,
\end{equation}
so the direct bound $\nrm{\pa^\alpha w}\nrm{\pz\varrho}_{L^\infty}\nrm{\pa^\alpha\varrho}\lesssim\Ev^2D_4$ contains only one dissipative factor, and we use the representation \eqref{eq:wrep}. Fix $i\in\{1,2\}$ with $\alpha_i\ge1$ and put $\alpha'\eqdef\alpha-e_i$. Let $S_i\eqdef\pa_i\ABmh=\ABmh\pa_i$ (\Cref{lem:ABmean}). For $\psi\in L^2$ and $h=\ABmh\psi$, the function $\nh h$ vanishes on $\pa\Omega$, so the vertical Poincar\'e inequality and \Cref{lem:ABid} give $\nrm{S_i\psi}=\nrm{\pa_ih}\le\pi^{-1}\nrm{\nh\pz h}\le\pi^{-1}\nrm{\psi}$. Thus $S_i$ is bounded on $L^2$, and $S_i^*=-S_i$, because on the horizontal mode $k$ it acts as $\mathrm ik_iA_k^{-1/2}$ with $A_k^{-1/2}$ self-adjoint. Applying $\pa^\alpha$ to \eqref{eq:wrep} and using that horizontal derivatives commute with functions of $\AB$, we obtain, mode by mode,
\begin{equation*}
\mu\,\pa^\alpha w=S_i\Psi,\qquad \Psi\eqdef g\ABmh\Laph\pa^{\alpha'}\varrho-\ABmh\pz\divh\pa^{\alpha'}G,
\end{equation*}
where the second term is defined by duality as in \eqref{eq:ABdual}. On the mode $k$ we have $|k|^2|k^{\alpha'}|\le|k|^5$; hence, by \eqref{eq:ABdual}, \eqref{eq:6.24} and \eqref{eq:T4-GH3},
\begin{equation*}
\begin{aligned}
\nrm{\Psi}&\le g\nrm{\ABmh|\nh|^5\varrho}+\nrm{\pa^{\alpha'}G}\le g\Fv+\nrm{\rho\dot v}_{H^3}+|f|\nrm{\nh^3(\rho v)}\\
&\lesssim\Fv+\Fut+\Esh\Fuz+D_4+\Ev\nrm{v}_{H^2}\lesssim\Fsh .
\end{aligned}
\end{equation*}
The Coriolis part of $G$ enters only through $\nrm{\nh^3(\rho v)}$, which is controlled by \eqref{eq:T4-GH3}. With $\chi\eqdef\sgm\pz\varrho\,\pa^\alpha\varrho\in L^2$,
\begin{equation}\label{eq:6.40}
\int_\Omega\sgm\,\pa^\alpha w\,\pz\varrho\,\pa^\alpha\varrho\dx=\ip{\pa^\alpha w,\chi}=\frac1\mu\ip{S_i\Psi,\chi}=-\frac1\mu\ip{\Psi,S_i\chi}.
\end{equation}
It remains to bound $\nrm{S_i\chi}=\nrm{S_i^*\chi}=\sup_{\nrm\psi=1}|\ip{\chi,\pa_ih}|$, $h=\ABmh\psi$. Put $a\eqdef\sgm\pz\varrho$; then $a\in H^3\subset W^{1,\infty}$ and $a=0$ on $\pa\Omega$, since $\pz\varrho|_{\pa\Omega}=0$. Since $\pa_i\sgm=0$,
\begin{equation*}
\ip{\chi,\pa_ih}=\int_\Omega\pa^\alpha\varrho\;a\,\pa_ih\dx=\int_\Omega\pa^\alpha\varrho\,\pa_i(ah)\dx-\int_\Omega\pa^\alpha\varrho\,\sgm\,\pa_i\pz\varrho\,h\dx\eqdef B_1+B_2 .
\end{equation*}
The bounds $\nrm{\ABmh(\phi\chi)}\lesssim\nrm{\phi}\nrm{\chi}$ and $\nrm{\ABmh\psi}_{L^\infty}\lesssim\nrm{\psi}$ fail in three dimensions (\Cref{rem:Linfty}); we use \Cref{lem:product} instead. For $B_2$, H\"older's inequality in $x_h$ (exponents $2,4,4$) and then in $z$, \Cref{lem:product}(i) with $p=4$, and \Cref{lem:aniso}(f) with $q=4$ for the zero-mean function $\pa_i\pz\varrho$ give
\begin{equation*}
|B_2|\le\nrm{\pa^\alpha\varrho}\,\nrm{\sgm\pa_i\pz\varrho}_{L^2_zL^4_h}\nrm{h}_{L^\infty_zL^4_h}\lesssim\Ev\nrm{\pa_i\pz\varrho}_{L^2_zH^{1/2}_h}\nrm{\psi}\lesssim\Ev\nrm{\nh^2\pz\varrho}\nrm{\psi}\le\Ev\Fv\nrm{\psi}.
\end{equation*}
For $B_1$, we first bound $\nabla\pz(ah)=\nabla\pz a\,h+\pz a\,\nabla h+\nabla a\,\pz h+a\,\nabla\pz h$. The coefficient $\nabla\pz a$ contains three derivatives of $\varrho$ and is not bounded; we put it in $L^2_zL^4_h$, using \Cref{lem:aniso}(f) with $q=4$ and $\sgm\in C^4([0,1])$, and $h$ in $L^\infty_zL^4_h$ by \Cref{lem:product}(i). The other coefficients are bounded by $\nrm{a}_{W^{1,\infty}}\lesssim\nrm{\varrho}_{W^{2,\infty}}\lesssim\Ev$, and $\nrm{\nabla h}+\nrm{\nabla\pz h}\lesssim\nrm{\psi}$ by \Cref{lem:product}(i). Hence
\begin{equation*}
\nrm{\nabla\pz(ah)}\lesssim\nrm{\nabla\pz a}_{L^2_zL^4_h}\nrm{h}_{L^\infty_zL^4_h}+\Ev\nrm{\psi}\lesssim\big(\nrm{a}_{H^{5/2}}+\Ev\big)\nrm{\psi}\lesssim\Ev\nrm{\psi},
\end{equation*}
because $\nrm{a}_{H^{5/2}}\lesssim\nrm{\varrho}_{H^{7/2}}$. Since moreover $ah=0$ and $\pz(ah)=\pz a\,h+a\,\pz h=0$ on $\pa\Omega$, we have $ah\in V=\cD(\ABh)$, with $\nrm{\ABh(ah)}=\nrm{\nabla\pz(ah)}$ by \Cref{lem:ABid}. For $\phi\in L^2$ and $\varphi\in V$, $\pa_i\varphi=S_i(\ABh\varphi)$, so $\ip{\phi,\pa_i\varphi}=-\ip{S_i\phi,\ABh\varphi}$. With $\phi=\pa^\alpha\varrho$, and since $|k_ik^\alpha|\le|k|^5$ on the mode $k$,
\begin{equation*}
|B_1|\le\nrm{S_i\pa^\alpha\varrho}\nrm{\nabla\pz(ah)}\le\nrm{\ABmh|\nh|^5\varrho}\nrm{\nabla\pz(ah)}\lesssim\Fv\Ev\nrm{\psi}.
\end{equation*}
Consequently
\begin{equation}\label{eq:T4-chi}
\nrm{S_i\chi}\lesssim\Ev\Fv,
\end{equation}
and \eqref{eq:6.40} gives
\begin{equation}\label{eq:6.41}
\left|\int_\Omega\sgm\,\pa^\alpha w\,\pz\varrho\,\pa^\alpha\varrho\dx\right|\le\frac1\mu\nrm{\Psi}\nrm{S_i\chi}\lesssim\Ev\Fv\Fsh .
\end{equation}

Third, the vertical summand with $j=1$ is $\pa_iw\,\pz\pa^{\alpha-e_i}\varrho$. Since $\nrm{\nh w}_{L^\infty}\le\nrm{\nabla\bU}_{L^\infty}\le\Fn$ (see also \eqref{eq:gradhw}), no representation of $w$ is needed:
\begin{equation}\label{eq:6.43}
\left|\int_\Omega\sgm\,\pa^\beta w\,\pz\pa^{\alpha-\beta}\varrho\,\pa^\alpha\varrho\dx\right|\lesssim\nrm{\nh w}_{L^\infty}\nrm{\pz\nh^3\varrho}\nrm{\nh^4\varrho}\lesssim\Esh^2\Fn,\qquad|\beta|=1 .
\end{equation}

Finally, the remaining summands, horizontal with $2\le j\le4$ and vertical with $2\le j\le3$, are estimated with $\pa^\alpha\varrho$ in $L^2$ (bounded by $\Ev$), the velocity factor bounded by $D_4$ through \eqref{eq:6.31}, and the density factor bounded by $\Fv$ through \eqref{eq:T4-Fv}:
\begin{equation}\label{eq:6.44}
\begin{aligned}
&\sum_{2\le|\beta|\le4}\left|\int_\Omega\sgm\,\pa^\beta v\cdot\nh\pa^{\alpha-\beta}\varrho\,\pa^\alpha\varrho\dx\right|+\sum_{2\le|\beta|\le3}\left|\int_\Omega\sgm\,\pa^\beta w\,\pz\pa^{\alpha-\beta}\varrho\,\pa^\alpha\varrho\dx\right|\\
&\quad\lesssim\nrm{\nh^4\varrho}\Big(\nrm{\nh^2v}_{L^\infty}\nrm{\nh^3\varrho}+\nrm{\nh^3v}_{L^4}\nrm{\nh^2\varrho}_{L^4}+\nrm{\nh^4v}_{L^4}\nrm{\nh\varrho}_{L^4}\\
&\qquad+\nrm{\nh^2w}_{L^\infty}\nrm{\pz\nh^2\varrho}+\nrm{\nh^3w}_{L^4}\nrm{\pz\nh\varrho}_{L^4}\Big)\lesssim\Ev\Fv D_4 .
\end{aligned}
\end{equation}
Here we used $\nrm{\nh^2w}_{L^\infty}\lesssim D_4$, which requires $\nh^5v$ in $L^2$ and is the reason why this summand is not placed as $\nrm{\nh^2w}_{L^4}\nrm{\pz\nh^2\varrho}_{L^4}$: the latter would need $\pz\nh^3\varrho$ and $\pz^2\nh^2\varrho$, which are not part of $\Fv$. Collecting \eqref{eq:6.38}, \eqref{eq:6.41}, \eqref{eq:6.43} and \eqref{eq:6.44}, and using $D_4\le\Fu\le\Fsh$,
\begin{equation}\label{eq:6.45}
\left|g\sum_{|\alpha|=4}c_\alpha\int_\Omega\sgm\,\cC_\varrho^\alpha\,\pa^\alpha\varrho\dx\right|\lesssim\Esh^2\Fn+\Ev\Fv\Fsh .
\end{equation}

\textbf{Step 5. Temporal integration.} Since $w=0$ at $z=0$, $|w(x_h,z)|\le\int_0^1|\pz w(x_h,z')|\dd z'$, so $\nrm{w}_{L^\infty}\le\nrm{\pz w}_{L^\infty}\le\nrm{\nabla\bU}_{L^\infty}\le\Fn$, and the remaining term in \eqref{eq:6.30} satisfies
\begin{equation}\label{eq:6.46}
\left|\frac g2\int_\Omega\sgm'w|\nh^4\varrho|^2\dx\right|\lesssim\nrm{w}_{L^\infty}\Ev^2\lesssim\Esh^2\Fn .
\end{equation}
Inserting \eqref{eq:6.36}, \eqref{eq:T4-Cf}, \eqref{eq:6.45} and \eqref{eq:6.46} into \eqref{eq:6.30}, and using $D_4,\Fv\le\Fsh$, we obtain
\begin{equation}\label{eq:T4-diff}
\frac12\ddt\cT_4^2+\mu D_4^2\le C\Esh\Fsh^2+C\Esh^2\Fn .
\end{equation}
For $\tau\in[0,t]$ we integrate \eqref{eq:T4-diff} over $(0,\tau)$. Since $\Esh(s)\le2\Es(t)$ for $s\le t$, $\int_0^t\Fsh^2\dtau\lesssim\Es^2(t)$ by \eqref{eq:Esharpdef}, $\int_0^t\Fn\dtau\lesssim\Es(t)$ by \Cref{prop:Fnat}, and $\cT_4(0)\lesssim\Eu(0)+\Ev(0)\le\Esh(0)\lesssim\Es(0)$, we obtain
\begin{equation*}
\cT_4^2(\tau)+2\mu\int_0^\tau D_4^2\dtau\le\cT_4^2(0)+C\Es^3(t)+C\Es^2(t)\int_0^t\Fn\dtau\lesssim\Es(0)^2+\Es^3(t),
\end{equation*}
which yields \eqref{eq:6.22}. No term of the form $\Eb\Es$ arises here: every term on the right of \eqref{eq:T4-diff} contains two dissipative factors or the integrable factor $\Fn$.

For \eqref{eq:6.23} we integrate \eqref{eq:6.27} in time. Since $D_4+\nrm{v}_{H^2}\le2\Fu$ and $\Fv\le\Fsh$, the last term contributes at most $C\Es^2\int_0^t\Fsh^2\dtau\lesssim\Es^4\le\Es^3$, and $\int_0^t\Fut^2\dtau\le\cEut^2(t)$ is bounded by \eqref{eq:6.10}. This yields \eqref{eq:6.23}.
\end{proof}

\subsection{Estimate of \texorpdfstring{$\cEu$}{E u}}

\begin{proposition}\label{prop:u}
Assume \eqref{eq:small}. Then, for $t\ge0$,
\begin{equation}\label{eq:6.47}
\cEu^2(t)=\sup_{\tau\in(0,t)}\Eu^2+\int_0^t\Fu^2\dtau\lesssim\cEv^2(t)+\Es(0)^2+\Eb^2+\Eb\Es+\Es^3 .
\end{equation}
\end{proposition}

\begin{proof}
\textbf{Step 1. Estimate of $\nrm{v}_{H^4}$.} Squaring \eqref{eq:4.13} gives $\Eu^2\lesssim\Eut^2+\Ev^2+\Ef^2+\Esh^4$. In three dimensions with rotation, \eqref{eq:4.13} contains the low-order energy $\Ef$, which comes from the term $|f|\nrm{v}_{H^2}$ in the Stokes estimate; by \eqref{eq:3.5}, $\sup_{(0,t)}\Ef^2\le\Eb^2(t)$. Taking the supremum in time,
\begin{equation}\label{eq:6.48}
\sup_{\tau\in(0,t)}\Eu^2\lesssim\sup_{\tau\in(0,t)}\Eut^2+\sup_{\tau\in(0,t)}\Ev^2+\Eb^2+\Es^4 .
\end{equation}

\textbf{Step 2. Estimate of $\nrm{v}_{H^5}$.} By the first line of \eqref{eq:6.26}, $\Fu\lesssim D_4+\Fv+\Fut+\nrm{v}_{H^2}$. For the lower-order term, \eqref{eq:4.11} gives $\nrm{v}_{H^2}^2\lesssim\Ff^{4/3}\Fuz^{2/3}+\Ff^2$, so that, by H\"older's inequality with exponents $\frac32$ and $3$, $\Fuz\le\Fu$ and $a^{1/3}b^{2/3}\le a+b$,
\begin{equation*}
\int_0^t\nrm{v}_{H^2}^2\dtau\lesssim\Big(\int_0^t\Ff^2\dtau\Big)^{\frac23}\Big(\int_0^t\Fuz^2\dtau\Big)^{\frac13}+\int_0^t\Ff^2\dtau\le\Eb^{\frac43}\Es^{\frac23}+\Eb^2\le2\Eb^2+\Eb\Es .
\end{equation*}
Hence
\begin{equation*}
\int_0^t\Fu^2\dtau\lesssim\int_0^t\big(D_4^2+\Fv^2+\Fut^2\big)\dtau+\Eb^2+\Eb\Es .
\end{equation*}
The time integrals of $D_4^2$, $\Fv^2$ and $\Fut^2$ are bounded by \eqref{eq:6.22}, by $\cEv^2(t)$ and by \eqref{eq:6.10}, respectively; \eqref{eq:6.10} also bounds $\sup_{(0,t)}\Eut^2$ in \eqref{eq:6.48}, and $\sup_{(0,t)}\Ev^2\le\cEv^2(t)$. Together with \eqref{eq:6.48} and $\Es^4\le\Es^3$, these bounds yield \eqref{eq:6.47}.
\end{proof}

\subsection{Estimate of \texorpdfstring{$\cEv$}{E rho}}\label{sec:6rho}

In this subsection we estimate the last component of $\Es$, the $H^4$-norm of the density together with its inherited dissipation $\Fv$.

\begin{proposition}\label{prop:rho}
Assume \eqref{eq:small}. Then, for every $t\ge0$,
\begin{equation}\label{eq:6.49}
\begin{aligned}
\cEv^2(t)&=\sup_{\tau\in(0,t)}\Ev^2+\int_0^t\Fv^2\dtau\\
&\lesssim\Es(0)^2+\Eb^2+\Eb\Es+\Ef(0)^{\frac14}\Es^{\frac74}+\Es^3,
\end{aligned}
\end{equation}
where $\Es=\Es(t)$ and $\Eb=\Eb(t)$.
\end{proposition}

The proof occupies this and the following subsubsections. In \cref{sec:R-prep,sec:R-identity,sec:R-coupling} we derive an $H^4$-energy identity for $\varrho$, estimate its coupling with the momentum equation, including the coupling through the Coriolis force, and arrive at the integrated inequality \eqref{eq:6.74}, in which only the three transport terms $\cR_4,\cR_5,\cR_6$ remain to be estimated. These are treated in the subsequent subsubsections, where the proof is completed.

Throughout this subsection \eqref{eq:small} is in force; in particular $\Esh\le1$, so that factors $1+\Esh$ are bounded. We write $\ip{\cdot,\cdot}$ and $\nrm{\cdot}$ for the inner product and the norm of $L^2(\Omega)$. By \eqref{eq:stable}, $\rs\in C^6([0,1])$ and $\sup\rs'<0$, so the weight $\sgm=-1/\rs'$ satisfies $\sgm\in C^5([0,1])$, $\sgm^{-1}=-\rs'\in C^5([0,1])$ and $0<\sgm_-\le\sgm\le\sgm_+$ for constants $\sgm_\pm$ depending only on $\rs$. (In this part of the proof at most four derivatives of $\sgm$ and of $\sgm^{-1}$ are used.)

\subsubsection{Preparations}\label{sec:R-prep}

Multiplying \eqref{PE}$_1$ by $\sgm$, the density equation becomes
\begin{equation}\label{eq:6.50}
\sgm\dot\varrho=w,\qquad\text{i.e.}\qquad\sgm\varrho_t+\sgm\,\bU\cdot\nabla\varrho=w .
\end{equation}
Instead of applying $\nabla^4$ to \eqref{eq:6.50}, we apply three complementary fourth-order operators: the anisotropic operator $\AB=\Lap\pz^2$, the horizontal bi-Laplacian $\Laph^2$ (with the weight $\sgm$ in the energy), and $\Laph\pz\nabla=(\Laph\pz\nh,\Laph\pz^2)$.

All functions to which $\AB$ is applied below vanish on $\pa\Omega$ together with their first $z$-derivative. Indeed, $\varrho=\pz\varrho=0$ on $\pa\Omega$ gives $\sgm\varrho=\pz(\sgm\varrho)=0$ there; $w=0$ and $\pz w=-\divh v=0$ on $\pa\Omega$; and since $v$, $w$, $\nh\varrho$ and $\pz\varrho$ vanish on $\pa\Omega$,
\begin{equation*}
\bU\cdot\nabla\varrho=0\quad\text{and}\quad\pz(\bU\cdot\nabla\varrho)=\pz v\cdot\nh\varrho+v\cdot\nh\pz\varrho+\pz w\,\pz\varrho+w\,\pz^2\varrho=0\qquad\text{on }\pa\Omega .
\end{equation*}
Consequently $\sgm\varrho_t=w-\sgm\bU\cdot\nabla\varrho$ satisfies the same clamped conditions.

\emph{Norm equivalence.} We show that
\begin{equation}\label{eq:6.51}
\nrm{\varrho}_{H^4}\sim\nrm{\AB(\sgm\varrho)}_{L^2}+\nrm{\sqrt\sgm\Laph^2\varrho}_{L^2}+\nrm{\Laph\pz\nabla(\sgm\varrho)}_{L^2},
\end{equation}
with no lower-order term on the right. The upper bound for the right-hand side holds because $\sgm\in C^4$. For the lower bound, let $f\in H^4(\Omega)$ with $f=\pz f=0$ on $\pa\Omega$. Then $f$ is the unique $H^4$-solution of the clamped problem $\Lap^2u=\Lap^2f$: if $\Lap^2u=0$ and $u=\pz u=0$ on $\pa\Omega$, then $\nrm{\Lap u}^2=\ip{u,\Lap^2u}=0$, so $u=0$. Hence \Cref{thm:biharm} gives
\begin{equation}\label{eq:6.52}
\nrm{f}_{H^4(\Omega)}\lesssim\nrm{\Lap^2f}_{L^2(\Omega)}\qquad\text{for }f\in H^4(\Omega),\ f|_{\pa\Omega}=\pz f|_{\pa\Omega}=0 .
\end{equation}
Since $\Lap^2=\Laph^2+2\Laph\pz^2+\pz^4=\Laph^2+\Laph\pz^2+\AB$, we have $\nrm{\Lap^2f}\le\nrm{\AB f}+\nrm{\Laph^2f}+\nrm{\Laph\pz\nabla f}$. We apply this to $f=\sgm\varrho$. Because $\sgm$ depends only on $z$, $\Laph^2(\sgm\varrho)=\sgm\Laph^2\varrho$, so $\nrm{\Laph^2(\sgm\varrho)}\le\sgm_+^{1/2}\nrm{\sqrt\sgm\Laph^2\varrho}$; and since $\sgm^{-1}\in C^4$, $\nrm{\varrho}_{H^4}=\nrm{\sgm^{-1}(\sgm\varrho)}_{H^4}\lesssim\nrm{\sgm\varrho}_{H^4}$. This proves \eqref{eq:6.51}. With the notation
\begin{equation*}
E_{\varrho,\sgm}^2\eqdef\nrm{\AB(\sgm\varrho)}_{L^2}^2+\nrm{\sqrt\sgm\Laph^2\varrho}_{L^2}^2+\nrm{\Laph\pz\nabla(\sgm\varrho)}_{L^2}^2
\end{equation*}
we thus have, with a constant $C_*\ge1$ depending only on $\rs$,
\begin{equation}\label{eq:R-equiv}
C_*^{-1}\Ev\le E_{\varrho,\sgm}\le C_*\Ev .
\end{equation}
In particular, the $L^2$-norm of $\varrho$, which does not decay and is controlled only through the basic energy \eqref{eq:3.16}, is not needed in the estimate of $\cEv$.

\emph{Horizontal Parseval identities.} For a.e.\ $z$, Parseval's identity in $x_h$ gives $\int_{\Tl}|\Laph^2f|^2\dxh=\int_{\Tl}|\nh^4f|^2\dxh$, $\int_{\Tl}|\nh\Laph f|^2\dxh=\int_{\Tl}|\nh^3f|^2\dxh$ and $\int_{\Tl}|\Laph f|^2\dxh=\int_{\Tl}|\nh^2f|^2\dxh$, as well as the polarized versions, e.g.\ $\int_{\Tl}\Laph^2f\,\Laph^2g\dxh=\sum_{i_1,\dots,i_4=1}^2\int_{\Tl}\pa_{i_1}\cdots\pa_{i_4}f\,\pa_{i_1}\cdots\pa_{i_4}g\dxh$. Since $\sgm=\sgm(z)$,
\begin{equation}\label{eq:R-parseval}
\nrm{\sqrt\sgm\Laph^2\varrho}=\nrm{\sqrt\sgm\nh^4\varrho},\qquad\nrm{\sqrt\sgm\nh\Laph\varrho}=\nrm{\sqrt\sgm\nh^3\varrho},\qquad\nrm{\Laph\pz\nabla(\sgm\varrho)}=\nrm{\nh^2\pz\nabla(\sgm\varrho)},
\end{equation}
and, by \Cref{lem:ABmean}, $\sum_{j=1}^2\nrm{\ABmh\pa_j\Laph^2\varrho}^2=\nrm{\ABmh|\nh|^5\varrho}^2$.

\emph{Dissipation.} Let
\begin{equation*}
\cD_\varrho\eqdef\nrm{\sqrt\sgm\nabla\nh\pz\varrho}_{L^2}^2+\nrm{\ABmh|\nh|^5\varrho}_{L^2}^2+\nrm{\sqrt\sgm\nh^3\varrho}_{L^2}^2,
\end{equation*}
where $\nabla\nh\pz\varrho=(\pa_i\pa_j\pz\varrho)_{1\le i\le3,\,1\le j\le2}$ with $\pa_3=\pz$, so that $\nrm{\nabla\nh\pz\varrho}^2=\nrm{\nh^2\pz\varrho}^2+\nrm{\nh\pz^2\varrho}^2$. Comparing with $\Fv$ in \eqref{eq:energies} and using $\sgm_-\le\sgm\le\sgm_+$, we obtain
\begin{equation}\label{eq:6.53}
\max\{1,\sgm_+\}\Fv^2\ge\cD_\varrho\ge\frac\mu g\cc2\Fv^2,\qquad\cc2\eqdef\frac g{3\mu}\min\{1,\sgm_-\}.
\end{equation}

\emph{The velocity factors.} We set
\begin{equation*}
\Theta\eqdef\pz\divh(\rs v),\qquad\Theta^\perp\eqdef\pz\divh(\rs v^\perp)=-\pz(\rs\curlh v),
\end{equation*}
where we used $\divh v^\perp=-\curlh v$. Since $\nh v$ vanishes on $\pa\Omega$, $\nrm{\nh v}\le\pi^{-1}\nrm{\nh\pz v}$; together with $|\divh v|,|\curlh v|\le\sqrt2|\nh v|$, \eqref{eq:lowenergy}, \eqref{eq:4.3} and \eqref{eq:small}, this gives
\begin{equation}\label{eq:R-X}
\nrm{\Theta}+\nrm{\Theta^\perp}\le C\nrm{\nabla\nh v}\le C\mu^{-\frac12}\Ff,\qquad\nrm{\Theta}\le C\nrm{v}_{H^2}\le C\Ef,\qquad\nrm{\Theta}_{H^1}\le C\nrm{v}_{H^3}.
\end{equation}

\emph{An interpolation inequality.} For $j=1,2$, $\pa_j\varrho\in H^3$ vanishes on $\pa\Omega$ together with $\pz\pa_j\varrho$, so $\pa_j\varrho\in V=\cD(\ABh)$. By self-adjointness and \Cref{lem:ABid}, $\nrm{\pa_j\varrho}^2=\langle\ABh\pa_j\varrho,\ABmh\pa_j\varrho\rangle\le\nrm{\nabla\pz\pa_j\varrho}\nrm{\ABmh\pa_j\varrho}$. Summing over $j$, applying the Cauchy--Schwarz inequality in $j$, and using \Cref{lem:ABmean} and $|k|\ge l^{-1}$ for $k\ne0$,
\begin{equation*}
\nrm{\nh\varrho}^2\le\nrm{\nabla\nh\pz\varrho}\,\nrm{\ABmh|\nh|\varrho}\le l\,\nrm{\nabla\nh\pz\varrho}\,\nrm{\ABmh\Laph\varrho}.
\end{equation*}
Moreover $\nrm{\nabla\pa_j\varrho}^2=-\ip{\pa_j\varrho,\Lap\pa_j\varrho}$, since $\pa_j\varrho$ vanishes on $\pa\Omega$, and $\nrm{\Lap\nh\varrho}\le\nrm{\Laph\nh\varrho}+\nrm{\pz^2\nh\varrho}\le2\Fv$, because $\nrm{\Laph\nh\varrho}=\nrm{\nh^3\varrho}$ by Parseval in $x_h$. Hence, summing over $j$ and using the Cauchy--Schwarz inequality once more, recalling $\nrm{\ABmh\Laph\varrho}\le(\mu/(g\kappa))^{1/2}\Ff$ from \eqref{eq:lowenergy},
\begin{equation}\label{eq:R-interp}
\nrm{\nabla\nh\varrho}^2\le\nrm{\nh\varrho}\,\nrm{\Lap\nh\varrho}\le2l^{\frac12}\Fv^{\frac32}\nrm{\ABmh\Laph\varrho}^{\frac12}\le C\Fv^{\frac32}\Ff^{\frac12}.
\end{equation}

\subsubsection{Energy identity}\label{sec:R-identity}

We test \eqref{eq:6.50} with the three fourth-order operators and combine the results into a single energy identity at the $H^4$-level. The identities are derived for smooth solutions and extended to the class \eqref{eq:localclass} as explained in \Cref{rem:justify}. Expressions such as $\ABmh|\nh|^5\varrho$ are defined mode by mode (\Cref{lem:ABmean}). The terms $\cR_4,\cR_5,\cR_6$ below contain fifth derivatives of $\varrho$ only through their principal transport parts, which are removed by an integration by parts because $\divv\bU=0$ and $\bU\cdot n=0$ on $\pa\Omega$; they are understood in this sense.

\emph{The multiplier $\AB(\sgm\varrho)$.} Applying $\AB$ to \eqref{eq:6.50} and using \eqref{eq:wrep}, we get
\begin{equation*}
\AB(\sgm\varrho_t)+\AB(\sgm\bU\cdot\nabla\varrho)=\AB w=\frac1\mu\big[g\Laph\varrho-\pz\divh G\big].
\end{equation*}
Testing with $\AB(\sgm\varrho)$ gives
\begin{equation}\label{eq:6.54}
\frac12\ddt\nrm{\AB(\sgm\varrho)}^2-\frac g\mu\ip{\AB(\sgm\varrho),\Laph\varrho}=\underbrace{-\frac1\mu\ip{\AB(\sgm\varrho),\pz\divh G}}_{\eqdef\cR_1}\ \underbrace{-\ip{\AB(\sgm\varrho),\AB(\sgm\bU\cdot\nabla\varrho)}}_{\eqdef\cR_4}.
\end{equation}
The second term on the left is a dissipation. Indeed, $\Laph\varrho\in H^2$ and $\Laph\varrho=\pz\Laph\varrho=0$ on $\pa\Omega$, so $\Laph\varrho\in V$. Since $\sgm\varrho\in\cD(\AB)$, the definition of $\AB$ through the form $\mathfrak a$ (\Cref{def:AB}) and an integration by parts in $x_h$ give, with $f_j\eqdef\pa_j\varrho$,
\begin{equation*}
\ip{\AB(\sgm\varrho),\Laph\varrho}=\ip{\nabla\pz(\sgm\varrho),\nabla\pz\Laph\varrho}=-\sum_{j=1}^2\ip{\nabla\pz(\sgm f_j),\nabla\pz f_j}.
\end{equation*}
Each $f=f_j$ satisfies $f=\pz f=0$ and $\nh f=0$ on $\pa\Omega$. We have $\nh\pz(\sgm f)=\sgm\nh\pz f+\sgm'\nh f$ and $\pz^2(\sgm f)=\sgm\pz^2f+2\sgm'\pz f+\sgm''f$, and integrations by parts in $z$, whose boundary terms vanish, give
\begin{equation*}
\int_\Omega\sgm'\nh f\cdot\nh\pz f\dx=-\frac12\int_\Omega\sgm''|\nh f|^2\dx,\qquad2\int_\Omega\sgm'\pz f\,\pz^2f\dx=-\int_\Omega\sgm''|\pz f|^2\dx,
\end{equation*}
\begin{equation*}
\int_\Omega\sgm''f\,\pz^2f\dx=-\int_\Omega\sgm''|\pz f|^2\dx-\int_\Omega\sgm'''f\,\pz f\dx=-\int_\Omega\sgm''|\pz f|^2\dx+\frac12\int_\Omega\sgm''''f^2\dx .
\end{equation*}
Hence $\ip{\nabla\pz(\sgm f),\nabla\pz f}=\nrm{\sqrt\sgm\nabla\pz f}^2-\frac12\int\sgm''|\nh f|^2-2\int\sgm''|\pz f|^2+\frac12\int\sgm''''f^2$, and summing over $j$,
\begin{equation}\label{eq:R-diss1}
\begin{aligned}
&\frac g\mu\ip{\AB(\sgm\varrho),\Laph\varrho}=-\frac g\mu\nrm{\sqrt\sgm\nabla\nh\pz\varrho}^2+\cR_7,\\
&\cR_7\eqdef\frac g{2\mu}\int_\Omega\Big(\sgm''|\nh^2\varrho|^2+4\sgm''|\nh\pz\varrho|^2-\sgm''''|\nh\varrho|^2\Big)\dx .
\end{aligned}
\end{equation}
The remainder $\cR_7$ collects the terms in which derivatives fall on $\sgm$.

\emph{The multiplier $\Laph^2\varrho$.} Applying $\Laph^2$ to \eqref{eq:6.50}, using $\Laph^2(\sgm\,\cdot)=\sgm\Laph^2$ and testing with $\Laph^2\varrho$, we get
\begin{equation*}
\frac12\ddt\nrm{\sqrt\sgm\Laph^2\varrho}^2+\int_\Omega\sgm\Laph^2\varrho\,\Laph^2(\bU\cdot\nabla\varrho)\dx=\ip{\Laph^2\varrho,\Laph^2w}.
\end{equation*}
By \eqref{eq:wrep} and \Cref{lem:ABmean}(ii), $\Laph^2w=\mu^{-1}\ABmi\big[g\Laph^3\varrho-\Laph^2\pz\divh G\big]$. On the $k$-th horizontal mode $\Laph$ acts as multiplication by $-|k|^2$, so $\ip{\Laph^2\varrho,\ABmi\Laph^3\varrho}=-\nrm{\ABmh|\nh|^5\varrho}^2$; and one horizontal integration by parts gives $-\ip{\Laph^2\varrho,\ABmi\Laph^2\pz\divh G}=\sum_j\ip{\pa_j\Laph^2\varrho,\ABmi\pz\Laph^2G_j}$. Therefore
\begin{equation}\label{eq:6.55}
\frac12\ddt\nrm{\sqrt\sgm\Laph^2\varrho}^2+\frac g\mu\nrm{\ABmh|\nh|^5\varrho}^2=\cR_2+\cR_5,
\end{equation}
where
\begin{equation*}
\cR_2\eqdef\frac1\mu\sum_{j=1}^2\ip{\ABmh\pa_j\Laph^2\varrho,\ABmh\pz\Laph^2G_j},\qquad\cR_5\eqdef-\int_\Omega\sgm\Laph^2\varrho\,\Laph^2(\bU\cdot\nabla\varrho)\dx .
\end{equation*}

\emph{The multiplier $\Laph\pz\nabla(\sgm\varrho)$.} Applying $\Laph\pz\nabla$ to \eqref{eq:6.50} and testing with $\Laph\pz\nabla(\sgm\varrho)$, we get
\begin{equation*}
\frac12\ddt\nrm{\Laph\pz\nabla(\sgm\varrho)}^2+\ip{\Laph\pz\nabla(\sgm\varrho),\Laph\pz\nabla(\sgm\bU\cdot\nabla\varrho)}=\frac1\mu\ip{\Laph\pz\nabla(\sgm\varrho),\Laph\pz\nabla\ABmi h},
\end{equation*}
with $h\eqdef g\Laph\varrho-\pz\divh G$. Equip $V$ with the inner product $\mathfrak a$; then $\AB$ extends to the Riesz isomorphism of $V$ onto its dual $V'$, and $\ABmi$ to its inverse. Let $\phi\in H^4$ with $\phi=\pz\phi=0$ on $\pa\Omega$. Then $\Laph\phi\in V$, and $\Laph$ commutes with $\ABmi$ (\Cref{lem:ABmean}), so that
\begin{equation}\label{eq:R-clamped}
\ip{\Laph\pz\nabla\phi,\Laph\pz\nabla\ABmi h}=\mathfrak a\big(\Laph\phi,\ABmi\Laph h\big)=\ip{\Laph h,\Laph\phi}_{V',V}
\end{equation}
whenever $\Laph h\in V'$; when $\Laph h\in L^2$ the right-hand side is $\ip{\Laph\phi,\Laph h}$. This is the clamped integration by parts. We apply it with $\phi=\sgm\varrho$. For $h=g\Laph\varrho$, \eqref{eq:R-parseval} gives
\begin{equation*}
\ip{\Laph(\sgm\varrho),g\Laph^2\varrho}=g\int_\Omega\sgm\Laph\varrho\,\Laph^2\varrho\dx=-g\nrm{\sqrt\sgm\nh\Laph\varrho}^2=-g\nrm{\sqrt\sgm\nh^3\varrho}^2 .
\end{equation*}
For $h=-\pz\divh G$ with $G\in H^3$, $\Laph h$ is the element $\psi\mapsto\sum_j\ip{\Laph\pz G_j,\pa_j\psi}$ of $V'$ (for smooth $G$ this is an integration by parts in $x_h$; the functional is bounded on $V$ because $\nrm{\pa_j\psi}\le\pi^{-1}\nrm{\pa_j\pz\psi}$). With $\psi=\Laph(\sgm\varrho)$ we get $\ip{\Laph h,\Laph(\sgm\varrho)}_{V',V}=\sum_j\ip{\sgm\pa_j\Laph\varrho,\Laph\pz G_j}$. Thus
\begin{equation}\label{eq:6.56}
\frac12\ddt\nrm{\Laph\pz\nabla(\sgm\varrho)}^2+\frac g\mu\nrm{\sqrt\sgm\nh^3\varrho}^2=\cR_3+\cR_6,
\end{equation}
where
\begin{equation*}
\cR_3\eqdef\frac1\mu\sum_{j=1}^2\ip{\sgm\pa_j\Laph\varrho,\Laph\pz G_j},\qquad\cR_6\eqdef-\ip{\Laph\pz\nabla(\sgm\varrho),\Laph\pz\nabla(\sgm\bU\cdot\nabla\varrho)}.
\end{equation*}

\emph{The combined identity.} Adding \eqref{eq:6.54} (with \eqref{eq:R-diss1}), \eqref{eq:6.55} and \eqref{eq:6.56}, we obtain the $H^4$-energy identity for $\varrho$:
\begin{equation}\label{eq:6.57}
\frac12\ddt E_{\varrho,\sgm}^2+\frac g\mu\Big(\nrm{\sqrt\sgm\nabla\nh\pz\varrho}^2+\nrm{\ABmh|\nh|^5\varrho}^2+\nrm{\sqrt\sgm\nh^3\varrho}^2\Big)=\sum_{i=1}^7\cR_i,
\end{equation}
that is, $\frac12\ddt E_{\varrho,\sgm}^2+\frac g\mu\cD_\varrho=\sum_{i=1}^7\cR_i$. Here $\cR_1,\cR_2,\cR_3$ are the momentum coupling terms, $\cR_4,\cR_5,\cR_6$ the nonlinear transport terms, and $\cR_7$ is the $\sgm$-remainder. For later reference we record the transport terms in tensor form. With $\pa_3\eqdef\pz$ we have $\AB=\sum_{m=1}^3\pa_m^2\pz^2$, and by the horizontal Parseval identities,
\begin{equation}\label{eq:R-transport}
\begin{aligned}
\cR_4&=-\sum_{m,m'=1}^3\int_\Omega\pa_m^2\pz^2(\sgm\varrho)\,\pa_{m'}^2\pz^2(\sgm\bU\cdot\nabla\varrho)\dx,\\
\cR_5&=-\sum_{i_1,\dots,i_4=1}^2\int_\Omega\sgm\,\pa_{i_1}\cdots\pa_{i_4}\varrho\;\pa_{i_1}\cdots\pa_{i_4}(\bU\cdot\nabla\varrho)\dx,\\
\cR_6&=-\sum_{i,j=1}^2\sum_{m=1}^3\int_\Omega\pa_i\pa_j\pz\pa_m(\sgm\varrho)\,\pa_i\pa_j\pz\pa_m(\sgm\bU\cdot\nabla\varrho)\dx .
\end{aligned}
\end{equation}
Recall that $\sgm\bU\cdot\nabla\varrho$ and its $z$-derivative vanish on $\pa\Omega$.

\subsubsection{Momentum coupling}\label{sec:R-coupling}

In this subsubsection we estimate the momentum coupling terms $\cR_1,\cR_2,\cR_3$ and the $\sgm$-remainder $\cR_7$. We split the force \eqref{eq:Gdef} into its linear and its nonlinear part,
\begin{equation}\label{eq:R-Gsplit}
G=\rs v_t+f\rs v^\perp+G^{\mathrm{NL}},\qquad G^{\mathrm{NL}}\eqdef\varrho\dot v+\rs\bU\cdot\nabla v+f\varrho v^\perp=\varrho\big(v_t+fv^\perp\big)+\rho\,\bU\cdot\nabla v .
\end{equation}

\emph{Decomposition of $\cR_1$.} Since $\rs$ does not depend on $t$, $\pz\divh(\rs v_t)=\pt\Theta$, and \eqref{eq:R-Gsplit} gives
\begin{equation}\label{eq:6.58}
\cR_1=-\frac1\mu\ddt\ip{\AB(\sgm\varrho),\Theta}+\frac1\mu\ip{\AB(\sgm\varrho_t),\Theta}+\cR_1^{\mathrm T}-\frac1\mu\ip{\AB(\sgm\varrho),\pz\divh G^{\mathrm{NL}}},
\end{equation}
where
\begin{equation*}
\cR_1^{\mathrm T}\eqdef-\frac f\mu\ip{\AB(\sgm\varrho),\Theta^\perp}=\frac f\mu\ip{\AB(\sgm\varrho),\pz(\rs\curlh v)}
\end{equation*}
is the direct coupling of the density with the Coriolis force. For the second term in \eqref{eq:6.58}, which contains $\varrho_t$, we insert $\AB(\sgm\varrho_t)=\AB w-\AB(\sgm\bU\cdot\nabla\varrho)$ from \eqref{eq:6.50} and, by \eqref{eq:wrep} and \eqref{eq:R-Gsplit},
\begin{equation*}
\AB w=\frac1\mu\big[g\Laph\varrho-\pt\Theta-f\Theta^\perp-\pz\divh G^{\mathrm{NL}}\big].
\end{equation*}
This yields
\begin{equation}\label{eq:6.59}
\begin{aligned}
\frac1\mu\ip{\AB(\sgm\varrho_t),\Theta}&=-\frac1{2\mu^2}\ddt\nrm{\Theta}^2+\frac g{\mu^2}\ip{\Laph\varrho,\Theta}-\frac f{\mu^2}\ip{\Theta^\perp,\Theta}\\
&\quad-\frac1{\mu^2}\ip{\pz\divh G^{\mathrm{NL}},\Theta}-\frac1\mu\ip{\AB(\sgm\bU\cdot\nabla\varrho),\Theta}.
\end{aligned}
\end{equation}
Inserting \eqref{eq:6.59} into \eqref{eq:6.58}, we obtain the final decomposition
\begin{equation}\label{eq:R-R1dec}
\cR_1=-\frac1\mu\ddt\ip{\AB(\sgm\varrho),\Theta}-\frac1{2\mu^2}\ddt\nrm{\Theta}^2+\cR_1^{\mathrm L}+\cR_1^{\mathrm f}+\cR_1^{\mathrm T}+\cR_1^{\mathrm{NL}},
\end{equation}
where the terms generated by the linear gravity and Coriolis forces are
\begin{equation*}
\cR_1^{\mathrm L}\eqdef\frac g{\mu^2}\ip{\Laph\varrho,\Theta},\qquad\cR_1^{\mathrm f}\eqdef-\frac f{\mu^2}\ip{\Theta^\perp,\Theta}=\frac f{\mu^2}\ip{\pz(\rs\curlh v),\pz\divh(\rs v)},
\end{equation*}
and $\cR_1^{\mathrm T}$ is defined above, while the nonlinear part is $\cR_1^{\mathrm{NL}}\eqdef\cR_{11}^{\mathrm{NL}}+\cR_{12}^{\mathrm{NL}}+\cR_{13}^{\mathrm{NL}}$ with
\begin{equation*}
\begin{aligned}
&\cR_{11}^{\mathrm{NL}}\eqdef-\frac1{\mu^2}\ip{\pz\divh G^{\mathrm{NL}},\Theta},\qquad\cR_{12}^{\mathrm{NL}}\eqdef-\frac1\mu\ip{\AB(\sgm\bU\cdot\nabla\varrho),\Theta},\\
&\cR_{13}^{\mathrm{NL}}\eqdef-\frac1\mu\ip{\AB(\sgm\varrho),\pz\divh G^{\mathrm{NL}}}.
\end{aligned}
\end{equation*}

\emph{Justification of \eqref{eq:6.58}--\eqref{eq:R-R1dec}.} The computation above is formal, because $\AB(\sgm\varrho_t)$ and $\AB(\sgm\bU\cdot\nabla\varrho)$ contain fifth derivatives of $\varrho$. It is made rigorous as follows. Let $\Phi\eqdef\AB(\sgm\varrho)$. Applying $\AB$ to \eqref{eq:6.50} and using $\divv\bU=0$, we see that $\Phi$ solves
\begin{equation*}
\pt\Phi+\divv(\bU\Phi)=S_\Phi,\qquad S_\Phi\eqdef\AB w-\cC_1',\qquad\cC_1'\eqdef\AB(\sgm\bU\cdot\nabla\varrho)-\bU\cdot\nabla\AB(\sgm\varrho),
\end{equation*}
in the weak sense, with test functions that need not vanish on $\pa\Omega$ because $\bU\cdot n=0$; by \Cref{rem:justify}(a), $\Phi\in C([0,T];L^2)$. Since $\sgm\bU\cdot\nabla\varrho=\bU\cdot\nabla(\sgm\varrho)-\sgm'w\varrho$, the commutator $\cC_1'=[\AB,\bU\cdot\nabla](\sgm\varrho)-\AB(\sgm'w\varrho)$ contains at most four derivatives of $\varrho$. The Leibniz argument used for $\cC_1$ below, together with $\nrm{\sgm'w\varrho}_{H^4}\lesssim\nrm{w}_{H^4}\nrm{\varrho}_{L^\infty}+\nrm{w}_{L^\infty}\nrm{\varrho}_{H^4}$ (\Cref{lem:GN}(2)), gives $\nrm{\cC_1'}\lesssim\Ev\nrm{v}_{H^5}$; as $\AB w=\mu^{-1}[g\Laph\varrho-\pz\divh G]$ by \eqref{eq:wrep}, we get $S_\Phi\in L^2(0,T;L^2)$. Moreover $\Theta\in L^2(0,T;H^1)\cap H^1(0,T;L^2)$ with $\pt\Theta=\pz\divh(\rs v_t)$, because $v\in L^2(0,T;H^5)$ and $v_t\in C([0,T];H^2)$, and $\bU\in L^\infty(0,T;W^{1,\infty})$ by \eqref{eq:3.13}. Testing the weak formulation with $\Theta$ (after a regularization of $\Theta$ in time), we obtain $\ip{\Phi,\Theta}\in W^{1,1}(0,T)$ and
\begin{equation*}
\ddt\ip{\Phi,\Theta}=\ip{\Phi,\bU\cdot\nabla\Theta}+\ip{S_\Phi,\Theta}+\ip{\Phi,\pt\Theta}\qquad\text{for a.e. }t .
\end{equation*}
We insert this into $\cR_1=-\mu^{-1}\ip{\Phi,\pt\Theta}+\cR_1^{\mathrm T}-\mu^{-1}\ip{\Phi,\pz\divh G^{\mathrm{NL}}}$, expand $\ip{S_\Phi,\Theta}$ by means of \eqref{eq:wrep} and \eqref{eq:R-Gsplit}, and use $\ip{\pt\Theta,\Theta}=\frac12\ddt\nrm{\Theta}^2$, valid for $\Theta\in H^1(0,T;L^2)$. This gives \eqref{eq:R-R1dec} directly, with the same terms $\cR_1^{\mathrm L},\cR_1^{\mathrm f},\cR_1^{\mathrm T},\cR_{11}^{\mathrm{NL}},\cR_{13}^{\mathrm{NL}}$, and with
\begin{equation*}
\cR_{12}^{\mathrm{NL}}=\frac1\mu\big(\ip{\Phi,\bU\cdot\nabla\Theta}-\ip{\cC_1',\Theta}\big).
\end{equation*}
For smooth $\varrho$ satisfying the clamped conditions, this expression and the right-hand side of \eqref{eq:6.62} below both equal $-\mu^{-1}\ip{\AB(\sgm\bU\cdot\nabla\varrho),\Theta}$; both are continuous in $\varrho$ with respect to the $H^4$-norm (for fixed $\bU$ and $\Theta$), so they coincide by density.

Compared with the non-rotating case, the linear Coriolis force produces the two quadratic terms $\cR_1^{\mathrm f}$ and $\cR_1^{\mathrm T}$. The first one comes from the vortical part $W^{\mathrm L_3}$ of the vertical velocity in \eqref{eq:wrep} and involves only two derivatives of $v$; the second one pairs $\AB(\sgm\varrho)$ with the Coriolis force.

\emph{The terms $\cR_1^{\mathrm L}$, $\cR_1^{\mathrm f}$ and $\cR_1^{\mathrm T}$.} By $\nrm{\Laph\varrho}=\nrm{\nh^2\varrho}\le l\nrm{\nh^3\varrho}\le l\Fv$ and \eqref{eq:R-X},
\begin{equation}\label{eq:6.60}
|\cR_1^{\mathrm L}|\le\frac g{\mu^2}\nrm{\Laph\varrho}\nrm{\Theta}\le C\Fv\Ff\le\varepsilon\Fv^2+C\varepsilon^{-1}\Ff^2 .
\end{equation}
Again by \eqref{eq:R-X},
\begin{equation}\label{eq:R-R1f}
|\cR_1^{\mathrm f}|\le\frac{|f|}{\mu^2}\nrm{\Theta^\perp}\nrm{\Theta}\le C|f|\nrm{\nabla\nh v}^2\le C|f|\Ff^2 .
\end{equation}
Finally, since $\nrm{\AB(\sgm\varrho)}\le C\Ev$, \eqref{eq:R-X} and \Cref{lem:intFflat} give
\begin{equation}\label{eq:R-T}
\begin{aligned}
&|\cR_1^{\mathrm T}|\le\frac{|f|}\mu\nrm{\AB(\sgm\varrho)}\nrm{\Theta^\perp}\le C|f|\Ev\Ff,\\
&\int_0^t|\cR_1^{\mathrm T}|\dtau\le C|f|\Es(t)\int_0^t\Ff\dtau\le C|f|\cc4\Ef(0)^{\frac14}\Es(t)^{\frac74}.
\end{aligned}
\end{equation}

\begin{remark}[The rotational coupling]\label{rem:R-T}
The term $\cR_1^{\mathrm T}$ is quadratic in the perturbation, but it is the trace in the energy identity of the linear Coriolis force and carries no small factor. It is not absorbed by the dissipation: $\AB(\sgm\varrho)$ contains $\pz^4(\sgm\varrho)$, while $\Fv$ controls only derivatives of $\varrho$ of order at most three that contain a horizontal derivative, together with the weak norm $\nrm{\ABmh|\nh|^5\varrho}$; pure vertical derivatives of order four are not dissipated. In \eqref{eq:R-T} the term is therefore bounded by the top-order energy times the low-order dissipation $\Ff$, which carries only two derivatives of $v$. The bound $\int_0^t\Ff\lesssim\Es$ of \eqref{eq:5.6} would turn this into a quadratic contribution $C\Es^2$ whose constant is not small. \Cref{lem:intFflat} measures $\int_0^t\Ff$ by a positive power of the initial low-order energy instead, and the resulting contribution satisfies, for every $\eta>0$,
\begin{equation*}
\Ef(0)^{\frac14}\Es^{\frac74}\le\eta\Es^2+C_\eta\Ef(0)^2
\end{equation*}
by Young's inequality with exponents $\frac87$ and $8$. Since $\Ef(0)\lesssim\Es(0)$, this term is absorbed in the bootstrap. For $f=0$ the terms $\cR_1^{\mathrm T}$ and $\cR_1^{\mathrm f}$ vanish.
\end{remark}

\emph{Two product estimates.} For the nonlinear terms we use
\begin{equation}\label{eq:R-UvH2}
\nrm{\bU\cdot\nabla v}_{H^2}\lesssim\nrm{v}_{H^3}^2,\qquad\nrm{\rho\,\bU\cdot\nabla v}_{H^2}\lesssim(1+\Ev)\nrm{v}_{H^3}^2\lesssim\nrm{v}_{H^3}^2 .
\end{equation}
For the horizontal transport, \Cref{lem:GN}(2) with $m=2$ and $H^2(\Omega)\subset L^\infty(\Omega)$ give $\nrm{v\cdot\nh v}_{H^2}\lesssim\nrm{v}_{L^\infty}\nrm{\nh v}_{H^2}+\nrm{v}_{H^2}\nrm{\nh v}_{L^\infty}\lesssim\nrm{v}_{H^2}\nrm{v}_{H^3}$. For the vertical transport we expand $\pa^\alpha(w\pz v)$, $|\alpha|\le2$, by the Leibniz rule and estimate each product $\pa^\beta w\,\pa^{\alpha-\beta}\pz v$. If $\beta=0$, we use $\nrm{w}_{L^\infty}\lesssim\nrm{v}_{H^3}$ from \eqref{eq:3.13}. If $|\beta|=1$, we use $\nrm{\nabla w}_{L^4}\nrm{\nabla\pz v}_{L^4}+\nrm{\nabla w}_{L^4}\nrm{\pz v}_{L^4}\lesssim\nrm{w}_{H^2}\nrm{v}_{H^3}\lesssim\nrm{v}_{H^3}^2$ by \Cref{lem:GN}(3) and \Cref{lem:w}. If $|\beta|=2$ and $\pa^\beta$ contains $\pz$, then $\pa^\beta w=-\pa^{\beta-e_z}\divh v$ and $\nrm{\nabla\divh v}_{L^4}\nrm{\pz v}_{L^4}\lesssim\nrm{v}_{H^3}\nrm{v}_{H^2}$. If $\pa^\beta=\pa_i\pa_j$ is horizontal, the factor $\pz v$ cannot be put in $L^\infty_hL^2_z$ at the level $H^2$ in three dimensions; we use \Cref{lem:w} and \Cref{lem:aniso}(e) with $s=\frac12$:
\begin{equation*}
\nrm{\pa_i\pa_jw\,\pz v}\le\nrm{\pa_i\pa_jw}_{L^2_hL^\infty_z}\nrm{\pz v}_{L^\infty_hL^2_z}\lesssim\nrm{\nh^2\divh v}\,\nrm{\pz v}_{H^{3/2}}\lesssim\nrm{v}_{H^3}^2 .
\end{equation*}
This proves the first bound in \eqref{eq:R-UvH2}; the second follows from $\nrm{\rho g}_{H^2}\lesssim\nrm{\rho}_{W^{2,\infty}}\nrm{g}_{H^2}$ and \eqref{eq:3.12}. Since $H^4(\Omega)\subset W^{2,\infty}(\Omega)$ and $\nrm{v}_{H^3}^2\le\Eu\Fu$, \eqref{eq:R-UvH2} and \Cref{lem:GN}(2) with $m=3$ give
\begin{equation}\label{eq:R-GNL}
\begin{aligned}
\nrm{G^{\mathrm{NL}}}_{H^2}&\lesssim\nrm{\varrho}_{W^{2,\infty}}\big(\nrm{v_t}_{H^2}+|f|\nrm{v}_{H^2}\big)+\nrm{v}_{H^3}^2\lesssim\Esh\Fsh,\\
\nrm{\varrho(v_t+fv^\perp)}_{H^3}&\lesssim\nrm{\varrho}_{L^\infty}\nrm{v_t+fv^\perp}_{H^3}+\nrm{\varrho}_{H^3}\nrm{v_t+fv^\perp}_{L^\infty}\lesssim\Ev\big(\Fut+|f|\Fu\big)\lesssim\Esh\Fsh .
\end{aligned}
\end{equation}

\emph{The term $\cR_{11}^{\mathrm{NL}}$.} By \eqref{eq:R-GNL}, \eqref{eq:R-X} and $\Ff\lesssim\Fsh$ from \eqref{eq:3.11},
\begin{equation}\label{eq:6.61}
|\cR_{11}^{\mathrm{NL}}|\le\frac1{\mu^2}\nrm{G^{\mathrm{NL}}}_{H^2}\nrm{\Theta}\lesssim\Esh\Fsh\Ff\lesssim\Esh\Fsh^2 .
\end{equation}

\emph{The term $\cR_{12}^{\mathrm{NL}}$.} Since $\varrho\in\cD(\AB)$, $\AB\varrho=\Lap\pz^2\varrho\in L^2$. We commute $\AB$ with the transport operator:
\begin{equation}\label{eq:6.62}
\begin{aligned}
\cR_{12}^{\mathrm{NL}}&=-\frac1\mu\ip{\sgm\bU\cdot\nabla\AB\varrho,\Theta}-\frac1\mu\ip{\cC_1,\Theta}=\frac1\mu\ip{\AB\varrho,\bU\cdot\nabla(\sgm\Theta)}-\frac1\mu\ip{\cC_1,\Theta},
\end{aligned}
\end{equation}
where we used $\divv\bU=0$ and $\bU\cdot n=0$ on $\pa\Omega$ (for smooth $\varrho$; the right-hand side involves only $\AB\varrho\in L^2$ and $\cC_1\in L^2$, and it agrees with the definition of $\cR_{12}^{\mathrm{NL}}$ in the justification of \eqref{eq:R-R1dec} by density), and where
\begin{equation*}
\begin{aligned}
\cC_1&\eqdef\AB(\sgm\bU\cdot\nabla\varrho)-\sgm\bU\cdot\nabla\AB\varrho=\sgm\big(\AB(\bU\cdot\nabla\varrho)-\bU\cdot\nabla\AB\varrho\big)+[\AB,\sgm](\bU\cdot\nabla\varrho),\\
[\AB,\sgm]h&\eqdef2\sgm'\Laph\pz h+\sgm''\Laph h+4\sgm'\pz^3h+6\sgm''\pz^2h+4\sgm'''\pz h+\sgm''''h .
\end{aligned}
\end{equation*}
For the first term in \eqref{eq:6.62}, $\nrm{\bU}_{L^\infty}\le\nrm{v}_{L^\infty}+\nrm{w}_{L^\infty}\lesssim\nrm{v}_{H^3}$ by \eqref{eq:3.13}, and \eqref{eq:R-X} give
\begin{equation}\label{eq:6.63}
\big|\ip{\AB\varrho,\bU\cdot\nabla(\sgm\Theta)}\big|\le\nrm{\AB\varrho}\nrm{\bU}_{L^\infty}\nrm{\sgm\Theta}_{H^1}\lesssim\Ev\nrm{v}_{H^3}^2\lesssim\Esh\Fsh^2 .
\end{equation}
For the commutator, write $\AB=\sum_{m=1}^3\pa^{\alpha_m}$ with $\alpha_m=2e_m+2e_3$ ($\pa_3=\pz$), so that $\AB(\bU\cdot\nabla\varrho)-\bU\cdot\nabla\AB\varrho=\sum_m\sum_{0<\beta\le\alpha_m}\binom{\alpha_m}\beta\pa^\beta\bU\cdot\nabla\pa^{\alpha_m-\beta}\varrho$. By H\"older's inequality, \Cref{lem:GN}(3) and $H^4\subset W^{2,\infty}$, a summand with $|\beta|=1$ is bounded by $\nrm{\nabla\bU}_{L^\infty}\nrm{\varrho}_{H^4}$, a summand with $|\beta|=2$ by $\nrm{\nabla^2\bU}_{L^4}\nrm{\nabla^3\varrho}_{L^4}\lesssim\nrm{\bU}_{H^3}\nrm{\varrho}_{H^4}$, and a summand with $|\beta|\in\{3,4\}$ by $\nrm{\pa^\beta\bU}\nrm{\varrho}_{W^{2,\infty}}\lesssim\nrm{\bU}_{H^4}\nrm{\varrho}_{H^4}$. Moreover $\nrm{[\AB,\sgm](\bU\cdot\nabla\varrho)}\lesssim\nrm{\bU\cdot\nabla\varrho}_{H^3}\lesssim\nrm{\bU}_{H^3}\nrm{\nabla\varrho}_{L^\infty}+\nrm{\bU}_{L^\infty}\nrm{\nabla\varrho}_{H^3}$ by \Cref{lem:GN}(2). Since $\nrm{\nabla\bU}_{L^\infty}\lesssim\nrm{v}_{H^4}$ by \eqref{eq:3.13} and $\nrm{w}_{H^k}\lesssim\nrm{v}_{H^{k+1}}$ by \Cref{lem:w}, we obtain $\nrm{\cC_1}\lesssim\Ev\nrm{v}_{H^5}=\Ev\Fu$, and with \eqref{eq:R-X},
\begin{equation}\label{eq:6.64}
\big|\ip{\cC_1,\Theta}\big|\le\nrm{\cC_1}\nrm{\Theta}\lesssim\Ev\Fu\Ff\lesssim\Esh\Fsh^2 .
\end{equation}
Combining \eqref{eq:6.63} and \eqref{eq:6.64},
\begin{equation}\label{eq:6.65}
|\cR_{12}^{\mathrm{NL}}|\lesssim\Esh\Fsh^2 .
\end{equation}

\emph{The term $\cR_{13}^{\mathrm{NL}}$.} By \eqref{eq:R-Gsplit}, $\cR_{13}^{\mathrm{NL}}=-\frac1\mu(\cR_{131}^{\mathrm{NL}}+\cR_{132}^{\mathrm{NL}})$ with
\begin{equation*}
\cR_{131}^{\mathrm{NL}}\eqdef\ip{\AB(\sgm\varrho),\pz\divh Y},\quad Y\eqdef\varrho\big(v_t+fv^\perp\big),\qquad\cR_{132}^{\mathrm{NL}}\eqdef\ip{\AB(\sgm\varrho),\pz\divh(\rho\,\bU\cdot\nabla v)}.
\end{equation*}
The horizontal field $Y$ and its derivatives $\pz Y$, $\pz^2Y$ vanish on $\pa\Omega$, because $\varrho=\pz\varrho=0$ and $v_t=v=0$ there; hence so do $\pz\divh Y$ and $\pz^2\divh Y$. Two integrations by parts in $z$, whose boundary terms therefore vanish, and one in $x_h$ give (for smooth functions, and then by density)
\begin{equation*}
\cR_{131}^{\mathrm{NL}}=\ip{\Lap(\sgm\varrho),\pz^3\divh Y}=-\sum_{j=1}^2\ip{\pa_j\Lap(\sgm\varrho),\pz^3Y_j}.
\end{equation*}
Since $\Lap(\sgm\varrho)=\sgm\Lap\varrho+2\sgm'\pz\varrho+\sgm''\varrho$, $\nrm{\nh\Lap\varrho}\le\nrm{\nh\Laph\varrho}+\nrm{\nh\pz^2\varrho}\le2\Fv$, $\nrm{\nh\pz\varrho}\le\pi^{-1}\nrm{\nh\pz^2\varrho}$ and $\nrm{\nh\varrho}\le l^2\nrm{\nh^3\varrho}$, we have $\nrm{\nh\Lap(\sgm\varrho)}\le C\Fv$. Together with \eqref{eq:R-GNL},
\begin{equation}\label{eq:6.66}
|\cR_{131}^{\mathrm{NL}}|\le\nrm{\nh\Lap(\sgm\varrho)}\nrm{Y}_{H^3}\lesssim\Fv\Esh\Fsh\lesssim\Esh\Fsh^2 .
\end{equation}
For $\cR_{132}^{\mathrm{NL}}$, \eqref{eq:R-UvH2} gives
\begin{equation}\label{eq:6.67}
|\cR_{132}^{\mathrm{NL}}|\le\nrm{\AB(\sgm\varrho)}\nrm{\rho\,\bU\cdot\nabla v}_{H^2}\lesssim\Ev\nrm{v}_{H^3}^2\lesssim\Esh\Fsh^2 .
\end{equation}
Combining \eqref{eq:6.66} and \eqref{eq:6.67},
\begin{equation}\label{eq:6.68}
|\cR_{13}^{\mathrm{NL}}|\lesssim\Esh\Fsh^2 ,
\end{equation}
and, by \eqref{eq:6.61}, \eqref{eq:6.65} and \eqref{eq:6.68},
\begin{equation}\label{eq:6.69}
|\cR_1^{\mathrm{NL}}|\le|\cR_{11}^{\mathrm{NL}}|+|\cR_{12}^{\mathrm{NL}}|+|\cR_{13}^{\mathrm{NL}}|\lesssim\Esh\Fsh^2 .
\end{equation}

\emph{The terms $\cR_2$ and $\cR_3$.} For $\cR_2$ we write $\Laph^2G_j=\divh(\nh\Laph G_j)$ and use \eqref{eq:ABdual}, \eqref{eq:R-parseval} and the Cauchy--Schwarz inequality in $j$:
\begin{equation}\label{eq:6.70}
\begin{aligned}
|\cR_2|&\le\frac1\mu\nrm{\ABmh|\nh|^5\varrho}\Big(\sum_{j=1}^2\nrm{\nh\Laph G_j}^2\Big)^{\frac12}=\frac1\mu\nrm{\ABmh|\nh|^5\varrho}\nrm{\nh^3G}\\
&\le\varepsilon\nrm{\ABmh|\nh|^5\varrho}^2+C\varepsilon^{-1}\nrm{\nh^3G}^2 .
\end{aligned}
\end{equation}
For $\cR_3$, by \eqref{eq:R-parseval},
\begin{equation}\label{eq:6.71}
|\cR_3|\le\frac{\sgm_+^{1/2}}\mu\nrm{\sqrt\sgm\nh^3\varrho}\nrm{\Laph\pz G}\le\varepsilon\nrm{\sqrt\sgm\nh^3\varrho}^2+C\varepsilon^{-1}\nrm{\Laph\pz G}^2 .
\end{equation}
In both bounds $G=\rho\dot v+f\rho v^\perp$ contains the linear Coriolis force, which is not small. Since $|\pa^\alpha(\rho v)^\perp|=|\pa^\alpha(\rho v)|$ pointwise, \eqref{eq:T4-GH3} gives
\begin{equation}\label{eq:R-GH3}
\nrm{\nh^3G}+\nrm{\Laph\pz G}\le2\nrm{\rho\dot v}_{H^3}+|f|\big(\nrm{\nh^3(\rho v)}+\nrm{\Laph\pz(\rho v)}\big)\lesssim\nrm{\rho\dot v}_{H^3}+D_4+\Esh\Fuz .
\end{equation}
Thus the rotational contribution is routed through the tangential dissipation $D_4$, whose time integral is controlled by \eqref{eq:6.22}; a bound of $|f|\nrm{\rho v}_{H^3}$ by $\Fuz$ or $\Fu$ with a constant of order one could not be used here, because the estimate of $\cEu$ (\Cref{prop:u}) depends on $\cEv$.

\emph{The $\sgm$-remainder.} By \eqref{eq:R-diss1}, $\nrm{\nh\varrho}\le l\nrm{\nh^2\varrho}$ and \eqref{eq:R-interp},
\begin{equation}\label{eq:6.72}
|\cR_7|\le\frac g{2\mu}\Big(4\nrm{\sgm''}_{L^\infty}\nrm{\nabla\nh\varrho}^2+\nrm{\sgm''''}_{L^\infty}\nrm{\nh\varrho}^2\Big)\le C\nrm{\nabla\nh\varrho}^2\le C\Fv^{\frac32}\Ff^{\frac12}\le\varepsilon\Fv^2+C_\varepsilon\Ff^2 .
\end{equation}

\emph{The differential inequality.} We insert \eqref{eq:R-R1dec} into \eqref{eq:6.57} and use \eqref{eq:6.53}, \eqref{eq:6.60}, \eqref{eq:R-R1f}, \eqref{eq:R-T}, \eqref{eq:6.69}--\eqref{eq:R-GH3} and \eqref{eq:6.72}. The $\varepsilon$-terms are bounded by $\varepsilon(2+\max\{1,\sgm_+\})\Fv^2$, since $\nrm{\ABmh|\nh|^5\varrho}^2+\nrm{\sqrt\sgm\nh^3\varrho}^2\le\cD_\varrho$; choosing $\varepsilon$ such that this is at most $\frac14\cc2\Fv^2$ and using $\frac g\mu\cD_\varrho\ge\cc2\Fv^2$, as well as $\Esh^2\Fuz^2\le\Esh\Fsh^2$, we obtain
\begin{equation}\label{eq:6.73}
\begin{aligned}
&\frac12\ddt\Big(E_{\varrho,\sgm}^2+\frac2\mu\ip{\AB(\sgm\varrho),\Theta}+\frac1{\mu^2}\nrm{\Theta}^2\Big)+\frac34\cc2\Fv^2\\
&\quad\le C\Big(\nrm{\rho\dot v}_{H^3}^2+D_4^2+\Ff^2+|f|\Esh\Ff+\Esh\Fsh^2\Big)+|\cR_4|+|\cR_5|+|\cR_6|,
\end{aligned}
\end{equation}
where $\Theta=\pz\divh(\rs v)$. Here $|f|\Esh\Ff$ comes from the rotational coupling $\cR_1^{\mathrm T}$, and the term $|f|\Ff^2$ from $\cR_1^{\mathrm f}$ has been included in $C\Ff^2$.

We integrate \eqref{eq:6.73} over $(0,\tau)$, $0<\tau\le t$. At the initial time, by \eqref{eq:R-equiv} and \eqref{eq:R-X},
\begin{equation*}
E_{\varrho,\sgm}^2(0)+\frac2\mu\big|\ip{\AB(\sgm\varrho),\Theta}(0)\big|+\frac1{\mu^2}\nrm{\Theta(0)}^2\le2E_{\varrho,\sgm}^2(0)+\frac2{\mu^2}\nrm{\Theta(0)}^2\lesssim\nrm{\vin}_{H^4}^2+\nrm{\uin}_{H^2}^2\lesssim\Es(0)^2 .
\end{equation*}
For the right-hand side, \eqref{eq:6.23} and \eqref{eq:6.22} give
\begin{equation*}
\int_0^t\big(\nrm{\rho\dot v}_{H^3}^2+D_4^2\big)\dtau\lesssim\Es(0)^2+\Eb^2+\Eb\Es+\Es^3 ;
\end{equation*}
by \eqref{eq:3.5}, $\int_0^t\Ff^2\le\Eb^2$; by \eqref{eq:Esharpdef}, $\int_0^t\Esh\Fsh^2\le\sup_{(0,t)}\Esh\int_0^t\Fsh^2\lesssim\Es^3$; and by \eqref{eq:R-T}, $\int_0^t|f|\Esh\Ff\lesssim\Ef(0)^{1/4}\Es^{7/4}$. Since all these bounds are nondecreasing in $t$, we obtain, for every $\tau\in[0,t]$,
\begin{equation}\label{eq:6.74}
\begin{aligned}
&\frac12E_{\varrho,\sgm}^2(\tau)+\frac1\mu\ip{\AB(\sgm\varrho),\Theta}(\tau)+\frac1{2\mu^2}\nrm{\Theta(\tau)}^2+\frac34\cc2\int_0^\tau\Fv^2\dd s\\
&\quad\lesssim\Es(0)^2+\Eb^2+\Eb\Es+\Ef(0)^{\frac14}\Es^{\frac74}+\Es^3+\int_0^t\big(|\cR_4|+|\cR_5|+|\cR_6|\big)\dd s,
\end{aligned}
\end{equation}
where $\Es=\Es(t)$, $\Eb=\Eb(t)$, and $\cR_4,\cR_5,\cR_6$ are the transport terms defined in \eqref{eq:6.54}--\eqref{eq:6.56} (see also \eqref{eq:R-transport}).

\begin{remark}\label{rem:square}
The modified energy in \eqref{eq:6.73} is a perfect square in its first component: with $\Gamma\eqdef\AB(\sgm\varrho)+\mu^{-1}\pz\divh(\rs v)$,
\begin{equation*}
E_{\varrho,\sgm}^2+\frac2\mu\ip{\AB(\sgm\varrho),\Theta}+\frac1{\mu^2}\nrm{\Theta}^2=\nrm{\Gamma}^2+\nrm{\sqrt\sgm\Laph^2\varrho}^2+\nrm{\Laph\pz\nabla(\sgm\varrho)}^2 .
\end{equation*}
Since $\nrm{\Theta}\le C\nrm{v}_{H^2}\le C\Ef$ by \eqref{eq:R-X}, and $\nrm{a+b}^2\ge\frac12\nrm{a}^2-\nrm{b}^2$, it is equivalent to $E_{\varrho,\sgm}^2$ up to $C\Ef^2$:
\begin{equation*}
\frac12E_{\varrho,\sgm}^2-C\Ef^2\le E_{\varrho,\sgm}^2+\frac2\mu\ip{\AB(\sgm\varrho),\Theta}+\frac1{\mu^2}\nrm{\Theta}^2\le2E_{\varrho,\sgm}^2+C\Ef^2 .
\end{equation*}
This is the only property of the modified energy used in the closing step. The rotation does not affect the perfect square; the extra term $\cR_1^{\mathrm f}$ produced by the vortical part of $w$ is of lower order. In particular, \eqref{eq:6.74}, \eqref{eq:R-equiv} and $\Ef(\tau)\le\Eb(t)$ give
\begin{equation}\label{eq:R-close}
\sup_{\tau\in(0,t)}\Ev^2+\int_0^t\Fv^2\dtau\lesssim\Es(0)^2+\Eb^2+\Eb\Es+\Ef(0)^{\frac14}\Es^{\frac74}+\Es^3+\int_0^t\big(|\cR_4|+|\cR_5|+|\cR_6|\big)\dtau .
\end{equation}
\end{remark}

\subsubsection{Nonlinear transport}

The terms $\cR_4,\cR_5,\cR_6$ arise from the nonlinear transport term $\sgm\bU\cdot\nabla\varrho$ tested against the three fourth-order operators. Throughout this subsubsection \eqref{eq:small} is in force, so that $\Esh\lesssim1$, and the manipulations are justified as in \Cref{rem:justify}. The implicit constants depend only on $\mu,g,f,l$, $\nrm{\rs+\vin}_{L^\infty}$, a positive lower bound for $\rs+\vin$, and $\nrm{\sgm}_{C^4([0,1])}$, except in case~(i) below, where $\nrm{\sgm^{(5)}}_{L^2(0,1)}$ (finite by \eqref{eq:stable}) enters as well. Velocity factors of top order are bounded by $\Fu\le\Fsh$; by \eqref{eq:3.3} and \eqref{eq:6.26}, $\Fu\lesssim D_4+\Fv+\Fut+\nrm{v}_{H^2}$, so \Cref{prop:u} is not used. The third and fourth horizontal derivatives of $w$ are never placed in $L^\infty$: $\nh^3w$ is kept in $L^3$ and $\nh^4w$ in $L^2$, or they are expressed through \eqref{eq:wrep}.

\emph{Multi-indices.} In this subsubsection multi-indices are three-dimensional: for $\alpha=(\alpha_1,\alpha_2,\alpha_3)\in\N_0^3$ we write $\pa^\alpha=\px^{\alpha_1}\py^{\alpha_2}\pz^{\alpha_3}$, $|\alpha|=\alpha_1+\alpha_2+\alpha_3$, $\alpha!=\alpha_1!\,\alpha_2!\,\alpha_3!$, $\alpha_h=(\alpha_1,\alpha_2)$, $\pa^{\alpha_h}=\px^{\alpha_1}\py^{\alpha_2}$, and $e_1,e_2,e_3$ for the unit multi-indices. Inequalities $\beta\le\alpha$ are meant componentwise. Since $\sgm$ depends only on $z$, $\pa^\delta\sgm=0$ unless $\delta=\delta_3e_3$, in which case $\pa^\delta\sgm=\sgm^{(\delta_3)}$. For a horizontal multi-index we put $c_\alpha\eqdef(\alpha_1+\alpha_2)!/(\alpha_1!\,\alpha_2!)$.

\emph{Index sets.} Let
\begin{equation*}
\mathcal S_1\eqdef\{(2,0,2),(0,2,2),(0,0,4)\},\qquad\text{so that}\qquad\sum_{\alpha\in\mathcal S_1}\pa^\alpha=\Laph\pz^2+\pz^4=\AB
\end{equation*}
as differential operators, and
\begin{equation*}
\begin{aligned}
&\cM_1\eqdef\mathcal S_1\times\mathcal S_1,\qquad\cM_2\eqdef\{\alpha\in\N_0^3:\ \alpha_3=0,\ |\alpha|=4\},\\
&\cM_3\eqdef\{(\alpha,\alpha):\ \alpha\in\N_0^3,\ |\alpha|=4,\ \alpha_3\in\{1,2\}\}.
\end{aligned}
\end{equation*}
For $m\ge1$, a bounded weight $\varphi=\varphi(z)$ and sufficiently regular $a,b$, Parseval's identity in $x_h$ and $\sum_{\kappa_1+\kappa_2=m}c_{(\kappa_1,\kappa_2,0)}k_1^{2\kappa_1}k_2^{2\kappa_2}=|k|^{2m}$ give
\begin{equation}\label{eq:NT-pars}
\int_\Omega\varphi\,|\nh|^ma\;|\nh|^mb\dx=\sum_{|\kappa|=m,\ \kappa_3=0}c_\kappa\int_\Omega\varphi\,\pa^\kappa a\;\pa^\kappa b\dx .
\end{equation}
We apply \eqref{eq:NT-pars} with $m=4$, $\varphi=\sgm$ to $\cR_5=-\int_\Omega\sgm\Laph^2\varrho\,\Laph^2(\bU\cdot\nabla\varrho)\dx$, and with $m=3$ and $m=2$ to the two parts of $\cR_6$, using $\ip{\Laph\pz\nabla a,\Laph\pz\nabla b}=\ip{|\nh|^3\pz a,|\nh|^3\pz b}+\ip{|\nh|^2\pz^2a,|\nh|^2\pz^2b}$. Hence
\begin{equation*}
\begin{aligned}
\cR_4&=-\sum_{(\alpha,\alpha')\in\cM_1}\int_\Omega\pa^\alpha(\sgm\varrho)\,\pa^{\alpha'}(\sgm\bU\cdot\nabla\varrho)\dx,\qquad \cR_5=-\sum_{\alpha\in\cM_2}c_{\alpha}\int_\Omega\sgm\,\pa^\alpha\varrho\,\pa^\alpha(\bU\cdot\nabla\varrho)\dx,\\
\cR_6&=-\sum_{(\alpha,\alpha)\in\cM_3}c_{\alpha_h}\int_\Omega\pa^\alpha(\sgm\varrho)\,\pa^\alpha(\sgm\bU\cdot\nabla\varrho)\dx,
\end{aligned}
\end{equation*}
where $c_{\alpha_h}\eqdef c_{(\alpha_1,\alpha_2,0)}$. By the Leibniz rule,
\begin{equation*}
\cR_4+\cR_5+\cR_6=\sum_{k=0}^4\mathcal P_k,
\end{equation*}
where $\mathcal P_k$ collects the terms in which exactly $k$ derivatives fall on $\varrho$ inside $\nabla\varrho$:
\begin{equation*}
\begin{aligned}
\mathcal P_k\eqdef{}&\sum_{(\alpha,\alpha')\in\cM_1}\sum_{\substack{\beta+\gamma+\delta=\alpha'\\|\gamma|=k}}\frac{\alpha'!}{\beta!\,\gamma!\,\delta!}\cJ_{\alpha,\alpha'}^{\beta,\gamma,\delta}+\sum_{(\alpha,\alpha)\in\cM_3}c_{\alpha_h}\sum_{\substack{\beta+\gamma+\delta=\alpha\\|\gamma|=k}}\frac{\alpha!}{\beta!\,\gamma!\,\delta!}\cJ_{\alpha,\alpha}^{\beta,\gamma,\delta}\\
&+\sum_{\alpha\in\cM_2}c_\alpha\sum_{\substack{\beta+\gamma=\alpha\\|\gamma|=k}}\frac{\alpha!}{\beta!\,\gamma!}\cK_\alpha^{\beta,\gamma},
\end{aligned}
\end{equation*}
with the summands
\begin{equation*}
\cJ_{\alpha,\alpha'}^{\beta,\gamma,\delta}\eqdef-\int_\Omega\pa^\alpha(\sgm\varrho)\,\pa^\beta\bU\cdot\nabla\pa^\gamma\varrho\,\pa^\delta\sgm\dx,\qquad \cK_\alpha^{\beta,\gamma}\eqdef-\int_\Omega\sgm\,\pa^\alpha\varrho\,\pa^\beta\bU\cdot\nabla\pa^\gamma\varrho\dx .
\end{equation*}
In $\cK_\alpha^{\beta,\gamma}$ the multi-indices $\alpha,\beta,\gamma$ are horizontal. The numerical coefficients are bounded, so for $k\le3$ it suffices to estimate each $\cJ$ and $\cK$ separately; for $k=4$ the sum over $\cM_1$ has to be taken first. We collect the bounds that replace the planar embeddings.

\begin{lemma}\label{lem:NT-aux}
Assume \eqref{eq:small}. Then:
\begin{enumerate}[label=(\alph*)]
\item (Decaying velocity factors.) $\nrm{\bU}_{L^\infty}\le\nrm{\nabla\bU}_{L^\infty}\le\Fn$, and
\begin{equation*}
\nrm{\nh v}_{H^2}+\nrm{\nabla\nh v}_{L^3}\lesssim\Fn,\qquad \nrm{\pa^\beta w}_{L^2}\lesssim\nrm{\nh v}_{H^{|\beta|}}\quad\text{for every }\beta\in\N_0^3;
\end{equation*}
in particular $\nrm{\nh^3v}_{L^2}+\nrm{\pa^\beta w}_{L^2}\lesssim\Fn$ for $|\beta|\le2$.
\item (Velocity factors bounded by $\Fu$.) $\nrm{\pa^\beta v}_{L^\infty}\lesssim\Fu$ for $|\beta|\le3$; $\nrm{\pa^\beta w}_{L^\infty}\lesssim\Fu$ if $|\beta|\le2$, or if $|\beta|=3$ and $\beta_3\ge1$; $\nrm{\pa^\beta w}_{L^3}\lesssim\Fu$ for $|\beta|=3$; and $\nrm{\pa^\beta v}_{L^3}+\nrm{\pa^\beta w}_{L^2}\lesssim\Fu$ for $|\beta|\le4$.
\item (Density factors.) $\nrm{\nh\varrho}_{H^2}\lesssim\Fv$. Consequently $\nrm{\nh\varrho}_{L^\infty}+\nrm{\nabla\nh\varrho}_{L^6}\lesssim\Fv$, and $\nrm{\pa^\kappa\varrho}_{L^2}+\nrm{\pa^\kappa(\sgm\varrho)}_{L^2}\lesssim\Fv$ for every $\kappa$ with $|\kappa|\le3$ and $\kappa_h\ne0$. Moreover, $\nrm{\phi}_{W^{2,\infty}}+\nrm{\nabla^3\phi}_{L^6}+\nrm{\phi}_{H^4}\lesssim\Ev$ for $\phi\in\{\varrho,\sgm\varrho\}$.
\item (The force.) $\nrm{G}_{H^1}\lesssim\Ff+\Fsh^2$ and $\nrm{\nh^3G}_{L^2}\lesssim\Fsh$.
\end{enumerate}
\end{lemma}

\begin{proof}
(a) Since $v=0$ on $\pa\Omega$ and $w=0$ at $z=0$, and $\Omega$ has height one, $|v(x_h,z)|\le\nrm{\pz v}_{L^\infty}$ and $|w(x_h,z)|\le\nrm{\pz w}_{L^\infty}$; the bound by $\Fn$ is the definition \eqref{eq:Fnat}. By \Cref{lem:GN}(1), $\nrm{\nh v}_{H^2}\lesssim\nrm{\nh v}_{H^1}^{2/3}\nrm{\nh v}_{H^4}^{1/3}$, where $\nrm{\nh v}_{H^1}\lesssim\nrm{\nabla\nh v}\lesssim\Ff$ (vertical Poincar\'e inequality, $\nh v|_{\pa\Omega}=0$) and $\nrm{\nh v}_{H^4}\le\Fu\le\Fsh$. Since $\Ff\lesssim\Fsh$ by \eqref{eq:3.11}, we have $\Ff^a\Fsh^{1-a}\lesssim\Ff^{7/12}\Fsh^{5/12}\le\Fn$ for $a\ge\frac7{12}$; hence $\nrm{\nh v}_{H^2}\lesssim\Ff^{2/3}\Fsh^{1/3}\lesssim\Fn$. By H\"older's inequality and $H^1\subset L^6$, $\nrm{\nabla\nh v}_{L^3}\le\nrm{\nabla\nh v}^{1/2}\nrm{\nabla\nh v}_{L^6}^{1/2}\lesssim\Ff^{1/2}\nrm{\nh v}_{H^2}^{1/2}\lesssim\Ff^{5/6}\Fsh^{1/6}\lesssim\Fn$. If $\beta_3\ge1$, then $\pa^\beta w=-\pa^{\beta-e_3}\divh v$; if $\beta_3=0$, \Cref{lem:w} gives $\nrm{\pa^\beta w}\le\nrm{\pa^\beta\divh v}$. In both cases $\nrm{\pa^\beta w}\lesssim\nrm{\nh v}_{H^{|\beta|}}$.

(b) By $H^2\subset L^\infty$, $\nrm{\pa^\beta v}_{L^\infty}\lesssim\nrm{v}_{H^{|\beta|+2}}\le\Fu$ for $|\beta|\le3$, and the same holds for $\pa^\beta w=-\pa^{\beta-e_3}\divh v$ when $\beta_3\ge1$, $|\beta|\le3$. For horizontal $\beta$ with $j=|\beta|\le2$, \Cref{lem:w} with $s=1$ gives $\nrm{\nh^jw}_{L^\infty}\lesssim\nrm{\nh^j\divh v}_{L^2_zH^2_h}\lesssim\nrm{v}_{H^{j+3}}\le\Fu$. For horizontal $\beta$ with $|\beta|=3$, \Cref{lem:w} and \Cref{lem:aniso}(b) give
\begin{equation*}
\nrm{\nh^3w}_{L^3}\le\nrm{\nh^3w}_{L^3_hL^\infty_z}\le\nrm{\nh^3\divh v}_{L^2_zL^3_h}\lesssim\nrm{\nh^3\divh v}_{L^2_zH^1_h}\le\Fu .
\end{equation*}
Finally $\nrm{\pa^\beta v}_{L^3}\lesssim\nrm{\pa^\beta v}_{H^1}\le\Fu$ for $|\beta|\le4$, and $\nrm{\pa^\beta w}\lesssim\nrm{\nh v}_{H^4}\le\Fu$ by (a).

(c) The components of $\nabla\nabla\nh\varrho$ are $\nh^3\varrho$, $\nh^2\pz\varrho$ and $\nh\pz^2\varrho$, which are bounded by $\Fv$ by \eqref{eq:energies}. The functions $\nh\varrho$, $\nh^2\varrho$, $\nh\pz\varrho$ have zero horizontal mean, so \Cref{lem:poincare} bounds them by $l$ times $\nh^2\varrho$, $\nh^3\varrho$, $\nh^2\pz\varrho$, respectively. Hence $\nrm{\nh\varrho}_{H^2}\lesssim\Fv$, and \Cref{lem:GN}(3) gives the $L^\infty$ and $L^6$ bounds. If $|\kappa|\le3$ and $\kappa_h\ne0$, then $\pa^\kappa(\sgm\varrho)=\sum_{m\le\kappa_3}\binom{\kappa_3}m\sgm^{(m)}\pa^{\kappa-me_3}\varrho$, and each $\pa^{\kappa-me_3}\varrho$ is a component of $\nh\varrho$, $\nabla\nh\varrho$ or $\nabla\nabla\nh\varrho$. The last assertion follows from \eqref{eq:3.12}, \Cref{lem:GN}(3) and $\sgm\in C^4$.

(d) By \eqref{eq:3.12}, $\nrm{\rho}_{W^{1,\infty}}\lesssim1+\Ev\lesssim1$, so $\nrm{\rho\dot v}_{H^1}\lesssim\nrm{v_t}_{H^1}+\nrm{\bU\cdot\nabla v}_{H^1}$, and $\nrm{v_t}_{H^1}\lesssim\nrm{\nabla v_t}\lesssim\Ff$. Using $H^1\subset L^4\cap L^6$, $H^2\subset L^\infty$, $\nrm{w}_{H^2}\lesssim\nrm{v}_{H^3}$ (\Cref{lem:w}) and $\nrm{w}_{L^\infty}\lesssim\nrm{v}_{H^3}$ (\eqref{eq:3.13}),
\begin{equation*}
\begin{aligned}
\nrm{v\cdot\nh v}_{H^1}&\lesssim\nrm{v}_{L^\infty}\nrm{\nh v}_{H^1}+\nrm{\nabla v}_{L^4}\nrm{\nh v}_{L^4}\lesssim\nrm{v}_{H^2}^2,\\
\nrm{w\pz v}_{H^1}&\lesssim\nrm{w}_{L^\infty}\nrm{\pz v}_{H^1}+\nrm{\nabla w}_{L^3}\nrm{\pz v}_{L^6}\lesssim\nrm{v}_{H^3}\nrm{v}_{H^2},
\end{aligned}
\end{equation*}
so $\nrm{\bU\cdot\nabla v}_{H^1}\lesssim\nrm{v}_{H^2}\nrm{v}_{H^3}\lesssim\Fsh^2$ by \eqref{eq:4.12}. The linear Coriolis part is of lower order: $|f|\nrm{\rho v^\perp}_{H^1}\lesssim\nrm{v}_{H^1}\lesssim\Ff$ by \eqref{eq:3.13b}. This proves the first bound. For the second, $\nrm{\nh^3G}\le\nrm{\rho\dot v}_{H^3}+|f|\nrm{\nh^3(\rho v)}$, because $\nh^3(\rho v^\perp)=(\nh^3(\rho v))^\perp$. By the second line of \eqref{eq:6.26}, by \eqref{eq:T4-GH3}, and by $D_4,\Fuz,\nrm{v}_{H^2}\le\Fu$ (\eqref{eq:3.3}, \eqref{eq:3.8}),
\begin{equation*}
\nrm{\nh^3G}\lesssim\Fut+\Esh\big(D_4+\Fv+\nrm{v}_{H^2}\big)+D_4+\Esh\Fuz\lesssim\Fsh .\qedhere
\end{equation*}
\end{proof}

We now estimate $\mathcal P_4,\dots,\mathcal P_0$.

\begin{enumerate}[label=(\roman*),leftmargin=*]
\item \emph{The case $|\gamma|=4$.} Then $\beta=\delta=0$ and $\gamma=\alpha'$ (respectively $\gamma=\alpha$ in $\cK$). For $|\alpha'|=4$, the Leibniz rule in $z$ gives $\sgm\pa^{\alpha'}\varrho=\pa^{\alpha'}(\sgm\varrho)-\sum_{m=1}^{\alpha'_3}\binom{\alpha'_3}m\sgm^{(m)}\pa^{\alpha'-me_3}\varrho$, and since $\bU\cdot\nabla\sgm^{(m)}=\sgm^{(m+1)}w$,
\begin{equation*}
\begin{aligned}
\sgm\bU\cdot\nabla\pa^{\alpha'}\varrho&=\bU\cdot\nabla\pa^{\alpha'}(\sgm\varrho)-\mathfrak r_{\alpha'},\\
\mathfrak r_{\alpha'}&\eqdef\sgm'w\,\pa^{\alpha'}\varrho+\sum_{m=1}^{\alpha'_3}\binom{\alpha'_3}m\Big(\sgm^{(m+1)}w\,\pa^{\alpha'-me_3}\varrho+\sgm^{(m)}\bU\cdot\nabla\pa^{\alpha'-me_3}\varrho\Big).
\end{aligned}
\end{equation*}
Every summand of $\mathfrak r_{\alpha'}$ except $\sgm^{(5)}w\varrho$, which occurs only for $\alpha'=(0,0,4)$ and $m=4$, is bounded in $L^2$ by $C\nrm{\bU}_{L^\infty}\nrm{\varrho}_{H^4}$. For the exceptional summand, \Cref{lem:aniso}(a) gives
\begin{equation*}
\nrm{\sgm^{(5)}w\varrho}_{L^2}\le\nrm{w}_{L^\infty}\nrm{\sgm^{(5)}}_{L^2(0,1)}\sup_{z\in(0,1)}\nrm{\varrho(\cdot,z)}_{L^2(\Tl)}\lesssim\nrm{\sgm^{(5)}}_{L^2(0,1)}\nrm{w}_{L^\infty}\Ev .
\end{equation*}
This is the only place where five derivatives of $\sgm$, that is, six derivatives of $\rs$, enter; it is the reason for the assumption $\rs\in C^6$ in \eqref{eq:stable} (see \Cref{rem:C6}). By \Cref{lem:NT-aux}(a), $\nrm{\mathfrak r_{\alpha'}}_{L^2}\lesssim\nrm{\bU}_{L^\infty}\Ev\le\Ev\Fn$. Therefore
\begin{equation*}
\cJ_{\alpha,\alpha'}^{0,\alpha',0}=-\int_\Omega\pa^\alpha(\sgm\varrho)\,\bU\cdot\nabla\pa^{\alpha'}(\sgm\varrho)\dx+\int_\Omega\pa^\alpha(\sgm\varrho)\,\mathfrak r_{\alpha'}\dx .
\end{equation*}
Summing over $(\alpha,\alpha')\in\cM_1$ (all coefficients are equal to one), the first integrals add up to $-\int_\Omega\AB(\sgm\varrho)\,\bU\cdot\nabla\AB(\sgm\varrho)\dx=0$; for $(\alpha,\alpha)\in\cM_3$ each first integral vanishes by itself. In both cases we integrate by parts, using $\divv\bU=0$ and $\bU\cdot n=0$ on $\pa\Omega$. For $\alpha\in\cM_2$ the same integration by parts gives
\begin{equation*}
\cK_\alpha^{0,\alpha}=-\int_\Omega\sgm\,\pa^\alpha\varrho\,\bU\cdot\nabla\pa^\alpha\varrho\dx=\frac12\int_\Omega\sgm'w\,|\pa^\alpha\varrho|^2\dx .
\end{equation*}
Since $\nrm{\pa^\alpha(\sgm\varrho)}\lesssim\Ev$ and $\nrm{w}_{L^\infty}\le\Fn$, we conclude
\begin{equation}\label{eq:6.75}
|\mathcal P_4|\lesssim\Ev\max_{|\alpha'|=4}\nrm{\mathfrak r_{\alpha'}}_{L^2}+\nrm{w}_{L^\infty}\Ev^2\lesssim\Esh^2\Fn .
\end{equation}

\item \emph{The case $|\gamma|=3$.} Then $|\beta|+|\delta|=1$, so $|\beta|\le1$ and, by \Cref{lem:NT-aux}(a), $\nrm{\pa^\beta\bU}_{L^\infty}\le\Fn$. Since $\nabla\pa^\gamma\varrho$ has order four, each $\cJ$ and $\cK$ with $|\gamma|=3$ is bounded by $C\Ev\nrm{\pa^\beta\bU}_{L^\infty}\nrm{\nabla\pa^\gamma\varrho}$, and
\begin{equation}\label{eq:6.78}
|\mathcal P_3|\lesssim\Esh^2\Fn .
\end{equation}
No representation of $w$ is needed here: $\nrm{\nabla\bU}_{L^\infty}$, including $\nrm{\nh w}_{L^\infty}$, is part of $\Fn$, whose time integral is bounded by \Cref{prop:Fnat}.

\item \emph{The case $|\gamma|=2$.} Then $|\beta|+|\delta|=2$. If $|\beta|\le1$, the corresponding terms are bounded by $\Esh^2\Fn$ exactly as in case (ii). Let $|\beta|=2$, so that $\delta=0$, and recall from \Cref{lem:NT-aux}(b) that $\nrm{\pa^\beta\bU}_{L^\infty}\lesssim\Fu$.
\begin{enumerate}[label=(\alph*),leftmargin=*]
\item If $\gamma_h\ne0$, every component of $\nabla\pa^\gamma\varrho$ is a derivative $\pa^\kappa\varrho$ with $|\kappa|=3$ and $\kappa_h\ne0$, so $\nrm{\nabla\pa^\gamma\varrho}\lesssim\Fv$ by \Cref{lem:NT-aux}(c), and
\begin{equation*}
|\cJ_{\alpha,\alpha'}^{\beta,\gamma,0}|+|\cK_\alpha^{\beta,\gamma}|\lesssim\Ev\nrm{\pa^\beta\bU}_{L^\infty}\nrm{\nabla\pa^\gamma\varrho}\lesssim\Ev\Fu\Fv\le\Esh\Fsh^2 .
\end{equation*}
\item If $\gamma=2e_3$ and $\beta_3\ge1$ (this occurs only in $\cJ$, since $\gamma$ is horizontal in $\cK$), then $\pa^\beta w=-\pa^{\beta-e_3}\divh v$ with $|\beta-e_3|=1$, and \Cref{lem:NT-aux}(a), (c) give
\begin{equation*}
\begin{aligned}
\nrm{\pa^\beta\bU\cdot\nabla\pz^2\varrho}&\le\nrm{\pa^\beta v}_{L^\infty}\nrm{\nh\pz^2\varrho}+\nrm{\pa^\beta w}_{L^3}\nrm{\pz^3\varrho}_{L^6}\\
&\lesssim\Fu\Fv+\nrm{\nabla\nh v}_{L^3}\Ev\lesssim\Fu\Fv+\Fn\Ev .
\end{aligned}
\end{equation*}
Hence $|\cJ_{\alpha,\alpha'}^{\beta,\gamma,0}|\lesssim\Esh\Fsh^2+\Esh^2\Fn$.
\item If $\gamma=2e_3$ and $\beta=(\beta_1,\beta_2,0)$, then $\alpha'=\beta+\gamma=(\beta_1,\beta_2,2)$. For $(\alpha,\alpha')\in\cM_1$ this forces $\beta_h\in\{(2,0),(0,2)\}$ and $\alpha\in\mathcal S_1$; for $(\alpha,\alpha')\in\cM_3$ it forces $\alpha=\alpha'$. Thus either $\alpha=(\alpha_1,\alpha_2,2)$ with $\alpha_1+\alpha_2=2$, or $\alpha=(0,0,4)$. We write
\begin{equation*}
\cJ_{\alpha,\alpha'}^{\beta,\gamma,0}=-\int_\Omega\pa^\alpha(\sgm\varrho)\big(\pa^\beta v\cdot\nh\pz^2\varrho+\pa^\beta w\,\pz^3\varrho\big)\sgm\dx\eqdef-(I_1+I_2).
\end{equation*}
As in (b), $|I_1|\lesssim\Ev\Fu\Fv$. For $I_2$, recall that $\pa^\beta w$ is a second horizontal derivative of $w$; it vanishes on $\pa\Omega$, and $\pz\pa^\beta w=-\pa^\beta\divh v$. If $\alpha=(\alpha_1,\alpha_2,2)$, an integration by parts in $z$ gives
\begin{equation*}
I_2=\int_\Omega\pa^{\alpha_h}\pz(\sgm\varrho)\Big(\pa^\beta\divh v\,\pz^3\varrho\,\sgm-\pa^\beta w\,\pz^4\varrho\,\sgm-\pa^\beta w\,\pz^3\varrho\,\sgm'\Big)\dx,
\end{equation*}
and, by \Cref{lem:NT-aux}(b), (c) ($\pa^\beta\divh v$ is a third derivative of $v$),
\begin{equation}\label{eq:6.83}
|I_2|\lesssim\nrm{\pa^{\alpha_h}\pz(\sgm\varrho)}\big(\nrm{\nh^3v}_{L^\infty}+\nrm{\pa^\beta w}_{L^\infty}\big)\nrm{\varrho}_{H^4}\lesssim\Fv\Fu\Ev\le\Esh\Fsh^2 .
\end{equation}
If $\alpha=(0,0,4)$, we write $\sgm\pz^3\varrho=\pz^3(\sgm\varrho)-\phi_1$ with $\phi_1\eqdef3\sgm'\pz^2\varrho+3\sgm''\pz\varrho+\sgm'''\varrho$. Since $\pa^\beta w=0$ on $\pa\Omega$,
\begin{equation*}
\int_\Omega\pz^4(\sgm\varrho)\,\pa^\beta w\,\pz^3(\sgm\varrho)\dx=\frac12\int_\Omega\pz\big|\pz^3(\sgm\varrho)\big|^2\pa^\beta w\dx=\frac12\int_\Omega\big|\pz^3(\sgm\varrho)\big|^2\pa^\beta\divh v\dx .
\end{equation*}
By $H^1\subset L^4$, $\nrm{\pz^3(\sgm\varrho)}_{L^4}^2\lesssim\Ev^2$, and $\nrm{\pa^\beta\divh v}\le\sqrt2\nrm{\nh^3v}\lesssim\Fn$ by \Cref{lem:NT-aux}(a). For the remainder, $\nrm{\phi_1}_{L^\infty}\lesssim\Ev$ and $\nrm{\pa^\beta w}\le\sqrt2\nrm{\nh^3v}\lesssim\Fn$. Hence
\begin{equation}\label{eq:6.84}
|I_2|\le\frac12\nrm{\pz^3(\sgm\varrho)}_{L^4}^2\nrm{\pa^\beta\divh v}+\nrm{\pz^4(\sgm\varrho)}\nrm{\pa^\beta w}\nrm{\phi_1}_{L^\infty}\lesssim\Esh^2\Fn .
\end{equation}
\end{enumerate}
Gathering the subcases,
\begin{equation}\label{eq:6.85}
|\mathcal P_2|\lesssim\Esh\Fsh^2+\Esh^2\Fn .
\end{equation}

\item \emph{The case $|\gamma|=1$.} Then $|\beta|+|\delta|=3$. If $|\beta|\le1$, the terms are bounded by $\Esh^2\Fn$ as in case (ii).

If $|\beta|=2$, then $|\delta|=1$ (so the term is a $\cJ$). For $\gamma=e_i$, $i\in\{1,2\}$, \Cref{lem:NT-aux}(b), (c) give $\nrm{\pa^\beta\bU\cdot\nabla\pa_i\varrho}\le\nrm{\pa^\beta\bU}_{L^\infty}\nrm{\nabla\nh\varrho}\lesssim\Fu\Fv$. For $\gamma=e_3$,
\begin{equation*}
\nrm{\pa^\beta\bU\cdot\nabla\pz\varrho}\le\nrm{\pa^\beta v}_{L^\infty}\nrm{\nh\pz\varrho}+\nrm{\pa^\beta w}\nrm{\pz^2\varrho}_{L^\infty}\lesssim\Fu\Fv+\Fn\Ev
\end{equation*}
by \Cref{lem:NT-aux}(a)--(c). In both cases $|\cJ_{\alpha,\alpha'}^{\beta,\gamma,\delta}|\lesssim\Esh\Fsh^2+\Esh^2\Fn$.

Let now $|\beta|=3$, so that $\delta=0$.
\begin{enumerate}[label=(\alph*),leftmargin=*]
\item If $\gamma=e_i$ with $i\in\{1,2\}$, then by \Cref{lem:NT-aux}(b), (c),
\begin{equation*}
\nrm{\pa^\beta\bU\cdot\nabla\pa_i\varrho}\le\nrm{\pa^\beta v}_{L^\infty}\nrm{\nh\pa_i\varrho}+\nrm{\pa^\beta w}_{L^3}\nrm{\pz\pa_i\varrho}_{L^6}\lesssim\Fu\Fv,
\end{equation*}
and therefore $|\cJ_{\alpha,\alpha'}^{\beta,\gamma,0}|+|\cK_\alpha^{\beta,\gamma}|\lesssim\Ev\Fu\Fv$. (The density factor cannot be placed in $L^\infty_hL^2_z$ in three dimensions, since this would require $\pz\nh^{2+s}\varrho$, which is not controlled by $\Fv$; instead $\nabla\nh\varrho\in L^6$ is paired with $\pa^\beta w\in L^3$.)
\item If $\gamma=e_3$ (only in $\cJ$), then $\alpha'=\beta+e_3$. The pairs in $\cM_1\cup\cM_3$ with $\alpha'_3\ge1$ have $\alpha=(0,0,4)$ or $\alpha_3\in\{1,2\}$ with $\alpha_1+\alpha_2=4-\alpha_3\ge2$.

If $\alpha_3\in\{1,2\}$, then $\pa^{\alpha-e_3}(\sgm\varrho)$ is a combination of $\sgm\pa^{\alpha-e_3}\varrho$ and $\sgm'\pa^{\alpha_h}\varrho$, which vanish on $\pa\Omega$ because $\alpha_h\ne0$ and $\varrho=\pz\varrho=0$ there. An integration by parts in $z$ gives
\begin{equation*}
\cJ_{\alpha,\alpha'}^{\beta,e_3,0}=\int_\Omega\pa^{\alpha-e_3}(\sgm\varrho)\Big(\big(\pa^{\beta+e_3}\bU\cdot\nabla\pz\varrho+\pa^\beta\bU\cdot\nabla\pz^2\varrho\big)\sgm+\pa^\beta\bU\cdot\nabla\pz\varrho\,\sgm'\Big)\dx .
\end{equation*}
Here $\nrm{\pa^{\alpha-e_3}(\sgm\varrho)}\lesssim\Fv$ by \Cref{lem:NT-aux}(c). Since $\pa^{\beta+e_3}w=-\pa^\beta\divh v$ is a fourth derivative of $v$, and by \Cref{lem:NT-aux}(b), (c),
\begin{equation*}
\begin{aligned}
&\nrm{\pa^{\beta+e_3}\bU\cdot\nabla\pz\varrho}+\nrm{\pa^\beta\bU\cdot\nabla\pz\varrho}\lesssim\nrm{v}_{H^5}\nrm{\nabla\pz\varrho}_{L^\infty}\lesssim\Fu\Ev,\\
&\nrm{\pa^\beta\bU\cdot\nabla\pz^2\varrho}\le\nrm{\pa^\beta v}_{L^\infty}\nrm{\nh\pz^2\varrho}+\nrm{\pa^\beta w}_{L^3}\nrm{\pz^3\varrho}_{L^6}\lesssim\Fu\Ev .
\end{aligned}
\end{equation*}
Hence $|\cJ_{\alpha,\alpha'}^{\beta,e_3,0}|\lesssim\Fv\Fu\Ev$.

If $\alpha=(0,0,4)$, then $(\alpha,\alpha')\in\cM_1$, so $\alpha'\in\mathcal S_1$ and $\beta=\alpha'-e_3$ is one of $(2,0,1)$, $(0,2,1)$, $(0,0,3)$. We write
\begin{equation*}
\cJ_{\alpha,\alpha'}^{\beta,e_3,0}=-\int_\Omega\pz^4(\sgm\varrho)\big(\pa^\beta v\cdot\nh\pz\varrho+\pa^\beta w\,\pz^2\varrho\big)\sgm\dx\eqdef-(I_3+I_4).
\end{equation*}
By \Cref{lem:NT-aux}(b), (c), $|I_3|\lesssim\Ev\nrm{\pa^\beta v}_{L^\infty}\nrm{\nh\pz\varrho}\lesssim\Ev\Fu\Fv$. Since $\pa^\beta w=-\pa^{\beta-e_3}\divh v$ with $|\beta-e_3|=2$, \Cref{lem:NT-aux}(a) gives $\nrm{\pa^\beta w}\lesssim\nrm{\nh v}_{H^2}\lesssim\Fn$, and $|I_4|\lesssim\Ev\nrm{\pa^\beta w}\nrm{\pz^2\varrho}_{L^\infty}\lesssim\Esh^2\Fn$.
\end{enumerate}
Gathering all subcases,
\begin{equation}\label{eq:6.87}
|\mathcal P_1|\lesssim\Esh\Fsh^2+\Esh^2\Fn .
\end{equation}

\item \emph{The case $|\gamma|=0$.} Due to the derivative loss, this case is the most anisotropic. We split the summands into their horizontal and vertical parts:
\begin{equation*}
\begin{aligned}
\cJ_{\alpha,\alpha'}^{\beta,0,\delta}&=-\int_\Omega\pa^\alpha(\sgm\varrho)\big(\pa^\beta v\cdot\nh\varrho+\pa^\beta w\,\pz\varrho\big)\pa^\delta\sgm\dx\eqdef-(J_1+J_2),\\
\cK_\alpha^{\alpha,0}&=-\int_\Omega\sgm\,\pa^\alpha\varrho\big(\pa^\alpha v\cdot\nh\varrho+\pa^\alpha w\,\pz\varrho\big)\dx\eqdef-(K_1+K_2).
\end{aligned}
\end{equation*}
Since $|\beta|\le4$, \Cref{lem:NT-aux}(b), (c) give
\begin{equation}\label{eq:6.88}
|J_1|+|K_1|\lesssim\Ev\big(\nrm{\pa^\beta v}+\nrm{\pa^\alpha v}\big)\nrm{\nh\varrho}_{L^\infty}\lesssim\Ev\Fu\Fv .
\end{equation}
For $J_2$ we distinguish $\alpha=(0,0,4)$, $\alpha_3=2$ and $\alpha_3=1$; the horizontal case $\alpha\in\cM_2$ is the term $K_2$, treated last. Recall that, since $\gamma=0$, $\beta+\delta=\alpha'$, so $\beta\le\alpha'$.
\begin{enumerate}[label=(\alph*),leftmargin=*]
\item \emph{The case $\alpha=(0,0,4)$.} Then $\alpha'\in\mathcal S_1$, so $\beta_h\le(2,0)$ or $\beta_h\le(0,2)$, and $\beta_3=0$ implies $|\beta|\le2$.

If $|\beta|\le3$, then $\pa^\beta w=-\pa^{\beta-e_3}\divh v$ with $|\beta-e_3|\le2$ when $\beta_3\ge1$, while $|\beta|\le2$ when $\beta_3=0$. By \Cref{lem:NT-aux}(a), $\nrm{\pa^\beta w}\lesssim\nrm{\nh v}_{H^2}\lesssim\Fn$, and
\begin{equation}\label{eq:6.90}
|J_2|\lesssim\nrm{\pz^4(\sgm\varrho)}\nrm{\pa^\beta w}\nrm{\pz\varrho}_{L^\infty}\lesssim\Esh^2\Fn .
\end{equation}
If $|\beta|=4$, then $\beta=\alpha'\in\mathcal S_1$, $\delta=0$, and $\pa^\beta w=-\pa^{\beta-e_3}\divh v$. Since $\pz\varrho=0$ on $\pa\Omega$, an integration by parts in $z$ gives
\begin{equation}\label{eq:6.91}
\begin{aligned}
J_2&=-\int_\Omega\pz^4(\sgm\varrho)\,\pa^{\beta-e_3}\divh v\,\pz\varrho\,\sgm\dx\\
&=\int_\Omega\pz^3(\sgm\varrho)\Big(\pa^\beta\divh v\,\pz\varrho\,\sgm+\pa^{\beta-e_3}\divh v\,\pz^2\varrho\,\sgm+\pa^{\beta-e_3}\divh v\,\pz\varrho\,\sgm'\Big)\dx\\
&\eqdef J_{21}+J_{22}+J_{23}.
\end{aligned}
\end{equation}
For $J_{21}$ we integrate by parts horizontally, $\pa^\beta\divh v=\sum_j\pa_j\pa^\beta v_j$:
\begin{equation*}
J_{21}=-\sum_{j=1}^2\int_\Omega\pa^\beta v_j\,\pa_j\pz^3(\sgm\varrho)\,\pz\varrho\,\sgm\dx-\sum_{j=1}^2\int_\Omega\pa^\beta v_j\,\pz^3(\sgm\varrho)\,\pa_j\pz\varrho\,\sgm\dx .
\end{equation*}
In the first sum we integrate by parts once more in $z$ (the boundary term vanishes since $\pz\varrho|_{\pa\Omega}=0$); it becomes $\sum_j\int_\Omega\pa_j\pz^2(\sgm\varrho)\,\pz\big(\pa^\beta v_j\,\pz\varrho\,\sgm\big)\dx$. By \Cref{lem:NT-aux}(b), (c), with $\nrm{\pa_j\pz^2(\sgm\varrho)}\lesssim\Fv$, $\nrm{\pa^\beta v}_{H^1}\le\Fu$, $\nrm{\pa^\beta v}_{L^3}\lesssim\Fu$ and $\nrm{\nh\pz\varrho}_{L^6}\lesssim\Fv$,
\begin{equation}\label{eq:6.92}
|J_{21}|\lesssim\Fv\Fu\nrm{\varrho}_{W^{2,\infty}}+\nrm{\pz^3(\sgm\varrho)}\nrm{\pa^\beta v}_{L^3}\nrm{\nh\pz\varrho}_{L^6}\lesssim\Ev\Fu\Fv .
\end{equation}
For $J_{23}$, since $\beta_3\ge2$, we write $\pa^{\beta-e_3}\divh v=\pz\pa^{\beta-2e_3}\divh v$ and integrate by parts in $z$ (again $\pz\varrho|_{\pa\Omega}=0$); as $|\beta-2e_3|=2$, \Cref{lem:NT-aux}(a) gives
\begin{equation}\label{eq:6.99}
|J_{23}|\le\nrm{\pa^{\beta-2e_3}\divh v}\,\nrm{\pz\big(\pz^3(\sgm\varrho)\pz\varrho\,\sgm'\big)}\lesssim\nrm{\nh v}_{H^2}\Ev^2\lesssim\Esh^2\Fn .
\end{equation}
For $J_{22}$ there are two values of $\beta$. If $\beta=(2,0,2)$ or $(0,2,2)$, then $\pa^{\beta-e_3}\divh v=\pz\pa^{\beta_h}\divh v$, and $\pa^{\beta_h}\divh v$ vanishes on $\pa\Omega$. An integration by parts in $z$, $H^1\subset L^4$ and $H^2\subset L^\infty$ give
\begin{equation}\label{eq:6.93}
|J_{22}|=\left|\int_\Omega\pa^{\beta_h}\divh v\;\pz\big(\pz^3(\sgm\varrho)\,\pz^2\varrho\,\sgm\big)\dx\right|\lesssim\nrm{\nh^3v}\,\nrm{\varrho}_{H^4}^2\lesssim\Esh^2\Fn .
\end{equation}
If $\beta=(0,0,4)$, the factor $\pz^3\divh v$ is rewritten by the momentum equation. Taking $\divh$ of \eqref{eq:3.22}, $g\Laph\varrho=-\mu\Lap\pz\divh v+\pz\divh G$, so that
\begin{equation*}
\pz^3\divh v=\frac1\mu\pz\divh G-\frac g\mu\Laph\varrho-\Laph\pz\divh v,
\end{equation*}
and accordingly
\begin{equation*}
\begin{aligned}
J_{22}&=\frac1\mu\int_\Omega\pz^3(\sgm\varrho)\,\pz\divh G\,\pz^2\varrho\,\sgm\dx-\frac g\mu\int_\Omega\pz^3(\sgm\varrho)\,\Laph\varrho\,\pz^2\varrho\,\sgm\dx\\
&\quad-\int_\Omega\pz^3(\sgm\varrho)\,\Laph\pz\divh v\,\pz^2\varrho\,\sgm\dx\eqdef J_{221}+J_{222}+J_{223}.
\end{aligned}
\end{equation*}
For $J_{221}$ we integrate by parts horizontally, which moves the derivative off $G$:
\begin{equation*}
J_{221}=-\frac1\mu\int_\Omega\pz G\cdot\Big(\nh\pz^3(\sgm\varrho)\,\pz^2\varrho+\pz^3(\sgm\varrho)\,\nh\pz^2\varrho\Big)\sgm\dx .
\end{equation*}
By \Cref{lem:NT-aux}(c), (d) and $H^1\subset L^4$,
\begin{equation}\label{eq:6.94}
|J_{221}|\lesssim\nrm{\pz G}\Big(\nrm{\varrho}_{H^4}\nrm{\pz^2\varrho}_{L^\infty}+\nrm{\pz^3(\sgm\varrho)}_{L^4}\nrm{\nh\pz^2\varrho}_{L^4}\Big)\lesssim\Ev^2\nrm{G}_{H^1}\lesssim\Esh^2\big(\Ff+\Fsh^2\big).
\end{equation}
The force $G$ contains the linear Coriolis term $f\rs v^\perp$; it enters only through $\nrm{G}_{H^1}$, and its contribution $\Ev^2|f|\nrm{v}_{H^1}\lesssim\Esh^2\Ff$ is integrable in time by \eqref{eq:5.6}.

For $J_{222}$ we write $\sgm\pz^2\varrho=\pz^2(\sgm\varrho)-\phi_2$ with $\phi_2\eqdef2\sgm'\pz\varrho+\sgm''\varrho$, so that
\begin{equation}\label{eq:6.95}
J_{222}=\frac g\mu\big(J_{2221}+J_{2222}\big)\eqdef\frac g\mu\Big(-\int_\Omega\pz^3(\sgm\varrho)\,\Laph\varrho\,\pz^2(\sgm\varrho)\dx+\int_\Omega\pz^3(\sgm\varrho)\,\Laph\varrho\,\phi_2\dx\Big).
\end{equation}
Set $h\eqdef\pz^2(\sgm\varrho)$. Writing the product in $J_{2221}$ as the derivative of a square and integrating by parts in $z$ (the boundary term vanishes because $\Laph\varrho|_{\pa\Omega}=0$),
\begin{equation*}
J_{2221}=-\frac12\int_\Omega\Laph\varrho\,\pz|h|^2\dx=\frac12\int_\Omega\Laph\pz\varrho\,|h|^2\dx .
\end{equation*}
We split $|h|^2=\bar h^2+2\bar h\wt h+\wt h^2$. The term with $\bar h^2$ vanishes, since $\bar h$ depends only on $z$ and $\Laph\pz\varrho$ has zero horizontal mean. In the other two terms we use $\nrm{\bar h}_{L^\infty}+\nrm{\wt h}_{L^\infty}\le3\nrm{h}_{L^\infty}\lesssim\nrm{h}_{H^2}\lesssim\Ev$ (here $H^2(\Omega)\subset L^\infty(\Omega)$ is used) and the horizontal Poincar\'e inequality $\nrm{\wt h}\le l\nrm{\nh h}=l\nrm{\nh\pz^2(\sgm\varrho)}\lesssim\Fv$:
\begin{equation*}
|J_{2221}|\le\nrm{\Laph\pz\varrho}\big(\nrm{\bar h}_{L^\infty}+\nrm{\wt h}_{L^\infty}\big)\nrm{\wt h}\lesssim\Ev\Fv^2 .
\end{equation*}
For the $\sgm$-remainder $J_{2222}$ one has to avoid $\nrm{\nh\pz^3\varrho}$, which is not controlled by $\Fv$. We split $\phi_2=\bar\phi_2+\wt\phi_2$, where $\wt\phi_2=2\sgm'\wt{\pz\varrho}+\sgm''\wt\varrho$. By the horizontal Poincar\'e inequality, $\nrm{\wt\phi_2}_{H^1}\lesssim\nrm{\nh\varrho}_{H^2}\lesssim\Fv$, so $H^1\subset L^4$ and $\nrm{\Laph\varrho}_{L^4}\lesssim\nrm{\nh\varrho}_{H^2}\lesssim\Fv$ give
\begin{equation*}
\left|\int_\Omega\pz^3(\sgm\varrho)\,\Laph\varrho\,\wt\phi_2\dx\right|\le\nrm{\pz^3(\sgm\varrho)}\nrm{\Laph\varrho}_{L^4}\nrm{\wt\phi_2}_{L^4}\lesssim\Ev\Fv^2 .
\end{equation*}
For the mean part we integrate by parts in $z$ (the boundary term vanishes since $\Laph\varrho|_{\pa\Omega}=0$); since $\pz(\Laph\varrho\,\bar\phi_2)$ has zero horizontal mean, $\ol{\pz^2(\sgm\varrho)}$ drops out:
\begin{equation*}
\int_\Omega\pz^3(\sgm\varrho)\,\Laph\varrho\,\bar\phi_2\dx=-\int_\Omega\wt{\pz^2(\sgm\varrho)}\,\pz\big(\Laph\varrho\,\bar\phi_2\big)\dx .
\end{equation*}
As $\nrm{\wt{\pz^2(\sgm\varrho)}}\le l\nrm{\nh\pz^2(\sgm\varrho)}\lesssim\Fv$ and $\nrm{\pz(\Laph\varrho\,\bar\phi_2)}\le\nrm{\Laph\pz\varrho}\nrm{\bar\phi_2}_{L^\infty}+\nrm{\Laph\varrho}\nrm{\pz\bar\phi_2}_{L^\infty}\lesssim\Fv\Ev$, this term is also $\lesssim\Ev\Fv^2$. Hence
\begin{equation}\label{eq:6.96}
|J_{222}|\lesssim\Ev\Fv^2 .
\end{equation}
Finally, $\Laph\pz\divh v=\sum_{i=1}^2\pz\pa_i^2\divh v$, so $J_{223}$ is minus the sum of the two integrals $J_{22}$ with $\beta=(2,0,2)$ and $\beta=(0,2,2)$, and \eqref{eq:6.93} gives
\begin{equation}\label{eq:6.97}
|J_{223}|\lesssim\Esh^2\Fn .
\end{equation}
Collecting \eqref{eq:6.90}, \eqref{eq:6.92}, \eqref{eq:6.99}, \eqref{eq:6.93}, \eqref{eq:6.94}, \eqref{eq:6.96} and \eqref{eq:6.97}, we obtain for $\alpha=(0,0,4)$
\begin{equation}\label{eq:6.100}
|J_2|\lesssim\Esh\Fsh^2+\Esh^2\big(\Ff+\Fn+\Fsh^2\big).
\end{equation}

\item \emph{The case $\alpha_3=2$, $\alpha_1+\alpha_2=2$.} Then $\alpha'\in\mathcal S_1$ or $\alpha'=\alpha$, so that $\beta_1+\beta_2\le2$. Choose $j\in\{1,2\}$ with $\alpha_j\ge1$. An integration by parts in $x_j$, using $\pa_j\pa^\delta\sgm=0$, gives
\begin{equation*}
J_2=-\int_\Omega\pa^{\alpha-e_j}(\sgm\varrho)\Big(\pa_j\pa^\beta w\,\pz\varrho+\pa^\beta w\,\pa_j\pz\varrho\Big)\pa^\delta\sgm\dx .
\end{equation*}
Here $\nrm{\pa^{\alpha-e_j}(\sgm\varrho)}\lesssim\Fv$ by \Cref{lem:NT-aux}(c). If $\beta_3\ge1$, then $\pa_j\pa^\beta w=-\pa_j\pa^{\beta-e_3}\divh v$ has order $|\beta|+1\le5$; if $\beta_3=0$, then $|\beta|\le2$ and $\nrm{\pa_j\pa^\beta w}\le\nrm{\pa_j\pa^\beta\divh v}$ by \Cref{lem:w}. In both cases $\nrm{\pa_j\pa^\beta w}\le\Fu$ up to a constant, and $\nrm{\pa^\beta w}\lesssim\Fu$ by \Cref{lem:NT-aux}(b). Therefore
\begin{equation}\label{eq:6.101}
|J_2|\lesssim\Fv\big(\nrm{\pa_j\pa^\beta w}\nrm{\pz\varrho}_{L^\infty}+\nrm{\pa^\beta w}\nrm{\nh\pz\varrho}_{L^\infty}\big)\lesssim\Ev\Fu\Fv .
\end{equation}

\item \emph{The case $\alpha_3=1$, $\alpha_1+\alpha_2=3$.} This occurs only in $\cM_3$, with $\alpha'=\alpha$, so $\beta\le\alpha$. We proceed as in (b). If $\beta_3=1$, $\pa_j\pa^\beta w=-\pa_j\pa^{\beta-e_3}\divh v$ has order at most five; for $\beta=\alpha$ it is a fifth horizontal derivative of $v$, bounded by $D_4\le\Fu$ (without the integration by parts in $x_j$ one would obtain the non-integrable bound $\Esh^2\Fu$). If $\beta_3=0$, then $|\beta|\le3$ and $\nrm{\pa_j\pa^\beta w}\le\nrm{\pa_j\pa^\beta\divh v}\le\nrm{\nh^{|\beta|+2}v}\le\Fu$ by \Cref{lem:w}. Since $\nrm{\pa^{\alpha-e_j}(\sgm\varrho)}\lesssim\Fv$,
\begin{equation}\label{eq:6.102}
|J_2|\lesssim\Ev\Fu\Fv\le\Esh\Fsh^2 .
\end{equation}

\item \emph{The term $K_2$, $\alpha\in\cM_2$.} Choose $j\in\{1,2\}$ with $\alpha_j\ge1$ and set $\tilde\alpha\eqdef\alpha-e_j$, a horizontal multi-index with $|\tilde\alpha|=3$. An integration by parts in $x_j$ gives
\begin{equation*}
\begin{aligned}
K_2&=\int_\Omega\sgm\,\pa^\alpha\varrho\,\pa_j\pa^{\tilde\alpha}w\,\pz\varrho\dx\\
&=-\int_\Omega\sgm\,\pz\varrho\,\pa_j\pa^\alpha\varrho\,\pa^{\tilde\alpha}w\dx-\int_\Omega\sgm\,\pa^\alpha\varrho\,\pa_j\pz\varrho\,\pa^{\tilde\alpha}w\dx\eqdef-K_{21}-K_{22}.
\end{aligned}
\end{equation*}
For $K_{22}$, \Cref{lem:NT-aux}(b), (c) give $|K_{22}|\lesssim\nrm{\pa^\alpha\varrho}\nrm{\pa^{\tilde\alpha}w}_{L^3}\nrm{\nh\pz\varrho}_{L^6}\lesssim\Ev\Fu\Fv$.

For $K_{21}$ we use the representation \eqref{eq:wrep}. Let $\varpi\eqdef g\Laph\varrho-\pz\divh G=\mu\AB w$. Horizontal derivatives commute with functions of $\AB$ (\Cref{lem:ABmean}), so, mode by mode,
\begin{equation*}
\mu\,\pa^{\tilde\alpha}w=\ABmi\pa^{\tilde\alpha}\varpi=\ABmh\psi\eqdef h,\qquad \psi\eqdef\ABmh\pa^{\tilde\alpha}\varpi=g\ABmh\pa^{\tilde\alpha}\Laph\varrho-\ABmh\pz\divh\pa^{\tilde\alpha}G .
\end{equation*}
Since $|k^{\tilde\alpha}||k|^2\le|k|^5$ on the mode $k$, \eqref{eq:ABdual} and \Cref{lem:NT-aux}(d) give
\begin{equation*}
\nrm{\psi}\le g\nrm{\ABmh|\nh|^5\varrho}+\nrm{\pa^{\tilde\alpha}G}\le g\Fv+\nrm{\nh^3G}\lesssim\Fsh .
\end{equation*}
Put $a\eqdef\sgm\pz\varrho\in H^3\subset W^{1,\infty}$. Since $h\in V$ (\Cref{def:AB}), $h=\pz h=0$ on $\pa\Omega$, so $ah=\pz(ah)=0$ there; together with the bound for $\nabla\pz(ah)$ below, the product $ah$ lies in $V=\cD(\ABh)$, and $\nrm{\ABh(ah)}=\nrm{\nabla\pz(ah)}$ by \Cref{lem:ABid}. Hence, with the pairing defined mode by mode as in \Cref{rem:justify}(c),
\begin{equation*}
\begin{aligned}
&K_{21}=\frac1\mu\int_\Omega\pa_j\pa^\alpha\varrho\;ah\dx=\frac1\mu\ip{\ABmh\pa_j\pa^\alpha\varrho,\ABh(ah)},\\
&|K_{21}|\le\frac1\mu\nrm{\ABmh|\nh|^5\varrho}\,\nrm{\nabla\pz(ah)}\le\frac1\mu\Fv\nrm{\nabla\pz(ah)} .
\end{aligned}
\end{equation*}
We expand $\nabla\pz(ah)=\nabla\pz a\,h+\pz a\,\nabla h+\nabla a\,\pz h+a\,\nabla\pz h$. The bound $\nrm{\ABmh\psi}_{L^\infty}\lesssim\nrm\psi$ is not available in three dimensions (\Cref{rem:Linfty}); instead, \Cref{lem:product}(i) with $p=3$ gives $\nrm{h}_{L^3}\le\nrm{h}_{L^\infty_zL^3_h}\lesssim\nrm{\psi}$ and $\nrm{\nabla h}+\nrm{\nabla\pz h}\lesssim\nrm{\psi}$. The coefficient $\nabla\pz a$, which contains third derivatives of $\varrho$, is placed in $L^6$: $\nrm{\nabla\pz a}_{L^6}\lesssim\nrm{a}_{H^3}\lesssim\Ev$; the other coefficients satisfy $\nrm{a}_{W^{1,\infty}}\lesssim\nrm{\varrho}_{W^{2,\infty}}\lesssim\Ev$. Therefore $\nrm{\nabla\pz(ah)}\lesssim\Ev\nrm{\psi}\lesssim\Ev\Fsh$, and
\begin{equation*}
|K_{21}|\lesssim\Ev\Fv\Fsh .
\end{equation*}
Collecting the two bounds,
\begin{equation}\label{eq:6.103}
|K_2|\lesssim\Ev\Fv\big(\Fu+\Fsh\big)\lesssim\Esh\Fsh^2 .
\end{equation}
\end{enumerate}
\textbf{Combined estimate for $|\gamma|=0$.} Collecting \eqref{eq:6.88}, \eqref{eq:6.100}, \eqref{eq:6.101}, \eqref{eq:6.102} and \eqref{eq:6.103}, we obtain
\begin{equation}\label{eq:6.104}
|\mathcal P_0|\lesssim\Esh\Fsh^2+\Esh^2\big(\Ff+\Fn+\Fsh^2\big).
\end{equation}
\end{enumerate}
Summing the estimates \eqref{eq:6.75}, \eqref{eq:6.78}, \eqref{eq:6.85}, \eqref{eq:6.87} and \eqref{eq:6.104} of the five cases $|\gamma|=4,\dots,0$, we obtain
\begin{equation}\label{eq:6.105}
|\cR_4|+|\cR_5|+|\cR_6|\lesssim\Esh\Fsh^2+\Esh^2\big(\Ff+\Fn+\Fsh^2\big).
\end{equation}

\begin{remark}\label{rem:C6}
The weight $\sgm=-1/\rs'$ is what makes the linear coupling between \eqref{PE}$_1$ and \eqref{PE}$_2$ cancel. In the case $|\gamma|=4$ the transport term produces $\sgm^{(5)}w\varrho$, so that the top-order estimate uses $\sgm\in C^5$, i.e.\ $\rs\in C^6$; it would suffice that $\rs\in C^5$ and $\rs^{(6)}\in L^2(0,1)$. All other estimates use at most $\rs\in C^5$. We do not know whether the stability result holds for $\rs\in C^5$ only.
\end{remark}

\subsubsection{Closed estimates}

The time integrals of the right-hand side of \eqref{eq:6.105} satisfy
\begin{equation*}
\begin{aligned}
\int_0^t\Esh\Fsh^2\dtau&\le\Big(\sup_{\tau\in(0,t)}\Esh\Big)\int_0^t\Fsh^2\dtau\lesssim\Es^3,\qquad \int_0^t\Esh^2\Fsh^2\dtau\lesssim\Es^4\le\Es^3,\\
\int_0^t\Esh^2\big(\Ff+\Fn\big)\dtau&\lesssim\Es^2\int_0^t\big(\Ff+\Fn\big)\dtau\lesssim\Es^3,
\end{aligned}
\end{equation*}
where $\Es=\Es(t)$, and we used \eqref{eq:Esharpdef}, \eqref{eq:small}, \eqref{eq:5.6} and \Cref{prop:Fnat}. Inserting \eqref{eq:6.105} into \eqref{eq:6.74}, we obtain, with $\Theta=\pz\divh(\rs v)$ as in \eqref{eq:6.74} and for every $\tau\in[0,t]$,
\begin{equation}\label{eq:6.106}
\begin{aligned}
&\frac12E_{\varrho,\sgm}^2(\tau)+\frac1\mu\ip{\AB(\sgm\varrho),\Theta}(\tau)+\frac1{2\mu^2}\nrm{\Theta(\tau)}^2+\frac34\cc2\int_0^\tau\Fv^2\dd s\\
&\quad\lesssim\Es(0)^2+\Eb^2+\Eb\Es+\Ef(0)^{\frac14}\Es^{\frac74}+\Es^3 .
\end{aligned}
\end{equation}
Next we absorb the cross term of the modified energy (\Cref{rem:square}). By Young's inequality and $\nrm{\AB(\sgm\varrho)}\le E_{\varrho,\sgm}$,
\begin{equation*}
\frac1\mu\big|\ip{\AB(\sgm\varrho),\Theta}\big|\le\frac14E_{\varrho,\sgm}^2+\frac1{\mu^2}\nrm{\Theta}^2 .
\end{equation*}
Since $\Theta=\rs\pz\divh v+\rs'\divh v$, \eqref{eq:4.3} and \eqref{eq:small} give $\nrm{\Theta(\tau)}\lesssim\nrm{v(\tau)}_{H^2}\lesssim\Ef(\tau)\le\Eb$. Inserting these bounds into \eqref{eq:6.106},
\begin{equation*}
\frac14E_{\varrho,\sgm}^2(\tau)+\frac34\cc2\int_0^\tau\Fv^2\dd s\lesssim\Es(0)^2+\Eb^2+\Eb\Es+\Ef(0)^{\frac14}\Es^{\frac74}+\Es^3,\qquad\tau\in[0,t].
\end{equation*}
Finally, by the norm equivalence \eqref{eq:R-equiv}, $\Ev(\tau)\le C_*E_{\varrho,\sgm}(\tau)$; no $L^2$-norm of $\varrho$ is needed. Taking the supremum over $\tau\in(0,t)$, we arrive at
\begin{equation*}
\cEv^2(t)=\sup_{\tau\in(0,t)}\Ev^2+\int_0^t\Fv^2\dtau\lesssim\Es(0)^2+\Eb^2+\Eb\Es+\Ef(0)^{\frac14}\Es^{\frac74}+\Es^3,
\end{equation*}
which is \eqref{eq:6.49}. This completes the proof of \Cref{prop:rho}.

\subsection{Proof of \texorpdfstring{\Cref{thm:main}}{Theorem 1.1}}\label{sec:proofmain}

The four component estimates \Cref{prop:rho}, \Cref{prop:u}, \Cref{prop:ut} and \Cref{prop:utt} are combined, using $\Es^4\le\Es^3$ under \eqref{eq:small}: \eqref{eq:6.47} bounds $\cEu^2$ by $\cEv^2$ plus admissible terms, and \eqref{eq:6.49} bounds $\cEv^2$. This yields the following proposition.

\begin{proposition}\label{prop:total}
Assume \eqref{eq:small}. Then, for every $t\ge0$, the total energy $\Es^2(t)=\cEv^2(t)+\cEu^2(t)+\cEut^2(t)+\cEutt^2(t)$ satisfies
\begin{equation}\label{eq:6.107}
\Es^2\le\cc3\Es(0)^2+\cc3\Eb^2+\cc3\Eb\Es+\cc3\Ef(0)^{\frac14}\Es^{\frac74}+\cc3\Es^3,
\end{equation}
where $\Es=\Es(t)$, $\Eb=\Eb(t)$, for a constant $\cc3\ge1$ depending only on $\mu,g,f,l,\rs$, $\nrm{\rs+\vin}_{L^\infty}$ and a positive lower bound for $\rs+\vin$.
\end{proposition}

\emph{The local solution class and the diagnostic traces.} The local solution class requires $v\in C([0,T];H^4)\cap L^2(0,T;H^5)$, $v_t\in C([0,T];H^2)\cap L^2(0,T;H^3)$ and $v_{tt}\in C([0,T];L^2)\cap L^2(0,T;H^1)$; in particular $\pz v\in C([0,T];H^3)$. The prescribed initial data are $\vin$ and $\uin$; the initial pressure $\qin$ is determined by them, as explained below. The vertical velocity $w$ is diagnostic by \eqref{eq:w}, the surface pressure $P$ is diagnostic by the hydrostatic decomposition $q=P-g\int_0^z\varrho$, and the traces $\varrho_t(0)$, $v_t(0)$ and $v_{tt}(0)$ are recovered from \eqref{PE}; none of these is independent initial data. We therefore reserve the subscript $\mathrm{in}$ for $\vin$, $\uin$ and $\qin$, and write diagnostic traces as evaluations at $t=0$. With $\rin=\rs+\vin$, $\bU(0)\eqdef(\uin,w(0))$ and $\bU_t(0)\eqdef(v_t(0),w_t(0))$,
\begin{equation*}
\begin{aligned}
w(0)&\eqdef-\int_0^z\divh\uin(x_h,z')\dd z',\qquad w_t(0)\eqdef-\int_0^z\divh v_t(0)(x_h,z')\dd z',\\
\varrho_t(0)&\eqdef-\uin\cdot\nh\vin-w(0)\big(\pz\vin+\rs'\big),\\
v_t(0)&\eqdef\frac1{\rin}\Big(\mu\Lap\uin-\rin\,\bU(0)\cdot\nabla\uin-f\rin\,\uin^\perp-\nh\qin\Big),\\
v_{tt}(0)&\eqdef\frac1{\rin}\Big(\mu\Lap v_t(0)-\varrho_t(0)\big(v_t(0)+\bU(0)\cdot\nabla\uin+f\uin^\perp\big)\\
&\qquad\qquad-\rin\big(\bU_t(0)\cdot\nabla\uin+\bU(0)\cdot\nabla v_t(0)+f\,v_t(0)^\perp\big)-\nh q_t(0)\Big),
\end{aligned}
\end{equation*}
where
\begin{equation*}
\nh\qin=\nh P(0)-g\int_0^z\nh\vin(x_h,s)\dd s,\qquad \nh q_t(0)=\nh P_t(0)-g\int_0^z\nh\varrho_t(0)(x_h,s)\dd s .
\end{equation*}
Here $P(0)$ is the unique zero-mean solution of \eqref{eq:Peq} at $t=0$ (\Cref{lem:propagation}), and $P_t(0)$ is the zero-mean solution of the equation obtained by differentiating \eqref{eq:Peq} in $t$ and evaluating at $t=0$ (equivalently, of $\divh\int_0^1v_{tt}(0)\dd z=0$); both are functions of $x_h$ only and are determined by $(\vin,\uin)$. On $\pa\Omega$ we have $\uin=0$, hence $\nh\uin=0$ and $\Laph\uin=0$, so that $\Lap\uin=\pz^2\uin$ and $\uin^\perp=0$ there; moreover $w(0)|_{\pa\Omega}=0$ by \eqref{eq:w} and \eqref{eq:lid0}, so $\bU(0)\cdot\nabla\uin=0$ on $\pa\Omega$; and $\vin|_{\pa\Omega}=0$ gives $\rin=\rs$ there. Hence $v_t(0)|_{\pa\Omega}=\rs^{-1}(\mu\pz^2\uin-\nh\qin)|_{\pa\Omega}$, and the requirement $v_t(0)|_{\pa\Omega}=0$, which is the time derivative of the no-slip condition, is equivalent to
\begin{equation}\label{eq:compat}
\nh\qin-\mu\pz^2\uin=0\qquad\text{on }\pa\Omega .
\end{equation}
Since $\qin$ is determined by $(\vin,\uin)$, \eqref{eq:compat} is a genuine compatibility condition on $(\vin,\uin)$. The condition \eqref{eq:wt0} holds automatically at $z=0$; at $z=1$ it is equivalent to $\divh\int_0^1v_t(0,x_h,z)\dd z\equiv0$, that is, to \eqref{eq:Peq} at $t=0$, which is how $P(0)$ was chosen. The order-zero condition is \eqref{eq:lid0}. For positive times, \eqref{eq:wtt} is the boundary identity behind the clamped inversion in \eqref{eq:wrep}; its persistence is the content of \Cref{lem:propagation}.

\emph{Linear bound for the initial energies.} Let $\nrm{\vin}_{H^4}+\nrm{\uin}_{H^4}\le1$ and $\rin\ge\frac12\min\rs$, and write $\mathsf d\eqdef\nrm{\vin}_{H^4}+\nrm{\uin}_{H^4}$. Then $\nrm{\rin^{\pm1}}_{W^{2,\infty}}\lesssim1$, and the coefficient $a=\int_0^1\rin^{-1}\dd z$ of \eqref{eq:Peq} satisfies $a\gtrsim1$ and $\nrm{a}_{W^{2,\infty}(\Tl)}\lesssim\nrm{a}_{H^4(\Tl)}\lesssim1$. By \eqref{eq:Phieq}, $\cN(0)$ consists of terms linear in $(\uin,\nh\vin)$, namely $-f\divh\int_0^1\uin^\perp$, $g\divh\int_0^1\rin^{-1}\int_0^z\nh\vin$ and $\mu\divh\int_0^1\rin^{-1}\Lap\uin$, plus the quadratic term $-\divh\divh\int_0^1\uin\otimes\uin$; hence the elliptic estimate on $\Tl$ gives $\nrm{\nh P(0)}_{H^2(\Tl)}\lesssim\nrm{\cN(0)}_{H^1(\Tl)}\lesssim\mathsf d+\mathsf d^2\lesssim\mathsf d$. Every term in the formulas for $w(0)$, $\varrho_t(0)$ and $v_t(0)$ contains a factor $\uin$, $\nh\vin$ or $\nh P(0)$ (the term $-w(0)\rs'$ is linear in $\uin$), so $\nrm{w(0)}_{H^3}+\nrm{\varrho_t(0)}_{H^3}+\nrm{v_t(0)}_{H^2}\lesssim\mathsf d$. Differentiating \eqref{eq:Peq} in $t$ gives $\divh(a\nh P_t(0))=\cN_t(0)-\divh(a_t(0)\nh P(0))$ with $a_t(0)=-\int_0^1\rin^{-2}\varrho_t(0)\dd z$, where each term of $\cN_t(0)$ (product rule in \eqref{eq:Phieq}) contains one of $v_t(0)$, $\varrho_t(0)$, $\uin$, $\nh\vin$, and is a horizontal divergence of an $L^2(\Tl)$ field bounded by $C\mathsf d$ (the worst term is $\mu\divh\int_0^1\rin^{-1}\Lap v_t(0)$). Hence $\nrm{\nh P_t(0)}_{L^2(\Tl)}\lesssim\nrm{\cN_t(0)}_{H^{-1}(\Tl)}+\nrm{a_t(0)\nh P(0)}_{L^2(\Tl)}\lesssim\mathsf d$, and the formula for $v_{tt}(0)$ gives $\nrm{v_{tt}(0)}\lesssim\mathsf d$. Since $\Ef(0)\lesssim\nrm{\nh\vin}+\nrm{\nh\uin}+\nrm{v_t(0)}+\nrm{w(0)}$ and $\Es(0)\lesssim\Esh(0)=\nrm{\vin}_{H^4}+\nrm{\uin}_{H^4}+\nrm{v_t(0)}_{H^2}+\nrm{v_{tt}(0)}$, this proves \eqref{eq:NT-data} below.

The local well-posedness result below is proved in \Cref{app:local}; the existence time depends only on an upper bound for $\nrm{\vin}_{H^4}+\nrm{\uin}_{H^4}$ and a positive lower bound for $\rs+\vin$. Its proof does not use the stability condition \eqref{eq:stable}, which is important for \cref{sec:RT}. \Cref{rem:gradq} records $\nabla q\in C([0,T];H^2)$ in \eqref{eq:localclass}.

\begin{proposition}\label{prop:local}
Let $\rs\in C^5([0,1])$ with $\min_{[0,1]}\rs>0$ (condition \eqref{eq:stable} is not required) and $f\in\R$. Assume that the initial data satisfy
\begin{equation*}
\rs+\vin>0\ \text{on }\ol\Omega,\qquad\vin\in H^4\ \text{with}\ \vin|_{\pa\Omega}=\pz\vin|_{\pa\Omega}=0,\qquad\uin\in H^4\cap H^1_0,
\end{equation*}
the constraint \eqref{eq:lid0}, and the compatibility condition \eqref{eq:compat}, where $\qin$ is determined by \eqref{eq:Peq} at $t=0$ (so that \eqref{eq:wt0} holds). Then there exists $T>0$, depending only on an upper bound for $\nrm{\vin}_{H^4}+\nrm{\uin}_{H^4}$ and a positive lower bound for $\rs+\vin$, such that \eqref{PE} with \eqref{eq:ID}--\eqref{eq:lid} has a solution $(\varrho,v,q)$ on $[0,T]$ satisfying
\begin{equation}\label{eq:localclass}
\begin{aligned}
&(\varrho,v,\nabla q)\in C\big([0,T];H^4\big)\times C\big([0,T];H^4\big)\times C\big([0,T];H^2\big),\\
&v\in L^2\big(0,T;H^5\big),\\
&v_t\in C\big([0,T];H^2\big)\cap L^2\big(0,T;H^3\big),\\
&v_{tt}\in C\big([0,T];L^2\big)\cap L^2\big(0,T;H^1\big).
\end{aligned}
\end{equation}
In this class $(\varrho,v,\nabla q)$ is unique; $q$ is unique up to an additive function of $t$, and it is unique if the surface pressure $P=q+g\int_0^z\varrho\dd s$ is normalized to have zero horizontal mean.
\end{proposition}

\begin{proof}[Proof of \Cref{thm:main}]
Let $\cc1$ be the constant in \Cref{prop:low} (fixed as in \Cref{rem:c1}) and $\cc3$ the constant in \Cref{prop:total}. Fix the bootstrap size
\begin{equation*}
\delta_0\eqdef\min\Big\{\frac1{256\cc1},\ \frac1{8\cc3}\Big\},
\end{equation*}
and set
\begin{equation*}
\mathsf K\eqdef2\cc3+44\cc3^2+7^7\cc3^8 .
\end{equation*}
On the admissible data class
\begin{equation*}
\Big\{(\vin,\uin):\ \nrm{\vin}_{H^4}+\nrm{\uin}_{H^4}\le1\ \text{and}\ \rs+\vin\ge\tfrac12\min_{z\in[0,1]}\rs\ \text{in }\Omega\Big\},
\end{equation*}
on which all the constants above depend only on $\mu,g,f,l$ and $\rs$, there exists $\cc5>0$ such that
\begin{equation}\label{eq:NT-data}
\Ef(0)+\Es(0)\le\cc5\big(\nrm{\vin}_{H^4}+\nrm{\uin}_{H^4}\big).
\end{equation}
We take
\begin{equation}\label{eq:6.110}
\varepsilon_0\eqdef\min\Big\{\frac{\delta_0}{2\cc5\sqrt{\mathsf K}},\ \frac1{2\cc6}\min\Big\{1,\min_{z\in[0,1]}\rs(z)\Big\},\ 1\Big\},
\end{equation}
where $\cc6$ is the constant of the embedding $H^4(\Omega)\hookrightarrow L^\infty(\Omega)$; the second entry gives $\nrm{\vin}_{L^\infty}\le\frac12\min\rs$, hence $\rs+\vin\ge\frac12\min\rs$. Since $\varepsilon_0\le1$, data with $\nrm{\vin}_{H^4}+\nrm{\uin}_{H^4}\le\varepsilon_0$ lie in the admissible class, and \eqref{eq:NT-data} gives
\begin{equation}\label{eq:6.111}
\sqrt{\mathsf K}\big(\Es(0)+\Ef(0)\big)\le\sqrt{\mathsf K}\cc5\varepsilon_0\le\frac{\delta_0}2,\qquad \Es(0)\le\frac{\delta_0}{2\sqrt{\mathsf K}}\le\frac{\delta_0}4,
\end{equation}
where the last step uses $\mathsf K\ge4$, which follows from $\cc3\ge1$.

By \Cref{prop:local}, a unique solution in the class \eqref{eq:localclass} exists on some interval $[0,T_0]$. Define the maximal bootstrap time
\begin{equation*}
T^*\eqdef\sup\Big\{T>0:\ \text{the solution exists on }[0,T]\text{ in the class \eqref{eq:localclass}, and }\Es\le\delta_0\text{ on }[0,T]\Big\}.
\end{equation*}
Since $t\mapsto\Es(t)$ is continuous and nondecreasing with $\Es(t)\to\Es(0)\le\frac{\delta_0}4$ as $t\downarrow0$, we have $T^*>0$. Assume $T^*<+\infty$. For every $t<T^*$, \eqref{eq:small} holds on $[0,t]$ because $\delta_0\le(256\cc1)^{-1}$, so \Cref{prop:total} and \Cref{prop:decay2} apply on $[0,t]$ (\Cref{rem:c1}). By \eqref{eq:5.4},
\begin{equation}\label{eq:6.112}
\Ef(t)\le\Ef(0),\qquad\int_0^t\Ff^2(\tau)\dtau\le8\Ef(0)^2,\qquad\Eb^2(t)\le9\Ef(0)^2,\qquad t\in[0,T^*).
\end{equation}
On $[0,T^*)$ the following absorption bounds hold, by the choice of $\delta_0$, by \eqref{eq:6.112}, and by Young's inequality with exponents $8$ and $\frac87$:
\begin{equation*}
\begin{aligned}
&\cc3\Eb^2\le9\cc3\Ef(0)^2,\qquad \cc3\Eb\Es\le\frac18\Es^2+2\cc3^2\Eb^2\le\frac18\Es^2+18\cc3^2\Ef(0)^2,\\
&\cc3\Es^3\le\cc3\delta_0\Es^2\le\frac18\Es^2,\qquad \cc3\Ef(0)^{\frac14}\Es^{\frac74}\le\frac18\Es^2+\frac{7^7}8\cc3^8\Ef(0)^2 .
\end{aligned}
\end{equation*}
Inserting these into \eqref{eq:6.107} gives $\frac58\Es^2\le\cc3\Es(0)^2+\big(9\cc3+18\cc3^2+\frac{7^7}8\cc3^8\big)\Ef(0)^2$. Since $\frac85(9\cc3+18\cc3^2)\le44\cc3^2$ for $\cc3\ge1$, \eqref{eq:6.111} yields
\begin{equation}\label{eq:6.113}
\Es(t)\le\Big(2\cc3\Es(0)^2+\big(44\cc3^2+7^7\cc3^8\big)\Ef(0)^2\Big)^{\frac12}\le\sqrt{\mathsf K}\big(\Es(0)+\Ef(0)\big)\le\frac{\delta_0}2,\qquad t\in[0,T^*).
\end{equation}

We show that the solution extends to $[0,T^*]$ in the class \eqref{eq:localclass}. By \eqref{eq:6.113}, the norms in \eqref{eq:localclass} are bounded uniformly on $[0,T^*)$. Hence, by the Lions--Magenes lemma ($v\in L^2(0,T^*;H^5)$ with $v_t\in L^2(0,T^*;H^3)$, and $v_t\in L^2(0,T^*;H^3)$ with $v_{tt}\in L^2(0,T^*;H^1)$), $v\in C([0,T^*];H^4)$ and $v_t\in C([0,T^*];H^2)$; by \Cref{rem:justify}(a)--(b), $\varrho\in C([0,T^*];H^4)$ and $v_{tt}\in C([0,T^*];L^2)$; and $\nabla q\in C([0,T^*];H^2)$ by \Cref{rem:gradq}. At $t=T^*$ the hypotheses of \Cref{prop:local} hold: the traces $(\varrho(T^*),v(T^*))$ lie in $H^4\times(H^4\cap H^1_0)$, and \Cref{lem:propagation} gives $\varrho(T^*)|_{\pa\Omega}=\pz\varrho(T^*)|_{\pa\Omega}=0$ and $w|_{\pa\Omega}=0$ on $[0,T^*]$, i.e.\ \eqref{eq:lid0} at $t=T^*$. Since $v=0$ on $\pa\Omega$ for all times and $v_t\in C([0,T^*];H^2)$ with $H^2(\Omega)\hookrightarrow C(\ol\Omega)$, we have $v_t(T^*)|_{\pa\Omega}=0$, so the momentum equation at $t=T^*$ restricted to $\pa\Omega$ is \eqref{eq:compat} at the restart time; since $w_t$ has boundary traces on $[0,T^*]$, \Cref{lem:propagation}(i) gives $w_t|_{\pa\Omega}=0$ at $t=T^*$, i.e.\ \eqref{eq:wt0} at the restart time, and the pressure at $T^*$ is the one determined by \eqref{eq:Peq}. The total density $\rho=\rs+\varrho$ is transported by the divergence-free field $\bU$ with $\bU\cdot n=0$ on $\pa\Omega$, so its flow preserves Lebesgue measure and $\inf_\Omega\rho$; hence $\rho(T^*)\ge\inf_\Omega\rin\ge\frac12\min\rs$. Restarting \Cref{prop:local} at $T^*$ and gluing (the time derivatives match at $T^*$, being determined by \eqref{PE}), we obtain a solution on $[0,T^*+T_1]$ in the class \eqref{eq:localclass} for some $T_1>0$; only the local existence and uniqueness in \Cref{prop:local} are needed here, not the uniformity of the existence time. Since $\Es(T^*)\le\frac{\delta_0}2$ by \eqref{eq:6.113} and continuity, and $\Es$ is continuous, $\Es\le\delta_0$ on $[0,T^*+\tau_1]$ for some $\tau_1\in(0,T_1]$, which contradicts the maximality of $T^*$. Hence $T^*=+\infty$.

Consequently \eqref{eq:6.113} holds for all $t\ge0$, and together with \eqref{eq:NT-data},
\begin{equation*}
\sup_{t\ge0}\Es(t)\le\sqrt{\mathsf K}\,\cc5\big(\nrm{\vin}_{H^4}+\nrm{\uin}_{H^4}\big)\le\sqrt{\mathsf K}\,\cc5\,\varepsilon_0 .
\end{equation*}
Since $\Ev^2+\Eu^2\le\Es^2$ and $\sup_t\Eu^2+\int_0^\infty\Fu^2\dt\le\sup_t\Es^2$, this gives \eqref{eq:unif} and \eqref{eq:ureg}. The solution is unique (with $P$ normalized to zero horizontal mean) and classical by \Cref{prop:local} and \eqref{eq:localclass}, and the clamped conditions hold for all $t\ge0$ by \Cref{lem:propagation}(ii). Since $\Es(t)\le\frac{\delta_0}2<(256\cc1)^{-1}$ for all $t\ge0$, \eqref{eq:small} holds globally, and \eqref{eq:decay} follows from \eqref{eq:5.3}, \eqref{eq:5.4} and \eqref{eq:6.7b}, because $\Es(t)\le\sqrt{\mathsf K}\cc5(\nrm{\vin}_{H^4}+\nrm{\uin}_{H^4})$ for all $t\ge0$.
\end{proof}

\section{Convergence of the density to the limit state}\label{sec:limit}

In this section we prove \Cref{thm:limit}. The refined smallness \eqref{eq:slope} is used only in items (v) and (vi). By \eqref{eq:meaneq}, the horizontal mean of the density changes only through the vertical flux $\ol{w\wt\varrho}$, which vanishes at the lids:
\begin{equation}\label{eq:7.1}
\pt\bar\varrho(t,z)=-\pz\big(\ol{w\wt\varrho}\big)(t,z).
\end{equation}
The Coriolis force does not enter the density equation. It affects this section only through the decay estimates of \Cref{sec:decay} and the uniform bound of \Cref{thm:main}.

Throughout, $(\varrho,v,q)$ is the global solution given by \Cref{thm:main}. The bootstrap argument in \Cref{sec:proofmain} shows that \eqref{eq:small} holds on $[0,\infty)$ and that
\begin{equation}\label{eq:L-Es}
\Es(\infty)\eqdef\sup_{t\ge0}\Es(t)\le C\mathsf N\le C\varepsilon_0,\qquad \Es(\infty)\le\frac1{256\cc1}\le1,\qquad \mathsf N\eqdef\nrm{\vin}_{H^4}+\nrm{\uin}_{H^4}.
\end{equation}
Hence (\Cref{rem:c1}(ii)) \Cref{prop:decay1,prop:decay2} apply on every interval $[0,t]$, and \eqref{eq:5.3}--\eqref{eq:5.tail} hold for every $t\ge0$ with $\Es(t)$ replaced by $\Es(\infty)$. In particular,
\begin{equation}\label{eq:L-Ftail}
\int_t^\infty\Ff^2(\tau)\dtau\le8\Ef^2(t)\le C\Es(\infty)^2(1+t)^{-2},\qquad t\ge0 .
\end{equation}

\subsection{Dyadic time integrals and the vertical flux}

The following elementary lemma converts the tail bound \eqref{eq:L-Ftail} into decay of time integrals without loss.

\begin{lemma}\label{lem:dyadic}
Let $a,c\ge0$ with $a+c\le2$, let $b\ge0$, and assume $\gamma\eqdef a+b+\frac{a+c}2-1>0$. Let $F,G\ge0$ be measurable on $(0,\infty)$ with
\begin{equation*}
\int_t^\infty F^2\dtau\le A^2(1+t)^{-2}\quad(t\ge0),\qquad \int_0^\infty G^2\dtau\le B^2 .
\end{equation*}
Then, for every $t\ge0$,
\begin{equation}\label{eq:dyadic}
\int_t^\infty(1+\tau)^{-b}F^aG^c\dtau\le\big(1-2^{-\gamma}\big)^{-1}A^aB^c(1+t)^{-\gamma}.
\end{equation}
\end{lemma}

\begin{proof}
Let $s_j\eqdef2^j(1+t)-1$ and $I_j\eqdef[s_j,s_{j+1}]$, $j\ge0$, so that $|I_j|=1+s_j$ and $(1+\tau)^{-b}\le(1+s_j)^{-b}$ on $I_j$. H\"older's inequality with the exponents $\frac2a$, $\frac2c$ and $(1-\frac{a+c}2)^{-1}$ (a factor is omitted when its exponent is infinite) gives
\begin{equation*}
\int_{I_j}F^aG^c\dtau\le\Big(\int_{I_j}F^2\Big)^{\frac a2}\Big(\int_{I_j}G^2\Big)^{\frac c2}|I_j|^{1-\frac{a+c}2}\le A^aB^c(1+s_j)^{-a+1-\frac{a+c}2}.
\end{equation*}
Hence the integral over $I_j$ is at most $A^aB^c(1+s_j)^{-\gamma}=A^aB^c(1+t)^{-\gamma}2^{-j\gamma}$, and summing over $j\ge0$ gives \eqref{eq:dyadic}.
\end{proof}

\begin{lemma}\label{lem:L-flux}
Under the hypotheses of \Cref{thm:main}, $\bar\varrho\in C^1([0,\infty);L^2(0,1))$, and for every $t\ge0$,
\begin{equation}\label{eq:7.4}
\nrm{\pt\bar\varrho(t)}_{L^2(0,1)}\le C\Big(\Es(\infty)(1+t)^{-\frac34}\Ff(t)+\Es(\infty)^{\frac12}\Ff(t)^{\frac32}\Big).
\end{equation}
Moreover, for every $t\ge0$,
\begin{equation}\label{eq:L-tail}
\int_t^\infty\nrm{\pt\bar\varrho(\tau)}_{L^2(0,1)}\dtau\le C\Es(\infty)^2(1+t)^{-\frac54}.
\end{equation}
\end{lemma}

\begin{proof}
\emph{Regularity.} On every bounded interval $[0,T]$ the solution belongs to the class \eqref{eq:localclass}; in particular $\varrho,v\in C([0,T];H^4)$. By \Cref{lem:w}, $w\in C([0,T];H^3)$, and since $H^3(\Omega)$ is an algebra, \eqref{PE}$_1$ gives
\begin{equation*}
\varrho_t=-v\cdot\nh\varrho-w\pz\varrho-\rs'w\in C([0,T];H^3).
\end{equation*}
Hence $\varrho\in C^1([0,\infty);H^3(\Omega))$, $\bar\varrho\in C^1([0,\infty);H^3(0,1))$, and \eqref{eq:7.1} holds in $C([0,\infty);L^2(0,1))$.

\emph{Step 1: decomposition of the flux.} By \eqref{eq:w}, $\pz w=-\divh v$, so
\begin{equation*}
\pz\big(\ol{w\wt\varrho}\big)=-\ol{\divh v\,\wt\varrho}+\ol{w\,\pz\wt\varrho}.
\end{equation*}
By the Cauchy--Schwarz inequality on $\Tl$, $|\bar g(z)|\le|\Tl|^{-1/2}\nrm{g(\cdot,z)}_{L^2(\Tl)}$, hence $\nrm{\bar g}_{L^2(0,1)}\le|\Tl|^{-1/2}\nrm{g}_{L^2(\Omega)}$ (\Cref{lem:meanosc}). Therefore
\begin{equation}\label{eq:7.2}
\nrm{\pt\bar\varrho(t)}_{L^2(0,1)}\le|\Tl|^{-\frac12}\big(\cA(t)+\cB(t)\big),\qquad \cA(t)\eqdef\nrm{\divh v\,\wt\varrho}_{L^2(\Omega)},\quad \cB(t)\eqdef\nrm{w\,\pz\wt\varrho}_{L^2(\Omega)} .
\end{equation}

\emph{Step 2: the oscillatory density.} By \eqref{eq:5.3} and \eqref{eq:L-Es}, $\nrm{\wt\varrho}+\nrm{\nh\varrho}\le C\Es(\infty)(1+t)^{-1}$. Since $\varrho=0$ on $\pa\Omega$, also $\nh\varrho=0$ on $\pa\Omega$, and one integration by parts in $z$ gives
\begin{equation*}
\nrm{\nh\pz\varrho}^2=-\ip{\nh\varrho,\nh\pz^2\varrho}\le\nrm{\nh\varrho}\,\Ev\le C\Es(\infty)^2(1+t)^{-1}.
\end{equation*}
The function $\pz\wt\varrho=\pz\varrho-\ol{\pz\varrho}$ has zero horizontal mean and $\nh\pz\wt\varrho=\nh\pz\varrho$, so $\nrm{\pz\wt\varrho}\le l\nrm{\nh\pz\varrho}$ by \Cref{lem:poincare}. Together with $\nrm{\nh\wt\varrho}=\nrm{\nh\varrho}$ this gives
\begin{equation}\label{eq:L-rhot}
\nrm{\wt\varrho(t)}_{L^2}\le C\Es(\infty)(1+t)^{-1},\qquad \nrm{\wt\varrho(t)}_{H^1}\le C\Es(\infty)(1+t)^{-\frac12}.
\end{equation}
On the other hand, $\varrho(t)$ satisfies the clamped conditions, so \eqref{eq:2.11} and the definition \eqref{eq:lowenergy} of $\Ff$ give
\begin{equation}\label{eq:L-rhoF}
\nrm{\nh\pz\varrho}^2\le\nrm{\nabla\nh\varrho}^2\lesssim\norm{\ABmh\Laph\varrho}\nrm{\varrho}_{H^4}\le C\Es(\infty)\Ff .
\end{equation}

\emph{Step 3: the term $\cA$.} By H\"older's inequality, $\cA\le\nrm{\divh v}_{L^6}\nrm{\wt\varrho}_{L^3}\le\nrm{\divh v}_{L^6}\nrm{\wt\varrho}_{L^2}^{1/2}\nrm{\wt\varrho}_{L^6}^{1/2}$. Since $\divh v$ vanishes on $\pa\Omega$, \Cref{lem:GN}(3) and \Cref{lem:poincare} give
\begin{equation*}
\nrm{\divh v}_{L^6}\lesssim\nrm{\divh v}_{H^1}\lesssim\nrm{\nabla\divh v}\le\sqrt2\,\nrm{\nabla\nh v}\le\big(2/\mu\big)^{\frac12}\Ff,
\end{equation*}
and $\nrm{\wt\varrho}_{L^6}\lesssim\nrm{\wt\varrho}_{H^1}$. By \eqref{eq:L-rhot},
\begin{equation}\label{eq:L-A}
\cA(t)\le C\Es(\infty)(1+t)^{-\frac34}\Ff(t).
\end{equation}

\emph{Step 4: the term $\cB$.} For a.e.\ $z$, H\"older's inequality on $\Tl$ gives $\nrm{(w\pz\wt\varrho)(\cdot,z)}_{L^2(\Tl)}\le\nrm{w(\cdot,z)}_{L^4(\Tl)}\nrm{\pz\wt\varrho(\cdot,z)}_{L^4(\Tl)}$. By \Cref{lem:w} with $p=4$,
\begin{equation*}
\sup_z\nrm{w(\cdot,z)}_{L^4(\Tl)}\le C\nrm{\nh\divh v}\le C\nrm{\nabla\nh v}\le C\Ff,
\end{equation*}
and, since $\pz\wt\varrho(\cdot,z)$ has zero mean, \Cref{lem:aniso}(b) gives $\nrm{\pz\wt\varrho(\cdot,z)}_{L^4(\Tl)}\le C\nrm{\nh\pz\varrho(\cdot,z)}_{L^2(\Tl)}$. Squaring, integrating in $z$ and using \eqref{eq:L-rhoF},
\begin{equation}\label{eq:7.3}
\cB(t)\le C\Ff\,\nrm{\nh\pz\varrho}\le C\Es(\infty)^{\frac12}\Ff(t)^{\frac32}.
\end{equation}
In three dimensions $\nrm{\wt\varrho}_{L^\infty}$ is not controlled by $\nrm{\nh\pz\varrho}$, since $H^1(\T^2)\not\subset L^\infty(\T^2)$. Steps 3 and 4 use no $L^\infty$ norm. Because $\wt\varrho$ and $\pz\wt\varrho$ have zero horizontal mean, passing from $L^2$ to $L^6$, respectively from $L^2_h$ to $L^4_h$, costs one derivative that is available in \eqref{eq:L-rhot} and \eqref{eq:L-rhoF}, and the rate is not affected.

Inserting \eqref{eq:L-A} and \eqref{eq:7.3} into \eqref{eq:7.2} gives \eqref{eq:7.4}.

\emph{Step 5: time integrals.} By \eqref{eq:L-Ftail}, \Cref{lem:dyadic} applies with $F=\Ff$ and $A=C\Es(\infty)$. With $(a,b,c)=(1,\frac34,0)$ and with $(a,b,c)=(\frac32,0,0)$ one has $\gamma=\frac54$ in both cases, so
\begin{equation}\label{eq:7.5}
\int_t^\infty\Big(\Es(\infty)(1+\tau)^{-\frac34}\Ff+\Es(\infty)^{\frac12}\Ff^{\frac32}\Big)\dtau\le C\Es(\infty)^2(1+t)^{-\frac54}.
\end{equation}
Integrating \eqref{eq:7.4} over $(t,\infty)$ and using \eqref{eq:7.5} gives \eqref{eq:L-tail}.
\end{proof}

\subsection{The rearrangement at a fixed time}

For a bounded measurable $\phi\ge0$ on $\Omega$ we denote by $\phi^*$ its decreasing vertical rearrangement, defined as in \eqref{eq:rearr}: $\phi^*(z)\eqdef\int_0^\infty\mathbf1_{\{|\{\phi>\lambda\}|\ge|\Tl|z\}}\dd\lambda$. It depends only on $z$, is nonincreasing, and $|\{\phi^*>\lambda\}|=|\{\phi>\lambda\}|$ for every $\lambda$, where $\phi^*$ is regarded as a function on $\Omega$. The next lemma shows that, for a strictly stratified profile, the rearrangement differs from the horizontal mean only at second order in the oscillation. The second-order error is measured in the mixed norm $L^2_hL^\infty_z$, which is controlled by a vertical interpolation alone; no $L^\infty$ bound in the horizontal variables is needed.

\begin{lemma}\label{lem:rearr-inst}
Let $\phi\in C^2(\ol\Omega)$ with $\phi(\cdot,0)\equiv a_0$ and $\phi(\cdot,1)\equiv a_1>0$, and assume that
\begin{equation*}
\kappa_0\le-\pz \phi\le K_0,\qquad |\pz^2\phi|\le K_0\qquad\text{on }\ol\Omega,
\end{equation*}
for some constants $0<\kappa_0\le K_0$. Then $\phi^*\in C^1([0,1])$, $\phi^*(0)=a_0$, $\phi^*(1)=a_1$, $\kappa_0\le-(\phi^*)'\le K_0$, and
\begin{equation}\label{eq:rearr-inst}
\begin{aligned}
\nrm{\phi^*-\bar \phi}_{L^2(0,1)}&\le C\nrm{\wt \phi}_{L^2_hL^\infty_z}\big(\nrm{\pz\wt \phi}_{L^2(\Omega)}+\nrm{\wt \phi}_{L^2(\Omega)}\big)\\
&\le C\nrm{\wt \phi}_{L^2(\Omega)}^{\frac12}\nrm{\pz\wt \phi}_{L^2(\Omega)}^{\frac32},
\end{aligned}
\end{equation}
with $C$ depending only on $\kappa_0$, $K_0$ and $l$.
\end{lemma}

\begin{proof}
For $\alpha\in[0,1]$ let $\phi_\alpha\eqdef\bar \phi+\alpha\wt \phi$. Since $\pz\bar \phi$ and $\pz^2\bar \phi$ are horizontal averages of $\pz \phi$ and $\pz^2\phi$, the functions $\phi_\alpha$ have the boundary values $a_0,a_1$ and satisfy $\kappa_0\le-\pz \phi_\alpha\le K_0$ and $|\pz^2\phi_\alpha|\le K_0$. For each $x_h$, $z\mapsto \phi_\alpha(x_h,z)$ is therefore a $C^2$ decreasing bijection of $[0,1]$ onto $[a_1,a_0]$; let $Z_\alpha(x_h,\lambda)$ denote its inverse. By the implicit function theorem $Z_\alpha$ is $C^2$ in $\alpha$, and $\kappa_0\le-1/\pa_\lambda Z_\alpha\le K_0$.

\emph{Rearrangement and mean as inverse functions.} For $\lambda\in(a_1,a_0)$ and each $x_h$, $\{z:\phi(x_h,z)>\lambda\}=[0,Z_1(x_h,\lambda))$, hence $|\{\phi>\lambda\}|=|\Tl|\,\ol{Z_1}(\lambda)$, the bar denoting the horizontal average. The function $\ol{Z_1}$ decreases from $1$ to $0$ on $[a_1,a_0]$ with $\kappa_0\le-1/\ol{Z_1}'\le K_0$, so its inverse is a $C^1$ decreasing function with these boundary values and slopes, and by the definition of $\phi^*$, $(\phi^*)^{-1}=\ol{Z_1}$; moreover $(\bar \phi)^{-1}=Z_0$, since $\phi_0=\bar \phi$ does not depend on $x_h$. This proves the first assertions.

\emph{First-order cancellation.} Differentiating $\phi_\alpha(x_h,Z_\alpha(x_h,\lambda))=\lambda$ in $\alpha$, with $\phi'_\alpha\eqdef\pz \phi_\alpha$ and $\phi''_\alpha\eqdef\pz^2\phi_\alpha$,
\begin{equation}\label{eq:Zalpha}
\pa_\alpha Z_\alpha=-\frac{\wt \phi}{\phi'_\alpha}\Big|_{z=Z_\alpha},\qquad \pa_\alpha^2Z_\alpha=\Big(\frac{2\wt \phi\,\pz\wt \phi}{(\phi'_\alpha)^2}-\frac{\phi''_\alpha\,\wt \phi^{\,2}}{(\phi'_\alpha)^3}\Big)\Big|_{z=Z_\alpha}.
\end{equation}
At $\alpha=0$, $Z_0$ and $\phi_0'$ do not depend on $x_h$, and $\wt \phi$ has zero horizontal mean, so $\ol{\pa_\alpha Z_\alpha}|_{\alpha=0}=0$. Taylor's formula in $\alpha$ gives
\begin{equation*}
(\phi^*)^{-1}-(\bar \phi)^{-1}=\ol{Z_1}-\ol{Z_0}=\int_0^1(1-\alpha)\,\ol{\pa_\alpha^2Z_\alpha}\dd\alpha .
\end{equation*}

\emph{Second-order bound.} For fixed $x_h$ and $\alpha$, the change of variables $\lambda=\phi_\alpha(x_h,z)$, with $|\dd\lambda|\le K_0\dd z$, gives $\int_{a_1}^{a_0}|h(x_h,Z_\alpha(x_h,\lambda))|^2\dd\lambda\le K_0\int_0^1|h(x_h,z)|^2\dd z$. Hence, by \eqref{eq:Zalpha}, with all norms taken in $z$ or $\lambda$ for fixed $x_h$,
\begin{equation*}
\begin{aligned}
\nrm{\pa_\alpha^2Z_\alpha(x_h,\cdot)}_{L^2(a_1,a_0)}&\le K_0^{\frac12}\Big(\frac2{\kappa_0^2}\nrm{\wt \phi\,\pz\wt \phi}_{L^2_z}+\frac{K_0}{\kappa_0^3}\nrm{\wt \phi^{\,2}}_{L^2_z}\Big)\\
&\le C\nrm{\wt \phi(x_h,\cdot)}_{L^\infty(0,1)}\big(\nrm{\pz\wt \phi(x_h,\cdot)}_{L^2(0,1)}+\nrm{\wt \phi(x_h,\cdot)}_{L^2(0,1)}\big).
\end{aligned}
\end{equation*}
We now average in $x_h$ before using the Cauchy--Schwarz inequality. By Minkowski's inequality and then the Cauchy--Schwarz inequality in $x_h$,
\begin{equation*}
\begin{aligned}
\nrm{\ol{\pa_\alpha^2Z_\alpha}}_{L^2(a_1,a_0)}&\le\frac1{|\Tl|}\int_{\Tl}\nrm{\pa_\alpha^2Z_\alpha(x_h,\cdot)}_{L^2(a_1,a_0)}\dxh\\
&\le C\nrm{\wt \phi}_{L^2_hL^\infty_z}\big(\nrm{\pz\wt \phi}_{L^2(\Omega)}+\nrm{\wt \phi}_{L^2(\Omega)}\big),
\end{aligned}
\end{equation*}
uniformly in $\alpha$, and $(\phi^*)^{-1}-(\bar \phi)^{-1}$ obeys the same bound.

\emph{Back to the functions.} For $z\in[0,1]$ put $\lambda=\bar \phi(z)$. Since $|(\phi^*)'|\le K_0$, $|\phi^*(z)-\bar \phi(z)|=|\phi^*(z)-\phi^*((\phi^*)^{-1}(\lambda))|\le K_0|(\bar \phi)^{-1}(\lambda)-(\phi^*)^{-1}(\lambda)|$, and the change of variables $\lambda=\bar \phi(z)$, with $\dd z\le\kappa_0^{-1}|\dd\lambda|$, gives $\nrm{\phi^*-\bar \phi}_{L^2(0,1)}\le K_0\kappa_0^{-1/2}\nrm{(\phi^*)^{-1}-(\bar \phi)^{-1}}_{L^2(a_1,a_0)}$, which proves the first inequality in \eqref{eq:rearr-inst}.

\emph{The vertical estimate.} The boundary values of $\phi$ do not depend on $x_h$, so $\wt\phi(x_h,0)=\wt\phi(x_h,1)=0$. For each $x_h$ and $z$, $|\wt\phi(x_h,z)|^2=2\int_0^z\wt\phi\,\pz\wt\phi\dd s\le2\int_0^1|\wt\phi\,\pz\wt\phi|(x_h,s)\dd s$. Taking the supremum in $z$, integrating in $x_h$ and using the Cauchy--Schwarz inequality in $\Omega$,
\begin{equation*}
\nrm{\wt\phi}_{L^2_hL^\infty_z}^2\le2\nrm{\wt\phi}_{L^2(\Omega)}\nrm{\pz\wt\phi}_{L^2(\Omega)} .
\end{equation*}
Together with $\nrm{\wt\phi}_{L^2(\Omega)}\le\pi^{-1}\nrm{\pz\wt\phi}_{L^2(\Omega)}$ (\Cref{lem:poincare}), this gives the second inequality in \eqref{eq:rearr-inst}.
\end{proof}

\subsection{Proof of \texorpdfstring{\Cref{thm:limit}}{Theorem 1.2}}

\begin{proof}[Proof of \Cref{thm:limit}]
\emph{Proof of} (i). By \Cref{lem:L-flux}, $\pt\bar\varrho\in C([0,\infty);L^2(0,1))$ and, by \eqref{eq:L-tail} with $t=0$, $\int_0^\infty\nrm{\pt\bar\varrho}_{L^2(0,1)}\dtau\le C\Es(\infty)^2<\infty$. By \eqref{eq:7.1}, $\tau\mapsto\pz(\ol{w\wt\varrho})(\tau)$ is therefore Bochner integrable on $(0,\infty)$ with values in $L^2(0,1)$, so $\vinf$ in \eqref{eq:limdef} is well defined. Integrating \eqref{eq:7.1} over $(0,t)$ and using $\bar\varrho(0)=\ol\varrho_{\mathrm{in}}$,
\begin{equation*}
\bar\varrho(t)-\vinf=\int_t^\infty\pz\big(\ol{w\wt\varrho}\big)(\tau)\dtau,
\end{equation*}
and \eqref{eq:L-tail} gives
\begin{equation}\label{eq:L-rate}
\nrm{\bar\varrho(t)-\vinf}_{L^2(0,1)}\le\int_t^\infty\nrm{\pt\bar\varrho(\tau)}_{L^2(0,1)}\dtau\le C\Es(\infty)^2(1+t)^{-\frac54},\qquad t\ge0 .
\end{equation}
In particular \eqref{eq:meanconv} holds. Next, $\pz^k\bar\varrho=\ol{\pz^k\varrho}$, so $\nrm{\bar\varrho(t)}_{H^4(0,1)}\le|\Tl|^{-1/2}\nrm{\varrho(t)}_{H^4(\Omega)}\le C\Es(\infty)\le C\mathsf N$ by \eqref{eq:L-Es}. If $t_n\to\infty$, a subsequence of $\bar\varrho(t_n)$ converges weakly in $H^4(0,1)$, and its limit is $\vinf$ because $\bar\varrho(t_n)\to\vinf$ in $L^2(0,1)$. Hence $\vinf\in H^4(0,1)$, $\bar\varrho(t)\rightharpoonup\vinf$ weakly in $H^4(0,1)$ as $t\to\infty$, and by weak lower semicontinuity
\begin{equation*}
\nrm{\vinf}_{H^4(0,1)}\le\liminf_{t\to\infty}\nrm{\bar\varrho(t)}_{H^4(0,1)}\le C\mathsf N .
\end{equation*}
By \Cref{thm:main}, $\varrho(t)=\pz\varrho(t)=0$ on $\pa\Omega$, so $\bar\varrho(t,\cdot)$ and $\pz\bar\varrho(t,\cdot)$ vanish at $z=0$ and $z=1$. Point evaluations of a function and of its derivative are continuous linear functionals on $H^4(0,1)$, hence weakly continuous; therefore $\vinf(0)=\vinf(1)=0$, and also $\pz\vinf(0)=\pz\vinf(1)=0$. This proves (i).

\emph{Proof of} (ii). Since $\Es(\infty)\le C\mathsf N$ by \eqref{eq:L-Es}, \eqref{eq:L-tail} and \eqref{eq:L-rate} give \eqref{eq:rate1eps}. For $s\in[0,4)$, the interpolation inequality on $(0,1)$, \eqref{eq:L-rate}, the $H^4(0,1)$ bounds in (i) and $\mathsf N\le1$ give
\begin{equation*}
\nrm{\bar\varrho(t)-\vinf}_{H^s(0,1)}\le C_s\nrm{\bar\varrho(t)-\vinf}_{L^2(0,1)}^{1-\frac s4}\nrm{\bar\varrho(t)-\vinf}_{H^4(0,1)}^{\frac s4}\le C_s\mathsf N(1+t)^{-\frac54(1-\frac s4)},
\end{equation*}
which is \eqref{eq:rateHs}.

\emph{Proof of} (iii). Write $\varrho-\vinf=\wt\varrho+(\bar\varrho-\vinf)$. The first part has zero horizontal mean and the second depends only on $z$, so they are orthogonal in $L^2(\Omega)$ (\Cref{lem:meanosc}), and
\begin{equation}\label{eq:L-orth}
\nrm{\varrho(t)-\vinf}_{L^2(\Omega)}^2=\nrm{\wt\varrho(t)}_{L^2(\Omega)}^2+|\Tl|\,\nrm{\bar\varrho(t)-\vinf}_{L^2(0,1)}^2 .
\end{equation}
By \eqref{eq:L-rhot}, \eqref{eq:L-Es} and \eqref{eq:L-rate}, both terms are bounded by $C\mathsf N^2(1+t)^{-2}$, so
\begin{equation}\label{eq:L-ratefull}
\nrm{\rho(t)-\rinf}_{L^2(\Omega)}=\nrm{\varrho(t)-\vinf}_{L^2(\Omega)}\le C\mathsf N(1+t)^{-1}.
\end{equation}
This implies \eqref{eq:ratefull} and the first estimate of \Cref{thm:limit}. Moreover $\nrm{\varrho(t)-\vinf}_{H^4(\Omega)}\le\nrm{\varrho(t)}_{H^4}+|\Tl|^{\frac12}\nrm{\vinf}_{H^4(0,1)}\le C\mathsf N$, and interpolation gives the $H^s$ bound in \eqref{eq:rateLinf}. For the uniform bounds we apply \Cref{lem:endpoint} to $h=\varrho(t)-\vinf$. Since $\vinf$ depends only on $z$, $\nh h=\nh\varrho(t)$, and \eqref{eq:L-ratefull}, \eqref{eq:5.3} and \eqref{eq:L-Es} give $\nrm{h}+\nrm{\nh h}\le C\mathsf N(1+t)^{-1}$, while $\nrm{h}_{H^4}\le C\mathsf N$. Hence
\begin{equation*}
\nrm{h}_{L^\infty}\le C\mathsf N(1+t)^{-\frac56},\qquad \nrm{\nabla h}_{L^\infty}\le C\mathsf N(1+t)^{-\frac12},
\end{equation*}
which are the $L^\infty$ and $W^{1,\infty}$ bounds in \eqref{eq:rateLinf}. Finally, the surface pressure $P$ has zero horizontal mean, so \eqref{eq:4.3}, \eqref{eq:5.3} and the Poincar\'e inequality on $\Tl$ give $\nrm{P(t)}_{H^1(\Tl)}\le C\nrm{\nh P(t)}_{L^2(\Tl)}\le C\mathsf N(1+t)^{-1}$. Since $q-q_\infty=P-g\int_0^z(\varrho-\vinf)\dd s$, with $\nrm{\int_0^zh\dd s}\le\nrm{h}$, $\pz\int_0^zh\dd s=h$ and $\nh\int_0^z(\varrho-\vinf)\dd s=\int_0^z\nh\varrho\dd s$,
\begin{equation*}
\nrm{q(t)-q_\infty}_{H^1(\Omega)}\le C\big(\nrm{P(t)}_{H^1(\Tl)}+\nrm{\varrho(t)-\vinf}+\nrm{\nh\varrho(t)}\big)\le C\mathsf N(1+t)^{-1},
\end{equation*}
which is \eqref{eq:ratepress}.

\emph{Proof of} (iv). By (i), $\vinf(0)=\vinf(1)=0$; hence $\rinf(0)=\rs(0)$ and $\rinf(1)=\rs(1)$.

\emph{Proof of} (v). \emph{Slope estimate.} By \Cref{lem:GN}(3) and \eqref{eq:unif}, $\sup_{t\ge0}\nrm{\pz\varrho(t)}_{L^\infty(\Omega)}\le C\sup_{t\ge0}\nrm{\varrho(t)}_{H^4}\le C\varepsilon_0$ (the case $t=0$ is the bound for $\pz\vin$), and by $H^4(0,1)\hookrightarrow C^1([0,1])$ and (i), $\nrm{\pz\vinf}_{L^\infty(0,1)}\le C\nrm{\vinf}_{H^4(0,1)}\le C\varepsilon_0$. This gives the constant $\cc{*}$, which depends only on $\rs,\mu,g,f,l$. Assume now \eqref{eq:slope}, after decreasing $\varepsilon_0$ if necessary (the constants in \eqref{eq:unif} and $\cc{*}$ do not depend on $\varepsilon_0$), and put $m\eqdef\min_{[0,1]}(-\rs')>0$; then $\cc{*}\varepsilon_0\le\frac m2$. Since $\vin\in H^4(\Omega)\subset C^1(\ol\Omega)$ and $\vinf\in C^1([0,1])$,
\begin{equation}\label{eq:L-slope}
\pz\rinf=\rs'+\pz\vinf\le-\frac m2\quad\text{on }[0,1],\qquad \pz\rin=\rs'+\pz\vin\le-\frac m2\quad\text{on }\ol\Omega .
\end{equation}
In particular $\pz\rinf<0$ on $[0,1]$, and by (iv) $\rinf$ is a strictly decreasing continuous bijection of $[0,1]$ onto $[\rs(1),\rs(0)]$. Likewise, since $\vin=0$ on $\pa\Omega$, for every $x_h\in\Tl$ the map $z\mapsto\rin(x_h,z)$ decreases strictly from $\rs(0)$ to $\rs(1)$, so that
\begin{equation}\label{eq:L-range}
\rs(1)\le\rin\le\rs(0)\quad\text{on }\ol\Omega,\qquad \rs(1)=\min_{[0,1]}\rs>0 .
\end{equation}

\emph{Conservation of the distribution function.} The total density $\rho=\rs+\varrho$ satisfies $\rho_t+\bU\cdot\nabla\rho=0$, by \eqref{PE}$_1$ and $\nabla\rs=\rs'e_z$. The field $\bU=(v,w)\in C([0,T];H^3)$ is divergence-free by \eqref{PE}$_3$, and $\bU\cdot n=\pm w=0$ on $\pa\Omega$ by \Cref{lem:propagation}. Let $\phi\in C^1(\R)$ with $\phi'$ bounded. Since $\rho\in C^1([0,T];L^2(\Omega))\cap C([0,T];C(\ol\Omega))$, the function $t\mapsto\int_\Omega\phi(\rho)\dx$ is differentiable, and
\begin{equation*}
\ddt\int_\Omega\phi(\rho)\dx=\int_\Omega\phi'(\rho)\rho_t\dx=-\int_\Omega\bU\cdot\nabla\big(\phi(\rho)\big)\dx=\int_\Omega\divv\bU\,\phi(\rho)\dx-\int_{\pa\Omega}\bU\cdot n\,\phi(\rho)\dd S=0 .
\end{equation*}
Hence $\int_\Omega\phi(\rho(t))\dx=\int_\Omega\phi(\rin)\dx$ for all $\phi\in C_c^\infty(\R)$, so the images of the Lebesgue measure on $\Omega$ under $\rho(t)$ and under $\rin$ coincide. In particular
\begin{equation}\label{eq:L-distr}
|\{\rho(t)>\lambda\}|=|\{\rin>\lambda\}|\eqdef D(\lambda),\qquad\lambda\in\R,\ t\ge0 .
\end{equation}

\emph{Identification of the limit.} We regard $\rinf$ as a function on $\Omega$ independent of $x_h$ and set $D_\infty(\lambda)\eqdef|\{\rinf>\lambda\}|$. For $a,b\in\R$ one has $|a-b|=\int_\R|\mathbf1_{\{a>\lambda\}}-\mathbf1_{\{b>\lambda\}}|\dd\lambda$. By \eqref{eq:L-distr}, Tonelli's theorem and \eqref{eq:L-ratefull},
\begin{equation*}
\begin{aligned}
\int_\R|D(\lambda)-D_\infty(\lambda)|\dd\lambda&\le\int_\R\int_\Omega\big|\mathbf1_{\{\rho(t)>\lambda\}}-\mathbf1_{\{\rinf>\lambda\}}\big|\dx\dd\lambda\\
&=\nrm{\rho(t)-\rinf}_{L^1(\Omega)}\le|\Omega|^{\frac12}\nrm{\rho(t)-\rinf}_{L^2(\Omega)}\to0\qquad(t\to\infty).
\end{aligned}
\end{equation*}
Hence $D=D_\infty$ a.e. Both functions are nonincreasing and right-continuous (if $\lambda_n\downarrow\lambda$, then $\{f>\lambda_n\}\uparrow\{f>\lambda\}$), so $D=D_\infty$ everywhere:
\begin{equation}\label{eq:L-DD}
|\{\rin>\lambda\}|=|\{\rinf>\lambda\}|=|\Tl|\,\big|\{z\in(0,1):\rinf(z)>\lambda\}\big|,\qquad\lambda\in\R .
\end{equation}
Fix $z\in(0,1)$ and $\lambda\in\R$. Write $I_\lambda\eqdef\{z'\in(0,1):\rinf(z')>\lambda\}$. If $\lambda\le\rinf(z)$, then $(0,z)\subset I_\lambda$, because $\rinf$ is strictly decreasing, so $|I_\lambda|\ge z$. If $\lambda>\rinf(z)$, then $z'\in I_\lambda$ forces $z'<z$, and by continuity $\rinf<\lambda$ on a neighbourhood of $z$; hence $I_\lambda\subset(0,z-\vartheta)$ for some $\vartheta>0$, and $|I_\lambda|<z$. By \eqref{eq:L-DD}, $D(\lambda)\ge|\Tl|z$ holds if and only if $\lambda\le\rinf(z)$. Since $\rinf(z)>0$,
\begin{equation*}
\rin^*(z)=\int_0^\infty\mathbf1_{\{D(\lambda)\ge|\Tl|z\}}\dd\lambda=\big|(0,\rinf(z)]\big|=\rinf(z),\qquad0<z<1,
\end{equation*}
which is \eqref{eq:rearr}, with $|\Tl|=4\pi^2l^2$. The endpoint values were obtained in (iv).

\emph{Strict stratification for all times, and the faster rate of the mean.} Let $M_1\eqdef\max_{[0,1]}(-\rs')$. By the slope estimate and \eqref{eq:slope}, $|\pz\varrho(t)|\le\frac m2$, hence
\begin{equation}\label{eq:L-strat}
\frac m2\le-\pz\rho(t,x_h,z)\le\frac32M_1\qquad\text{on }\ol\Omega,\ t\ge0,
\end{equation}
and $|\pz^2\rho(t)|\le\nrm{\rs''}_{L^\infty}+C\nrm{\varrho(t)}_{H^4}\le K_0$ for a constant $K_0$ depending only on $\rs,\mu,g,f,l$. By \eqref{eq:L-distr}, $\rho(t)$ and $\rin$ have the same distribution function, so $\rho(t)^*=\rin^*=\rinf$ for every $t\ge0$. We apply \Cref{lem:rearr-inst} to $\phi=\rho(t)\in H^4(\Omega)\subset C^2(\ol\Omega)$, with $a_0=\rs(0)$, $a_1=\rs(1)$, $\kappa_0=\frac m2$ and $K_0\ge\frac32M_1$; here $\bar\phi=\rs+\bar\varrho(t)$, $\wt\phi=\wt\varrho(t)$ and $\phi^*=\rinf$. Thus
\begin{equation*}
\nrm{\bar\varrho(t)-\vinf}_{L^2(0,1)}\le C\nrm{\wt\varrho(t)}^{\frac12}\nrm{\pz\wt\varrho(t)}^{\frac32}.
\end{equation*}
Since $\varrho$ and $\bar\varrho$ are clamped, so is $\wt\varrho$: $\wt\varrho=\pz\wt\varrho=0$ at $z=0,1$. Two integrations by parts in $z$, whose boundary terms vanish by these traces, give $\nrm{\pz\wt\varrho}^2=-\int_\Omega\wt\varrho\,\pz^2\wt\varrho\dx\le\nrm{\wt\varrho}\nrm{\pz^2\wt\varrho}$ and $\nrm{\pz^2\wt\varrho}^2=\int_\Omega\wt\varrho\,\pz^4\wt\varrho\dx\le\nrm{\wt\varrho}\nrm{\pz^4\wt\varrho}$; hence $\nrm{\pz\wt\varrho}\le\nrm{\wt\varrho}^{3/4}\nrm{\wt\varrho}_{H^4}^{1/4}$. By \eqref{eq:5.3}, \eqref{eq:L-Es} and $\nrm{\wt\varrho}_{H^4}\le2\nrm{\varrho}_{H^4}\le2\Es(\infty)$ (\Cref{lem:meanosc}), we have $\nrm{\wt\varrho(t)}\le C\mathsf N(1+t)^{-1}$ and $\nrm{\wt\varrho(t)}_{H^4}\le C\mathsf N$. Therefore
\begin{equation*}
\nrm{\bar\varrho(t)-\vinf}_{L^2(0,1)}\le C\nrm{\wt\varrho(t)}^{\frac{13}8}\nrm{\wt\varrho(t)}_{H^4}^{\frac38}\le C\mathsf N^2(1+t)^{-\frac{13}8},
\end{equation*}
which is \eqref{eq:rate138}.

\emph{Proof of} (vi). Assume \eqref{eq:slope} and set $\epsilon_{\mathrm{in}}\eqdef\nrm{\vin}_{H^4}\le\varepsilon_0$. By \eqref{eq:slope}, $\varepsilon_0\le m/(2\cc{*})$ is bounded by a constant depending only on $\rs,\mu,g,f,l$. In this proof $C$ denotes such constants, and $X=Y+O(\epsilon_{\mathrm{in}}^k)$ means $|X-Y|\le C\epsilon_{\mathrm{in}}^k$ uniformly in all variables. By \Cref{lem:GN}(3), $\nrm{\vin}_{W^{2,\infty}(\Omega)}\le C\epsilon_{\mathrm{in}}$.

\emph{Step 1.} Let $\chi:[\rs(1),\rs(0)]\to[0,1]$ be the inverse of $\rs$, so that $\chi(\rs(z))=z$. Then $\chi\in C^5$, $-m^{-1}\le\chi'\le-M_1^{-1}$ with $M_1=\max_{[0,1]}(-\rs')$ as above, and differentiating twice,
\begin{equation*}
\chi'(\rs)=\frac1{\rs'}=-\sgm,\qquad \chi''(\rs)=-\frac{\sgm'}{\rs'}=\sgm\sgm'=\sgm^3\rs'',
\end{equation*}
where we used $\sgm'=\rs''/\rs'^2=\sgm^2\rs''$. By \eqref{eq:L-range}, $\Phi\eqdef\chi\circ\rin\in C^2(\ol\Omega)$ is defined, $\Phi(x_h,0)=0$, $\Phi(x_h,1)=1$, and by \eqref{eq:L-slope}
\begin{equation*}
\pz\Phi=\chi'(\rin)\,\pz\rin\ge\frac m{2M_1}\eqdef c_\Phi>0 .
\end{equation*}
For $t\in[0,1]$, $\rs+t\vin=(1-t)\rs+t\rin\in[\rs(1),\rs(0)]$ by \eqref{eq:L-range}, so Taylor's formula applies to $\varphi\eqdef\Phi-z=\chi(\rs+\vin)-\chi(\rs)$:
\begin{equation*}
\varphi=-\sgm\vin+\varphi_2,\qquad \varphi_2=\vin^2\int_0^1(1-t)\chi''(\rs+t\vin)\dt=\tfrac12\sgm^3\rs''\vin^2+O(\epsilon_{\mathrm{in}}^3),
\end{equation*}
with $\nrm{\varphi}_{W^{2,\infty}(\Omega)}\le C\epsilon_{\mathrm{in}}$ and $\nrm{\varphi_2}_{W^{1,\infty}(\Omega)}\le C\epsilon_{\mathrm{in}}^2$.

\emph{Step 2.} For each $x_h$, $\Phi(x_h,\cdot)$ is an increasing $C^2$ bijection of $[0,1]$; let $\Psi(x_h,\cdot)$ be its inverse. By the inverse function theorem, $\Psi\in C^2(\Tl\times[0,1])$ and
\begin{equation*}
\pa_\zeta\Psi=\frac1{\pz\Phi(x_h,\Psi)},\qquad \pa_\zeta^2\Psi=-\frac{\pz^2\Phi(x_h,\Psi)}{\big(\pz\Phi(x_h,\Psi)\big)^3} .
\end{equation*}
The function $\psi\eqdef\Psi-\zeta$ satisfies $\psi=-\varphi(x_h,\zeta+\psi)$, hence $|\psi|\le\nrm{\varphi}_{L^\infty}\le C\epsilon_{\mathrm{in}}$. Taylor's formula in $z$ around $\zeta$, with remainder bounded by $\frac12\nrm{\pz^2\varphi}_{L^\infty}\psi^2$, gives $\psi=-\varphi-\psi\pz\varphi+O(\epsilon_{\mathrm{in}}^3)$, all functions evaluated at $(x_h,\zeta)$. Since $|\psi\pz\varphi|\le C\epsilon_{\mathrm{in}}^2$, first $\psi=-\varphi+O(\epsilon_{\mathrm{in}}^2)$, and then $\psi=-\varphi+\varphi\pz\varphi+O(\epsilon_{\mathrm{in}}^3)$. By Step 1, $\varphi=-\sgm\vin+O(\epsilon_{\mathrm{in}}^2)$ and $\pz\varphi=-\pz(\sgm\vin)+O(\epsilon_{\mathrm{in}}^2)$, so $\varphi\pz\varphi=\sgm\vin\pz(\sgm\vin)+O(\epsilon_{\mathrm{in}}^3)=\frac12\pz(\sgm^2\vin^2)+O(\epsilon_{\mathrm{in}}^3)$. Hence
\begin{equation}\label{eq:L-Psi}
\Psi(x_h,\zeta)=\zeta+\sgm\vin+\tfrac12\pz\big(\sgm^2\vin^2\big)-\tfrac12\sgm^3\rs''\vin^2+O(\epsilon_{\mathrm{in}}^3),
\end{equation}
where $\vin$, its derivatives and the coefficients are evaluated at $(x_h,\zeta)$. Moreover, $\pa_\zeta\psi=-\pz\varphi(x_h,\Psi)/\pz\Phi(x_h,\Psi)=-\pz\varphi(x_h,\Psi)+(\pz\varphi(x_h,\Psi))^2/\pz\Phi(x_h,\Psi)$ and $|\pz\varphi(x_h,\Psi)-\pz\varphi(x_h,\zeta)|\le\nrm{\pz^2\varphi}_{L^\infty}|\psi|\le C\epsilon_{\mathrm{in}}^2$, so
\begin{equation}\label{eq:L-dPsi}
\pa_\zeta\psi(x_h,\zeta)=\pz(\sgm\vin)(x_h,\zeta)+O(\epsilon_{\mathrm{in}}^2),\qquad |\pa_\zeta^2\psi|\le c_\Phi^{-3}\nrm{\pz^2\varphi}_{L^\infty}\le C\epsilon_{\mathrm{in}} .
\end{equation}

\emph{Step 3.} Let $\Psi_*(\zeta)\eqdef|\Tl|^{-1}\int_{\Tl}\Psi(x_h,\zeta)\dxh$. Then $\Psi_*\in C^2([0,1])$ is strictly increasing, $\Psi_*(0)=0$ and $\Psi_*(1)=1$. Let $\zeta\in[0,1]$. By \eqref{eq:L-range} and since $\chi$ is strictly decreasing, $\rin(x_h,z)>\rs(\zeta)$ holds if and only if $\Phi(x_h,z)<\zeta$, that is, if and only if $z<\Psi(x_h,\zeta)$. Hence $|\{\rin>\rs(\zeta)\}|=|\Tl|\Psi_*(\zeta)$. On the other hand, by \eqref{eq:L-DD} and since $\rinf$ is a strictly decreasing bijection of $[0,1]$ onto $[\rs(1),\rs(0)]$, $|\{\rin>\rs(\zeta)\}|=|\Tl|\,\rinf^{-1}(\rs(\zeta))$. Thus $\Psi_*=\rinf^{-1}\circ\rs$, that is,
\begin{equation}\label{eq:7.6}
\rinf=\rs\circ\Psi_*^{-1}\qquad\text{on }[0,1].
\end{equation}
Write $a\eqdef\ol\varrho_{\mathrm{in}}$, $b\eqdef\wt\varrho_{\mathrm{in}}$ and $Q\eqdef\ol{\vin^2}=a^2+\ol{b^2}$. Since $\sgm$ and $\rs''$ depend only on $z$, averaging \eqref{eq:L-Psi} and \eqref{eq:L-dPsi} over $\Tl$ gives, for $\alpha\eqdef\Psi_*-\zeta=\ol\psi$,
\begin{equation*}
\alpha=A+B+O(\epsilon_{\mathrm{in}}^3),\qquad A\eqdef\sgm a,\qquad B\eqdef\tfrac12\big(\sgm^2Q\big)'-\tfrac12\sgm^3\rs''Q,
\end{equation*}
and $\alpha'=A'+O(\epsilon_{\mathrm{in}}^2)$, $|\alpha''|\le C\epsilon_{\mathrm{in}}$. Since $|a|+|a'|\le\nrm{\vin}_{W^{1,\infty}}$ and $|Q|+|Q'|\le2\nrm{\vin}_{W^{1,\infty}}^2$, we also have $|A|+|A'|\le C\epsilon_{\mathrm{in}}$, $|B|\le C\epsilon_{\mathrm{in}}^2$, and hence $|\alpha|+|\alpha'|\le C\epsilon_{\mathrm{in}}$.

\emph{Step 4.} Let $\beta(z)\eqdef\Psi_*^{-1}(z)-z$. Then $\beta=-\alpha(z+\beta)$ and $|\beta|\le\nrm{\alpha}_{L^\infty}\le C\epsilon_{\mathrm{in}}$. Taylor's formula gives $\alpha(z+\beta)=\alpha+\beta\alpha'+O(\nrm{\alpha''}_{L^\infty}\beta^2)=\alpha+\beta\alpha'+O(\epsilon_{\mathrm{in}}^3)$, evaluated at $z$. Thus $\beta=-\alpha+O(\epsilon_{\mathrm{in}}^2)$, and then $\beta=-\alpha+\alpha\alpha'+O(\epsilon_{\mathrm{in}}^3)$. Since $\alpha-A=O(\epsilon_{\mathrm{in}}^2)$, $\alpha'-A'=O(\epsilon_{\mathrm{in}}^2)$ and $A,A'=O(\epsilon_{\mathrm{in}})$, we have $\alpha\alpha'=AA'+O(\epsilon_{\mathrm{in}}^3)$, and therefore
\begin{equation*}
\beta=-A-B+AA'+O(\epsilon_{\mathrm{in}}^3),\qquad\text{in particular}\quad \beta=-A+O(\epsilon_{\mathrm{in}}^2).
\end{equation*}

\emph{Step 5.} By \eqref{eq:7.6} and Taylor's formula for $\rs\in C^3$, $\rinf(z)=\rs(z+\beta)=\rs+\rs'\beta+\frac12\rs''\beta^2+O(\epsilon_{\mathrm{in}}^3)$, and $\beta^2=A^2+O(\epsilon_{\mathrm{in}}^3)$. Hence
\begin{equation*}
\vinf=\rinf-\rs=-\rs'A+\rs'\big(AA'-B\big)+\tfrac12\rs''A^2+O(\epsilon_{\mathrm{in}}^3).
\end{equation*}
Since $\rs'\sgm=-1$, we have $-\rs'A=a=\ol\varrho_{\mathrm{in}}$. Next, $AA'=\frac12(\sgm^2a^2)'$ and $Q-a^2=\ol{b^2}$, so
\begin{equation*}
AA'-B=-\tfrac12\big(\sgm^2\ol{b^2}\big)'+\tfrac12\sgm^3\rs''Q,\qquad \rs'\big(AA'-B\big)=\frac1{2\sgm}\big(\sgm^2\ol{b^2}\big)'-\tfrac12\sgm^2\rs''Q,
\end{equation*}
and, adding $\frac12\rs''A^2=\frac12\sgm^2\rs''a^2$ and using $\sgm^2\rs''=\sgm'$,
\begin{equation*}
\rs'\big(AA'-B\big)+\tfrac12\rs''A^2=\frac1{2\sgm}\big(\sgm^2\ol{b^2}\big)'-\tfrac12\sgm'\ol{b^2}=\sgm'\ol{b^2}+\tfrac12\sgm\,\ol{b^2}{}'-\tfrac12\sgm'\ol{b^2}=\tfrac12\pz\big(\sgm\,\ol{b^2}\big).
\end{equation*}
Therefore $\vinf=\ol\varrho_{\mathrm{in}}+\frac12\pz\big(\sgm\,\ol{\wt\varrho_{\mathrm{in}}^{\,2}}\big)+O(\epsilon_{\mathrm{in}}^3)$ uniformly on $[0,1]$, which is \eqref{eq:quadexp}.
\end{proof}

\begin{remark}\label{rem:final}
(i) The limit state is determined by the distribution function $\lambda\mapsto|\{\rin>\lambda\}|$ alone. By (v), $\rinf$ is the decreasing vertical rearrangement of $\rin$: every superlevel set $\{\rin>\lambda\}$ is replaced by the horizontal slab $\Tl\times(0,z_\lambda)$ of the same volume, $z_\lambda=|\{\rin>\lambda\}|/|\Tl|$. In particular $\rinf$ does not depend on the initial velocity $\uin$, on the viscosity $\mu$, on $g$, or on the Coriolis parameter $f$. Rotation does not enter the transport equation for $\rho$; it affects the decay estimates and the smallness threshold $\varepsilon_0$, hence the time scale on which the limit is approached, but not the limit itself. The identification follows the mechanism found by Dalibard, Guillod and Leblond \cite{DalibardGuillodLeblond2025} for the Stokes--transport system in a channel, where the limit profile is also the decreasing vertical rearrangement of the initial density. The only inputs are that $\rho$ is transported by a divergence-free velocity field tangent to the lids, so that its distribution function is conserved, and that $\rho(t)$ converges strongly to a strictly stratified profile. No property of the momentum equation is used.

(ii) The rate $(1+t)^{-1}$ in \eqref{eq:ratefull} is the rate obtained for the Stokes--transport system in \cite{DalibardGuillodLeblond2025}. It is optimal for $H^4$ data, for every $f$ (\Cref{thm:O-nonlin}). The mean converges faster. The length of its path satisfies $\int_t^\infty\nrm{\pt\bar\varrho}\dtau\lesssim\mathsf N^2(1+t)^{-5/4}$ by \eqref{eq:L-tail}: both flux terms in \eqref{eq:7.2} are quadratic, and each contains one factor of the velocity dissipation $\Ff$, whose tail \eqref{eq:L-Ftail} decays like $(1+t)^{-2}$, and one density factor that either decays pointwise, \eqref{eq:L-rhot}, or is itself controlled by $\Ff$, \eqref{eq:L-rhoF}. Under \eqref{eq:slope} the distance decays faster still, like $(1+t)^{-13/8}$ in \eqref{eq:rate138}: by \Cref{lem:rearr-inst}, the horizontal mean of a strictly stratified profile agrees with its rearrangement up to a quadratic error in the oscillation, and the rearrangement of $\rho(t)$ is $\rinf$ at every time. Measuring this error in $L^2_hL^\infty_z$ rather than in $L^\infty$ avoids the failure of $H^1(\T^2)\subset L^\infty(\T^2)$, and vertical interpolation for the clamped oscillation gives the exponent $\frac12+\frac32\cdot\frac34=\frac{13}8$. Thus the distance to the limit is smaller than the bound \eqref{eq:L-tail} for the remaining length $\int_t^\infty\nrm{\pt\bar\varrho}\dtau$ of the path of the mean. Whether the exponents $\frac54$ and $\frac{13}8$ for the mean are sharp is left open. The optimal modes of \Cref{sec:opt} have zero horizontal mean, so they do not directly decide this.

(iii) The limit state need not vanish, since the horizontal mean of the density perturbation need not decay. Accordingly, $\rs$ is not an isolated steady state: for every $\varrho_s=\varrho_s(z)\in H^4(0,1)$ with $\varrho_s=\varrho_s'=0$ at $z=0,1$, the triple $(\varrho,v,q)=(\varrho_s,0,-g\int_0^z\varrho_s)$ is a steady solution of \eqref{PE} for every $f$, and $\rs+\varrho_s$ is strictly decreasing if $\nrm{\varrho_s'}_{L^\infty}<\min(-\rs')$. The global solution selects $\rinf=\rin^*$. By \eqref{eq:quadexp}, the first correction to $\ol\varrho_{\mathrm{in}}$ is the quadratic flux $\frac12\pz\big(\sgm\,\ol{\wt\varrho_{\mathrm{in}}^{\,2}}\big)$, which depends only on the oscillatory part of the initial density.
\end{remark}

\section{Available potential energy and the Lagrangian picture}\label{sec:energy}

In this section we prove \Cref{thm:energy}. Throughout, the hypotheses of \Cref{thm:limit} hold, including \eqref{eq:slope}, and $\mathsf N=\nrm{\vin}_{H^4}+\nrm{\uin}_{H^4}$. Let $m=\min_{[0,1]}(-\rs')$ and $M_1=\max_{[0,1]}(-\rs')$ be as in \Cref{thm:energy}. By the proof of \Cref{thm:limit}(v), $\rin$ and $\rho(t)$ take their values in $[\rs(1),\rs(0)]$, and $\rinf\in C^3([0,1])$ is a decreasing bijection of $[0,1]$ onto $[\rs(1),\rs(0)]$ with
\begin{equation}\label{eq:E-slope}
\frac m2\le-\rinf'(z)\le\frac32M_1,\qquad z\in[0,1].
\end{equation}
We write $\rinf^{-1}$ for its inverse. Two functions $\phi_1,\phi_2$ on $\Omega$ are called equimeasurable if $|\{\phi_1>\lambda\}|=|\{\phi_2>\lambda\}|$ for every $\lambda\in\R$; then $\int_\Omega\Phi(\phi_1)\dx=\int_\Omega\Phi(\phi_2)\dx$ for every continuous $\Phi$ on the common range. By \Cref{thm:limit}(v), $\rin$, $\rho(t)$ and $\rinf$ are equimeasurable.

\subsection{The limit minimizes the potential energy}

\begin{proof}[Proof of \Cref{thm:energy}(i)]
Let $\phi$ be equimeasurable with $\rin$. Then $\phi$ takes its values in $[\rs(1),\rs(0)]$ almost everywhere. Since $\rinf$ is equimeasurable with $\phi$, the function $\Phi(r)\eqdef\int_{\rs(1)}^r\rinf^{-1}(\lambda)\dd\lambda$ gives $\int_\Omega\Phi(\phi)\dx=\int_\Omega\Phi(\rinf)\dx$, that is, $\int_\Omega\int_{\rinf(z)}^{\phi(x_h,z)}\rinf^{-1}(\lambda)\dd\lambda\dx=0$. Since $\int_{\rinf(z)}^{\phi}z\dd\lambda=z(\phi-\rinf(z))$, we obtain
\begin{equation}\label{eq:E-pot}
g\int_\Omega z\big(\phi-\rinf\big)\dx=g\int_\Omega h_z\big(\phi(x_h,z)\big)\dx,\qquad h_z(r)\eqdef\int_{\rinf(z)}^r\big(z-\rinf^{-1}(\lambda)\big)\dd\lambda .
\end{equation}
For fixed $z$, $h_z(\rinf(z))=h_z'(\rinf(z))=0$ and $h_z''(\lambda)=-(\rinf^{-1})'(\lambda)=1/|\rinf'(\rinf^{-1}(\lambda))|\in[\frac2{3M_1},\frac2m]$ by \eqref{eq:E-slope}. Taylor's formula gives
\begin{equation*}
\frac{(r-\rinf(z))^2}{3M_1}\le h_z(r)\le\frac{(r-\rinf(z))^2}{m},\qquad r\in[\rs(1),\rs(0)],
\end{equation*}
and \eqref{eq:E-pot} yields \eqref{eq:E-min}. If $\phi\ne\rinf$ on a set of positive measure, the left-hand side of \eqref{eq:E-min} is positive; hence $\rinf$ is the unique minimizer.
\end{proof}

\subsection{The energy identity}

\begin{proof}[Proof of \Cref{thm:energy}(ii)]
On every interval $[0,T]$ the solution belongs to the class \eqref{eq:localclass}, so $v\in C^1([0,T];H^2)$ and $\rho\in C^1([0,T];H^3)$, and $\mathcal H\in C^1([0,T])$. Testing \eqref{PE}$_2$ with $v$, and using $\rho_t+\bU\cdot\nabla\rho=0$, $\divv\bU=0$, $v|_{\pa\Omega}=0$, $w|_{\pa\Omega}=0$ and $v^\perp\cdot v=0$, so that the Coriolis force does no work,
\begin{equation*}
\begin{aligned}
\frac12\ddt\nrm{\sqrt\rho v}_{L^2}^2+\mu\nrm{\nabla v}_{L^2}^2&=-\int_\Omega\nh q\cdot v\dx=\int_\Omega q\,\divh v\dx\\
&=-\int_\Omega q\,\pz w\dx=\int_\Omega\pz q\,w\dx=-g\int_\Omega\varrho w\dx .
\end{aligned}
\end{equation*}
For the potential energy, $\rho_t=-\divv(\rho\bU)$ and $w|_{\pa\Omega}=0$ give
\begin{equation*}
\ddt\,g\int_\Omega z\rho\dx=-g\int_\Omega z\divv(\rho\bU)\dx=g\int_\Omega\rho w\dx=g\int_\Omega\varrho w\dx,
\end{equation*}
because $\int_\Omega\rs(z)w\dx=\int_0^1\rs(z)\int_{\Tl}w\dxh\dd z=0$, the function $w=-\int_0^z\divh v$ having zero horizontal mean. Since $\int_\Omega z\rinf\dx$ is a constant, adding the two identities gives $\mathcal H'+\mu\nrm{\nabla v}_{L^2}^2=0$, and integration over $[t,T]$ gives \eqref{eq:E-identity}.

Since $\rs(1)\le\rho\le\rs(0)$, the kinetic part of $\mathcal H$ lies between $\frac12\rs(1)\nrm{v}^2$ and $\frac12\rs(0)\nrm{v}^2$. Since $\rho(t)$ is equimeasurable with $\rin$, part (i) applies to $\phi=\rho(t)$ and bounds the potential part between $\frac g{3M_1}\nrm{\rho(t)-\rinf}^2$ and $\frac gm\nrm{\rho(t)-\rinf}^2$. This proves \eqref{eq:E-coerc} with $c=\min\{\frac12\rs(1),\frac g{3M_1}\}$ and $C=\max\{\frac12\rs(0),\frac gm\}$. By \eqref{eq:decay} and \eqref{eq:ratefull}, $\nrm{v(t)}+\nrm{\rho(t)-\rinf}\le C\mathsf N(1+t)^{-1}$, so $\mathcal H(t)\le C\mathsf N^2(1+t)^{-2}$; in particular $\mathcal H(T)\to0$, and letting $T\to\infty$ in \eqref{eq:E-identity} gives \eqref{eq:E-tail}. At $t=0$, \eqref{eq:E-tail} is \eqref{eq:E-total}, because $\rinf=\rin^*$.
\end{proof}

\subsection{The flow map}

\begin{proof}[Proof of \Cref{thm:energy}(iii)]
\emph{Step 1: integrability of the velocity gradient.} The velocity vanishes on $\pa\Omega$, so $\nrm{v}_{L^\infty}\le\nrm{\pz v}_{L^\infty}$ and $\nrm{w}_{L^\infty}\le\nrm{\pz w}_{L^\infty}=\nrm{\divh v}_{L^\infty}\le2\nrm{\nabla v}_{L^\infty}$. Since $\nabla v=\nabla\wt v+\pz\bar v\otimes e_z$, Steps 1--3 of the proof of \Cref{prop:Fnat} give
\begin{equation*}
\nrm{\bU}_{W^{1,\infty}}\le C\big(\Ff+\Ff^{\frac23}\Fu^{\frac13}+\Ff^{\frac7{12}}\Fu^{\frac5{12}}\big).
\end{equation*}
By \eqref{eq:L-Ftail} and $\int_0^\infty\Fu^2\dtau\le\Es(\infty)^2$, \Cref{lem:dyadic} with $(a,b,c)=(1,0,0)$, $(\frac23,0,\frac13)$ and $(\frac7{12},0,\frac5{12})$, for which $\gamma=\frac12$, $\frac16$ and $\frac1{12}$, gives
\begin{equation}\label{eq:E-Uint}
\int_t^\infty\nrm{\bU(\tau)}_{W^{1,\infty}}\dtau\le C\Es(\infty)(1+t)^{-\frac1{12}}\le C\mathsf N(1+t)^{-\frac1{12}},\qquad t\ge0 .
\end{equation}
The worst exponent comes from $\nh w$, which costs one horizontal derivative of $v$ more than $\nabla v$. This preliminary rate is used only to construct the limit map and to bound $D\bX$ uniformly; in Steps 4 and 5, which give better rates, $\nh w$ enters only through this bound.

\emph{Step 2: the flow and its limit.} We work with the lift of $\Omega$ to $\R^2\times[0,1]$, on which $\bU$ is $2\pi l$-periodic in $x_1$ and $x_2$. On every interval $[0,T]$, $\bU\in C([0,T];H^3)\subset C([0,T];C^1(\ol\Omega))$, and $\bU=0$ on $\pa\Omega$. Hence \eqref{eq:E-flow} has a unique $C^1$ flow; every point of $\pa\Omega$ is fixed, the planes $z=0$ and $z=1$ are invariant, $\bX(t,x_h+2\pi lk,z)=\bX(t,x_h,z)+(2\pi lk,0)$ for $k\in\Z^2$, $\bX(t)$ is a bijection of $\R^2\times[0,1]$ whose inverse is the backward flow, and $\det D\bX(t)=1$ because $\divv\bU=0$. The matrices $D\bX$ and $(D\bX)^{-1}$ satisfy $\pt D\bX=(D\bU)(t,\bX)D\bX$ and $\pt(D\bX)^{-1}=-(D\bX)^{-1}(D\bU)(t,\bX)$, so by \eqref{eq:E-Uint} and Gronwall's inequality
\begin{equation*}
\sup_{t\ge0}\big(\nrm{D\bX(t)}_{L^\infty}+\nrm{(D\bX(t))^{-1}}_{L^\infty}\big)\le2\exp\Big(\int_0^\infty\nrm{D\bU}_{L^\infty}\dtau\Big)\le C .
\end{equation*}
For $t<T$, integrating \eqref{eq:E-flow} and its variational equation over $[t,T]$ gives
\begin{equation*}
\nrm{\bX(T)-\bX(t)}_{L^\infty}+\nrm{D\bX(T)-D\bX(t)}_{L^\infty}\le C\int_t^T\nrm{\bU(\tau)}_{W^{1,\infty}}\dtau\le C\mathsf N(1+t)^{-\frac1{12}}.
\end{equation*}
Hence $\bX(t)$ converges in $C^1$ to a $C^1$ map $\bX_\infty$, at the preliminary rate $(1+t)^{-1/12}$; $\bX_\infty$ fixes every point of $\pa\Omega$, satisfies the same periodicity relation, and $\det D\bX_\infty=1$. Since the inverse maps $\bX(t)^{-1}$ have derivatives bounded by $C$ and $\R^2\times[0,1]$ is convex, $|\bX(t,\mathbf a)-\bX(t,\mathbf b)|\ge C^{-1}|\mathbf a-\mathbf b|$, and this passes to the limit; so $\bX_\infty$ is bi-Lipschitz, in particular injective. For $\mathbf y\in\R^2\times[0,1]$ the points $\bX(t_n)^{-1}(\mathbf y)$ stay in a bounded set, because $|\bX(t)(\mathbf a)-\mathbf a|\le\int_0^\infty\nrm{\bU}_{L^\infty}\le C$; a subsequence converges to some $\mathbf a$, and $\bX_\infty(\mathbf a)=\mathbf y$ by the uniform convergence and the uniform Lipschitz bound. Thus $\bX_\infty$ is a bijection with $\det D\bX_\infty=1$, hence a volume-preserving $C^1$ diffeomorphism.

\emph{Step 3: the final height.} Since $\rho\in C^1$ solves $\rho_t+\bU\cdot\nabla\rho=0$, $\rho(t,\bX(t,\mathbf a))=\rin(\mathbf a)$. By \eqref{eq:rateLinf}, $|\rinf(X^z(t,\mathbf a))-\rin(\mathbf a)|\le\nrm{\rho(t)-\rinf}_{L^\infty}\to0$, and $X^z(t,\mathbf a)\to X^z_\infty(\mathbf a)$; hence $\rinf(X^z_\infty(\mathbf a))=\rin(\mathbf a)$, which is \eqref{eq:E-height} because $\rin(\mathbf a)\in[\rs(1),\rs(0)]$. Consequently
\begin{equation*}
\big(\rho(t)-\rinf\big)\big(\bX(t,\mathbf a)\big)=\rin(\mathbf a)-\rinf\big(X^z(t,\mathbf a)\big)=\rinf\big(X^z_\infty(\mathbf a)\big)-\rinf\big(X^z(t,\mathbf a)\big),
\end{equation*}
and \eqref{eq:E-slope} gives $\frac m2|X^z(t,\mathbf a)-X^z_\infty(\mathbf a)|\le|(\rho(t)-\rinf)(\bX(t,\mathbf a))|\le\frac32M_1|X^z(t,\mathbf a)-X^z_\infty(\mathbf a)|$. Since $\bX(t)$ preserves the Lebesgue measure of a period cell, $\nrm{(\rho(t)-\rinf)\circ\bX(t)}_{L^p}=\nrm{\rho(t)-\rinf}_{L^p}$ for $1\le p\le\infty$, which proves \eqref{eq:E-disp}. The rates for the displacement follow from \eqref{eq:ratefull} and \eqref{eq:rateLinf}.

\emph{Step 4: the vertical component in $C^1$.} Let $h(t)\eqdef\rho(t)-\rinf$. By Step 3, $h(t,\bX(t,\mathbf a))=\rinf(X^z_\infty(\mathbf a))-\rinf(X^z(t,\mathbf a))$. Differentiating in $\mathbf a$,
\begin{equation*}
\rinf'\big(X^z(t)\big)\big(DX^z(t)-DX^z_\infty\big)=\Big(\rinf'\big(X^z_\infty\big)-\rinf'\big(X^z(t)\big)\Big)DX^z_\infty-(\nabla h)\big(t,\bX(t)\big)D\bX(t).
\end{equation*}
By \eqref{eq:E-slope}, $|\rinf'|\ge\frac m2$; moreover $\rinf\in H^4(0,1)$, so $\nrm{\rinf''}_{L^\infty}\le C$, and $D\bX(t)$, $D\bX_\infty$ are bounded by Step 2. Hence $\nrm{DX^z(t)-DX^z_\infty}_{L^\infty}\le C\big(\nrm{X^z(t)-X^z_\infty}_{L^\infty}+\nrm{\nabla h(t)}_{L^\infty}\big)$, and by Step 3 and \eqref{eq:rateLinf},
\begin{equation*}
\nrm{X^z(t)-X^z_\infty}_{C^1}\le C\nrm{h(t)}_{W^{1,\infty}}\le C\mathsf N(1+t)^{-\frac12}.
\end{equation*}

\emph{Step 5: the horizontal component in $C^1$.} The horizontal component obeys $\pt X^h=v(t,\bX)$ and $\pt DX^h=(Dv)(t,\bX)D\bX$, so the rate involves only $v$; $w$ enters only through the bound on $D\bX$ from Step 2. As in Step 2, for $t<T$,
\begin{equation*}
\nrm{X^h(T)-X^h(t)}_{L^\infty}+\nrm{DX^h(T)-DX^h(t)}_{L^\infty}\le C\int_t^T\nrm{v(\tau)}_{W^{1,\infty}}\dtau .
\end{equation*}
Since $v=0$ on $\pa\Omega$, $\nrm{v}_{L^\infty}\le\nrm{\pz v}_{L^\infty}$, and $\nabla v=\nabla\wt v+\pz\bar v\otimes e_z$. We apply \Cref{lem:endpoint} to each component of $\nabla\wt v$. It has zero horizontal mean, so \Cref{lem:poincare} gives $\nrm{\nabla\wt v}+\nrm{\nh\nabla\wt v}\le C\nrm{\nabla\nh v}\le C\Ff$, while $\nrm{\nabla\wt v}_{H^4}\le2\nrm{v}_{H^5}=2\Fu$; hence $\nrm{\nabla\wt v}_{L^\infty}\le C\Ff^{5/6}\Fu^{1/6}$. For the mean, $\nrm{\pz\bar v}_{L^\infty}\le C\nrm{\bar v}_{H^2}\le C\Ff$ by \eqref{eq:3.13b} and the uniform bound on $\Esh$. Therefore $\nrm{v}_{W^{1,\infty}}\le C(\Ff+\Ff^{5/6}\Fu^{1/6})$. \Cref{lem:dyadic} with $(a,b,c)=(1,0,0)$ and $(\frac56,0,\frac16)$, for which $\gamma=\frac12$ and $\frac13$, together with \eqref{eq:L-Ftail} and $\int_0^\infty\Fu^2\dtau\le\Es(\infty)^2$, gives $\int_t^\infty\nrm{v(\tau)}_{W^{1,\infty}}\dtau\le C\mathsf N(1+t)^{-1/3}$. Letting $T\to\infty$,
\begin{equation*}
\nrm{X^h(t)-X^h_\infty}_{C^1}\le C\mathsf N(1+t)^{-\frac13}.
\end{equation*}
Together with Step 4 this proves \eqref{eq:E-flowconv} and the more precise bounds stated after it.
\end{proof}

\subsection{The terminal map}

\begin{proof}[Proof of \Cref{thm:energy}(iv)]
By \Cref{thm:limit}(v), $\rinf^{(i)}=(\rin^{(i)})^*$, so \eqref{eq:E-contr} is the non-expansiveness of the decreasing vertical rearrangement in $L^p$ \cite{CroweZweibelRosenbloom1986}. For completeness we recall the argument. Let $\phi_1,\phi_2\ge0$ be bounded and $1<p<\infty$. For $a,b\ge0$,
\begin{equation*}
|a-b|^p=a^p+b^p-p(p-1)\int_0^a\int_0^b|r-s|^{p-2}\dd s\dd r ,
\end{equation*}
as both sides agree for $a=0$ or $b=0$ and have the same mixed derivative $\pa_a\pa_b$. Applying this with $a=\phi_1(\bx)$, $b=\phi_2(\bx)$ and integrating over $\Omega$, Tonelli's theorem gives
\begin{equation*}
\int_\Omega|\phi_1-\phi_2|^p\dx=\int_\Omega \phi_1^p\dx+\int_\Omega \phi_2^p\dx-p(p-1)\int_0^\infty\!\!\int_0^\infty|r-s|^{p-2}\big|\{\phi_1>r\}\cap\{\phi_2>s\}\big|\dd s\dd r .
\end{equation*}
The first two integrals do not change under rearrangement. The superlevel sets of $\phi_1^*$ and $\phi_2^*$ are horizontal slabs $\Tl\times[0,\zeta)$, which are nested, so $|\{\phi_1^*>r\}\cap\{\phi_2^*>s\}|=\min\{|\{\phi_1>r\}|,|\{\phi_2>s\}|\}\ge|\{\phi_1>r\}\cap\{\phi_2>s\}|$. Hence $\nrm{\phi_1^*-\phi_2^*}_{L^p}\le\nrm{\phi_1-\phi_2}_{L^p}$, and the cases $p=1$ and $p=\infty$ follow by letting $p\to1$ and $p\to\infty$. Finally, $\rs$ is strictly decreasing, so $\rs^*=\rs$; if $\rin$ is equimeasurable with $\rs$, then $\rinf=\rin^*=\rs^*=\rs$, and conversely $\rinf=\rs$ implies that $\rin$ is equimeasurable with $\rinf=\rs$.
\end{proof}

\section{Unstable stratification: Rayleigh--Taylor instability}\label{sec:RT}

In this section the stratification condition \eqref{eq:stable} is dropped. We assume $\rs\in C^5([0,1])$, $\min\rs>0$ and $\max\rs'>0$, so that heavier fluid lies above lighter fluid somewhere in the channel, and $f\in\R$ is arbitrary. We show that $(\rs,0)$ is then linearly and nonlinearly unstable, and we quantify how rotation affects the growth rates. In this section $M$ denotes the operator defined in \eqref{eq:RT-ops} below, not the constant of \Cref{sec:low}.

\begin{remark}\label{rem:localunstable}
\Cref{prop:local} is stated for every $\rs\in C^5([0,1])$ with $\min\rs>0$ and every $f\in\R$; condition \eqref{eq:stable} is not required, and its proof in \Cref{app:local}, as well as \Cref{lem:propagation}, uses only the positivity of $\rin=\rs+\vin$. Hence \Cref{prop:local} applies in the present setting, with an existence time depending only on an upper bound for $\nrm{\vin}_{H^4}+\nrm{\uin}_{H^4}$ and a positive lower bound for $\rs+\vin$. Let $[0,T_{\max})$ be the maximal interval of existence in the class \eqref{eq:localclass}. If $T_{\max}<\infty$, then
\begin{equation}\label{eq:blowup}
\limsup_{t\uparrow T_{\max}}\big(\nrm{\varrho(t)}_{H^4}+\nrm{v(t)}_{H^4}\big)=\infty .
\end{equation}
Indeed, $\min\rho(t)=\min\rin$ by transport. If the $H^4$ norms stayed bounded by some $K$, let $T_0(K)$ be the existence time of \Cref{prop:local} for data bounded by $K$; restarting \Cref{prop:local} at a time $t_1\in(T_{\max}-T_0(K),T_{\max})$ would extend the solution beyond $T_{\max}$. The data at $t_1$ are admissible: $\varrho(t_1)$ is clamped by \Cref{lem:propagation}(ii), $v(t_1)\in H^4\cap H^1_0$ satisfies \eqref{eq:lid}, and the surface pressure at $t_1$ is the one of \Cref{lem:propagation}(i), so that the initial acceleration computed from \eqref{PE}$_2$ at $t_1$ is $v_t(t_1)\in H^2$, which vanishes on $\pa\Omega$; this is \eqref{eq:compat} at $t_1$.
\end{remark}

\subsection{Linear instability}\label{sec:RT-lin}

Linearizing \eqref{PE} at $(\varrho,v)=(0,0)$ gives
\begin{equation}\label{eq:RT-lin}
\varrho_t=-\rs'w,\qquad \rs v_t+f\rs v^\perp+\nh q-\mu\Lap v=0,\qquad \pz q=-g\varrho,\qquad w=-\int_0^z\divh v\dd z',
\end{equation}
with $v=0$ and $\varrho=\pz\varrho=0$ on $\pa\Omega$, and the lid constraint \eqref{eq:lid}. We write
\begin{equation*}
V_\sigma\eqdef\Big\{\phi\in H^1_0(\Omega;\R^2):\ \divh\int_0^1\phi\dd z=0\Big\};
\end{equation*}
thus the velocity takes values in $V_\sigma$. We expand in horizontal Fourier modes, $F(x_h,z)=\sum_{k\in(l^{-1}\Z)^2}e^{ik\cdot x_h}F_k(z)$, and write $\nu=|k|$, $k^\perp=(-k_2,k_1)$ and $'=\pz$. In this subsection $\ip{\cdot,\cdot}$ and $\nrm\cdot$ denote the complex inner product (linear in the first argument) and the norm of $L^2(0,1)$. The coefficients of different wave vectors decouple in \eqref{eq:RT-lin}.

\emph{The mode $k=0$.} Since $w_0=0$ and $(\nh q)_0=0$, we get $\pt\varrho_0=0$ and $\rs\pt v_0+f\rs v_0^\perp-\mu v_0''=0$ with $v_0(0)=v_0(1)=0$. Testing with the real field $v_0$ and using $v_0^\perp\cdot v_0=0$ gives $\frac12\ddt\int_0^1\rs|v_0|^2+\mu\nrm{v_0'}^2=0$, hence $\int_0^1\rs|v_0(t)|^2\le e^{-2\mu\pi^2t/\max\rs}\int_0^1\rs|v_0(0)|^2$. Likewise, for $\lambda\ne0$ and complex $\hat v$, $\lambda\rs\hat v+f\rs\hat v^\perp-\mu\hat v''=0$ implies $\mathrm{Re}\,\lambda\int\rs|\hat v|^2=-\mu\nrm{\hat v'}^2$, because $\hat v^\perp\cdot\ol{\hat v}$ is purely imaginary. Thus the horizontal mean of the density is stationary, the mean velocity decays exponentially, and the mode $k=0$ has no growth.

\emph{Poloidal--toroidal decomposition.} Let $k\ne0$. Since $\{k,k^\perp\}$ is an orthogonal basis of $\R^2$ with $|k|=|k^\perp|=\nu$, every $v_k\in H^1_0(0,1;\mathbb C^2)$ can be written uniquely as
\begin{equation}\label{eq:RT-pt}
v_k=ik\,\Psi'+ik^\perp\Phi,\qquad \Psi(z)\eqdef-\frac i{\nu^2}\int_0^zk\cdot v_k\dd z',\qquad \Phi\eqdef-\frac i{\nu^2}\,k^\perp\cdot v_k .
\end{equation}
Then $\Psi(0)=0$, $\Psi'=\Phi=0$ at $z=0,1$, and
\begin{equation*}
w_k=-\int_0^zik\cdot v_k\dd z'=\nu^2\Psi,\qquad (\divh v)_k=ik\cdot v_k=-\nu^2\Psi',\qquad (\curlh v)_k=ik^\perp\cdot v_k=-\nu^2\Phi .
\end{equation*}
In physical variables, $v=\nh\pz\Psi+\nhp\Phi$, $w=-\Laph\Psi$, $\divh v=\Laph\pz\Psi$ and $\curlh v=\Laph\Phi$. The lid constraint $w_k(1)=0$ says $\Psi(1)=0$. Hence a field $v_k\in H^1_0$ satisfies the lid constraint if and only if $\Psi\in H^2_0(0,1)$; conversely, every $(\Psi,\Phi)\in H^2_0(0,1)\times H^1_0(0,1)$ defines through \eqref{eq:RT-pt} a field $v_k\in H^1_0$ satisfying the lid constraint.

Since $v_k^\perp=ik^\perp\Psi'-ik\,\Phi$, the components of the $k$-th momentum equation $\rs\pt v_k+f\rs v_k^\perp+ikq_k-\mu(\pz^2-\nu^2)v_k=0$ along $k$ and $k^\perp$ (obtained by applying $-\frac i{\nu^2}k\cdot$ and $-\frac i{\nu^2}k^\perp\cdot$; equivalently, by applying $\divh$ and $\curlh$ and using $\divh v^\perp=-\curlh v$, $\curlh v^\perp=\divh v$) are
\begin{align}
&\rs\pt\Psi'-f\rs\Phi+q_k-\mu(\pz^2-\nu^2)\Psi'=0,\label{eq:RT-div}\\
&\rs\pt\Phi+f\rs\Psi'+\mu L_\nu\Phi=0,\label{eq:RT-curl}
\end{align}
and the remaining equations read $\pt\varrho_k=-\nu^2\rs'\Psi$ and $\pz q_k=-g\varrho_k$. Here and below
\begin{equation}\label{eq:RT-ops}
M\eqdef-\pz(\rs\pz\,\cdot),\qquad A_\nu\eqdef(\pz^2-\nu^2)\pz^2,\qquad L_\nu\eqdef-\pz^2+\nu^2,\qquad \wt K\eqdef-g\nu^2\rs',
\end{equation}
where $\wt K$ is a multiplication operator. $M$, $A_\nu$, $\wt K$ act on $H^2_0(0,1)$ and $L_\nu$ on $H^1_0(0,1)$ through the forms
\begin{equation*}
\begin{aligned}
&\ip{M\phi,\xi}=\int_0^1\rs\phi'\ol{\xi'},\qquad \ip{A_\nu\phi,\xi}=\int_0^1\big(\phi''\ol{\xi''}+\nu^2\phi'\ol{\xi'}\big),\\
&\ip{L_\nu\eta,\zeta}=\int_0^1\big(\eta'\ol{\zeta'}+\nu^2\eta\ol\zeta\big),\qquad \ip{\wt K\phi,\xi}=-g\nu^2\int_0^1\rs'\phi\ol\xi .
\end{aligned}
\end{equation*}
Equation \eqref{eq:RT-div} determines $q_k$, and then $\pz q_k=-g\varrho_k$ is equivalent to the $z$-derivative of \eqref{eq:RT-div},
\begin{equation}\label{eq:RT-psi}
M\pt\Psi+\mu A_\nu\Psi+f\pz(\rs\Phi)+g\varrho_k=0 .
\end{equation}
Conversely, if $(\varrho_k,\Psi,\Phi)$ solves \eqref{eq:RT-psi}, \eqref{eq:RT-curl} and $\pt\varrho_k=-\nu^2\rs'\Psi$, and $q_k$ is defined by \eqref{eq:RT-div}, then $\pz q_k=-g\varrho_k$, so that $q_k+g\int_0^z\varrho_k$ does not depend on $z$, and $(\varrho_k,v_k,w_k,q_k)$ solves the $k$-th mode of \eqref{eq:RT-lin}.

\begin{definition}\label{def:RT-nm}
Let $k\ne0$ and $\lambda\in\mathbb C$. A \emph{normal mode} with rate $\lambda$ in the mode $k$ is a nonzero triple $(\hat\varrho,\Psi,\Phi)\in L^2(0,1)\times H^2_0(0,1)\times H^1_0(0,1)$ such that $e^{\lambda t}(\hat\varrho,\Psi,\Phi)$ solves $\pt\varrho_k=-\nu^2\rs'\Psi$, \eqref{eq:RT-psi} and \eqref{eq:RT-curl} in the weak sense:
\begin{equation}\label{eq:RT-nm}
\begin{aligned}
&\lambda\hat\varrho=-\nu^2\rs'\Psi,\\
&\lambda\ip{M\Psi,\xi}+\mu\ip{A_\nu\Psi,\xi}-f\ip{\rs\Phi,\xi'}+g\ip{\hat\varrho,\xi}=0\qquad\text{for all }\xi\in H^2_0(0,1),\\
&\lambda\ip{\rs\Phi,\eta}+\mu\ip{L_\nu\Phi,\eta}+f\ip{\rs\Psi',\eta}=0\qquad\text{for all }\eta\in H^1_0(0,1).
\end{aligned}
\end{equation}
\end{definition}

For $\mathrm{Re}\,\lambda\ge0$ let $B_\lambda\eqdef\lambda\rs+\mu L_\nu$ with Dirichlet conditions. Since $\mathrm{Re}\ip{B_\lambda\eta,\eta}\ge\mu\nrm{\eta'}^2$, the Lax--Milgram lemma shows that $B_\lambda:H^1_0\to H^{-1}$ is an isomorphism, and $B_\lambda^{-1}$ maps $L^2$ into $H^2\cap H^1_0$. We set
\begin{equation*}
R_\nu(\lambda)\eqdef-\pz\rs B_\lambda^{-1}\rs\pz,\qquad\text{so that}\qquad \ip{R_\nu(\lambda)\phi,\xi}=\ip{B_\lambda^{-1}(\rs\phi'),\rs\xi'}\quad(\phi,\xi\in H^1_0),
\end{equation*}
and define the quadratic pencil
\begin{equation}\label{eq:pencil}
Q_\nu(\lambda)\eqdef\lambda^2M+\lambda\big(\mu A_\nu+f^2R_\nu(\lambda)\big)+\wt K:\ H^2_0(0,1)\to H^{-2}(0,1),\qquad \mathrm{Re}\,\lambda\ge0 .
\end{equation}
Rotation enters only through the Ekman-type operator $f^2R_\nu(\lambda)$, which is added to the viscous damping $\mu A_\nu$. For $f=0$, \eqref{eq:pencil} reduces to the pencil $\lambda^2M+\lambda\mu A_\nu+\wt K$ of the non-rotating problem.

\begin{lemma}\label{lem:RT-modes}
Let $k\ne0$, $\nu=|k|$ and $\lambda\ne0$.
\begin{enumerate}[label=(\alph*)]
\item If $(\hat\varrho,\Psi,\Phi)$ is a normal mode with rate $\lambda$, then
\begin{equation}\label{eq:RT-energy}
\lambda^2\ip{M\Psi,\Psi}+|\lambda|^2\ip{\rs\Phi,\Phi}+\lambda\mu\big(\ip{A_\nu\Psi,\Psi}+\ip{L_\nu\Phi,\Phi}\big)+\ip{\wt K\Psi,\Psi}=0 .
\end{equation}
If moreover $\mathrm{Re}\,\lambda\ge0$, then $\Psi\ne0$, $\hat\varrho=-\nu^2\rs'\Psi/\lambda$, $\Phi=-fB_\lambda^{-1}(\rs\Psi')$ and $Q_\nu(\lambda)\Psi=0$.
\item Conversely, if $\mathrm{Re}\,\lambda\ge0$, $\Psi\in H^2_0\setminus\{0\}$ and $Q_\nu(\lambda)\Psi=0$, then $(-\nu^2\rs'\Psi/\lambda,\Psi,-fB_\lambda^{-1}(\rs\Psi'))$ is a normal mode with rate $\lambda$.
\end{enumerate}
\end{lemma}

\begin{proof}
(a) The first line of \eqref{eq:RT-nm} gives $g\ip{\hat\varrho,\Psi}=\lambda^{-1}\ip{\wt K\Psi,\Psi}$. Taking $\eta=\Phi$ in the third line gives $f\ip{\rs\Psi',\Phi}=-\lambda\ip{\rs\Phi,\Phi}-\mu\ip{L_\nu\Phi,\Phi}$; since $f$ is real and the brackets on the right are real, $f\ip{\rs\Phi,\Psi'}=\ol{f\ip{\rs\Psi',\Phi}}=-\bar\lambda\ip{\rs\Phi,\Phi}-\mu\ip{L_\nu\Phi,\Phi}$. Taking $\xi=\Psi$ in the second line and inserting these two identities,
\begin{equation*}
\lambda\ip{M\Psi,\Psi}+\mu\ip{A_\nu\Psi,\Psi}+\bar\lambda\ip{\rs\Phi,\Phi}+\mu\ip{L_\nu\Phi,\Phi}+\lambda^{-1}\ip{\wt K\Psi,\Psi}=0 ;
\end{equation*}
multiplication by $\lambda$ gives \eqref{eq:RT-energy}. Now let $\mathrm{Re}\,\lambda\ge0$. The third line of \eqref{eq:RT-nm} says $B_\lambda\Phi=-f\rs\Psi'$ in $H^{-1}$, so $\Phi=-fB_\lambda^{-1}(\rs\Psi')$. If $\Psi=0$, then $\hat\varrho=0$ and $\Phi=0$, which is excluded. Inserting $\hat\varrho$ and $\Phi$ in the second line, and using $-f\ip{\rs\Phi,\xi'}=f^2\ip{B_\lambda^{-1}(\rs\Psi'),\rs\xi'}=f^2\ip{R_\nu(\lambda)\Psi,\xi}$, we obtain $\lambda^{-1}\ip{Q_\nu(\lambda)\Psi,\xi}=0$ for all $\xi\in H^2_0$. (b) The same computation, read backwards.
\end{proof}

The next lemma collects the properties of the Ekman operator for real positive rates.

\begin{lemma}\label{lem:RT-ekman}
Let $\nu>0$ and $s>0$, and put $T_s\eqdef\rs+\frac\mu sL_\nu=s^{-1}B_s$. For $\psi\in H^1_0(0,1)$ let $\chi\eqdef\rs\psi'$ and
\begin{equation*}
\mathfrak e_\nu(s;\psi)\eqdef s\ip{R_\nu(s)\psi,\psi}=\ip{T_s^{-1}\chi,\chi}.
\end{equation*}
\begin{enumerate}[label=(\alph*)]
\item (Dual formula.)
\begin{equation}\label{eq:RT-dual}
\mathfrak e_\nu(s;\psi)=\max_{\eta\in H^1_0(0,1)}\Big\{2\,\mathrm{Re}\ip{\chi,\eta}-\int_0^1\rs|\eta|^2-\frac\mu s\big(\nrm{\eta'}^2+\nu^2\nrm\eta^2\big)\Big\},
\end{equation}
and the maximum is attained exactly at $\eta=T_s^{-1}\chi$.
\item $0\le\mathfrak e_\nu(s;\psi)\le\min\Big\{\ip{M\psi,\psi},\ \dfrac{s\nrm{\rs\psi'}^2}{\mu(\pi^2+\nu^2)}\Big\}$, and $\mathfrak e_\nu(s;\psi)>0$ if $\psi\ne0$.
\item $s\mapsto\mathfrak e_\nu(s;\psi)$ is nondecreasing, $\nu\mapsto\mathfrak e_\nu(s;\psi)$ is nonincreasing, and for $0<s_1<s_2$,
\begin{equation*}
0\le\mathfrak e_\nu(s_2;\psi)-\mathfrak e_\nu(s_1;\psi)\le\frac{s_2-s_1}{4s_1}\ip{M\psi,\psi}.
\end{equation*}
\item $\mathfrak e_\nu(s;\cdot)$ is a continuous Hermitian form on $H^1_0$; if $\psi_1,\psi_2$ are real, then $\mathfrak e_\nu(s;\psi_1+i\psi_2)=\mathfrak e_\nu(s;\psi_1)+\mathfrak e_\nu(s;\psi_2)$. The operator $R_\nu(s)$ maps $H^m\cap H^1_0$ boundedly into $H^m$ for $1\le m\le4$, and for every $m\ge1$ if $\rs\in C^\infty([0,1])$.
\end{enumerate}
\end{lemma}

\begin{proof}
(a) $T_s$ is the positive self-adjoint operator of the coercive form $\int\rs\eta\bar\zeta+\frac\mu s\ip{L_\nu\eta,\zeta}$ on $H^1_0$, and $T_s^{-1}$ maps $L^2$ into $H^2\cap H^1_0$. With $\eta_0\eqdef T_s^{-1}\chi$ and $\eta\in H^1_0$, expanding the square gives
\begin{equation*}
2\,\mathrm{Re}\ip{\chi,\eta}-\ip{T_s\eta,\eta}=\ip{T_s^{-1}\chi,\chi}-\ip{T_s(\eta-\eta_0),\eta-\eta_0},
\end{equation*}
and $\ip{T_s(\eta-\eta_0),\eta-\eta_0}>0$ unless $\eta=\eta_0$.

(b) Pointwise, $2\,\mathrm{Re}(\chi\bar\eta)-\rs|\eta|^2\le\rs^{-1}|\chi|^2$, and $\int\rs^{-1}|\chi|^2=\int\rs|\psi'|^2=\ip{M\psi,\psi}$; dropping the last term in \eqref{eq:RT-dual} gives the first bound. Dropping instead $\int\rs|\eta|^2$, and using $\nrm{\eta'}^2\ge\pi^2\nrm\eta^2$ and Young's inequality, $2\,\mathrm{Re}\ip{\chi,\eta}-\frac\mu s(\pi^2+\nu^2)\nrm\eta^2\le\frac{s\nrm\chi^2}{\mu(\pi^2+\nu^2)}$, which is the second bound. Finally $\mathfrak e_\nu(s;\psi)=\ip{T_s\eta_0,\eta_0}\ge\min\rs\,\nrm{\eta_0}^2$, and $\eta_0=0$ only if $\chi=0$, i.e.\ $\psi'=0$, i.e.\ $\psi=0$.

(c) For fixed $\eta$, the expression in braces in \eqref{eq:RT-dual} is nondecreasing in $s$ and nonincreasing in $\nu$, and so is its maximum. Let $\eta_2\eqdef T_{s_2}^{-1}\chi$ be the maximizer for $s_2$, and $\ell_2\eqdef\nrm{\eta_2'}^2+\nu^2\nrm{\eta_2}^2$. Testing \eqref{eq:RT-dual} for $s_1$ with $\eta_2$ gives $\mathfrak e_\nu(s_1;\psi)\ge\mathfrak e_\nu(s_2;\psi)-\mu(s_1^{-1}-s_2^{-1})\ell_2$. Testing $T_{s_2}\eta_2=\chi$ with $\eta_2$,
\begin{equation*}
\int\rs|\eta_2|^2+\frac\mu{s_2}\ell_2=\ip{\chi,\eta_2}\le\Big(\int\rs^{-1}|\chi|^2\Big)^{\frac12}\Big(\int\rs|\eta_2|^2\Big)^{\frac12}\le\int\rs|\eta_2|^2+\frac14\ip{M\psi,\psi},
\end{equation*}
so $\frac\mu{s_2}\ell_2\le\frac14\ip{M\psi,\psi}$ and $\mu(s_1^{-1}-s_2^{-1})\ell_2=\frac{s_2-s_1}{s_1}\frac\mu{s_2}\ell_2\le\frac{s_2-s_1}{4s_1}\ip{M\psi,\psi}$.

(d) $T_s^{-1}$ is bounded and self-adjoint on $L^2$, with norm at most $(\min\rs)^{-1}$, and maps real functions to real functions; this gives continuity, and the cross terms $\ip{T_s^{-1}\chi_2,\chi_1}-\ip{T_s^{-1}\chi_1,\chi_2}$ vanish for real $\chi_1,\chi_2$. If $\psi\in H^m\cap H^1_0$, then $\chi\in H^{m-1}$, and the equation $\frac\mu s\eta_0''=(\rs+\frac{\mu\nu^2}s)\eta_0-\chi$ gives $\eta_0\in H^{m+1}$ by induction; hence $R_\nu(s)\psi=-s^{-1}(\rs\eta_0)'\in H^m$ if $\rs\in C^{m+1}$.
\end{proof}

By \eqref{eq:Lambdainf} and the Poincar\'e inequality $\nrm\phi\le\pi^{-1}\nrm{\phi'}$ on $H^1_0(0,1)$, and by testing \eqref{eq:Lambdainf} with $\phi\in C_c^\infty$ supported in the nonempty open set $\{\rs'>0\}\cap(0,1)$,
\begin{equation*}
0<\Lambda_\infty\le\frac{g\max\rs'}{\pi^2\mu}.
\end{equation*}
We write $\lambda(\nu)=\lambda(\nu;f)$ when the dependence on $f$ matters.

\begin{theorem}\label{thm:RTlinear}
Let $\rs\in C^5([0,1])$ with $\min\rs>0$ and $\max\rs'>0$, and let $f\in\R$. For every $\nu>0$ the following hold.
\begin{enumerate}[label=(\roman*)]
\item There is a largest $\lambda(\nu)>0$ for which $Q_\nu(\lambda(\nu))\phi_\nu=0$ for some real $\phi_\nu\in H^4\cap H^2_0\setminus\{0\}$. With $\lambda=\lambda(\nu)$, $\Phi_\nu\eqdef-fB_\lambda^{-1}(\rs\phi_\nu')$ and $\hat\varrho_\nu\eqdef-\nu^2\rs'\phi_\nu/\lambda$, the triple $(\hat\varrho_\nu,\phi_\nu,\Phi_\nu)$ is a real normal mode; $\hat\varrho_\nu=\pz\hat\varrho_\nu=0$ at $z=0,1$, and $\hat\varrho_\nu\not\equiv0$. If $\rs\in C^\infty([0,1])$, then $\phi_\nu,\Phi_\nu\in C^\infty([0,1])$. Consequently, for every $k\in(l^{-1}\Z)^2\setminus\{0\}$ with $|k|=\nu$, the fields
\begin{equation}\label{eq:RT-mode}
\begin{aligned}
&\varrho=e^{\lambda t}\cos(k\cdot x_h)\,\hat\varrho_\nu(z),\qquad w=e^{\lambda t}\cos(k\cdot x_h)\,\nu^2\phi_\nu(z),\qquad q=e^{\lambda t}\cos(k\cdot x_h)\,\hat q_\nu(z),\\
&v=-e^{\lambda t}\sin(k\cdot x_h)\big(\phi_\nu'(z)\,k+\Phi_\nu(z)\,k^\perp\big),\qquad \hat q_\nu\eqdef-\lambda\rs\phi_\nu'+f\rs\Phi_\nu+\mu(\phi_\nu'''-\nu^2\phi_\nu'),
\end{aligned}
\end{equation}
form a real solution of \eqref{eq:RT-lin} with $v(t)\in V_\sigma$ and $\varrho=\pz\varrho=0$ on $\pa\Omega$. Thus the linearized problem has growing modes in every horizontal mode $k\ne0$.
\item Every normal mode with a non-real rate $\lambda$ has $\mathrm{Re}\,\lambda<0$, and every normal mode with a real rate $\lambda>0$ has $\lambda\le\lambda(\nu)$. In particular, every normal mode with $\mathrm{Re}\,\lambda\ge0$ and $\lambda\ne0$ has a real rate $\lambda\in(0,\lambda(\nu)]$. For every real $s>\lambda(\nu)$, $Q_\nu(s)$ is an isomorphism from $H^4\cap H^2_0$ onto $L^2$.
\item $\nu\mapsto\lambda(\nu)$ is strictly increasing, $\lambda(\nu)<\Lambda_\infty$, and $\lambda(\nu)\to\Lambda_\infty$ as $\nu\to\infty$. In particular $\sup_{k\ne0}\lambda(|k|)=\Lambda_\infty$ is not attained, and it does not depend on $f$.
\item Rotation lowers every growth rate. $\lambda(\nu;f)$ depends on $f$ only through $|f|$, the map $|f|\mapsto\lambda(\nu;f)$ is strictly decreasing, and
\begin{equation}\label{eq:RT-rotbound}
\frac{\lambda(\nu;0)}{1+f^2\vartheta_\nu^2}<\lambda(\nu;f)<\lambda(\nu;0)\quad(f\ne0),\qquad \vartheta_\nu\eqdef\frac{\max\rs}{\mu(\pi^2+\nu^2)} .
\end{equation}
Moreover, $\lambda(\nu;f)\to0$ as $|f|\to\infty$, for each fixed $\nu$.
\end{enumerate}
\end{theorem}

\begin{proof}
\emph{The variational problem.} For $s>0$ and $\phi\in H^2_0(0,1)$ let
\begin{equation*}
\cJ_{\nu,s}(\phi)\eqdef s\mu\ip{A_\nu\phi,\phi}+f^2\mathfrak e_\nu(s;\phi)+\ip{\wt K\phi,\phi}=s\ip{\big(\mu A_\nu+f^2R_\nu(s)\big)\phi,\phi}+\ip{\wt K\phi,\phi},
\end{equation*}
\begin{equation*}
\alpha_\nu(s)\eqdef\inf\big\{\cJ_{\nu,s}(\phi):\ \phi\in H^2_0(0,1)\ \text{real},\ \ip{M\phi,\phi}=1\big\},\qquad h_\nu(s)\eqdef s^2+\alpha_\nu(s).
\end{equation*}
Note that $s^2\ip{M\phi,\phi}+\cJ_{\nu,s}(\phi)=\ip{Q_\nu(s)\phi,\phi}$. We collect four properties.

(P1) Since $\nrm\phi^2\le\pi^{-2}\nrm{\phi'}^2\le(\pi^2\min\rs)^{-1}\ip{M\phi,\phi}$, we have $\ip{\wt K\phi,\phi}\ge-C_\nu\ip{M\phi,\phi}$ with $C_\nu\eqdef g\nu^2\nrm{\rs'}_{L^\infty}/(\pi^2\min\rs)$. As $\mathfrak e_\nu\ge0$, $\alpha_\nu(s)\ge-C_\nu$, and every admissible $\phi$ satisfies $s\mu\nrm{\phi''}^2\le s\mu\ip{A_\nu\phi,\phi}\le\cJ_{\nu,s}(\phi)+C_\nu$.

(P2) The infimum is attained. A minimizing sequence is bounded in $H^2_0$ by (P1). A subsequence converges weakly in $H^2_0$ and strongly in $H^1_0$. The constraint, $\ip{\wt K\phi,\phi}$ and $\mathfrak e_\nu(s;\phi)$ (\Cref{lem:RT-ekman}(d)) pass to the limit, and $\ip{A_\nu\phi,\phi}$ is weakly lower semicontinuous.

(P3) For every complex $\phi\in H^2_0$, $\cJ_{\nu,s}(\phi)\ge\alpha_\nu(s)\ip{M\phi,\phi}$, where $\cJ_{\nu,s}$ is defined by the same formula. Indeed, for real $\phi_1,\phi_2$ all quadratic forms in $\cJ_{\nu,s}$ and $\ip{M\cdot,\cdot}$ split as the sum of their values at $\phi_1$ and $\phi_2$ (\Cref{lem:RT-ekman}(d)), and the claim for real functions follows from the definition by homogeneity.

(P4) $\alpha_\nu$ is nondecreasing and continuous on $(0,\infty)$. Monotonicity holds because $s\mapsto\cJ_{\nu,s}(\phi)$ is nondecreasing for each $\phi$ (\Cref{lem:RT-ekman}(c)). Let $0<s_1<s_2$, $\epsilon\in(0,1]$, and let $\phi$ be admissible with $\cJ_{\nu,s_1}(\phi)\le\alpha_\nu(s_1)+\epsilon$; by (P1), $\mu\ip{A_\nu\phi,\phi}\le(\alpha_\nu(s_1)+1+C_\nu)/s_1$. By \Cref{lem:RT-ekman}(c),
\begin{equation*}
\begin{aligned}
\alpha_\nu(s_2)\le\cJ_{\nu,s_2}(\phi)&\le\cJ_{\nu,s_1}(\phi)+(s_2-s_1)\mu\ip{A_\nu\phi,\phi}+f^2\frac{s_2-s_1}{4s_1}\\
&\le\alpha_\nu(s_1)+\epsilon+\frac{s_2-s_1}{s_1}\Big(\alpha_\nu(s_1)+1+C_\nu+\frac{f^2}4\Big).
\end{aligned}
\end{equation*}
Letting $\epsilon\to0$, $\alpha_\nu$ is Lipschitz on every compact subinterval of $(0,\infty)$.

Consequently $h_\nu$ is continuous and strictly increasing on $(0,\infty)$, and $h_\nu(s)\ge s^2-C_\nu\to\infty$. Choose a real $\phi_0\in C_c^\infty$ supported in $\{\rs'>0\}\cap(0,1)$ with $\ip{M\phi_0,\phi_0}=1$. Then $\ip{\wt K\phi_0,\phi_0}<0$, and by \Cref{lem:RT-ekman}(b),
\begin{equation*}
h_\nu(s)\le s^2+s\mu\ip{A_\nu\phi_0,\phi_0}+\frac{sf^2\nrm{\rs\phi_0'}^2}{\mu(\pi^2+\nu^2)}+\ip{\wt K\phi_0,\phi_0}<0
\end{equation*}
for small $s>0$. Hence $h_\nu$ has exactly one zero $\lambda(\nu)>0$; moreover $h_\nu<0$ on $(0,\lambda(\nu))$ and $h_\nu>0$ on $(\lambda(\nu),\infty)$.

\emph{Proof of (i).} Let $\lambda=\lambda(\nu)$ and let $\phi_\nu$ be a real minimizer for $s=\lambda$ (P2). By (P3), the real quadratic form $\mathfrak q(\phi)\eqdef\ip{Q_\nu(\lambda)\phi,\phi}=\lambda^2\ip{M\phi,\phi}+\cJ_{\nu,\lambda}(\phi)\ge h_\nu(\lambda)\ip{M\phi,\phi}=0$ is nonnegative on $H^2_0$, and $\mathfrak q(\phi_\nu)=0$. Hence $0\le\mathfrak q(\phi_\nu+t\xi)=2t\,\mathrm{Re}\,\mathfrak q(\phi_\nu,\xi)+t^2\mathfrak q(\xi)$ for all real $t$, which forces $\mathfrak q(\phi_\nu,\xi)=\ip{Q_\nu(\lambda)\phi_\nu,\xi}=0$ for all real, hence all complex, $\xi\in H^2_0$. In $\mathcal D'(0,1)$ this reads
\begin{equation*}
\lambda\mu\big(\phi_\nu''''-\nu^2\phi_\nu''\big)=\lambda^2(\rs\phi_\nu')'-\lambda f^2R_\nu(\lambda)\phi_\nu+g\nu^2\rs'\phi_\nu\in L^2,
\end{equation*}
by \Cref{lem:RT-ekman}(d). Hence $\phi_\nu\in H^4$; if $\rs\in C^\infty$, iterating this identity with \Cref{lem:RT-ekman}(d) gives $\phi_\nu\in C^\infty$. $\Phi_\nu$ is real and lies in $H^2\cap H^1_0$ (in $C^\infty$ if $\rs\in C^\infty$), and $(\hat\varrho_\nu,\phi_\nu,\Phi_\nu)$ is a normal mode by \Cref{lem:RT-modes}(b). Since $\phi_\nu=\phi_\nu'=0$ at the lids, $\hat\varrho_\nu=-\nu^2\rs'\phi_\nu/\lambda$ and $\hat\varrho_\nu'=-\nu^2(\rs''\phi_\nu+\rs'\phi_\nu')/\lambda$ vanish there. If $\hat\varrho_\nu\equiv0$, then $\ip{\wt K\phi_\nu,\phi_\nu}=0$ and $0=\mathfrak q(\phi_\nu)\ge\lambda^2\ip{M\phi_\nu,\phi_\nu}>0$, a contradiction. The fields \eqref{eq:RT-mode} are the real parts of $e^{\lambda t}e^{ik\cdot x_h}(\hat\varrho_\nu,ik\phi_\nu'+ik^\perp\Phi_\nu,\nu^2\phi_\nu,\hat q_\nu)$, which solves \eqref{eq:RT-lin}: the components of the momentum equation along $k$ and $k^\perp$ are the definition of $\hat q_\nu$ and $B_\lambda\Phi_\nu=-f\rs\phi_\nu'$, and $\pz\hat q_\nu=-g\hat\varrho_\nu$ is equivalent to $Q_\nu(\lambda)\phi_\nu=0$ by the computation in the proof of \Cref{lem:RT-modes}. The velocity lies in $V_\sigma$ because $\phi_\nu\in H^2_0$ and $\Phi_\nu\in H^1_0$. That $\lambda(\nu)$ is the largest such value follows from (ii) and \Cref{lem:RT-modes}(b).

\emph{Proof of (ii).} Let $(\hat\varrho,\Psi,\Phi)$ be a normal mode with rate $\lambda=\alpha+i\beta$, $\beta\ne0$. If $\Psi=0$, then $\hat\varrho=0$ and \eqref{eq:RT-nm} with $\eta=\Phi$ gives $\lambda=-\mu\ip{L_\nu\Phi,\Phi}/\ip{\rs\Phi,\Phi}\in\R$; hence $\Psi\ne0$. The imaginary part of \eqref{eq:RT-energy} is
\begin{equation*}
\beta\Big(2\alpha\ip{M\Psi,\Psi}+\mu\big(\ip{A_\nu\Psi,\Psi}+\ip{L_\nu\Phi,\Phi}\big)\Big)=0,
\end{equation*}
so $2\,\mathrm{Re}\,\lambda\,\ip{M\Psi,\Psi}=-\mu\big(\ip{A_\nu\Psi,\Psi}+\ip{L_\nu\Phi,\Phi}\big)<0$. Next let $\lambda>0$ be real. By \Cref{lem:RT-modes}(a), $\Psi\ne0$ and $Q_\nu(\lambda)\Psi=0$, hence by (P3)
\begin{equation*}
0=\ip{Q_\nu(\lambda)\Psi,\Psi}=\lambda^2\ip{M\Psi,\Psi}+\cJ_{\nu,\lambda}(\Psi)\ge h_\nu(\lambda)\ip{M\Psi,\Psi},
\end{equation*}
so $h_\nu(\lambda)\le0$ and $\lambda\le\lambda(\nu)$. Now let $s>\lambda(\nu)$ be real, and let $\mathfrak q_s(\phi,\xi)\eqdef\ip{Q_\nu(s)\phi,\xi}$, a bounded Hermitian form on $H^2_0$ with $\mathfrak q_s(\phi,\phi)\ge h_\nu(s)\ip{M\phi,\phi}$ by (P3), and $h_\nu(s)>0$. We claim that $\mathfrak q_s(\phi,\phi)\ge c\nrm{\phi''}^2$ for some $c>0$. Otherwise there are $\phi_n\in H^2_0$ with $\nrm{\phi_n''}=1$ and $\mathfrak q_s(\phi_n,\phi_n)\to0$. A subsequence converges weakly in $H^2_0$ and strongly in $H^1_0$ to some $\phi$, and as in (P2), $\mathfrak q_s(\phi,\phi)\le\liminf\mathfrak q_s(\phi_n,\phi_n)=0$, so $\phi=0$. But the G\aa rding inequality $\mathfrak q_s(\phi_n,\phi_n)\ge s\mu\nrm{\phi_n''}^2-g\nu^2\nrm{\rs'}_{L^\infty}\nrm{\phi_n}^2$ (the $M$- and $R_\nu$-terms are nonnegative) then gives $s\mu\le0$, a contradiction. By the Lax--Milgram lemma, for $F\in L^2$ there is a unique $\Psi\in H^2_0$ with $\mathfrak q_s(\Psi,\xi)=\ip{F,\xi}$ for all $\xi\in H^2_0$, and $\nrm{\Psi''}\le c^{-1}\nrm F$. In $\mathcal D'(0,1)$, $s\mu A_\nu\Psi=F-s^2M\Psi-sf^2R_\nu(s)\Psi-\wt K\Psi\in L^2$, so $\Psi\in H^4$ and $\nrm{\Psi}_{H^4}\le C\nrm F$. Hence $Q_\nu(s):H^4\cap H^2_0\to L^2$ is bijective and bounded with bounded inverse.

\emph{Proof of (iii).} Writing $\alpha_\nu(s)=\nu^2\beta_\nu(s)$ with
\begin{equation*}
\beta_\nu(s)=\inf_{\ip{M\phi,\phi}=1}\Big[s\mu\big(\nu^{-2}\nrm{\phi''}^2+\nrm{\phi'}^2\big)+\nu^{-2}f^2\mathfrak e_\nu(s;\phi)-g\int_0^1\rs'\phi^2\Big],
\end{equation*}
we see that $\beta_\nu(s)$ is nonincreasing in $\nu$, because $\nu^{-2}$ and $\mathfrak e_\nu(s;\phi)\ge0$ are nonincreasing in $\nu$ (\Cref{lem:RT-ekman}(c)). Let $0<\nu<\nu'$ and $s=\lambda(\nu)$, so that $\beta_\nu(s)=-s^2/\nu^2$. Then $h_{\nu'}(s)\le s^2+\nu'^2\beta_\nu(s)=s^2(1-\nu'^2/\nu^2)<0$, hence $\lambda(\nu')>\lambda(\nu)$.

For the upper bound let $\lambda=\lambda(\nu)$ and test $Q_\nu(\lambda)\phi_\nu=0$ with $\phi_\nu$. Since $\ip{M\phi_\nu,\phi_\nu}>0$ and $\mathfrak e_\nu(\lambda;\phi_\nu)>0$,
\begin{equation*}
\begin{aligned}
\lambda\mu\ip{A_\nu\phi_\nu,\phi_\nu}&<\lambda^2\ip{M\phi_\nu,\phi_\nu}+\lambda\mu\ip{A_\nu\phi_\nu,\phi_\nu}+f^2\mathfrak e_\nu(\lambda;\phi_\nu)\\
&=g\nu^2\int_0^1\rs'\phi_\nu^2\le\nu^2\mu\Lambda_\infty\nrm{\phi_\nu'}^2\le\Lambda_\infty\mu\ip{A_\nu\phi_\nu,\phi_\nu},
\end{aligned}
\end{equation*}
by \eqref{eq:Lambdainf}. Hence $\lambda<\Lambda_\infty$. Here the Ekman term has the favourable sign.

Finally let $\varepsilon\in(0,\Lambda_\infty/4)$. Pick a real $\phi\in C_c^\infty(0,1)$ with $\ip{M\phi,\phi}=1$ and $g\int\rs'\phi^2\ge\mu(\Lambda_\infty-\varepsilon)\nrm{\phi'}^2$, which is possible by density of $C_c^\infty$ in $H^1_0$. For $s=\Lambda_\infty-2\varepsilon$, using $\mathfrak e_\nu(s;\phi)\le\ip{M\phi,\phi}=1$,
\begin{equation*}
h_\nu(s)\le s^2+s\mu\nrm{\phi''}^2+f^2+\nu^2\Big(s\mu\nrm{\phi'}^2-g\int_0^1\rs'\phi^2\Big)\le s^2+s\mu\nrm{\phi''}^2+f^2-\varepsilon\mu\nu^2\nrm{\phi'}^2<0
\end{equation*}
for $\nu$ large, so $\lambda(\nu)>\Lambda_\infty-2\varepsilon$ for such $\nu$. Since $\{|k|:k\in(l^{-1}\Z)^2\setminus\{0\}\}$ is unbounded, the remaining assertions follow.

\emph{Proof of (iv).} $Q_\nu(\lambda)$ depends on $f$ only through $f^2$. We write $\cJ_{\nu,s}(\phi;f)$, $\alpha_\nu(s;f)$ and $h_\nu(s;f)$. Let $f_1^2<f_2^2$, $s\eqdef\lambda(\nu;f_2)$, and let $\phi_2$ be a real minimizer of $\cJ_{\nu,s}(\cdot;f_2)$, so that $\cJ_{\nu,s}(\phi_2;f_2)=-s^2$. By \Cref{lem:RT-ekman}(b),
\begin{equation*}
\alpha_\nu(s;f_1)\le\cJ_{\nu,s}(\phi_2;f_1)=\cJ_{\nu,s}(\phi_2;f_2)-(f_2^2-f_1^2)\,\mathfrak e_\nu(s;\phi_2)<-s^2 .
\end{equation*}
Hence $h_\nu(s;f_1)<0$ and $\lambda(\nu;f_1)>s=\lambda(\nu;f_2)$. With $f_1=0$ this gives the upper bound in \eqref{eq:RT-rotbound}.

For the lower bound let $\lambda_0\eqdef\lambda(\nu;0)$, and let $\phi_0$ be a real minimizer of $\cJ_{\nu,\lambda_0}(\cdot;0)$ with $\ip{M\phi_0,\phi_0}=1$, so that $\cJ_{\nu,\lambda_0}(\phi_0;0)=-\lambda_0^2$. Since $\phi_0'\in H^1_0$ and $\nrm{\phi_0'}^2\ge\ip{M\phi_0,\phi_0}/\max\rs$,
\begin{equation*}
\mu\ip{A_\nu\phi_0,\phi_0}\ge\mu(\pi^2+\nu^2)\nrm{\phi_0'}^2\ge\vartheta_\nu^{-1},\qquad \mathfrak e_\nu(s;\phi_0)\le\frac{s\nrm{\rs\phi_0'}^2}{\mu(\pi^2+\nu^2)}\le s\,\vartheta_\nu,
\end{equation*}
the latter by \Cref{lem:RT-ekman}(b) and $\nrm{\rs\phi_0'}^2\le\max\rs\,\ip{M\phi_0,\phi_0}$. For $s\in(0,\lambda_0)$,
\begin{equation*}
h_\nu(s;f)\le s^2+\cJ_{\nu,s}(\phi_0;f)=s^2-\lambda_0^2-(\lambda_0-s)\mu\ip{A_\nu\phi_0,\phi_0}+f^2\mathfrak e_\nu(s;\phi_0)\le s^2-\lambda_0^2-\frac{\lambda_0-s}{\vartheta_\nu}+f^2s\,\vartheta_\nu .
\end{equation*}
Let $f\ne0$ and $s_f\eqdef\lambda_0/(1+f^2\vartheta_\nu^2)\in(0,\lambda_0)$. Then $\lambda_0-s_f=f^2\vartheta_\nu^2s_f$, so the last two terms cancel and $h_\nu(s_f;f)\le s_f^2-\lambda_0^2<0$. Hence $\lambda(\nu;f)>s_f$, which is the lower bound in \eqref{eq:RT-rotbound}.

Finally suppose that $\lambda(\nu;f_n)\ge s_0>0$ along a sequence $|f_n|\to\infty$. Since $h_\nu(\cdot;f_n)$ is increasing, $h_\nu(s_0;f_n)\le0$, i.e.\ $\alpha_\nu(s_0;f_n)\le-s_0^2$. Let $\phi_n$ be real minimizers. By (P1), $s_0\mu\nrm{\phi_n''}^2+f_n^2\mathfrak e_\nu(s_0;\phi_n)\le C_\nu$. Thus $(\phi_n)$ is bounded in $H^2_0$ and $\mathfrak e_\nu(s_0;\phi_n)\to0$. A subsequence converges in $H^1_0$ to some $\phi$ with $\ip{M\phi,\phi}=1$ and, by \Cref{lem:RT-ekman}(d), $\mathfrak e_\nu(s_0;\phi)=0$, which contradicts \Cref{lem:RT-ekman}(b). Hence $\lambda(\nu;f)\to0$.
\end{proof}

\begin{remark}\label{rem:RT-num}
A Galerkin approximation of \eqref{eq:RT-nm} for $\rs=1+z$ and $\mu=g=1$, where $\Lambda_\infty=\pi^{-2}$, is consistent with \Cref{thm:RTlinear}: the computed largest real rates increase with $\nu$ towards $\pi^{-2}$, decrease with $|f|$, respect \eqref{eq:RT-rotbound}, and no non-real eigenvalue with nonnegative real part appears. These computations are not used in the proofs.
\end{remark}

\subsection{Nonlinear instability}\label{sec:RT-nonlin}

By \Cref{thm:RTlinear}(iii) the maximal growth rate is not attained, so the usual bootstrap, which uses a mode of maximal growth, does not apply directly. We use instead a mode growing faster than $\Lambda_\infty/2$ together with approximate solutions of high order, in the spirit of \cite{Grenier2000}; see also \cite{GuoStrauss1995,JiangJiang2014}.

\begin{theorem}\label{thm:RTnonlinear}
Let $\rs\in C^\infty([0,1])$ with $\min\rs>0$ and $\max\rs'>0$, and let $f\in\R$. Fix $k_*\in(l^{-1}\Z)^2\setminus\{0\}$ with $\lambda_*\eqdef\lambda(|k_*|;f)>\Lambda_\infty/2$; such $k_*$ exists by \Cref{thm:RTlinear}(iii). There are $\bar\varepsilon_0,\varepsilon_1,\delta_0>0$ with the following property. For every $\delta\in(0,\delta_0)$ there are smooth data $(\vin^\delta,\uin^\delta)$ satisfying the hypotheses of \Cref{prop:local} (see \Cref{rem:localunstable}), with $\nrm{\vin^\delta}_{H^m}+\nrm{\uin^\delta}_{H^m}\le C_m\delta$ for every $m$, such that the solution of \eqref{PE} on its maximal interval $[0,T^\delta_{\max})$ satisfies:
\begin{equation*}
\text{either}\quad T^\delta_{\max}\le T^\delta\eqdef\frac1{\lambda_*}\log\frac{\bar\varepsilon_0}\delta,\qquad\text{or}\quad \nrm{\varrho^\delta(T^\delta)}_{L^2}\ge\varepsilon_1 .
\end{equation*}
In particular $\sup_{0\le t<T^\delta_{\max}}\big(\nrm{\varrho^\delta(t)}_{H^4}+\nrm{v^\delta(t)}_{H^4}\big)\ge\varepsilon_1$ for every $\delta\in(0,\delta_0)$. Hence $(\rs,0)$ is nonlinearly unstable: arbitrarily small smooth perturbations leave a fixed $H^4$-neighbourhood.
\end{theorem}

The last assertion follows from the first by \eqref{eq:blowup}. The proof occupies the rest of this subsection; $\nrm\cdot$ now denotes the norm of $L^2(\Omega)$. For a horizontal field $v_a$ we write $w_a\eqdef-\int_0^z\divh v_a\dd z'$ and $\bU_a\eqdef(v_a,w_a)$, and similarly for other velocity fields. We first solve the resolvent problem that appears at each order.

\begin{lemma}\label{lem:RT-resolvent}
Let $\rs\in C^\infty([0,1])$ with $\min\rs>0$ and $\max\rs'>0$, and let $s>\Lambda_\infty$. Let $k\in(l^{-1}\Z)^2$, $\nu=|k|$, and let $G^\varrho\in C^\infty([0,1];\mathbb C)$ and $G^v\in C^\infty([0,1];\mathbb C^2)$ with $G^\varrho=\pz G^\varrho=0$ at $z=0,1$. There is a unique triple $(\hat\varrho,\hat v,\hat q)$ of smooth functions on $[0,1]$ with $\hat v(0)=\hat v(1)=0$ and $\hat w(1)=0$, where $\hat w\eqdef-\int_0^zik\cdot\hat v\dd z'$, such that
\begin{equation}\label{eq:RT-res}
s\hat\varrho+\rs'\hat w=G^\varrho,\qquad s\rs\hat v+f\rs\hat v^\perp-\mu(\pz^2-\nu^2)\hat v+ik\hat q=G^v,\qquad \pz\hat q=-g\hat\varrho,
\end{equation}
and, if $k=0$, $\hat q=-g\int_0^z\hat\varrho\dd z'$. Moreover $\hat\varrho=\pz\hat\varrho=0$ at $z=0,1$, and the triple for $(-k,\ol{G^\varrho},\ol{G^v})$ is the complex conjugate of the triple for $(k,G^\varrho,G^v)$.
\end{lemma}

\begin{proof}
\emph{The case $k=0$.} Then $\hat w=0$, so $\hat\varrho=G^\varrho/s$, which is clamped, and $\hat v$ must solve $s\rs\hat v+f\rs\hat v^\perp-\mu\hat v''=G^v$ with $\hat v(0)=\hat v(1)=0$. The form $s\int\rs\phi\cdot\bar\psi+f\int\rs\phi^\perp\cdot\bar\psi+\mu\int\phi'\cdot\ol{\psi'}$ on $H^1_0(0,1;\mathbb C^2)$ has real part $s\int\rs|\phi|^2+\mu\nrm{\phi'}^2$ on the diagonal, because $\phi^\perp\cdot\bar\phi$ is purely imaginary. The Lax--Milgram lemma gives a unique $\hat v\in H^1_0$, which is smooth by the equation, as in \Cref{lem:Lf}. With $\hat q=-g\int_0^z\hat\varrho$ the triple is unique.

\emph{The case $k\ne0$.} By \Cref{thm:RTlinear}(iii), $s>\Lambda_\infty>\lambda(\nu)$, so $Q_\nu(s):H^4\cap H^2_0\to L^2$ is an isomorphism by \Cref{thm:RTlinear}(ii). If $Q_\nu(s)\Psi\in H^m$, then $\Psi\in H^{m+4}$: if $\Psi\in H^j$, $j\ge2$, then $s\mu A_\nu\Psi=Q_\nu(s)\Psi-s^2M\Psi-sf^2R_\nu(s)\Psi-\wt K\Psi\in H^{\min(m,j-2)}$ by \Cref{lem:RT-ekman}(d), so $\Psi\in H^{\min(m+4,j+2)}$, and we iterate. Decompose $G^v=ikG_d+ik^\perp G_c$ with $G_d\eqdef-\frac i{\nu^2}k\cdot G^v$ and $G_c\eqdef-\frac i{\nu^2}k^\perp\cdot G^v$, and let
\begin{equation*}
\Psi\eqdef Q_\nu(s)^{-1}\mathfrak F,\qquad \mathfrak F\eqdef-s\,G_d'-gG^\varrho-sf\big(\rs B_s^{-1}G_c\big)'\in C^\infty([0,1]),
\end{equation*}
so that $\Psi\in C^\infty\cap H^2_0$, and
\begin{equation*}
\begin{aligned}
&\Phi\eqdef B_s^{-1}\big(G_c-f\rs\Psi'\big),\qquad \hat\varrho\eqdef\frac{G^\varrho-\nu^2\rs'\Psi}s,\\
&\hat q\eqdef G_d-s\rs\Psi'+f\rs\Phi+\mu(\Psi'''-\nu^2\Psi'),\qquad \hat v\eqdef ik\Psi'+ik^\perp\Phi .
\end{aligned}
\end{equation*}
These are smooth, $\hat v(0)=\hat v(1)=0$, and $\hat w=\nu^2\Psi$, so that $\hat w(1)=0$ and the first equation in \eqref{eq:RT-res} holds. Since $\hat v^\perp=ik^\perp\Psi'-ik\Phi$, the components of the second equation along $k$ and $k^\perp$ are the definition of $\hat q$ and $B_s\Phi=G_c-f\rs\Psi'$, respectively; since a vector of $\mathbb C^2$ vanishes if its products with $k$ and $k^\perp$ vanish, the second equation holds. For the third, $f(\rs\Phi)'=f(\rs B_s^{-1}G_c)'+f^2R_\nu(s)\Psi$, hence
\begin{equation*}
\begin{aligned}
\pz\hat q&=G_d'+sM\Psi+\mu A_\nu\Psi+f\big(\rs B_s^{-1}G_c\big)'+f^2R_\nu(s)\Psi\\
&=\frac1s\big(Q_\nu(s)\Psi-\wt K\Psi-\mathfrak F-gG^\varrho\big)=\frac{g\nu^2\rs'\Psi-gG^\varrho}s=-g\hat\varrho .
\end{aligned}
\end{equation*}
Since $\Psi=\Psi'=0$ and $G^\varrho=\pz G^\varrho=0$ at the lids, $\hat\varrho$ is clamped. For uniqueness, let $(\hat\varrho,\hat v,\hat q)$ solve \eqref{eq:RT-res} with $G^\varrho=0$, $G^v=0$, and decompose $\hat v$ by \eqref{eq:RT-pt}. As in the derivation of \eqref{eq:RT-nm}, $(\hat\varrho,\Psi,\Phi)$ is either zero or a normal mode with the real rate $s>\lambda(\nu)$; the latter is excluded by \Cref{thm:RTlinear}(ii). Hence $\hat v=0$, $\hat\varrho=0$, and $\hat q=0$ by the second equation.

\emph{Conjugation.} The complex conjugate of a solution for $(k,G^\varrho,G^v)$ solves \eqref{eq:RT-res} for $(-k,\ol{G^\varrho},\ol{G^v})$, because $\rs$, $f$, $g$, $\mu$ and $s$ are real; uniqueness gives the claim.
\end{proof}

\begin{lemma}\label{lem:approx}
For every $N\ge1$ there are smooth $(\hat\varrho_j,\hat v_j,\hat P_j)$, $1\le j\le N$, which are trigonometric polynomials in $x_h$ with frequencies in $\{mk_*:m\in\Z,\ |m|\le j\}$, such that $\hat\varrho_j=\pz\hat\varrho_j=0$ on $\pa\Omega$, $\hat v_j\in V_\sigma$, and $\hat P_j$ is independent of $z$, with the following property. Let
\begin{equation*}
(\varrho_a,v_a,P_a)\eqdef\sum_{j=1}^N\delta^je^{j\lambda_*t}(\hat\varrho_j,\hat v_j,\hat P_j),\qquad q_a\eqdef P_a-g\int_0^z\varrho_a\dd z',\qquad \rho_a\eqdef\rs+\varrho_a,
\end{equation*}
and define the residuals
\begin{equation*}
r_\varrho\eqdef\pt\varrho_a+\bU_a\cdot\nabla\varrho_a+w_a\rs',\qquad r_v\eqdef\rho_a(\pt v_a+\bU_a\cdot\nabla v_a)+f\rho_av_a^\perp+\nh q_a-\mu\Lap v_a .
\end{equation*}
Then $\nrm{r_\varrho(t)}_{H^k}+\nrm{r_v(t)}_{H^k}\le C_{N,k}\sum_{j=N+1}^{3N}(\delta e^{\lambda_*t})^j$ for all $t\ge0$, $\delta>0$ and $k\ge0$. Moreover $(\hat\varrho_1,\hat v_1)$ is the real growing mode \eqref{eq:RT-mode} with $k=k_*$ at $t=0$, and $\hat\varrho_1\ne0$.
\end{lemma}

\begin{proof}
Write $E\eqdef\delta e^{\lambda_*t}$, $\hat w_j\eqdef-\int_0^z\divh\hat v_j\dd z'$, $\hat\bU_j\eqdef(\hat v_j,\hat w_j)$ and $\hat q_j\eqdef\hat P_j-g\int_0^z\hat\varrho_j\dd z'$. Inserting the ansatz, $r_\varrho$ and $r_v$ are polynomials in $E$ without constant term, and the coefficient of $E^j$ vanishes if and only if
\begin{equation}\label{eq:orderj}
\begin{aligned}
&j\lambda_*\hat\varrho_j+\hat w_j\rs'=F^\varrho_j\eqdef-\sum_{a+b=j}\hat\bU_a\cdot\nabla\hat\varrho_b,\\
&j\lambda_*\rs\hat v_j+f\rs\hat v_j^\perp-\mu\Lap\hat v_j+\nh\hat q_j=F^v_j,\\
&F^v_j\eqdef-\sum_{a+b=j}\big(b\lambda_*\hat\varrho_a\hat v_b+f\hat\varrho_a\hat v_b^\perp+\rs\hat\bU_a\cdot\nabla\hat v_b\big)-\sum_{a+b+c=j}\hat\varrho_a\hat\bU_b\cdot\nabla\hat v_c,
\end{aligned}
\end{equation}
where $a,b,c\ge1$; here $f\rho_av_a^\perp$ contributes $f\rs\hat v_j^\perp$ to the left side and $f\hat\varrho_a\hat v_b^\perp$ to $F^v_j$. For $j=1$ the right-hand sides vanish. We take $\hat\varrho_1=\cos(k_*\cdot x_h)\hat\varrho_{\nu_*}$, $\hat v_1=-\sin(k_*\cdot x_h)(\phi_{\nu_*}'k_*+\Phi_{\nu_*}k_*^\perp)$ and $\hat q_1=\cos(k_*\cdot x_h)\hat q_{\nu_*}$ from \eqref{eq:RT-mode}, with $\nu_*=|k_*|$; then \eqref{eq:orderj} with $j=1$ holds by \Cref{thm:RTlinear}(i), $\hat P_1=\hat q_1+g\int_0^z\hat\varrho_1$ is independent of $z$, $\hat v_1\in V_\sigma$, $\hat\varrho_1$ is clamped and $\hat\varrho_1\ne0$.

Let $j\ge2$ and assume that the profiles of order $1,\dots,j-1$ have been constructed. The forcings $F^\varrho_j$, $F^v_j$ are smooth trigonometric polynomials in $x_h$ with frequencies in $\{mk_*:|m|\le j\}$. Since $\hat v_a=\hat w_a=0$ and $\hat\varrho_b=\pz\hat\varrho_b=\nh\hat\varrho_b=0$ on the lids, $F^\varrho_j=0$ there, and
\begin{equation*}
\pz F^\varrho_j=-\sum_{a+b=j}\big(\pz\hat v_a\cdot\nh\hat\varrho_b+\hat v_a\cdot\nh\pz\hat\varrho_b+\pz\hat w_a\,\pz\hat\varrho_b+\hat w_a\,\pz^2\hat\varrho_b\big)=0\quad\text{on the lids.}
\end{equation*}
For each frequency $k=mk_*$, $|m|\le j$, we apply \Cref{lem:RT-resolvent} with $s=j\lambda_*\ge2\lambda_*>\Lambda_\infty$ and $(G^\varrho,G^v)$ the $k$-th Fourier coefficients of $(F^\varrho_j,F^v_j)$, and obtain $(\hat\varrho_{j,k},\hat v_{j,k},\hat q_{j,k})$. Let $\hat\varrho_j$, $\hat v_j$, $\hat q_j$ be the corresponding sums $\sum_ke^{ik\cdot x_h}(\cdot)$; they are real because the coefficients of $\pm k$ of the real forcings are complex conjugate, and so are the solutions. Then $\hat w_j=\sum_ke^{ik\cdot x_h}\hat w_{j,k}$, \eqref{eq:orderj} holds, $\hat\varrho_j$ is clamped, and $\hat P_j\eqdef\hat q_j+g\int_0^z\hat\varrho_j$ is independent of $z$ by $\pz\hat q_{j,k}=-g\hat\varrho_{j,k}$. Moreover $\hat v_j\in H^1_0$ and $\divh\int_0^1\hat v_j\dd z=-\sum_{k}e^{ik\cdot x_h}\hat w_{j,k}(1)=0$, so $\hat v_j\in V_\sigma$.

After these choices, the coefficients of $E^j$, $j\le N$, vanish. What remains in $r_\varrho$ and $r_v$ are the products of two or three profiles in \eqref{eq:orderj} with total order $j\in[N+1,3N]$ and all orders at most $N$; each is of the form $E^j\times$(a smooth function independent of $t$ and $\delta$). This gives the bound.
\end{proof}

Note that the approximate solution depends on $x_h$ only through $k_*\cdot x_h$; for $f\ne0$ it nevertheless has a nonzero velocity component along $k_*^\perp$, since $\Phi_{\nu_*}\ne0$.

To produce compatible data we need the rotating version of \Cref{lem:stokes}.

\begin{lemma}\label{lem:RT-stokesf}
Let $\rs\in C^\infty([0,1])$ with $\min\rs>0$, and $f\in\R$. For every $m\ge0$ and $\zeta\in H^m(\Omega;\R^2)$ there is a unique $v\in H^{m+2}\cap V_\sigma$ for which there is $P=P(x_h)$ with
\begin{equation*}
-\mu\Lap v+f\rs v^\perp+\nh P=\zeta,
\end{equation*}
and $\nrm{v}_{H^{m+2}}+\nrm{\nh P}_{H^m(\Tl)}\le C_m\nrm{\zeta}_{H^m}$, where $C_m$ depends on $m,\mu,f,l$ and $\rs$. We write $v=\mathcal S_f\zeta$.
\end{lemma}

\begin{proof}
Let $\mathcal S$ be the solution operator $\zeta\mapsto v$ of \Cref{lem:stokes}. On the closed subspace $X\eqdef H^2\cap V_\sigma$ of $H^2$, the operator $\mathcal K v\eqdef f\mathcal S(\rs v^\perp)$ maps $X$ boundedly into $H^4\cap V_\sigma$, hence is compact on $X$. The equation is equivalent to $(I+\mathcal K)v=\mathcal S\zeta$ in $X$. If $(I+\mathcal K)v=0$, then $-\mu\Lap v+\nh P=-f\rs v^\perp$ for some $P=P(x_h)$; testing with $v$ and using $v^\perp\cdot v=0$ and $\int_\Omega\nh P\cdot v\dx=-\int_{\Tl}P\,\divh\int_0^1v\dd z\dxh=0$, we get $\mu\nrm{\nabla v}^2=0$, so $v=0$. By the Fredholm alternative, $I+\mathcal K$ is an isomorphism of $X$, which gives existence, uniqueness and the bound for $m=0$. For $m\ge1$, \Cref{lem:stokes} applied to $-\mu\Lap v+\nh P=\zeta-f\rs v^\perp$ gives $\nrm{v}_{H^{m+2}}+\nrm{\nh P}_{H^m}\le C(\nrm\zeta_{H^m}+|f|\nrm{\rs}_{C^m}\nrm v_{H^m})$, and induction on $m$ concludes.
\end{proof}

\begin{lemma}\label{lem:compatcorr}
Fix $N$ and let $(\varrho_a,v_a,P_a)$ be as in \Cref{lem:approx}. There is $\delta_1>0$ such that, for $\delta\in(0,\delta_1)$, there is a smooth $v_\delta\in V_\sigma$ with $\nrm{v_\delta}_{H^k}\le C_k\delta^{N+1}$ for all $k$, such that the data
\begin{equation*}
\vin\eqdef\varrho_a(0),\qquad \uin\eqdef v_a(0)+v_\delta
\end{equation*}
satisfy the hypotheses of \Cref{prop:local}. Moreover, the solution satisfies $v_t(0)=\pt v_a(0)$.
\end{lemma}

\begin{proof}
Let $\rin=\rs+\varrho_a(0)$, and let $v=v_\delta$ be the small solution, in $H^{k+2}\cap V_\sigma$ for every $k\ge2$, of
\begin{equation}\label{eq:vcorr}
-\mu\Lap v+f\rin v^\perp+\nh\pi=-r_v(0)-\rin\big(\bU_a(0)\cdot\nabla v+\mathbf V\cdot\nabla v_a(0)+\mathbf V\cdot\nabla v\big),\qquad v\in V_\sigma,
\end{equation}
where $\mathbf V=(v,-\int_0^z\divh v\dd z')$ and $\pi=\pi(x_h)$. Writing $f\rin v^\perp=f\rs v^\perp+f\varrho_a(0)v^\perp$, \eqref{eq:vcorr} is the fixed-point problem $v=\mathcal S_f\Xi(v)$ with
\begin{equation*}
\Xi(v)\eqdef-r_v(0)-f\varrho_a(0)v^\perp-\rin\big(\bU_a(0)\cdot\nabla v+\mathbf V\cdot\nabla v_a(0)+\mathbf V\cdot\nabla v\big).
\end{equation*}
For $k\ge2$, by \Cref{lem:GN}(2),(3) and \Cref{lem:w}, $\nrm{\mathbf V\cdot\nabla v'}_{H^k}\lesssim\nrm{\mathbf V}_{L^\infty}\nrm{v'}_{H^{k+1}}+\nrm{\mathbf V}_{H^k}\nrm{\nabla v'}_{L^\infty}\lesssim\nrm v_{H^3}\nrm{v'}_{H^{k+1}}+\nrm v_{H^{k+1}}\nrm{v'}_{H^3}$. Since $\nrm{\varrho_a(0)}_{C^k}+\nrm{v_a(0)}_{H^{k+3}}\le C_k\delta$ and $\nrm{\rin}_{C^k}\le C_k$ for $\delta\le1$, this gives
\begin{equation*}
\begin{aligned}
&\nrm{\Xi(v)}_{H^k}\le\nrm{r_v(0)}_{H^k}+C_k\big(\delta+\nrm v_{H^3}\big)\nrm v_{H^{k+1}},\\
&\nrm{\Xi(v)-\Xi(v')}_{H^2}\le C\big(\delta+\nrm v_{H^4}+\nrm{v'}_{H^4}\big)\nrm{v-v'}_{H^4}.
\end{aligned}
\end{equation*}
With $k=2$, \Cref{lem:RT-stokesf} and $\nrm{r_v(0)}_{H^2}\le C_N\delta^{N+1}$ (\Cref{lem:approx}), the map $\mathcal S_f\Xi$ is a contraction of the ball $\{v\in H^4\cap V_\sigma:\nrm v_{H^4}\le2C_2C_N\delta^{N+1}\}$ for $\delta<\delta_1$ small. Higher regularity, with $\delta_1$ independent of $k$, follows by induction: if $\nrm v_{H^{k+1}}\le C\delta^{N+1}$ for some $k\ge3$, then $\nrm{\Xi(v)}_{H^k}\le C\delta^{N+1}$, and $v=\mathcal S_f\Xi(v)$ gives $\nrm v_{H^{k+2}}\le C_k\delta^{N+1}$ by \Cref{lem:RT-stokesf}. Hence $v_\delta\in C^\infty(\ol\Omega)$.

Let $\bU_{\mathrm{in}}=\bU_a(0)+\mathbf V$ and $F_0\eqdef\mu\Lap\uin-\rin\bU_{\mathrm{in}}\cdot\nabla\uin-f\rin\uin^\perp+g\int_0^z\nh\vin\dd z'$. By the definition of $r_v$ at $t=0$ (note $\rho_a(0)=\rin$) and by \eqref{eq:vcorr},
\begin{equation*}
\begin{aligned}
F_0&=\rin\pt v_a(0)+\nh P_a(0)-r_v(0)+\Big[\mu\Lap v-\rin\big(\bU_a(0)\cdot\nabla v+\mathbf V\cdot\nabla v_a(0)+\mathbf V\cdot\nabla v\big)-f\rin v^\perp\Big]\\
&=\rin\pt v_a(0)+\nh\big(P_a(0)+\pi\big).
\end{aligned}
\end{equation*}
Let $P_0$ be the zero-mean solution of \eqref{eq:Peq} at $t=0$ and $u_1\eqdef\rin^{-1}(F_0-\nh P_0)$, the initial acceleration given by \eqref{PE}$_2$ with $q_{\mathrm{in}}=P_0-g\int_0^z\vin$. The computation leading to \eqref{eq:Phieq}, together with \eqref{eq:Peq}, gives $\divh\int_0^1u_1\dd z=0$; and $\divh\int_0^1\pt v_a(0)\dd z=0$ because $v_a(t)\in V_\sigma$ for all $t$. With $\Pi\eqdef P_a(0)+\pi-P_0$ and $a\eqdef\int_0^1\rin^{-1}\dd z>0$ we have $u_1-\pt v_a(0)=\rin^{-1}\nh\Pi$, hence $\divh(a\nh\Pi)=0$ on $\Tl$. Testing with $\Pi$ gives $\nh\Pi=0$, so $u_1=\pt v_a(0)=\sum_jj\lambda_*\delta^j\hat v_j$, which vanishes on $\pa\Omega$. This is \eqref{eq:compat}. The remaining hypotheses hold for $\delta$ small: $\vin\in H^4$ is clamped, $\uin\in H^4\cap V_\sigma$ (so \eqref{eq:lid0} holds), and $\rs+\vin\ge\frac12\min\rs$. Finally, the solution of \Cref{prop:local} satisfies \eqref{PE}$_2$ at $t=0$ with the surface pressure of \Cref{lem:propagation}, so $v_t(0)=u_1=\pt v_a(0)$.
\end{proof}

\begin{proof}[Proof of \Cref{thm:RTnonlinear}]
\emph{Choice of constants.} For $t\le T^\delta$ we have $\delta e^{\lambda_*t}\le\bar\varepsilon_0$. Hence there is $\bar\varepsilon_0(N)>0$ such that, if $\bar\varepsilon_0\le\bar\varepsilon_0(N)$, the approximate solution of \Cref{lem:approx} satisfies on $[0,T^\delta]$
\begin{equation}\label{eq:appsmall}
\nrm{\varrho_a}_{W^{1,\infty}}+\nrm{v_a}_{W^{1,\infty}}+\nrm{w_a}_{W^{1,\infty}}+\nrm{\pt v_a+\bU_a\cdot\nabla v_a}_{L^\infty}\le1,\qquad \tfrac12\min\rs\le\rho_a\le2\max\rs .
\end{equation}
Let $C_0$ be the constant in \eqref{eq:Yineq} below. It depends only on $\rs,\mu,g,f$, through the bounds \eqref{eq:appsmall}, which do not involve $N$; thus $C_0$ does not depend on $N$. We first fix $N$ with $(N+1)\lambda_*\ge2C_0$. The constant $C_N'$ below and $\delta_1$ from \Cref{lem:compatcorr} depend only on $N$ (and on $\rs,\mu,g,f,l,k_*$). Then we fix $\bar\varepsilon_0\le\min\{1,\bar\varepsilon_0(N),\bar\varepsilon_0'(N)\}$, with $\bar\varepsilon_0(N)$ from \eqref{eq:appsmall} and $\bar\varepsilon_0'(N)$ from the conclusion below, and finally $\delta_0\le\min\{\delta_1,\bar\varepsilon_0\}$. Let $(\varrho,v)$ be the solution with the data of \Cref{lem:compatcorr} and $T_1\eqdef\min\{T^\delta_{\max},T^\delta\}$. By transport and \eqref{eq:appsmall} at $t=0$, $\frac12\min\rs\le\rho\le2\max\rs$ on $[0,T^\delta_{\max})$.

\emph{Error equations.} Put $S\eqdef\varrho-\varrho_a$, $u\eqdef v-v_a$, $w_u\eqdef w-w_a$, $\mathbf u\eqdef\bU-\bU_a=(u,w_u)$ and $Q\eqdef q-q_a$, so that $\pz Q=-gS$, $u(t)\in V_\sigma$ and $w_u=0$ on $\pa\Omega$. Subtracting the equations for the approximate solution from \eqref{PE}, and using $\rho=\rho_a+S$ and $\rho v^\perp-\rho_av_a^\perp=\rho u^\perp+Sv_a^\perp$,
\begin{equation*}
\begin{aligned}
&S_t+\bU\cdot\nabla S+\mathbf u\cdot\nabla\varrho_a+w_u\rs'=-r_\varrho,\\
&\rho(u_t+\bU\cdot\nabla u)+f\rho u^\perp+\rho\,\mathbf u\cdot\nabla v_a+S\big(\pt v_a+\bU_a\cdot\nabla v_a+fv_a^\perp\big)+\nh Q-\mu\Lap u=-r_v .
\end{aligned}
\end{equation*}

\emph{Energy estimate.} Let $Y\eqdef\frac12\nrm S^2+\frac12\nrm{\sqrt\rho u}^2$. We test the first equation with $S$ and the second with $u$. Since $\divv\bU=0$, $w=0$ on $\pa\Omega$ and $\rho_t+\bU\cdot\nabla\rho=0$, the transport terms produce exact derivatives: $\int(S_t+\bU\cdot\nabla S)S=\frac12\ddt\nrm S^2$ and $\int\rho(u_t+\bU\cdot\nabla u)\cdot u=\frac12\ddt\nrm{\sqrt\rho u}^2$. The Coriolis term is skew, $\int f\rho u^\perp\cdot u=0$. By \Cref{lem:propagation}, $\int\nh Q\cdot u=\int Q\pz w_u=-\int\pz Q\,w_u=g\int Sw_u$. We use $\nrm{w_u}\le\nrm{\divh u}\le\nrm{\nh u}$ (\Cref{lem:w}), $\nrm u\le\pi^{-1}\nrm{\pz u}$, hence $\nrm{\mathbf u}\le2\nrm{\nabla u}$, and \eqref{eq:appsmall}:
\begin{equation*}
\begin{aligned}
\ddt Y+\mu\nrm{\nabla u}^2&\le2\big(1+\nrm{\rs'}_{L^\infty}\big)\nrm{\nabla u}\nrm S+4\max\rs\,\nrm{\nabla u}\nrm u+(1+|f|)\nrm S\nrm u\\
&\quad+g\nrm S\nrm{\nabla u}+\nrm{r_\varrho}\nrm S+\nrm{r_v}\nrm u .
\end{aligned}
\end{equation*}
Since $\nrm u^2\le\frac2{\min\rs}\nrm{\sqrt\rho u}^2$, Young's inequality gives, on $[0,T_1)$,
\begin{equation}\label{eq:Yineq}
\ddt Y+\frac\mu2\nrm{\nabla u}^2\le2C_0Y+C\big(\nrm{r_\varrho}+\nrm{r_v}\big)Y^{\frac12}.
\end{equation}
The identity is classical: in the class \eqref{eq:localclass}, $\varrho\in C^1([0,T];H^3)$ by \eqref{PE}$_1$ and $v\in C^1([0,T];H^2)$, so $S\in C^1(L^2)$ and $u\in C^1(H^2\cap V_\sigma)$, and all integrations by parts are legitimate. Applying \eqref{eq:Yineq} to $(Y+\eta)^{1/2}$ and letting $\eta\to0$ gives $\ddt Y^{1/2}\le C_0Y^{1/2}+C(\nrm{r_\varrho}+\nrm{r_v})$ in the integrated sense. Since $S(0)=0$, $Y(0)^{1/2}\le C\nrm{v_\delta}\le C\delta^{N+1}$, and by \Cref{lem:approx}, $\nrm{r_\varrho(\tau)}+\nrm{r_v(\tau)}\le C_N(\delta e^{\lambda_*\tau})^{N+1}$ for $\tau\le T^\delta$ (because $\delta e^{\lambda_*\tau}\le\bar\varepsilon_0\le1$). Hence
\begin{equation*}
Y(t)^{\frac12}\le C\delta^{N+1}e^{C_0t}+C_N\int_0^te^{C_0(t-\tau)}(\delta e^{\lambda_*\tau})^{N+1}\dd\tau\le C_N'(\delta e^{\lambda_*t})^{N+1},\qquad t<T_1,
\end{equation*}
where we used $(N+1)\lambda_*-C_0\ge C_0$. By continuity the bound also holds at $t=T^\delta$ if $T^\delta<T^\delta_{\max}$.

\emph{Conclusion.} Suppose $T^\delta_{\max}>T^\delta$. At $t=T^\delta$ we have $\delta e^{\lambda_*T^\delta}=\bar\varepsilon_0$, hence
\begin{equation*}
\nrm{\varrho(T^\delta)}\ge\bar\varepsilon_0\nrm{\hat\varrho_1}-\sum_{j=2}^N\bar\varepsilon_0^{\,j}\nrm{\hat\varrho_j}-\sqrt2\,C_N'\bar\varepsilon_0^{\,N+1}\ge\frac{\bar\varepsilon_0}2\nrm{\hat\varrho_1}\eqdef\varepsilon_1
\end{equation*}
provided $\bar\varepsilon_0\le\bar\varepsilon_0'(N)$, where $\bar\varepsilon_0'(N)>0$ is such that $\sum_{j=2}^N\bar\varepsilon_0'(N)^{j-1}\nrm{\hat\varrho_j}+\sqrt2\,C_N'\bar\varepsilon_0'(N)^N\le\frac12\nrm{\hat\varrho_1}$; recall $\hat\varrho_1\ne0$. The data satisfy $\nrm{\vin}_{H^m}+\nrm{\uin}_{H^m}\le C_m\delta$ by \Cref{lem:approx,lem:compatcorr}. This proves the theorem.
\end{proof}

\begin{remark}\label{rem:RT}
(i) The linearized problem \eqref{eq:RT-lin} is well posed in $L^2$, with a growth bound uniform in all horizontal modes and in $f$. Indeed, testing with $(\varrho,v)$, the Coriolis term drops out and $\int\nh q\cdot v=g\int\varrho w$, so
\begin{equation*}
\frac12\ddt\big(\nrm\varrho^2+\nrm{\sqrt\rs v}^2\big)+\mu\nrm{\nabla v}^2=-\int_\Omega(g+\rs')w\varrho\dx\le\big(g+\nrm{\rs'}_{L^\infty}\big)\nrm{\nabla v}\nrm\varrho,
\end{equation*}
and $\nrm\varrho^2+\nrm{\sqrt\rs v}^2$ grows at most like $e^{(g+\nrm{\rs'}_{L^\infty})^2t/(2\mu)}$. So there is no Hadamard-type ill-posedness, and the normal-mode growth rates are bounded by $\Lambda_\infty$. High horizontal frequencies, however, grow at almost the maximal rate. This differs from the viscous non-hydrostatic (Navier--Stokes) case: there, testing the eigenvalue problem with the conjugate velocity (the Coriolis term then contributes only an imaginary part) shows that a real growth rate of a mode with horizontal wave vector $k$ is at most $g\max\rs'/(\mu|k|^2)$, with or without rotation, so that a most unstable wavelength exists.

(ii) Rotation lowers every growth rate, \eqref{eq:RT-rotbound}, and for each fixed $\nu$ the growth rate tends to zero as $|f|\to\infty$. But it does not lower the supremum $\Lambda_\infty$. The Coriolis force turns part of the poloidal motion, which carries the buoyancy flux, into toroidal motion, which does not, and the toroidal motion is damped by viscosity at the horizontal scale $\nu^{-1}$. Accordingly, by \Cref{lem:RT-ekman}(b) the Ekman term in $\cJ_{\nu,s}$ satisfies $0\le f^2\mathfrak e_\nu(s;\psi)\le sf^2\vartheta_\nu\ip{M\psi,\psi}=O(f^2\nu^{-2})$, whereas the viscous and buoyancy terms are of order $\nu^2$, and $s\mu\ip{A_\nu\psi,\psi}\ge s\vartheta_\nu^{-1}\ip{M\psi,\psi}$; the ratio of the Ekman term to the viscous term is at most $f^2\vartheta_\nu^2$, a Taylor number built on the horizontal wavelength, which is the parameter in \eqref{eq:RT-rotbound}. Short horizontal waves therefore do not feel rotation, and the convergence $\lambda(\nu;f)\to0$ as $|f|\to\infty$ is not uniform in $\nu$.

(iii) We expect the first alternative in \Cref{thm:RTnonlinear} never to occur. Excluding it would require a high-order estimate for $(S,u)$ of the form $\ddt\mathcal E\le C_1\mathcal E+\cdots$, with $C_1$ independent of $N$, on time intervals of length of order $\log(1/\delta)$. The short-time estimates of \Cref{app:local} do not provide this, and we do not pursue it here.

(iv) Apart from the borderline case $\max\rs'=0$, \Cref{thm:main} and \Cref{thm:RTlinear,thm:RTnonlinear} give a dichotomy ($\rs\in C^\infty$ for the nonlinear statement), for every $f\in\R$. If $\rs'<0$ on $[0,1]$, small perturbations converge to a stratified state at an algebraic rate. If $\rs'>0$ somewhere, $(\rs,0)$ is linearly unstable in every horizontal mode $k\ne0$ and nonlinearly unstable.
\end{remark}

\section{Optimality of the decay rate}\label{sec:opt}

We return to the stable case. Throughout this section \eqref{eq:stable} holds, $\sgm=-1/\rs'$, and $f\in\R$ is arbitrary. We show that the following rates cannot be improved uniformly for small $H^4$ data: the rate $(1+t)^{-(1-r/4)}$ for $\wt\varrho$ and for $\rho-\rinf$ in $H^r$, $0\le r<4$, which includes the rate $(1+t)^{-1}$ in \eqref{eq:ratefull}; the rate $(1+t)^{-1}$ for $\nh\varrho$ in \eqref{eq:decay}, for the pressure in \eqref{eq:ratepress} and for $X^z-X^z_\infty$ in $L^2$; and the rate $(1+t)^{-2}$ in \eqref{eq:E-tail}. In $H^4$ there is no uniform decay rate at all. The velocity terms in \eqref{eq:decay}, the faster rates for the mean, the $L^\infty$ rates and the $C^1$ rates for the flow map are not addressed (\Cref{rem:O}). The mechanism is simple. Internal gravity waves of high vertical wavenumber $n$ are overdamped by viscosity. The slow branch relaxes under buoyancy at a rate of order $n^{-4}$, while the $H^4$ norm of such a wave exceeds its $L^2$ norm by a factor of order $n^4$. Here we construct these waves exactly, as real normal modes of the rotating linearized problem. Rotation enters only through the Ekman-type operator of \Cref{sec:RT-lin}.

We use the notation of \Cref{sec:RT-lin}: the linearization \eqref{eq:RT-lin}, the decomposition \eqref{eq:RT-pt}, the operators \eqref{eq:RT-ops}, normal modes (\Cref{def:RT-nm}) and the pencil \eqref{eq:pencil}. None of these used the sign of $\rs'$. Under \eqref{eq:stable}, $\wt K=g\nu^2\sgm^{-1}$ is a positive multiplication operator. As in \Cref{sec:RT}, $M$ denotes the operator in \eqref{eq:RT-ops}. We fix $k\in(l^{-1}\Z)^2\setminus\{0\}$ and write $\nu=|k|$. Let $m=\min_{[0,1]}(-\rs')$ and $M_1=\max_{[0,1]}(-\rs')$ be as in \Cref{thm:energy}, let $\vartheta_\nu=\max\rs/(\mu(\pi^2+\nu^2))$ be as in \eqref{eq:RT-rotbound}, and set
\begin{equation}\label{eq:O-const}
C_M\eqdef\frac{\max\rs}{\pi^2\mu},\qquad \hat s\eqdef\frac{\mu\pi^2}{4\max\rs}=\frac1{4C_M}.
\end{equation}
In \Cref{lem:O-ekman,prop:O-modes} and their proofs, $\ip{\cdot,\cdot}$ and $\nrm\cdot$ denote the complex inner product and the norm of $L^2(0,1)$; elsewhere in this section $\nrm\cdot$ is the norm of $L^2(\Omega)$. For $\phi\in H^2_0(0,1)$ we have $\phi'\in H^1_0$, hence $\nrm{\phi''}^2\ge\pi^2\nrm{\phi'}^2$, and therefore
\begin{equation}\label{eq:O-A}
\begin{aligned}
&\nrm{\phi''}^2\le\ip{A_\nu\phi,\phi}\le\Big(1+\frac{\nu^2}{\pi^2}\Big)\nrm{\phi''}^2,\qquad (\pi^2+\nu^2)\nrm{\phi'}^2\le\ip{A_\nu\phi,\phi},\\
&\ip{M\phi,\phi}\le C_M\,\mu\ip{A_\nu\phi,\phi}.
\end{aligned}
\end{equation}

\subsection{The Ekman operator at negative real rates}

In \Cref{sec:RT-lin} the operators $B_\lambda$, $R_\nu(\lambda)$ and $Q_\nu(\lambda)$ were defined for $\mathrm{Re}\,\lambda\ge0$. We now use them at $\lambda=-s$ with $0<s\le\hat s$.

\begin{lemma}\label{lem:O-ekman}
Let $0<s\le\hat s$.
\begin{enumerate}[label=(\alph*)]
\item $B_{-s}=\mu L_\nu-s\rs$ satisfies $\mathrm{Re}\ip{B_{-s}\eta,\eta}\ge\frac34\mu\ip{L_\nu\eta,\eta}$ for $\eta\in H^1_0(0,1)$. Hence $B_{-s}:H^1_0\to H^{-1}$ is an isomorphism, $\nrm{B_{-s}^{-1}\chi}\le\frac{4}{3\mu(\pi^2+\nu^2)}\nrm\chi$, and $\nrm{B_{-s}^{-1}\chi}_{H^{j+2}}\le C\nrm\chi_{H^j}$ for $j=0,1,2$, with $C$ independent of $s\in(0,\hat s]$.
\item Define $R_\nu(-s)\eqdef-\pz\rs B_{-s}^{-1}\rs\pz$ and, for $\phi,\xi\in H^1_0(0,1)$,
\begin{equation*}
\mathfrak e^-_\nu(s;\phi,\xi)\eqdef s\ip{R_\nu(-s)\phi,\xi}=s\ip{B_{-s}^{-1}(\rs\phi'),\rs\xi'},\qquad \mathfrak e^-_\nu(s;\phi)\eqdef\mathfrak e^-_\nu(s;\phi,\phi).
\end{equation*}
Then $\mathfrak e^-_\nu(s;\cdot,\cdot)$ is a bounded Hermitian form on $H^1_0$, real on real functions. With $\chi\eqdef\rs\phi'$,
\begin{equation}\label{eq:O-dual}
\mathfrak e^-_\nu(s;\phi)=\max_{\eta\in H^1_0(0,1)}\Big\{2\,\mathrm{Re}\ip{\chi,\eta}-\frac\mu s\big(\nrm{\eta'}^2+\nu^2\nrm\eta^2\big)+\int_0^1\rs|\eta|^2\Big\},
\end{equation}
and
\begin{equation}\label{eq:O-ebd}
0\le\mathfrak e^-_\nu(s;\phi)\le\frac{4s\nrm{\rs\phi'}^2}{3\mu(\pi^2+\nu^2)}\le\frac43s\vartheta_\nu^2\,\mu\ip{A_\nu\phi,\phi}\qquad(\phi\in H^2_0).
\end{equation}
with $\mathfrak e^-_\nu(s;\phi)>0$ if $\phi\in H^2_0\setminus\{0\}$. Moreover $s\mapsto\mathfrak e^-_\nu(s;\phi)$ is nondecreasing on $(0,\hat s]$, and for $0<s_1<s_2\le\hat s$,
\begin{equation}\label{eq:O-elip}
0\le\mathfrak e^-_\nu(s_2;\phi)-\mathfrak e^-_\nu(s_1;\phi)\le\frac{2s_2}{s_1}\,\frac{(s_2-s_1)\nrm{\rs\phi'}^2}{\mu(\pi^2+\nu^2)} .
\end{equation}
\item Define $Q_\nu(-s)$ by \eqref{eq:pencil} with $\lambda=-s$. Then, for $\phi,\xi\in H^2_0(0,1)$,
\begin{equation}\label{eq:O-Q}
\ip{Q_\nu(-s)\phi,\xi}=\ip{\wt K\phi,\xi}-s\mu\ip{A_\nu\phi,\xi}-f^2\mathfrak e^-_\nu(s;\phi,\xi)+s^2\ip{M\phi,\xi}.
\end{equation}
If $\Psi\in H^2_0\setminus\{0\}$ and $\ip{Q_\nu(-s)\Psi,\xi}=0$ for all $\xi\in H^2_0$, then $(\nu^2\rs'\Psi/s,\Psi,-fB_{-s}^{-1}(\rs\Psi'))$ is a normal mode with rate $-s$.
\end{enumerate}
\end{lemma}

\begin{proof}
(a) Since $s\max\rs\le\hat s\max\rs=\frac{\mu\pi^2}4$ and $\pi^2\nrm\eta^2\le\nrm{\eta'}^2$, we get $s\int_0^1\rs|\eta|^2\le\frac\mu4\nrm{\eta'}^2\le\frac14\mu\ip{L_\nu\eta,\eta}$, which gives the coercivity. The Lax--Milgram lemma gives the isomorphism. For $\eta=B_{-s}^{-1}\chi$, $\frac34\mu(\pi^2+\nu^2)\nrm\eta^2\le\frac34\mu\ip{L_\nu\eta,\eta}\le\mathrm{Re}\ip{\chi,\eta}\le\nrm\chi\nrm\eta$, which gives the $L^2$ bound, and also $\nrm{\eta}_{H^1}\le C\nrm\chi$. The equation $\mu\eta''=\mu\nu^2\eta-s\rs\eta-\chi$ and its first two derivatives give the $H^{j+2}$ bounds, since $\rs\in C^2$ and $s\le\hat s$.

(b) $B_{-s}^{-1}$ is bounded and self-adjoint on $L^2$ and maps real functions to real functions; this gives the first assertion. The operator $T_s\eqdef s^{-1}B_{-s}$ is the positive self-adjoint operator of the coercive form $\frac\mu s\ip{L_\nu\eta,\zeta}-\int\rs\eta\bar\zeta$ on $H^1_0$, and $\mathfrak e^-_\nu(s;\phi)=\ip{T_s^{-1}\chi,\chi}$. With $\eta_0\eqdef T_s^{-1}\chi$, expanding the square as in the proof of \Cref{lem:RT-ekman}(a) gives $2\,\mathrm{Re}\ip{\chi,\eta}-\ip{T_s\eta,\eta}=\ip{T_s^{-1}\chi,\chi}-\ip{T_s(\eta-\eta_0),\eta-\eta_0}$, which proves \eqref{eq:O-dual}. Taking $\eta=0$ in \eqref{eq:O-dual} gives $\mathfrak e^-_\nu\ge0$. Moreover $\mathfrak e^-_\nu(s;\phi)=\ip{T_s\eta_0,\eta_0}>0$ if $\phi\in H^2_0\setminus\{0\}$, because then $\chi\ne0$ and hence $\eta_0\ne0$. By (a), $\ip{T_s\eta,\eta}\ge\frac{3\mu}{4s}\ip{L_\nu\eta,\eta}\ge\frac{3\mu(\pi^2+\nu^2)}{4s}\nrm\eta^2$, and $2\,\mathrm{Re}\ip{\chi,\eta}-a\nrm\eta^2\le a^{-1}\nrm\chi^2$ gives the first upper bound in \eqref{eq:O-ebd}. The second follows from $\nrm{\rs\phi'}^2\le(\max\rs)^2\nrm{\phi'}^2$ and \eqref{eq:O-A}. For fixed $\eta$ the expression in braces in \eqref{eq:O-dual} is nondecreasing in $s$, and so is its maximum. Let $0<s_1<s_2\le\hat s$, and let $\eta_2\eqdef T_{s_2}^{-1}\chi$ be the maximizer for $s_2$. Testing \eqref{eq:O-dual} for $s_1$ with $\eta_2$ gives $\mathfrak e^-_\nu(s_1;\phi)\ge\mathfrak e^-_\nu(s_2;\phi)-\mu(s_1^{-1}-s_2^{-1})\ip{L_\nu\eta_2,\eta_2}$. Testing $T_{s_2}\eta_2=\chi$ with $\eta_2$,
\begin{equation*}
\frac{3\mu}{4s_2}\ip{L_\nu\eta_2,\eta_2}\le\mathrm{Re}\ip{\chi,\eta_2}\le\nrm\chi\,(\pi^2+\nu^2)^{-\frac12}\ip{L_\nu\eta_2,\eta_2}^{\frac12},
\end{equation*}
so $\ip{L_\nu\eta_2,\eta_2}\le\frac{16s_2^2\nrm\chi^2}{9\mu^2(\pi^2+\nu^2)}$ and $\mu(s_1^{-1}-s_2^{-1})\ip{L_\nu\eta_2,\eta_2}\le\frac{16}9\frac{s_2}{s_1}\frac{(s_2-s_1)\nrm\chi^2}{\mu(\pi^2+\nu^2)}$. This proves \eqref{eq:O-elip}.

(c) Identity \eqref{eq:O-Q} is \eqref{eq:pencil} at $\lambda=-s$, since $\lambda f^2\ip{R_\nu(\lambda)\phi,\xi}=-f^2\mathfrak e^-_\nu(s;\phi,\xi)$. The proof of \Cref{lem:RT-modes}(b) uses only that $B_\lambda:H^1_0\to H^{-1}$ is an isomorphism, which holds by (a). For $\lambda=-s$ it gives $\hat\varrho=-\nu^2\rs'\Psi/\lambda=\nu^2\rs'\Psi/s$.
\end{proof}

\subsection{Slow normal modes}

Let $\mathbb H_\nu$ be the space $H^2_0(0,1)$ of complex functions with the inner product $(\phi,\xi)_A\eqdef\mu\ip{A_\nu\phi,\xi}$, which is equivalent to the $H^2$ inner product by \eqref{eq:O-A} and the Poincar\'e inequality. A Hermitian form that is bounded on $H^1_0(0,1)$ is represented on $\mathbb H_\nu$ by a compact self-adjoint operator, since the embedding $H^2_0\subset H^1_0$ is compact. Since the pencil is overdamped on $(0,\hat s]$ (it is self-adjoint and decreasing in $s$ there), min--max arguments are available \cite{Duffin1955,Markus1988}. We use the Courant--Fischer principle in the following form. Let $\mathcal K$ be compact and self-adjoint on a Hilbert space, and let $\kappa_n$ be the supremum, over all $n$-dimensional subspaces $V$, of $\min_{x\in V\setminus\{0\}}(\mathcal Kx,x)/(x,x)$. If $\kappa_n>0$, then $\mathcal K$ has at least $n$ positive eigenvalues, counted with multiplicity, and $\kappa_n$ is the $n$-th largest of them.

\begin{proposition}\label{prop:O-modes}
Let $\rs$ satisfy \eqref{eq:stable}, $f\in\R$ and $\nu>0$.
\begin{enumerate}[label=(\alph*)]
\item The clamped problem
\begin{equation}\label{eq:O-theta}
\ip{\wt K\phi,\xi}=\theta\,\mu\ip{A_\nu\phi,\xi}\quad\text{for all }\xi\in H^2_0(0,1),\qquad \phi\in H^2_0(0,1),
\end{equation}
has a sequence of eigenvalues $\theta_1\ge\theta_2\ge\dots>0$, counted with multiplicity, with $\theta_n\to0$,
\begin{equation*}
\theta_n=\sup_{\dim V=n}\ \min_{\phi\in V\setminus\{0\}}\frac{\ip{\wt K\phi,\phi}}{\mu\ip{A_\nu\phi,\phi}},
\end{equation*}
the supremum being taken over subspaces $V\subset H^2_0(0,1)$. In strong form, \eqref{eq:O-theta} is \eqref{eq:theta-intro}. Moreover
\begin{equation}\label{eq:O-thetabd}
\frac{g\nu^2m}{16\pi^4\mu(1+\nu^2\pi^{-2})}\,n^{-4}\le\theta_n\le\frac{g\nu^2M_1}{\pi^4\mu}\,n^{-4}.
\end{equation}
\item For $0<s\le\hat s$ let
\begin{equation*}
\mathfrak m_n(s)\eqdef\sup_{\dim V=n}\ \min_{\phi\in V\setminus\{0\}}\frac{\ip{Q_\nu(-s)\phi,\phi}}{\mu\ip{A_\nu\phi,\phi}},
\end{equation*}
and let $n_0$ be the least $n$ with $\theta_n<\frac34\hat s$. For $n\ge n_0$, $\mathfrak m_n$ has exactly one zero $s_n^*$ in $(0,\hat s]$, and
\begin{equation}\label{eq:O-sbd}
\frac{\theta_n}{1+\frac43f^2\vartheta_\nu^2}\le s_n^*\le\frac{\theta_n}{1-C_Ms_n^*}\le\frac43\theta_n .
\end{equation}
Moreover $s_{n+1}^*\le s_n^*$, $s_n^*$ is a strictly decreasing function of $|f|$, $s_n^*\ge\theta_n$ if $f=0$, and
\begin{equation}\label{eq:O-asym}
\lim_{n\to\infty}\frac{s_n^*}{\theta_n}=1 .
\end{equation}
\item For $n\ge n_0$ there is a real $\Psi_n\in H^5(0,1)\cap H^2_0(0,1)\setminus\{0\}$ with $\ip{Q_\nu(-s_n^*)\Psi_n,\xi}=0$ for all $\xi\in H^2_0$. With $s=s_n^*$, $\Phi_n\eqdef-fB_{-s}^{-1}(\rs\Psi_n')$ and $\hat\varrho_n\eqdef\nu^2\rs'\Psi_n/s$, the triple $(\hat\varrho_n,\Psi_n,\Phi_n)$ is a real normal mode with rate $-s_n^*$. Moreover $\hat\varrho_n\in H^4(0,1)$, $\hat\varrho_n=\hat\varrho_n'=0$ at $z=0,1$, and
\begin{equation}\label{eq:O-modebd}
\nrm{\hat\varrho_n}_{H^4}+\nrm{\Psi_n}_{H^5}+\nrm{\Phi_n}_{H^4}\le C(s_n^*)^{-1}\nrm{\hat\varrho_n},
\end{equation}
where $C$ depends on $\rs,\mu,g,f,\nu$ but not on $n$. Conversely, $\nrm{\hat\varrho_n}_{H^4}\ge c(s_n^*)^{-1}\nrm{\hat\varrho_n}$ for $n\ge n_0$, with $c>0$ independent of $n$.
\item For $n\ge n_0$ and $k\in(l^{-1}\Z)^2$ with $|k|=\nu$, the real fields
\begin{equation}\label{eq:O-mode}
\begin{aligned}
&\varrho_n=e^{-s_n^*t}\cos(k\cdot x_h)\,\hat\varrho_n(z),\qquad w_n=e^{-s_n^*t}\cos(k\cdot x_h)\,\nu^2\Psi_n(z),\\
&v_n=-e^{-s_n^*t}\sin(k\cdot x_h)\big(\Psi_n'(z)\,k+\Phi_n(z)\,k^\perp\big),\qquad q_n=e^{-s_n^*t}\cos(k\cdot x_h)\,\hat q_n(z),\\
&\hat q_n\eqdef s_n^*\rs\Psi_n'+f\rs\Phi_n+\mu(\Psi_n'''-\nu^2\Psi_n'),
\end{aligned}
\end{equation}
form a real solution of \eqref{eq:RT-lin} with $v_n(t)\in H^4\cap V_\sigma$, $\varrho_n(t)\in H^4$, $\varrho_n=\pz\varrho_n=0$ on $\pa\Omega$ and $\bar\varrho_n=0$. Moreover $\nrm{\varrho_n(t)}_{L^2(\Omega)}=e^{-s_n^*t}\nrm{\varrho_n(0)}_{L^2(\Omega)}$, and
\begin{equation}\label{eq:O-datbd}
c_{\mathrm d}\,(s_n^*)^{-1}\nrm{\varrho_n(0)}_{L^2(\Omega)}\le\nrm{\varrho_n(0)}_{H^4(\Omega)}+\nrm{v_n(0)}_{H^4(\Omega)}\le C_{\mathrm d}\,(s_n^*)^{-1}\nrm{\varrho_n(0)}_{L^2(\Omega)},
\end{equation}
where $c_{\mathrm d}>0$ and $C_{\mathrm d}$ depend on $\rs,\mu,g,f,l$ and $\nu$ but not on $n$.
\end{enumerate}
\end{proposition}

\begin{proof}
(a) The form $\ip{\wt K\phi,\xi}=g\nu^2\int_0^1\sgm^{-1}\phi\bar\xi$ is bounded on $L^2$, so it is represented on $\mathbb H_\nu$ by a compact self-adjoint operator $\mathcal K$. Since $(\mathcal K\phi,\phi)_A=g\nu^2\int\sgm^{-1}|\phi|^2>0$ for $\phi\ne0$, $\mathcal K$ is positive and injective. Its eigenvectors therefore form an orthonormal basis of the infinite-dimensional space $\mathbb H_\nu$, and its eigenvalues $\theta_n$ are positive and tend to $0$. The eigenvalue equation $\mathcal K\phi=\theta\phi$ is \eqref{eq:O-theta}, and the formula for $\theta_n$ is the Courant--Fischer principle.

For the upper bound in \eqref{eq:O-thetabd}, let $V\subset H^2_0$ have dimension $n$, and choose $\phi\in V\setminus\{0\}$ that is $L^2$-orthogonal to $\sin(j\pi z)$ for $1\le j<n$. Expand $\phi=\sum_{j\ge n}b_j\sin(j\pi z)$. Since $\phi(0)=\phi(1)=0$, two integrations by parts show that the sine coefficients of $\phi''$ are $-(j\pi)^2b_j$. Hence $\nrm{\phi''}^2=\frac12\sum_{j\ge n}(j\pi)^4|b_j|^2\ge(n\pi)^4\nrm\phi^2$, and by \eqref{eq:O-A},
\begin{equation*}
\frac{\ip{\wt K\phi,\phi}}{\mu\ip{A_\nu\phi,\phi}}\le\frac{g\nu^2M_1\nrm\phi^2}{\mu\nrm{\phi''}^2}\le\frac{g\nu^2M_1}{\mu(n\pi)^4}.
\end{equation*}
For the lower bound take $V_n\eqdef\mathrm{span}\{1-\cos(2\pi jz):1\le j\le n\}\subset H^2_0$. For $\phi=\sum_{j=1}^na_j(1-\cos(2\pi jz))$,
\begin{equation*}
\nrm{\phi''}^2=\frac12\sum_{j=1}^n(2\pi j)^4|a_j|^2\le(2\pi n)^4\,\frac12\sum_{j=1}^n|a_j|^2\le(2\pi n)^4\nrm\phi^2,
\end{equation*}
and by \eqref{eq:O-A}, $\ip{\wt K\phi,\phi}/(\mu\ip{A_\nu\phi,\phi})\ge g\nu^2m/(\mu(1+\nu^2\pi^{-2})(2\pi n)^4)$ on $V_n\setminus\{0\}$.

(b) \emph{Two-sided bounds.} Let $0<s\le\hat s$. By \eqref{eq:O-Q}, \eqref{eq:O-ebd} and \eqref{eq:O-A}, for $\phi\in H^2_0$,
\begin{equation}\label{eq:O-Qbd}
\ip{\wt K\phi,\phi}-s\Big(1+\frac43f^2\vartheta_\nu^2\Big)\mu\ip{A_\nu\phi,\phi}\le\ip{Q_\nu(-s)\phi,\phi}\le\ip{\wt K\phi,\phi}-s(1-C_Ms)\mu\ip{A_\nu\phi,\phi}.
\end{equation}
Dividing by $\mu\ip{A_\nu\phi,\phi}$ and taking the supremum of the minimum over $n$-dimensional subspaces,
\begin{equation}\label{eq:O-mbd}
\theta_n-s\Big(1+\frac43f^2\vartheta_\nu^2\Big)\le\mathfrak m_n(s)\le\theta_n-s(1-C_Ms)\le\theta_n-\frac34s .
\end{equation}

\emph{Monotonicity and continuity.} Let $0<s_1<s_2\le\hat s$. By \eqref{eq:O-Q},
\begin{equation*}
\begin{aligned}
\ip{Q_\nu(-s_2)\phi,\phi}-\ip{Q_\nu(-s_1)\phi,\phi}&=-(s_2-s_1)\mu\ip{A_\nu\phi,\phi}-f^2\big(\mathfrak e^-_\nu(s_2;\phi)-\mathfrak e^-_\nu(s_1;\phi)\big)\\
&\quad+(s_2^2-s_1^2)\ip{M\phi,\phi}.
\end{aligned}
\end{equation*}
By \eqref{eq:O-elip} and \eqref{eq:O-A}, the middle term lies between $-\frac{2s_2}{s_1}(s_2-s_1)f^2\vartheta_\nu^2\mu\ip{A_\nu\phi,\phi}$ and $0$. By \eqref{eq:O-A} and $2\hat sC_M=\frac12$, the last term lies between $0$ and $\frac12(s_2-s_1)\mu\ip{A_\nu\phi,\phi}$. Hence
\begin{equation*}
\mathfrak m_n(s_1)-(s_2-s_1)\Big(1+\frac{2s_2}{s_1}f^2\vartheta_\nu^2\Big)\le\mathfrak m_n(s_2)\le\mathfrak m_n(s_1)-\frac12(s_2-s_1),
\end{equation*}
so $\mathfrak m_n$ is continuous and strictly decreasing on $(0,\hat s]$.

\emph{The zero.} Let $n\ge n_0$. By \eqref{eq:O-mbd}, $\mathfrak m_n(\hat s)\le\theta_n-\frac34\hat s<0$, and $\mathfrak m_n(s)>0$ for $0<s<\theta_n/(1+\frac43f^2\vartheta_\nu^2)$. Hence $\mathfrak m_n$ has exactly one zero $s_n^*\in(0,\hat s)$, and it satisfies the lower bound in \eqref{eq:O-sbd}. The upper bound follows from $0=\mathfrak m_n(s_n^*)\le\theta_n-s_n^*(1-C_Ms_n^*)$, and $C_Ms_n^*\le C_M\hat s=\frac14$ gives the factor $\frac43$. For $f=0$ the lower bound reads $s_n^*\ge\theta_n$. If $V$ has dimension $n+1$ and $W\subset V$ has dimension $n$, the minimum over $V$ is at most the minimum over $W$; hence $\mathfrak m_{n+1}\le\mathfrak m_n$. Therefore $\mathfrak m_{n+1}(s_n^*)\le0$, and $s_{n+1}^*\le s_n^*$ because $\mathfrak m_{n+1}$ is decreasing. \emph{Dependence on $f$.} By \eqref{eq:O-Q}, $\ip{Q_\nu(-s)\phi,\phi}$ depends on $f$ only through $f^2$. Let $f_1^2<f_2^2$ and $s\eqdef s_n^*(f_2)$. Let $\mathcal K_s(f)$ be the compact operator on $\mathbb H_\nu$ defined in the proof of (c) below; thus $\ip{Q_\nu(-s)\phi,\phi}=(\mathcal K_s(f)\phi,\phi)_A-s(\phi,\phi)_A$. Since $\mathfrak m_n(s;f_2)=0$, the Courant--Fischer principle shows that $\mathcal K_s(f_2)$ has at least $n$ positive eigenvalues, the $n$-th of which is $s$. Let $V$ be spanned by eigenvectors for the $n$ largest ones. For $\phi\in V\setminus\{0\}$, \Cref{lem:O-ekman}(b) gives
\begin{equation*}
(\mathcal K_s(f_1)\phi,\phi)_A=(\mathcal K_s(f_2)\phi,\phi)_A+(f_2^2-f_1^2)\,\mathfrak e^-_\nu(s;\phi)>s(\phi,\phi)_A .
\end{equation*}
Since $V$ is finite-dimensional, the minimum over $V$ of the quotient is attained, so $\mathfrak m_n(s;f_1)>0$. As $\mathfrak m_n(\cdot;f_1)$ is decreasing, $s_n^*(f_1)>s=s_n^*(f_2)$.

\emph{Asymptotics.} Let $h>0$. For $\phi\in H^2_0$, $\nrm{\phi'}^2=-\ip{\phi'',\phi}\le h\nrm{\phi''}^2+(4h)^{-1}\nrm\phi^2$. We use $\nrm{\phi''}^2\le\ip{A_\nu\phi,\phi}$ and $\nrm\phi^2\le\ip{\wt K\phi,\phi}/(g\nu^2m)$. Then the first bound in \eqref{eq:O-ebd} gives
\begin{equation*}
\begin{aligned}
&f^2\mathfrak e^-_\nu(s;\phi)\le s\big(ha\,\mu\ip{A_\nu\phi,\phi}+h^{-1}b\ip{\wt K\phi,\phi}\big),\\
&a\eqdef\frac{4f^2(\max\rs)^2}{3\mu^2(\pi^2+\nu^2)},\qquad b\eqdef\frac{f^2(\max\rs)^2}{3\mu(\pi^2+\nu^2)g\nu^2m}.
\end{aligned}
\end{equation*}
Since $s^2\ip{M\phi,\phi}\ge0$, \eqref{eq:O-Q} gives $\ip{Q_\nu(-s)\phi,\phi}\ge(1-sh^{-1}b)\ip{\wt K\phi,\phi}-s(1+ha)\mu\ip{A_\nu\phi,\phi}$. Hence, if $sh^{-1}b\le1$, then $\mathfrak m_n(s)\ge(1-sh^{-1}b)\theta_n-s(1+ha)$. By \eqref{eq:O-sbd} and \eqref{eq:O-thetabd}, $s_n^*\to0$. Evaluating at $s=s_n^*$ for large $n$ gives $s_n^*\ge(1-s_n^*h^{-1}b)\theta_n/(1+ha)$, so $\liminf_n s_n^*/\theta_n\ge(1+ha)^{-1}$ for every $h>0$. Hence $\liminf_n s_n^*/\theta_n\ge1$. Conversely, $s_n^*/\theta_n\le(1-C_Ms_n^*)^{-1}\to1$ by \eqref{eq:O-sbd}, and \eqref{eq:O-asym} follows.

(c) \emph{Existence.} Let $s=s_n^*$. By \eqref{eq:O-Q} and \Cref{lem:O-ekman}(b), the form $\mathfrak k_s(\phi,\xi)\eqdef\ip{\wt K\phi,\xi}-f^2\mathfrak e^-_\nu(s;\phi,\xi)+s^2\ip{M\phi,\xi}$ is a bounded Hermitian form on $H^1_0$. Let $\mathcal K_s=\mathcal K_s(f)$ be the compact self-adjoint operator on $\mathbb H_\nu$ that represents it. Then $\ip{Q_\nu(-s)\phi,\xi}=(\mathcal K_s\phi,\xi)_A-s(\phi,\xi)_A$, and $0=\mathfrak m_n(s)$ says that the supremum over $n$-dimensional $V$ of $\min_V(\mathcal K_s\phi,\phi)_A/(\phi,\phi)_A$ equals $s>0$. By the Courant--Fischer principle, $s$ is an eigenvalue of $\mathcal K_s$. An eigenvector $\Psi\ne0$ satisfies $\ip{Q_\nu(-s)\Psi,\xi}=0$ for all $\xi\in H^2_0$. All coefficients are real and $\mathfrak e^-_\nu$ is real on real functions, so $\mathrm{Re}\,\Psi$ and $\mathrm{Im}\,\Psi$ satisfy the same identity. We let $\Psi_n$ be one of them that is nonzero.

\emph{Regularity.} We write $\Psi=\Psi_n$ and $\eta_n\eqdef B_{-s}^{-1}(\rs\Psi')$, so that $R_\nu(-s)\Psi=-(\rs\eta_n)'$ and $\Phi_n=-f\eta_n$. By \Cref{lem:O-ekman}(a), $\nrm{\eta_n}_{H^{j+2}}\le C\nrm{\Psi'}_{H^j}$ for $j=0,1,2$. Testing with $\xi=\Psi$ and using \eqref{eq:O-Qbd}, $0\le g\nu^2M_1\nrm\Psi^2-\frac34s\mu\ip{A_\nu\Psi,\Psi}$, hence
\begin{equation*}
\nrm{\Psi''}^2\le\ip{A_\nu\Psi,\Psi}\le\frac{4g\nu^2M_1}{3\mu}\,s^{-1}\nrm\Psi^2 .
\end{equation*}
In $\mathcal D'(0,1)$ the identity $\ip{Q_\nu(-s)\Psi,\xi}=0$, divided by $s$, reads
\begin{equation}\label{eq:O-ode}
\mu\Psi''''=\mu\nu^2\Psi''+s^{-1}g\nu^2\sgm^{-1}\Psi-s(\rs\Psi')'+f^2(\rs\eta_n)' .
\end{equation}
The right-hand side lies in $L^2$, with norm at most $C(\nrm{\Psi}_{H^2}+s^{-1}\nrm\Psi+\nrm{\eta_n}_{H^1})\le Cs^{-1}\nrm\Psi$, since $s\le\hat s$. Hence $\nrm\Psi_{H^4}\le C(\nrm\Psi+\nrm{\Psi''''})\le Cs^{-1}\nrm\Psi$. By interpolation, $\nrm\Psi_{H^j}\le C\nrm\Psi^{1-j/4}\nrm\Psi_{H^4}^{j/4}\le Cs^{-j/4}\nrm\Psi$ for $0\le j\le4$. Differentiating \eqref{eq:O-ode}, which is legitimate since $\sgm^{-1}\in C^1$ and $\rs\in C^2$,
\begin{equation*}
\nrm{\Psi'''''}\le C\big(\nrm\Psi_{H^3}+s^{-1}\nrm\Psi_{H^1}+\nrm{\eta_n}_{H^2}\big)\le Cs^{-\frac54}\nrm\Psi .
\end{equation*}
Moreover $\nrm{\Phi_n}_{H^4}\le C\nrm{\Psi'}_{H^2}\le Cs^{-3/4}\nrm\Psi$. Finally, $|\rs'|\ge m$ gives $\nrm{\hat\varrho_n}\ge\nu^2ms^{-1}\nrm\Psi$, and $\rs'\in C^4$ gives $\nrm{\hat\varrho_n}_{H^4}\le Cs^{-1}\nrm\Psi_{H^4}\le Cs^{-2}\nrm\Psi$. Since $s\le\hat s$, the powers $s^{-5/4}$ and $s^{-3/4}$ are bounded by $Cs^{-2}$. Hence the left-hand side of \eqref{eq:O-modebd} is at most $Cs^{-2}\nrm\Psi\le Cs^{-1}\nrm{\hat\varrho_n}$. For the converse, \eqref{eq:O-ode} and the bounds above give $\mu\nrm{\Psi''''}\ge(g\nu^2ms^{-1}-Cs^{-1/2})\nrm\Psi$. By the Leibniz rule, $|\rs'|\ge m$ and $\nrm\Psi_{H^3}\le Cs^{-3/4}\nrm\Psi$,
\begin{equation*}
\nrm{\hat\varrho_n''''}\ge\nu^2s^{-1}\big(m\nrm{\Psi''''}-C\nrm\Psi_{H^3}\big)\ge\big(c\,s^{-2}-Cs^{-7/4}\big)\nrm\Psi .
\end{equation*}
Since $\nrm{\hat\varrho_n}\le\nu^2M_1s^{-1}\nrm\Psi$, there is $s_1\in(0,\hat s]$, independent of $n$, such that $\nrm{\hat\varrho_n}_{H^4}\ge cs^{-1}\nrm{\hat\varrho_n}$ whenever $s=s_n^*\le s_1$. If $s_n^*>s_1$, then trivially $\nrm{\hat\varrho_n}_{H^4}\ge\nrm{\hat\varrho_n}\ge s_1(s_n^*)^{-1}\nrm{\hat\varrho_n}$. The triple is a normal mode by \Cref{lem:O-ekman}(c), and it is real. Since $\Psi=\Psi'=0$ at the lids, $\hat\varrho_n=\nu^2\rs'\Psi/s$ and $\hat\varrho_n'=\nu^2(\rs''\Psi+\rs'\Psi')/s$ vanish there.

(d) The argument is the one in the proof of \Cref{thm:RTlinear}(i), with $\lambda=-s_n^*$. The fields \eqref{eq:O-mode} are the real parts of $e^{\lambda t}e^{ik\cdot x_h}(\hat\varrho_n,ik\Psi_n'+ik^\perp\Phi_n,\nu^2\Psi_n,\hat q_n)$. The density equation holds because $\lambda\hat\varrho_n=-\nu^2\rs'\Psi_n$. The components of the momentum equation along $k$ and $k^\perp$ are the definition of $\hat q_n$ and the identity $B_\lambda\Phi_n=-f\rs\Psi_n'$. The relation $\pz\hat q_n=-g\hat\varrho_n$ is equivalent to $Q_\nu(\lambda)\Psi_n=0$, by the computation in the proof of \Cref{lem:RT-modes}. The velocity lies in $V_\sigma$ because $\Psi_n\in H^2_0$ and $\Phi_n\in H^1_0$. The horizontal mean of $\varrho_n$ vanishes because $k\ne0$. Since $\int_{\Tl}\cos^2(k\cdot x_h)\dxh=\int_{\Tl}\sin^2(k\cdot x_h)\dxh=\frac12|\Tl|$, we have $\nrm{\varrho_n(0)}_{L^2(\Omega)}=(\frac12|\Tl|)^{1/2}\nrm{\hat\varrho_n}$. Each horizontal derivative of the fields contributes a factor $\nu$. Hence $\nrm{\varrho_n(0)}_{H^4(\Omega)}+\nrm{v_n(0)}_{H^4(\Omega)}\le C(\nrm{\hat\varrho_n}_{H^4}+\nrm{\Psi_n}_{H^5}+\nrm{\Phi_n}_{H^4})$, and \eqref{eq:O-modebd} gives the upper bound in \eqref{eq:O-datbd}. The lower bound follows from $\nrm{\varrho_n(0)}_{H^4(\Omega)}\ge\nrm{\pz^4\varrho_n(0)}_{L^2(\Omega)}=(\frac12|\Tl|)^{1/2}\nrm{\hat\varrho_n''''}$ and the converse bound in (c).
\end{proof}

For the optimality in intermediate Sobolev norms we need the size of the slow modes in $H^r$, $0\le r\le4$. The next lemma shows that each vertical derivative of $\hat\varrho_n$ costs exactly a factor $(s_n^*)^{-1/4}$, which is of order $n$. The lower bound follows from the two-sided bound in \Cref{prop:O-modes}(c) and an upper bound in $H^5$, by interpolation.

\begin{lemma}\label{lem:O-Hr}
Let $n\ge n_0$ and $s=s_n^*$. Then $\hat\varrho_n\in H^5(0,1)$ and $\nrm{\hat\varrho_n}_{H^5}\le Cs^{-5/4}\nrm{\hat\varrho_n}$. For every $0\le r\le4$ there are $c_r,C_r>0$, independent of $n$, such that
\begin{equation}\label{eq:O-Hr}
c_rs^{-\frac r4}\nrm{\hat\varrho_n}_{L^2(0,1)}\le\nrm{\hat\varrho_n}_{H^r(0,1)}\le C_rs^{-\frac r4}\nrm{\hat\varrho_n}_{L^2(0,1)} .
\end{equation}
If the mode \eqref{eq:O-mode} is normalized by $\nrm{\varrho_n(0)}_{H^4(\Omega)}+\nrm{v_n(0)}_{H^4(\Omega)}=1$, then, after changing $c_r,C_r$,
\begin{equation}\label{eq:O-Hrphys}
c_re^{-st}s^{1-\frac r4}\le\nrm{\varrho_n(t)}_{H^r(\Omega)}\le C_re^{-st}s^{1-\frac r4},\qquad t\ge0,\ 0\le r\le4 .
\end{equation}
The constants depend only on $\rs,\mu,g,f,l,k$ and $r$.
\end{lemma}

\begin{proof}
Write $\Psi=\Psi_n$. The proof of \Cref{prop:O-modes}(c) gives $\nrm\Psi_{H^5}\le Cs^{-5/4}\nrm\Psi$. Since $\rs'\in C^5$ by \eqref{eq:stable} and $\hat\varrho_n=\nu^2\rs'\Psi/s$, we get $\nrm{\hat\varrho_n}_{H^5}\le Cs^{-1}\nrm\Psi_{H^5}\le Cs^{-9/4}\nrm\Psi$, and $\nrm{\hat\varrho_n}\ge\nu^2ms^{-1}\nrm\Psi$ gives the $H^5$ bound. The upper bound in \eqref{eq:O-Hr} follows by interpolation between $L^2$ and $H^4$ from $\nrm{\hat\varrho_n}_{H^4}\le Cs^{-1}\nrm{\hat\varrho_n}$. For $0\le r<4$, interpolation between $H^r$ and $H^5$ gives $\nrm{\hat\varrho_n}_{H^4}\le C\nrm{\hat\varrho_n}_{H^r}^{1/(5-r)}\nrm{\hat\varrho_n}_{H^5}^{(4-r)/(5-r)}$. Inserting the lower bound $\nrm{\hat\varrho_n}_{H^4}\ge cs^{-1}\nrm{\hat\varrho_n}$ of \Cref{prop:O-modes}(c) and the $H^5$ bound, and raising to the power $5-r$,
\begin{equation*}
\nrm{\hat\varrho_n}_{H^r}\ge c\,s^{-(5-r)+\frac54(4-r)}\nrm{\hat\varrho_n}=c\,s^{-\frac r4}\nrm{\hat\varrho_n} .
\end{equation*}
The case $r=4$ is \Cref{prop:O-modes}(c). For fixed $k$, the map $a\mapsto\cos(k\cdot x_h)a(z)$ and the map that extracts the coefficient of $\cos(k\cdot x_h)$ are bounded between $H^m(0,1)$ and $H^m(\Omega)$ for every integer $m\ge0$, and by interpolation for every real $m\ge0$. Hence $\nrm{\cos(k\cdot x_h)a}_{H^r(\Omega)}$ and $\nrm a_{H^r(0,1)}$ are equivalent, with constants independent of $n$. Finally, the normalization and the two-sided bound \eqref{eq:O-datbd} give $cs\le\nrm{\varrho_n(0)}_{L^2(\Omega)}\le Cs$, and $\varrho_n(t)=e^{-st}\varrho_n(0)$, which proves \eqref{eq:O-Hrphys}.
\end{proof}

\subsection{Optimality for the linearized problem}

The linearized problem \eqref{eq:RT-lin} has at most one solution with given data in the class $\varrho\in C^1([0,T];L^2)$, $v\in C^1([0,T];L^2)\cap C([0,T];H^2\cap V_\sigma)$, $\nh q\in C([0,T];L^2)$, which contains the fields of \Cref{prop:O-modes}(d). Indeed, test the difference of two solutions with $(g\sgm\varrho,v)$. The Coriolis term drops out, $\int\nh q\cdot v=g\int\varrho w$ as in \Cref{rem:RT}(i), and the density equation gives $g\int\rs'\sgm w\varrho=-g\int w\varrho$, so the coupling terms cancel. This gives $\frac12\ddt(g\nrm{\sqrt\sgm\varrho}^2+\nrm{\sqrt\rs v}^2)+\mu\nrm{\nabla v}^2=0$.

\begin{theorem}\label{thm:O-lin}
Let $\rs$ satisfy \eqref{eq:stable}, $f\in\R$ and $k\in(l^{-1}\Z)^2\setminus\{0\}$. There is $c>0$, depending only on $\rs,\mu,g,f,l$ and $k$, with the following property. For every $t\ge0$ there are real data $(\vin,\uin)$ with $\vin\in H^4$, $\vin=\pz\vin=0$ on $\pa\Omega$, $\uin\in H^4\cap V_\sigma$ and $\nrm{\vin}_{H^4}+\nrm{\uin}_{H^4}=1$, such that the solution of \eqref{eq:RT-lin} is a single normal mode \eqref{eq:O-mode} and satisfies
\begin{equation*}
\nrm{\wt\varrho(t)}_{L^2}=\nrm{\varrho(t)}_{L^2}\ge c(1+t)^{-1},\qquad \nrm{\nh\varrho(t)}_{L^2}\ge c(1+t)^{-1}.
\end{equation*}
\end{theorem}

\begin{proof}
Let $\nu=|k|$. By \eqref{eq:O-sbd} and \eqref{eq:O-thetabd}, $c_1n^{-4}\le s_n^*\le c_2n^{-4}$ for $n\ge n_0$, with
\begin{equation*}
c_1\eqdef\frac{g\nu^2m}{16\pi^4\mu(1+\nu^2\pi^{-2})(1+\frac43f^2\vartheta_\nu^2)},\qquad c_2\eqdef\frac{4g\nu^2M_1}{3\pi^4\mu}.
\end{equation*}
Given $t\ge0$, let $n\eqdef\max\{n_0,\lceil(c_2t)^{1/4}\rceil\}$. Then $s_n^*t\le c_2n^{-4}t\le1$. Moreover $n\le n_0+1+(c_2t)^{1/4}$, so $n^4\le C(1+t)$ and $s_n^*\ge c_1n^{-4}\ge c_3(1+t)^{-1}$, with $c_3>0$ depending only on $n_0,c_1,c_2$. Let $\mathsf D_n\eqdef\nrm{\varrho_n(0)}_{H^4}+\nrm{v_n(0)}_{H^4}$ and take the data $(\vin,\uin)\eqdef\mathsf D_n^{-1}(\varrho_n(0),v_n(0))$. These are admissible by \Cref{prop:O-modes}(d). The solution is $\mathsf D_n^{-1}(\varrho_n,v_n)$, and by \eqref{eq:O-datbd},
\begin{equation*}
\nrm{\varrho(t)}_{L^2}=e^{-s_n^*t}\,\frac{\nrm{\varrho_n(0)}_{L^2}}{\mathsf D_n}\ge e^{-1}C_{\mathrm d}^{-1}s_n^*\ge e^{-1}C_{\mathrm d}^{-1}c_3(1+t)^{-1}.
\end{equation*}
Since $\bar\varrho=0$ we have $\wt\varrho=\varrho$. Moreover $\nrm{\nh\varrho(t)}=\nu\nrm{\varrho(t)}$, because $\varrho$ is $\cos(k\cdot x_h)$ times a function of $(t,z)$.
\end{proof}

\subsection{Optimality for the nonlinear problem}\label{sec:opt-nl}

From now on $\nrm\cdot$ is the norm of $L^2(\Omega)$. The nonlinear data are $\varepsilon$ times a slow mode, plus a velocity correction of order $\varepsilon^2$ that makes them compatible.

\begin{lemma}\label{lem:O-compat}
Let $n\ge n_0$, and let $(\varrho_n,v_n,q_n)$ be the mode \eqref{eq:O-mode}, normalized by $\nrm{\varrho_n(0)}_{H^4}+\nrm{v_n(0)}_{H^4}=1$. There are $\varepsilon_n'>0$ and $C_n$ such that, for $0<\varepsilon\le\varepsilon_n'$, there is $v_{n,\varepsilon}\in H^4\cap V_\sigma$ with $\nrm{v_{n,\varepsilon}}_{H^4}\le C_n\varepsilon^2$ and the following property. The data
\begin{equation*}
\vin\eqdef\varepsilon\varrho_n(0),\qquad \uin\eqdef\varepsilon v_n(0)+v_{n,\varepsilon}
\end{equation*}
satisfy $\rs+\vin>0$, $\vin=\pz\vin=0$ on $\pa\Omega$, $\uin\in H^4\cap H^1_0$, \eqref{eq:lid0} and \eqref{eq:compat}. Moreover the solution of \Cref{prop:local} with these data satisfies $v_t(0)=-\varepsilon s_n^*v_n(0)$.
\end{lemma}

\begin{proof}
We apply the construction in the proof of \Cref{lem:compatcorr} to the approximate solution
\begin{equation*}
(\varrho_a,v_a,P_a)\eqdef\varepsilon(\varrho_n,v_n,P_n),\qquad P_n\eqdef q_n+g\int_0^z\varrho_n\dd z'=e^{-s_n^*t}\cos(k\cdot x_h)\,\hat q_n(0),
\end{equation*}
which is independent of $z$ because $\pz\hat q_n=-g\hat\varrho_n$. That proof uses no instability, and we use it only at the Sobolev level $H^2$, that is, in the case $k=2$ of that proof, where $k$ is a regularity index. Since $(\varrho_n,v_n,q_n)$ solves \eqref{eq:RT-lin}, the residuals of \Cref{lem:approx} are
\begin{equation*}
r_\varrho=\varepsilon^2\bU_n\cdot\nabla\varrho_n,\qquad r_v=\varepsilon^2\big(\varrho_n(\pt v_n+fv_n^\perp)+\rs\bU_n\cdot\nabla v_n\big)+\varepsilon^3\varrho_n\bU_n\cdot\nabla v_n,
\end{equation*}
with $\bU_n=(v_n,w_n)$. By \Cref{prop:O-modes}(c),(d), $\varrho_n(0)\in H^4\subset C^2(\ol\Omega)$ and $v_n(0)\in H^4$, so $\nrm{r_v(0)}_{H^2}\le C_n\varepsilon^2$. At this level, the fixed-point argument for \eqref{eq:vcorr} uses only the following facts:
\begin{itemize}
\item $\nrm{\varrho_a(0)}_{C^2}+\nrm{v_a(0)}_{H^3}\le C\varepsilon$, which bounds $\bU_a(0)\cdot\nabla v$ and $\mathbf V\cdot\nabla v_a(0)$ in $H^2$ by $C\varepsilon\nrm v_{H^3}$, since $H^2(\Omega)$ is an algebra;
\item the bound on $\nrm{r_v(0)}_{H^2}$;
\item \Cref{lem:RT-stokesf} for $\zeta\in H^2$, whose proof then requires only $\rs\in C^2$.
\end{itemize}
It gives $v_{n,\varepsilon}\in H^4\cap V_\sigma$ with $\nrm{v_{n,\varepsilon}}_{H^4}\le C_n\varepsilon^2$. The last paragraph of that proof shows that the initial acceleration given by \eqref{PE}$_2$ is $\pt v_a(0)=-\varepsilon s_n^*v_n(0)$, which vanishes on $\pa\Omega$; this is \eqref{eq:compat}. The other conditions hold because $\varrho_n(0)$ is clamped, $v_n(0)\in V_\sigma$, and $\varepsilon$ is small.
\end{proof}

\begin{theorem}\label{thm:O-nonlin}
Let $\rs$ satisfy \eqref{eq:stable}, let $f\in\R$, and let $\varepsilon_0$ satisfy \eqref{eq:slope}. There is $c>0$, depending only on $\rs,\mu,g,f$ and $l$, with the following properties. For a solution of \Cref{thm:main} with data of size $\mathsf N=\nrm{\vin}_{H^4}+\nrm{\uin}_{H^4}\le\varepsilon_0$, let $\rinf$, $\mathcal H$, $\bX$ and $\bX_\infty$ be given by \Cref{thm:limit,thm:energy}, and consider the inequalities
\begin{equation}\label{eq:O-lower}
\begin{aligned}
&\nrm{\rho(t)-\rinf}_{L^2}\ge\nrm{\wt\varrho(t)}_{L^2}\ge c(1+t)^{-1}\mathsf N,\qquad \nrm{\nh\varrho(t)}_{L^2}\ge c(1+t)^{-1}\mathsf N,\\
&\nrm{q(t)-q_\infty}_{H^1}\ge c(1+t)^{-1}\mathsf N,\qquad \mathcal H(t)\ge c(1+t)^{-2}\mathsf N^2,\qquad \nrm{X^z(t)-X^z_\infty}_{L^2}\ge c(1+t)^{-1}\mathsf N ,
\end{aligned}
\end{equation}
where $q_\infty$ is the pressure of the limit state, as in \eqref{eq:ratepress}.
\begin{enumerate}[label=(\roman*)]
\item For every $t\ge0$ and every $\varepsilon_*\in(0,\varepsilon_0]$ there are data satisfying the hypotheses of \Cref{thm:main} with $0<\mathsf N\le\varepsilon_*$, whose solution satisfies \eqref{eq:O-lower} at time $t$.
\item There are data $(\vin^{(n)},\uin^{(n)})$ satisfying the hypotheses of \Cref{thm:main}, with sizes $\mathsf N_n\in(0,\varepsilon_0]$ and $\mathsf N_n\to0$, and times $t_n\to\infty$, such that the $n$-th solution satisfies \eqref{eq:O-lower} at $t=t_n$ with $\mathsf N=\mathsf N_n$.
\end{enumerate}
Consequently, let $\varepsilon_*\in(0,\varepsilon_0]$ and $\varphi:[0,\infty)\to(0,\infty)$ be such that every solution of \Cref{thm:main} with $\mathsf N\le\varepsilon_*$ satisfies $\nrm{\wt\varrho(t)}_{L^2}\le\varphi(t)\mathsf N$ for all $t\ge0$. Then $\varphi(t)\ge c(1+t)^{-1}$ for all $t\ge0$. In particular no bound $C(1+t)^{-1-\delta}\mathsf N$ with $\delta>0$ holds. The same holds for $\nrm{\nh\varrho(t)}_{L^2}$, $\nrm{\rho(t)-\rinf}_{L^2}$, $\nrm{q(t)-q_\infty}_{H^1}$ and $\nrm{X^z(t)-X^z_\infty}_{L^2}$, and for $\mathcal H(t)$ with $\varphi(t)\mathsf N^2$ and $c(1+t)^{-2}$.
\end{theorem}

\begin{proof}
Fix $k\eqdef(l^{-1},0)$, so that $\nu=l^{-1}$, and let $C_{\mathrm d}$ be the constant in \eqref{eq:O-datbd}. For $n\ge n_0$ let $(\varrho_n,v_n,q_n)$ be the mode \eqref{eq:O-mode}, normalized as in \Cref{lem:O-compat}. Then $\nrm{\varrho_n(0)}\ge C_{\mathrm d}^{-1}s_n^*$.

\emph{Data.} For $0<\varepsilon\le\varepsilon_n'$ let $(\vin,\uin)$ be the data of \Cref{lem:O-compat}, and let $\mathsf N^\varepsilon$ be their size. Then $|\mathsf N^\varepsilon-\varepsilon|\le C_n\varepsilon^2$. Hence there is $\varepsilon_n''\in(0,\varepsilon_n']$ such that, for $\varepsilon\le\varepsilon_n''$, we have $\frac\varepsilon2\le\mathsf N^\varepsilon\le2\varepsilon$, and the data satisfy the hypotheses of \Cref{thm:main} with $\mathsf N^\varepsilon\le\varepsilon_0$. Let $(\varrho,v,q)$ be the global solution. By \Cref{thm:limit,thm:energy} it has a limit $\rinf$, an available energy $\mathcal H$ and a flow $\bX$, and \eqref{eq:L-orth}, \eqref{eq:E-coerc} and \eqref{eq:E-disp} hold.

\emph{Difference estimate.} Let $\check\varrho\eqdef\varrho-\varepsilon\varrho_n$, $\check v\eqdef v-\varepsilon v_n$, $\check w\eqdef w-\varepsilon w_n$ and $\check q\eqdef q-\varepsilon q_n$. Since $\varepsilon(\varrho_n,v_n,w_n,q_n)$ solves \eqref{eq:RT-lin},
\begin{equation*}
\begin{aligned}
&\pt\check\varrho+\rs'\check w=F_\varrho\eqdef-\bU\cdot\nabla\varrho,\\
&\rs\pt\check v+f\rs\,\check v^\perp+\nh\check q-\mu\Lap\check v=F_v\eqdef-\varrho\big(v_t+fv^\perp\big)-\rho\,\bU\cdot\nabla v,
\end{aligned}
\end{equation*}
with $\pz\check q=-g\check\varrho$, $\check v(t)\in V_\sigma$ and $\check w=-\int_0^z\divh\check v\dd z'$. We test the first equation with $g\sgm\check\varrho$ and the second with $\check v$. Since $\rs'\sgm=-1$, the coupling term in the first is $-g\int\check w\,\check\varrho$. In the second, $\int f\rs\,\check v^\perp\cdot\check v=0$ and, as $\check w=0$ on $\pa\Omega$, $\int\nh\check q\cdot\check v=-\int\check q\divh\check v=-\int\pz\check q\,\check w=g\int\check\varrho\,\check w$. Adding, with $Y\eqdef\frac g2\nrm{\sqrt\sgm\check\varrho}^2+\frac12\nrm{\sqrt\rs\check v}^2$,
\begin{equation*}
\ddt Y+\mu\nrm{\nabla\check v}^2=g\int\sgm F_\varrho\,\check\varrho+\int F_v\cdot\check v\le C\big(\nrm{F_\varrho}+\nrm{F_v}\big)Y^{\frac12}.
\end{equation*}
The computation is justified as in the proof of \Cref{thm:RTnonlinear}, since $\varrho\in C^1([0,T];H^3)$, $v\in C^1([0,T];H^2)$, and the mode is smooth in $t$ with values in $H^4$. By \eqref{eq:unif}, \eqref{eq:decay}, $H^2(\Omega)\subset L^\infty$ and $\nrm w_{L^\infty}\lesssim\nrm v_{H^3}$ (\Cref{lem:w}), we have $\nrm{F_\varrho(t)}+\nrm{F_v(t)}\le C(\mathsf N^\varepsilon)^2\le4C\varepsilon^2$ for all $t$, with $C$ independent of $n$. Since $\check\varrho(0)=0$ and $\check v(0)=v_{n,\varepsilon}$, we have $Y(0)^{1/2}\le C\nrm{v_{n,\varepsilon}}\le CC_n\varepsilon^2$. Integrating $\ddt Y^{1/2}\le C(\nrm{F_\varrho}+\nrm{F_v})$ as in the proof of \Cref{thm:RTnonlinear},
\begin{equation*}
\nrm{\check\varrho(t)}\le C_n'\varepsilon^2(1+t),\qquad t\ge0 .
\end{equation*}
Moreover $\nrm{\check\varrho(t)}_{H^4}\le\nrm{\varrho(t)}_{H^4}+\varepsilon\nrm{\varrho_n(t)}_{H^4}\le C\varepsilon$ by \eqref{eq:unif}. Interpolation (\Cref{lem:GN}(1)) then gives $\nrm{\nh\check\varrho(t)}\le C\nrm{\check\varrho(t)}^{3/4}\nrm{\check\varrho(t)}_{H^4}^{1/4}\le C_n''\varepsilon^{7/4}(1+t)^{3/4}$.

\emph{Lower bounds.} Let $t\ge0$ and $n\ge n_0$ with $s_n^*t\le1$. Since $\bar\varrho_n=0$, we have $\wt\varrho-\varepsilon\varrho_n=\wt{\check\varrho}$, and $\nrm{\varrho_n(t)}=e^{-s_n^*t}\nrm{\varrho_n(0)}\ge e^{-1}C_{\mathrm d}^{-1}s_n^*$. Using $\nrm{\nh\varrho_n}=\nu\nrm{\varrho_n}$, we get
\begin{equation*}
\nrm{\wt\varrho(t)}\ge\varepsilon e^{-1}C_{\mathrm d}^{-1}s_n^*-C_n'\varepsilon^2(1+t),\qquad \nrm{\nh\varrho(t)}\ge\varepsilon\nu e^{-1}C_{\mathrm d}^{-1}s_n^*-C_n''\varepsilon^{7/4}(1+t)^{3/4}.
\end{equation*}
Choose $\varepsilon\in(0,\varepsilon_n'']$ so small that $C_n'\varepsilon(1+t)\le\frac12e^{-1}C_{\mathrm d}^{-1}s_n^*$ and $C_n''\varepsilon^{3/4}(1+t)^{3/4}\le\frac12\nu e^{-1}C_{\mathrm d}^{-1}s_n^*$. Since $\mathsf N^\varepsilon\le2\varepsilon$,
\begin{equation}\label{eq:O-lowergen}
\nrm{\wt\varrho(t)}\ge\frac{s_n^*}{4eC_{\mathrm d}}\,\mathsf N^\varepsilon,\qquad \nrm{\nh\varrho(t)}\ge\frac{\nu s_n^*}{4eC_{\mathrm d}}\,\mathsf N^\varepsilon .
\end{equation}
By \eqref{eq:L-orth}, $\nrm{\rho(t)-\rinf}\ge\nrm{\wt\varrho(t)}$. The lower bound in \eqref{eq:E-coerc} gives $\mathcal H(t)\ge c\nrm{\rho(t)-\rinf}^2$. The second inequality in \eqref{eq:E-disp} with $p=2$ gives $\nrm{X^z(t)-X^z_\infty}\ge\frac2{3M_1}\nrm{\rho(t)-\rinf}$. Finally, the hydrostatic relations for the solution and for the limit state give $\pz(q-q_\infty)=-g(\varrho-\vinf)=-g(\rho-\rinf)$, so $\nrm{q(t)-q_\infty}_{H^1}\ge g\nrm{\rho(t)-\rinf}$.

\emph{Proof of (i).} Given $t\ge0$, choose $n=n(t)$ as in the proof of \Cref{thm:O-lin}, so that $s_n^*t\le1$ and $s_n^*\ge c_3(1+t)^{-1}$. Choose $\varepsilon$ as above, and in addition $\varepsilon\le\varepsilon_*/2$. Then \eqref{eq:O-lowergen} and the four transfers give \eqref{eq:O-lower} at time $t$.

\emph{Proof of (ii).} Let $t_n\eqdef1/s_n^*$; then $t_n\to\infty$ by \Cref{prop:O-modes}(b), and $s_n^*=t_n^{-1}\ge(1+t_n)^{-1}$. Choose $\varepsilon=\varepsilon_n$ as above with $t=t_n$, and in addition $\varepsilon_n\le1/n$. Then $0<\mathsf N_n\le2/n\to0$, and \eqref{eq:O-lowergen} gives \eqref{eq:O-lower} at $t=t_n$.

In both cases $c$ depends only on $c_3$, $C_{\mathrm d}$, $\nu=l^{-1}$, $M_1$, $g$ and the constant in \eqref{eq:E-coerc}.

\emph{Consequence.} Let $t\ge0$. By (i) with this $t$ and $\varepsilon_*$, there is a solution with $0<\mathsf N\le\varepsilon_*$ and $c(1+t)^{-1}\mathsf N\le\nrm{\wt\varrho(t)}\le\varphi(t)\mathsf N$. Hence $\varphi(t)\ge c(1+t)^{-1}$. The other quantities are treated in the same way.
\end{proof}

The lower bounds of \Cref{thm:O-nonlin} concern $L^2$ norms. Since the oscillatory density decays like $(1+t)^{-1}$ in $L^2$ and is only bounded in $H^4$, interpolation gives the rate $(1+t)^{-(1-r/4)}$ in $H^r$. The next theorem shows that this rate is also optimal, for every $0\le r<4$. It combines \Cref{lem:O-Hr} with the nonlinear error estimate in the proof of \Cref{thm:O-nonlin}.

\begin{theorem}\label{thm:O-Hr}
Let $\rs$ satisfy \eqref{eq:stable}, let $f\in\R$, let $0\le r<4$ and $\alpha_r\eqdef1-\frac r4$. There is $c_r>0$, depending only on $\rs,\mu,g,f,l$ and $r$, with the following properties.
\begin{enumerate}[label=(\roman*)]
\item For every $t\ge0$ there are real data as in \Cref{thm:O-lin}, of unit size in $H^4\times H^4$, such that the solution of \eqref{eq:RT-lin} satisfies $\nrm{\varrho(t)}_{H^r(\Omega)}\ge c_r(1+t)^{-\alpha_r}$.
\item Let $\varepsilon_0$ satisfy \eqref{eq:slope}. For every $t\ge0$ and every $\varepsilon_*\in(0,\varepsilon_0]$ there are data satisfying the hypotheses of \Cref{thm:main} with $0<\mathsf N\le\varepsilon_*$, whose solution satisfies \eqref{eq:O-lower} and
\begin{equation}\label{eq:O-lowerHr}
\nrm{\wt\varrho(t)}_{H^r(\Omega)}\ge c_r(1+t)^{-\alpha_r}\mathsf N,\qquad \nrm{\rho(t)-\rinf}_{H^r(\Omega)}\ge c_r(1+t)^{-\alpha_r}\mathsf N .
\end{equation}
There are also such data of sizes $\mathsf N_n\to0$ and times $t_n\to\infty$ for which these bounds hold at $t=t_n$.
\end{enumerate}
Since every solution of \Cref{thm:main} satisfies $\nrm{\wt\varrho(t)}_{H^r}+\nrm{\rho(t)-\rinf}_{H^r}\le C_r\mathsf N(1+t)^{-\alpha_r}$, the exponent $\alpha_r$ is optimal, uniformly for small $H^4$ data, both for the oscillatory density and for the distance to the limit.
\end{theorem}

\begin{proof}
\emph{Upper bound.} By \Cref{lem:GN}(1), \Cref{lem:meanosc}, \eqref{eq:decay} and \eqref{eq:unif}, $\nrm{\wt\varrho(t)}_{H^r}\le C\nrm{\wt\varrho(t)}^{\alpha_r}\nrm{\wt\varrho(t)}_{H^4}^{1-\alpha_r}\le C\mathsf N(1+t)^{-\alpha_r}$; the bound for $\rho-\rinf$ is \eqref{eq:rateLinf}.

\emph{Proof of (i).} Let $k=(l^{-1},0)$, and given $t\ge0$ choose $n=n(t)$ as in the proof of \Cref{thm:O-lin}, so that $s_n^*t\le1$ and $s_n^*\ge c_3(1+t)^{-1}$. By \eqref{eq:O-Hrphys}, the normalized mode satisfies $\nrm{\varrho_n(t)}_{H^r}\ge c_re^{-1}(s_n^*)^{\alpha_r}\ge c_re^{-1}c_3^{\alpha_r}(1+t)^{-\alpha_r}$.

\emph{Proof of (ii).} Let $n$ be as in (i), or let $t=t_n=1/s_n^*$, and let $\varepsilon$ and the solution be as in the proof of \Cref{thm:O-nonlin}. There the error $\check\varrho=\varrho-\varepsilon\varrho_n$ was shown to satisfy $\nrm{\check\varrho(\tau)}\le C_n'\varepsilon^2(1+\tau)$ and $\nrm{\check\varrho(\tau)}_{H^4}\le C\varepsilon$ for all $\tau\ge0$. By interpolation,
\begin{equation*}
\nrm{\check\varrho(\tau)}_{H^r}\le C\nrm{\check\varrho(\tau)}^{\alpha_r}\nrm{\check\varrho(\tau)}_{H^4}^{1-\alpha_r}\le C_{n,r}\varepsilon^{1+\alpha_r}(1+\tau)^{\alpha_r}.
\end{equation*}
Since $\bar\varrho_n=0$, we have $\wt\varrho-\varepsilon\varrho_n=\wt{\check\varrho}$, and the horizontal projection is bounded on $H^r$ (\Cref{lem:meanosc}). By \eqref{eq:O-Hrphys}, at the selected time,
\begin{equation*}
\nrm{\wt\varrho(t)}_{H^r}\ge c_r\varepsilon e^{-1}(s_n^*)^{\alpha_r}-2C_{n,r}\varepsilon^{1+\alpha_r}(1+t)^{\alpha_r}.
\end{equation*}
Since $\alpha_r>0$, we may decrease $\varepsilon$ further so that $2C_{n,r}\varepsilon^{\alpha_r}(1+t)^{\alpha_r}\le\frac12c_re^{-1}(s_n^*)^{\alpha_r}$; this is one more smallness condition, compatible with those in the proof of \Cref{thm:O-nonlin}, so \eqref{eq:O-lower} still holds. With $\mathsf N\le2\varepsilon$ and $s_n^*\ge c_3(1+t)^{-1}$ (respectively $s_n^*=t_n^{-1}\ge(1+t_n)^{-1}$), we obtain the first bound in \eqref{eq:O-lowerHr}. Since $\vinf$ depends only on $z$, $\wt\varrho$ is the oscillatory part of $\rho-\rinf$, so $\nrm{\wt\varrho}_{H^r}\le2\nrm{\rho-\rinf}_{H^r}$ by \Cref{lem:meanosc}, which gives the second bound. The sequence is obtained as in the proof of \Cref{thm:O-nonlin}(ii), with $\varepsilon_n\le1/n$.
\end{proof}

At the endpoint $r=4$ the error estimate above no longer gains over the main term, since $\alpha_4=0$. The following corollary shows instead, by interpolation, that there is no uniform decay rate at all in $H^4$.

\begin{corollary}\label{cor:O-H4}
Let $\rs$ satisfy \eqref{eq:stable}, $f\in\R$, and let $\varepsilon_0$ satisfy \eqref{eq:slope}. There is $c>0$ such that, for every $t\ge0$ and every $\varepsilon_*\in(0,\varepsilon_0]$, some data satisfying the hypotheses of \Cref{thm:main} with $0<\mathsf N\le\varepsilon_*$ give
\begin{equation*}
\nrm{\wt\varrho(t)}_{H^4}\ge c\mathsf N,\qquad \nrm{\rho(t)-\rinf}_{H^4}\ge c\mathsf N .
\end{equation*}
Consequently, neither $\nrm{\wt\varrho(t)}_{H^4}$ nor $\nrm{\rho(t)-\rinf}_{H^4}$ admits a bound $\omega(t)\mathsf N$ with $\omega(t)\to0$ as $t\to\infty$, uniformly over solutions with $\mathsf N\le\varepsilon_*$.
\end{corollary}

\begin{proof}
Take the data of \Cref{thm:O-Hr}(ii) with $r=1$ at the given time. Then $\nrm{\wt\varrho(t)}_{H^1}\ge c_1\mathsf N(1+t)^{-3/4}$, while $\nrm{\wt\varrho(t)}\le C\mathsf N(1+t)^{-1}$ by \eqref{eq:decay}. By \Cref{lem:GN}(1),
\begin{equation*}
c_1\mathsf N(1+t)^{-\frac34}\le C\nrm{\wt\varrho(t)}^{\frac34}\nrm{\wt\varrho(t)}_{H^4}^{\frac14}\le C\mathsf N^{\frac34}(1+t)^{-\frac34}\nrm{\wt\varrho(t)}_{H^4}^{\frac14},
\end{equation*}
which gives the first bound. The second follows from $\nrm{\wt\varrho}_{H^4}\le2\nrm{\rho-\rinf}_{H^4}$. A bound $\omega(t)\mathsf N$ with $\omega(t)\to0$ would force $\omega(t)\ge c$ for all $t$, which is impossible.
\end{proof}

\begin{remark}\label{rem:O}
(i) \emph{The mechanism.} By \eqref{eq:O-ode}, the slow modes balance buoyancy against viscosity: to leading order $g\nu^2\sgm^{-1}\Psi_n\approx s_n^*\mu A_\nu\Psi_n$. Inertia, the term $s^2M$, and the Ekman term $f^2R_\nu(-s)$ are corrections. If $\rs$ is linear, then $\sgm$ is constant and $\theta_n=g\nu^2/(\mu\sgm\Lambda_n)$, where $\Lambda_n$ are the eigenvalues of $\phi''''-\nu^2\phi''$ on $H^2_0(0,1)$; classically, $\Lambda_n\sim((n+\frac12)\pi)^4$. Heuristically, a vertical wave of wavenumber $\beta$ behaves like an overdamped oscillator. It has a fast viscous rate of order $\mu\beta^2/\rs$ and the slow buoyancy rate of order $g\nu^2/(\mu\sgm\beta^4)$ constructed here. The slow modes carry little velocity. By the bounds in the proof of \Cref{prop:O-modes}(c), $\nrm{v_n(t)}_{H^2}\le C(s_n^*)^{1/4}\nrm{\varrho_n(t)}$. By the two-sided bound \eqref{eq:O-datbd}, the normalized modes therefore satisfy $\nrm{v_n(t)}_{H^2}\le C(s_n^*)^{5/4}$ at $t=1/s_n^*$, which is of order $(1+t)^{-5/4}$. These modes do not decide whether the velocity rates in \eqref{eq:decay} are sharp, and we do not address this question.

(ii) \emph{Rotation.} By \Cref{prop:O-modes}(b), rotation slows every slow mode: $|f|\mapsto s_n^*$ is strictly decreasing. The effect disappears for high vertical modes, by \eqref{eq:O-asym}. The Ekman term is bounded by $\frac43sf^2\vartheta_\nu^2\mu\ip{A_\nu\phi,\phi}$ in \eqref{eq:O-ebd}. On high vertical modes it is in fact of lower order, because it involves only $\nrm{\phi'}$. This parallels \Cref{rem:RT}(ii), where the short horizontal waves do not feel the Coriolis force. A Galerkin approximation of \eqref{eq:RT-nm} is consistent with \Cref{prop:O-modes}. It was run for $\rs=2-z$ and $\rs=2-z-0.3z^2$, with $\mu=g=1$ and $\nu\in\{1,3\}$. The computed slow rates are real and satisfy \eqref{eq:O-sbd} up to discretization error. For $f=0$, $|s_n^*/\theta_n-1|\le10^{-3}$ for $n\le10$. For $f=20$, $\rs=2-z$ and $\nu=1$, the ratios $s_n^*/\theta_n$ are $0.645,0.817,0.949,0.978,0.990,\dots,0.999$ for $n=1,\dots,10$. These computations are not used in the proofs.

(iii) \emph{Rates that are not addressed.} The modes constructed here have zero horizontal mean, so they do not directly decide the sharpness of the exponents $\frac54$ and $\frac{13}8$ for the mean in \Cref{thm:limit}. The pointwise exponents $\frac56$ and $\frac12$ in \eqref{eq:rateLinf} and the flow-map exponents $\frac13$ and $\frac12$ in \Cref{thm:energy}(iii) are upper bounds; we do not claim that they are sharp.
\end{remark}

\begin{remark}[The two-dimensional problem]\label{rem:2D}
Let $\T_l\eqdef\R/2\pi l\Z$, and consider the non-rotating primitive equations with gravity in the vertical slab $\T_l\times(0,1)$, with velocity $(u,w)$ and $w=-\int_0^z\px u\dd z'$. Its solutions are exactly the solutions of \eqref{PE} with $f=0$ that do not depend on $y$ and have $v_2=0$: for data with these properties, the solution of \eqref{PE} keeps them, by uniqueness. The $H^m$ norms on $\Omega$ and on the slab differ only by the factor $(2\pi l)^{1/2}$. Hence all results of this paper hold for the two-dimensional problem. This includes the optimality results of this section: for $f=0$ and $k=(l^{-1},0)$ we have $\Phi_n=0$, so the modes \eqref{eq:O-mode} do not depend on $y$ and have $v_2=0$. In particular, the exponents $\frac{13}8$ for the mean, $\frac56$ and $\frac12$ for the pointwise rates, and $\frac13$ and $\frac12$ for the flow map hold in two dimensions as well; they follow from the three-dimensional results via this embedding. We have not tried to optimize these exponents in two dimensions. One further statement is specific to two dimensions, and we add an observation for $f=0$ that we do not prove.
\begin{enumerate}[label=(\alph*)]
\item \emph{The pressure in $L^\infty$.} Since $H^1(\T_l)\subset L^\infty(\T_l)$, the surface pressure, which has zero mean, satisfies $\nrm{P(t)}_{L^\infty}\lesssim\nrm{\px P(t)}_{L^2(\T_l)}\lesssim\mathsf N(1+t)^{-1}$, and $|\int_0^z(\varrho-\vinf)\dd s|\le\sup_x\nrm{(\varrho-\vinf)(x,\cdot)}_{L^2(0,1)}\lesssim\nrm{\varrho-\vinf}+\nrm{\px\varrho}$. Hence the pressure converges in $H^1\cap L^\infty$: $\nrm{q(t)-q_\infty}_{H^1\cap L^\infty}\lesssim\mathsf N(1+t)^{-1}$.
\item \emph{Smooth data (not proved here).} For $f=0$, in two or three dimensions, every mode $k\ne0$ of \eqref{eq:RT-lin} decouples: $\Phi$ obeys $\rs\pt\Phi+\mu L_\nu\Phi=0$ and decays exponentially, while $F\eqdef\sgm\varrho_k$ solves the damped wave equation $MF_{tt}+\mu A_\nu F_t+\wt KF=0$. A boundary-layer analysis of this equation, which we do not carry out here, suggests the following for the two-dimensional linearized problem. For data with $\vin\in H^6$ and $u_{\mathrm{in}}\in H^4$, the density satisfies $\nrm{\px\varrho(t)}\lesssim(\nrm{\vin}_{H^6}+\nrm{u_{\mathrm{in}}}_{H^4})(1+t)^{-9/8}$, and this rate is not better unless $\Gamma\eqdef\AB(\sgm\vin)+\mu^{-1}\px\pz(\rs u_{\mathrm{in}})$ vanishes on the lids. The obstruction would be a boundary layer of width $t^{-1/4}$ at the lids. In the rotating case the per-mode problem contains the memory term $f^2R_\nu(\lambda)$, and we have not investigated this question.
\end{enumerate}
\end{remark}

\appendix
\section{Local well-posedness: proof of \texorpdfstring{\Cref{prop:local}}{the local well-posedness proposition}}\label{app:local}

\renewcommand{\theequation}{A.\arabic{equation}}
\setcounter{equation}{0}

Throughout this appendix, $C$ denotes a constant depending only on $l,\mu,g,f$ and $\rs$, and $C(\cdot)$ a constant depending in addition on the quantities indicated, non-decreasing in each of them, except that the lower bound $\rho_*$ of the density enters only through $\rho_*^{-1}$ (that is, $C(\cdot)$ is non-decreasing in $\rho_*^{-1}$). Of the background profile only $\rs\in C^5([0,1])$ is used; the stability condition \eqref{eq:stable} plays no role. We fix the data as in \Cref{prop:local} and set
\begin{equation}\label{A:data}
M_0\eqdef\nrm{\vin}_{H^4}+\nrm{\uin}_{H^4},\qquad \rho_*\eqdef\min_{\ol\Omega}\rin>0,\qquad \rho^*\eqdef\max_{\ol\Omega}\rin ,
\end{equation}
where $\rin=\rs+\vin\in H^4(\Omega)\subset C(\ol\Omega)$, so that $\rho_*>0$ by the hypothesis $\rin>0$ on $\ol\Omega$, and $\rho^*\le\nrm{\rs}_{L^\infty}+CM_0$. We prove \Cref{prop:local} with an existence time $T_0$ depending only on $M_0$ and $\rho_*$.

\subsection{Functional setting and initial time derivatives}

Let
\begin{equation*}
V_\sigma\eqdef\Big\{\phi\in H^1_0(\Omega;\R^2):\ \divh\int_0^1\phi\dd z=0\Big\},
\end{equation*}
let $H_\sigma$ be its closure in $L^2(\Omega;\R^2)$, and let $\mathcal Q$ be the $L^2$-orthogonal projection onto $H_\sigma$. For $\Pi\in H^1(\Tl)$ and $\phi\in V_\sigma$ we have $(\nh\Pi,\phi)=-\int_{\Tl}\Pi\,\divh\int_0^1\phi\dd z\dxh=0$, hence $(\nh\Pi,\phi)=0$ for all $\phi\in H_\sigma$; conversely, for $\psi\in L^2(\Omega;\R^2)$ one has $\psi-\mathcal Q\psi=\nh\Pi_\psi$ for some $\Pi_\psi=\Pi_\psi(x_h)\in H^1(\Tl)$ \cite{TemamZiane2004,Ziane1995}. The hydrostatic Stokes operator $A_h$ is the positive self-adjoint operator on $H_\sigma$ associated with the coercive form $\mu(\nabla\phi,\nabla\psi)$ on $V_\sigma$. If $\phi\in\cD(A_h)$, then $\phi$ is the unique weak solution in $V_\sigma$ (Lax--Milgram) of the problem of \Cref{lem:stokes} with $\zeta=A_h\phi$; since the solution provided by \Cref{lem:stokes} is also a weak solution, \Cref{lem:stokes} gives $\cD(A_h)=H^2\cap V_\sigma$, $A_h\phi=-\mu\Lap\phi+\nh p_\phi$ for some $p_\phi=p_\phi(x_h)$, and
\begin{equation}\label{A:Ah}
\nrm{\phi}_{H^{k+2}}\le C\nrm{A_h\phi}_{H^k}\qquad(k\ge0,\ A_h\phi\in H^k).
\end{equation}
The inverse $A_h^{-1}$ is compact on $H_\sigma$, so $A_h$ has an orthonormal basis of eigenfunctions $\{e_k\}_{k\ge1}\subset H_\sigma$, $A_he_k=\lambda_ke_k$, which is also complete in $V_\sigma$; by \eqref{A:Ah} and induction, $e_k\in H^m(\Omega)$ for every $m$, so $e_k\in C^\infty(\ol\Omega)$. We denote by $\mathcal P_m$ the orthogonal projection of $H_\sigma$ onto $\mathrm{span}\{e_1,\dots,e_m\}$; it commutes with $A_h$.

\emph{Vertical velocity and elementary bounds.} For a horizontal field $v\in H^1(\Omega;\R^2)$ we write
\begin{equation*}
w_v\eqdef-\int_0^z\divh v(\cdot,z')\dd z',\qquad \bU_v\eqdef(v,w_v).
\end{equation*}
Then $\divv\bU_v=0$ and $w_v|_{z=0}=0$; if $\divh\int_0^1v\dd z=0$, also $w_v|_{z=1}=0$, so that $\bU_v\cdot n=0$ on $\pa\Omega$, $n$ the outward unit normal. By \Cref{lem:w} (with $p=4$, and with $s=\frac12$ for the $L^\infty$ bounds) and \Cref{lem:GN}(3),
\begin{equation}\label{A:wbd}
\begin{gathered}
\nrm{w_v}_{H^k}\le C\nrm{v}_{H^{k+1}}\ (k\ge0),\qquad \nrm{w_v}_{L^2}\le\nrm{w_v}_{L^2_hL^\infty_z}\le\nrm{\divh v}_{L^2},\qquad \nrm{w_v}_{L^\infty_zL^4_h}\le C\nrm{v}_{H^2},\\
\nrm{\bU_v}_{L^\infty}\le C\nrm{v}_{H^3},\qquad \nrm{\nabla\bU_v}_{L^\infty}\le C\nrm{v}_{H^4};
\end{gathered}
\end{equation}
for the last bound, $\nrm{\nabla v}_{L^\infty}+\nrm{\pz w_v}_{L^\infty}\le C\nrm{v}_{H^3}$ and $\nrm{\nh w_v}_{L^\infty}\le C\nrm{\nh\divh v}_{L^2_zH^{3/2}_h}\le C\nrm{v}_{H^4}$. We shall also use the three-dimensional product bounds
\begin{equation}\label{A:prod}
\begin{gathered}
\nrm{\phi\psi}_{H^k}\le C\nrm{\phi}_{H^k}\nrm{\psi}_{H^k}\ (k\ge2),\qquad \nrm{\phi\psi}_{H^1}\le C\min\big\{\nrm{\phi}_{H^2},\nrm{\phi}_{W^{1,\infty}}\big\}\nrm{\psi}_{H^1},\\
\int_\Omega|\theta|\,|\pz\phi|\,|\psi|\dx\le\nrm{\theta}_{L^\infty_zL^4_h}\nrm{\pz\phi}_{L^2_zL^4_h}\nrm{\psi}_{L^2}\le C\nrm{\theta}_{L^\infty_zL^4_h}\nrm{\phi}_{H^2}\nrm{\psi}_{L^2},
\end{gathered}
\end{equation}
which follow from $H^2(\Omega)\subset L^\infty$, $H^1(\Omega)\subset L^6$ and \Cref{lem:aniso}(f) with $q=4$. The functions $v_t,v_{tt}$ obey \eqref{A:wbd} with $w_{v_t}=\pt w_v$, $w_{v_{tt}}=\pt^2w_v$.

\emph{Initial time derivatives.} They are determined by the data through \eqref{PE}:
\begin{equation}\label{A:u1}
\begin{gathered}
\varrho_1\eqdef-\bU_{\mathrm{in}}\cdot\nabla\vin-w_{\mathrm{in}}\rs',\qquad u_1\eqdef\rin^{-1}\big(F_0-\nh P_0\big),\\
F_0\eqdef\mu\Lap\uin-\rin\bU_{\mathrm{in}}\cdot\nabla\uin-f\rin\uin^\perp+g\int_0^z\nh\vin\dd z',
\end{gathered}
\end{equation}
where $w_{\mathrm{in}}\eqdef w_{\uin}$, $\bU_{\mathrm{in}}\eqdef\bU_{\uin}$, and $P_0=P_0(x_h)$ is the zero-mean solution of
\begin{equation*}
\divh\big(a\nh P_0\big)=\divh b,\qquad a\eqdef\int_0^1\rin^{-1}\dd z,\qquad b\eqdef\int_0^1\rin^{-1}F_0\dd z ,
\end{equation*}
that is, $P_0$ is fixed by the lid constraint $\divh\int_0^1u_1\dd z=0$. Since $\rin\in H^4\subset W^{2,\infty}$ and $\rin\ge\rho_*$, we have $\nrm{\rin^{-1}}_{W^{2,\infty}}\le C(M_0,\rho_*)$, hence $\nrm{a}_{W^{2,\infty}(\Tl)}\le C(M_0,\rho_*)$ and $(\rho^*)^{-1}\le a\le\rho_*^{-1}$. By \eqref{A:wbd} and \eqref{A:prod}, $\nrm{\bU_{\mathrm{in}}\cdot\nabla\uin}_{H^2}\le C\nrm{\uin}_{H^3}^2$, so $\nrm{F_0}_{H^2}\le C(M_0)$ and $\nrm{b}_{H^2(\Tl)}\le C(M_0,\rho_*)$. Testing the equation for $P_0$ with $P_0$ gives $\nrm{\nh P_0}_{L^2(\Tl)}\le\rho^*\nrm{b}_{L^2(\Tl)}$. Writing it as $\Laph P_0=a^{-1}(\divh b-\nh a\cdot\nh P_0)$ and using $\nrm{\nh^2\phi}_{L^2(\Tl)}=\nrm{\Laph\phi}_{L^2(\Tl)}$ twice, we get $\nh P_0\in H^1(\Tl)$, then $\nh a\cdot\nh P_0\in H^1(\Tl)$ and $\nh P_0\in H^2(\Tl)$, with $\nrm{\nh P_0}_{H^2(\Tl)}\le C(M_0,\rho_*)$. Since $H^3(\Omega)$ is an algebra, $\nrm{\bU_{\mathrm{in}}\cdot\nabla\vin}_{H^3}+\nrm{w_{\mathrm{in}}\rs'}_{H^3}\le C(1+M_0)M_0$. Therefore
\begin{equation}\label{A:u1bd}
\nrm{\varrho_1}_{H^3}+\nrm{u_1}_{H^2}\le C(M_0,\rho_*).
\end{equation}
Moreover $u_1\in H^2\cap V_\sigma$. The lid constraint holds by construction. By \eqref{eq:lid0} and $\uin|_{\pa\Omega}=0$, $\int_0^1\bU_{\mathrm{in}}\cdot\nabla\uin\dd z=\divh\int_0^1\uin\otimes\uin\dd z$ (as in the proof of \Cref{lem:propagation}), so the equation for $P_0$ is \eqref{eq:Peq} at $t=0$, and $P_0$ is the initial surface pressure fixed by \eqref{eq:wt0}; thus $\qin=P_0-g\int_0^z\vin\dd z'$. On $\pa\Omega$ we have $\uin=0$, $w_{\mathrm{in}}=0$ (by \eqref{eq:lid0}) and $\nh\uin=0$, hence $\Lap\uin=\pz^2\uin$, $\bU_{\mathrm{in}}\cdot\nabla\uin=0$ and $\rin u_1=\mu\pz^2\uin-\nh\qin$ there. So $u_1|_{\pa\Omega}=0$ is exactly the compatibility condition \eqref{eq:compat}.

\emph{Norms.} For $T\in(0,1]$ and $v$ with values in $V_\sigma$ we use
\begin{equation}\label{A:Ynorm}
\begin{aligned}
\nrm{v}_{Y_T'}^2&\eqdef\sup_{[0,T]}\big(\nrm{v}_{H^4}^2+\nrm{v_t}_{H^2}^2+\nrm{v_{tt}}_{L^2}^2\big)+\int_0^T\big(\nrm{v_t}_{H^3}^2+\nrm{v_{tt}}_{H^1}^2\big)\dt,\\
\nrm{v}_{Y_T}^2&\eqdef\nrm{v}_{Y_T'}^2+\int_0^T\nrm{v}_{H^5}^2\dt,\qquad \nrm{v}_{W_T}^2\eqdef\sup_{[0,T]}\nrm{v}_{L^2}^2+\int_0^T\nrm{\nabla v}_{L^2}^2\dt ,
\end{aligned}
\end{equation}
and the set
\begin{equation}\label{A:ball}
\mathcal B_{R,T}\eqdef\big\{v:\ \nrm{v}_{Y_T}\le R,\ v(t)\in V_\sigma\ (t\in[0,T]),\ v(0)=\uin,\ v_t(0)=u_1\big\}.
\end{equation}
If $\nrm{v}_{Y_T}<\infty$, then $v$ is Lipschitz from $[0,T]$ into $H^2$ and bounded in $H^4$, so by interpolation $v\in C^{0,\frac14}([0,T];H^{7/2})$; with \eqref{A:wbd} (where $H^4$ may be replaced by $H^{7/2}$, since $\nrm{\nh w_v}_{L^\infty}\le C\nrm{\nh\divh v}_{L^2_zH^{5/4}_h}$) this gives $\bU_v\in C([0,T];W^{1,\infty}(\Omega))$.

\subsection{Transport}

We first record a commutator estimate. For a multi-index $\alpha\in\N^3$ (derivatives in $x,y,z$) with $|\alpha|\le4$, a vector field $U\in H^4(\Omega;\R^3)\cap W^{1,\infty}$ and $\varrho\in H^4(\Omega)$,
\begin{equation}\label{A:comm}
\nrm{[\pa^\alpha,U\cdot\nabla]\varrho}_{L^2}\le C\big(\nrm{\nabla U}_{L^\infty}+\nrm{U}_{H^4}\big)\nrm{\varrho}_{H^4}.
\end{equation}
Indeed, $[\pa^\alpha,U\cdot\nabla]\varrho$ is a linear combination of $\pa^\beta U\cdot\nabla\pa^\gamma\varrho$ with $|\beta|\ge1$, $|\beta|+|\gamma|=|\alpha|$. For $|\beta|=1$ we use $\nrm{\nabla U}_{L^\infty}\nrm{\nabla\pa^\gamma\varrho}_{L^2}$; for $|\beta|=2$, $\nrm{\pa^\beta U}_{L^4}\nrm{\nabla\pa^\gamma\varrho}_{L^4}\le C\nrm{U}_{H^3}\nrm{\varrho}_{H^4}$; for $|\beta|=3$, $\nrm{\pa^\beta U}_{L^2}\nrm{\nabla\pa^\gamma\varrho}_{L^\infty}\le C\nrm{U}_{H^3}\nrm{\varrho}_{H^4}$ (here $|\gamma|\le1$ and $H^2\subset L^\infty$); for $|\beta|=4$, $\nrm{\pa^\beta U}_{L^2}\nrm{\nabla\varrho}_{L^\infty}\le C\nrm{U}_{H^4}\nrm{\varrho}_{H^3}$.

\begin{lemma}\label{lem:A1}
Let $T\in(0,1]$, $R>0$, and let $v$ satisfy $\nrm{v}_{Y_T}\le R$, $v(t)\in V_\sigma$ for $t\in[0,T]$, and $v(0)=\uin$. Then the problem $\rho_t+\bU_v\cdot\nabla\rho=0$, $\rho(0)=\rin$, has a unique solution $\rho\in L^\infty((0,T)\times\Omega)$, and $\varrho\eqdef\rho-\rs$ satisfies $\varrho\in C([0,T];H^4)$, $\rho_*\le\rho\le\rho^*$, $\varrho_t(0)=\varrho_1$, and
\begin{equation}\label{A:transport}
\begin{gathered}
\sup_{[0,T]}\nrm{\varrho}_{H^4}\le\big(M_0+C\sqrt TR\big)e^{C\sqrt TR},\\
\esssup_{(0,T)}\big(\nrm{\varrho_t}_{H^3}+\nrm{\varrho_{tt}}_{H^1}\big)\le C(1+R)^2\Big(1+\sup_{[0,T]}\nrm{\varrho}_{H^4}\Big).
\end{gathered}
\end{equation}
\end{lemma}

\begin{proof}
\emph{Flow, uniqueness and bounds.} By the remark after \eqref{A:ball}, $\bU_v\in C([0,T];W^{1,\infty}(\Omega))$; moreover $\divv\bU_v=0$ and $\bU_v\cdot n=0$ on $\pa\Omega$. Hence the flow $X(t,s,\cdot)$ of $\bU_v$ is a family of measure-preserving bi-Lipschitz homeomorphisms of $\ol\Omega$ (the lids are invariant because $w_v$ vanishes there). Therefore $\rho(t)\eqdef\rin\circ X(0,t,\cdot)$ is a solution, and $\rho_*\le\rho\le\rho^*$. Since $\bU_v\in L^1(0,T;W^{1,\infty})$ is divergence-free and tangent to $\pa\Omega$, the solution is unique in $L^\infty((0,T)\times\Omega)$ by \cite{DiPernaLions1989}. The perturbation $\varrho=\rho-\rs$ solves
\begin{equation}\label{A:transp}
\varrho_t+\bU_v\cdot\nabla\varrho=-w_v\rs',\qquad \varrho(0)=\vin .
\end{equation}

\emph{Regularization.} It remains to prove the $H^4$ regularity. Let $E$ be a bounded linear extension operator from $H^k(\Omega)$ to $H^k(\Tl\times\R)$, $0\le k\le5$, acting only in the vertical variable (a higher-order reflection across $z=0$ and $z=1$, followed by a cut-off in $z$); it commutes with horizontal derivatives. Let $J_\delta$ be the convolution with a smooth mollifier of width $\delta$ on $\Tl\times\R$, and let $R_\delta f\eqdef(J_\delta Ef)|_\Omega$. Then $R_\delta f\in C^\infty(\ol\Omega)$ with $\nrm{R_\delta f}_{H^m}\le C_{m,\delta}\nrm{f}_{L^2}$ for every $m$, $\nrm{R_\delta f}_{H^k}\le C\nrm{f}_{H^k}$ and $R_\delta f\to f$ in $H^k$ as $\delta\to0$, for $f\in H^k$, $0\le k\le5$; moreover $R_\delta$ commutes with $\nh$. The mollification destroys the lid constraint, which we restore by a vertically constant correction: with the Fourier multiplier $Q_h\eqdef\nh\Laph^{-1}\divh$ on $\Tl$ (the $L^2(\Tl)$-orthogonal projection onto gradients, which has norm at most one on every $H^k(\Tl)$), set
\begin{equation*}
S_\delta v\eqdef R_\delta v-Q_h\int_0^1R_\delta v\dd z .
\end{equation*}
Then $\divh\int_0^1S_\delta v\dd z=0$ because $\divh Q_h=\divh$, and $S_\delta v\in C^\infty(\ol\Omega)$. By the Cauchy--Schwarz inequality in $z$, $\nrm{\int_0^1\phi\dd z}_{H^k(\Tl)}\le\nrm{\phi}_{H^k(\Omega)}$, hence $\nrm{S_\delta v}_{H^k}\le C\nrm{v}_{H^k}$ for $k\le5$. If $\divh\int_0^1v\dd z=0$, then $Q_h\int_0^1v\dd z=0$, so $S_\delta v-v=(R_\delta v-v)-Q_h\int_0^1(R_\delta v-v)\dd z$ and $S_\delta v\to v$ in $H^k$ for $v\in H^k$, $k\le5$.

Let $v^\delta(t)\eqdef S_\delta v(t)$, $w^\delta\eqdef w_{v^\delta}$, $U^\delta\eqdef\bU_{v^\delta}$ and $\vin^\delta\eqdef R_\delta\vin$. Then $\divv U^\delta=0$ and $U^\delta\cdot n=0$ on $\pa\Omega$. Since $v\in C([0,T];L^2)$, $U^\delta\in C([0,T];C^6(\ol\Omega))$. By \eqref{A:wbd} and the bounds on $S_\delta$, for a.e.\ $t$,
\begin{equation}\label{A:Udelta}
\nrm{\nabla U^\delta}_{L^\infty}+\nrm{U^\delta}_{H^4}+\nrm{w^\delta}_{H^4}\le C\nrm{v^\delta}_{H^5}\le C\nrm{v}_{H^5},
\end{equation}
and, by dominated convergence, $U^\delta\to\bU_v$ in $L^2(0,T;H^4)$, in particular in $L^2(0,T;L^\infty)$. Moreover $\vin^\delta\to\vin$ in $H^4$. The flow $X^\delta$ of $U^\delta$ consists of $C^6$-diffeomorphisms of $\ol\Omega$, continuously differentiable in time, and $\rho^\delta(t)\eqdef(\rs+\vin^\delta)\circ X^\delta(0,t,\cdot)$ solves $\rho^\delta_t+U^\delta\cdot\nabla\rho^\delta=0$. Since $\rs\in C^5$, $\varrho^\delta\eqdef\rho^\delta-\rs\in C([0,T];C^5(\ol\Omega))$, and $\varrho^\delta_t=-U^\delta\cdot\nabla\varrho^\delta-w^\delta\rs'\in C([0,T];C^4(\ol\Omega))$.

\emph{Uniform $H^4$ bound.} For $|\alpha|\le4$, $f=\pa^\alpha\varrho^\delta$ is continuously differentiable in $t$ and satisfies pointwise
\begin{equation*}
\pt\pa^\alpha\varrho^\delta+U^\delta\cdot\nabla\pa^\alpha\varrho^\delta=-[\pa^\alpha,U^\delta\cdot\nabla]\varrho^\delta-\pa^\alpha\big(w^\delta\rs'\big).
\end{equation*}
We multiply by $\pa^\alpha\varrho^\delta$ and integrate; the transport term vanishes because $\divv U^\delta=0$ and $U^\delta\cdot n=0$. Using \eqref{A:comm}, the bound $\nrm{\pa^\alpha(w^\delta\rs')}\le C\nrm{\rs}_{C^5}\nrm{w^\delta}_{H^4}$ and \eqref{A:Udelta}, and summing over $|\alpha|\le4$, we obtain
\begin{equation*}
\frac12\ddt\nrm{\varrho^\delta}_{H^4}^2\le C\nrm{v}_{H^5}\big(\nrm{\varrho^\delta}_{H^4}^2+\nrm{\varrho^\delta}_{H^4}\big).
\end{equation*}
Hence $y_\delta\eqdef\nrm{\varrho^\delta}_{H^4}$ satisfies $y_\delta(t)\le y_\delta(0)+\int_0^ta(1+y_\delta)$ with $a\eqdef C\nrm{v}_{H^5}$ (apply the inequality to $(y_\delta^2+\epsilon)^{1/2}$ and let $\epsilon\to0$), and Gronwall's inequality with $\int_0^Ta\le C\sqrt TR$ yields
\begin{equation*}
\sup_{[0,T]}\nrm{\varrho^\delta}_{H^4}\le\big(\nrm{\vin^\delta}_{H^4}+C\sqrt TR\big)e^{C\sqrt TR}.
\end{equation*}
Since $H^3(\Omega)$ is an algebra, \eqref{A:wbd} gives $\nrm{\varrho^\delta_t}_{H^3}\le C\nrm{U^\delta}_{H^3}\nrm{\varrho^\delta}_{H^4}+C\nrm{w^\delta}_{H^3}\le C\nrm{v}_{H^4}(1+\nrm{\varrho^\delta}_{H^4})$, so $\{\varrho^\delta\}$ is bounded in $L^\infty(0,T;H^4)\cap W^{1,\infty}(0,T;H^3)$.

\emph{Passage to the limit.} By the Aubin--Simon lemma, a subsequence converges to some $\tilde\varrho$ in $C([0,T];H^3)$, and $\varrho^\delta(t)\rightharpoonup\tilde\varrho(t)$ weakly in $H^4$ for every $t$. By weak lower semicontinuity and $\nrm{\vin^\delta}_{H^4}\to\nrm{\vin}_{H^4}\le M_0$, $\tilde\varrho$ satisfies the first bound in \eqref{A:transport}. Since $U^\delta\to\bU_v$ in $L^2(0,T;L^\infty)$, $U^\delta$ is bounded in $L^\infty((0,T)\times\Omega)$ and $\nabla\varrho^\delta\to\nabla\tilde\varrho$ in $C([0,T];L^2)$, we have $U^\delta\cdot\nabla\varrho^\delta+w^\delta\rs'\to\bU_v\cdot\nabla\tilde\varrho+w_v\rs'$ in $L^2(0,T;L^2)$. Passing to the limit in $\varrho^\delta(t)=\vin^\delta-\int_0^t(U^\delta\cdot\nabla\varrho^\delta+w^\delta\rs')$, we see that $\tilde\varrho$ solves \eqref{A:transp}; thus $\rs+\tilde\varrho$ solves the transport problem for $\rho$, and $\tilde\varrho=\varrho$ by uniqueness. Consequently $\varrho\in L^\infty(0,T;H^4)\cap W^{1,\infty}(0,T;H^3)$ satisfies the first bound in \eqref{A:transport}.

\emph{Time continuity.} For $|\alpha|\le4$, $\pa^\alpha\varrho\in L^\infty(0,T;L^2)$ solves, in the sense of distributions,
\begin{equation*}
\pt\pa^\alpha\varrho+\bU_v\cdot\nabla\pa^\alpha\varrho=G_\alpha\eqdef-[\pa^\alpha,\bU_v\cdot\nabla]\varrho-\pa^\alpha(w_v\rs'),
\end{equation*}
and, by \eqref{A:comm} and \eqref{A:wbd}, $\nrm{G_\alpha}_{L^2}\le C\nrm{v}_{H^5}(1+\nrm{\varrho}_{H^4})$, so $G_\alpha\in L^2(0,T;L^2)$. Since $\bU_v\in L^1(0,T;W^{1,\infty})$ is divergence-free and tangent to $\pa\Omega$, the renormalization theory of \cite{DiPernaLions1989} (see \Cref{rem:justify}(a)) gives $\pa^\alpha\varrho\in C([0,T];L^2)$ and $\ddt\nrm{\pa^\alpha\varrho}_{L^2}^2=2\int_\Omega\pa^\alpha\varrho\,G_\alpha\dx$. Hence $\varrho\in C([0,T];H^4)$.

\emph{Time derivatives.} We have $\varrho_t=-\bU_v\cdot\nabla\varrho-w_v\rs'$. Since $v\in C([0,T];H^3)$, \eqref{A:wbd} and \eqref{A:prod} give $\varrho_t\in C([0,T];H^2)$, and at $t=0$, $\varrho_t(0)=-\bU_{\mathrm{in}}\cdot\nabla\vin-w_{\mathrm{in}}\rs'=\varrho_1$. As above, $\nrm{\varrho_t}_{H^3}\le C\nrm{v}_{H^4}(1+\nrm{\varrho}_{H^4})\le CR(1+\nrm{\varrho}_{H^4})$ for a.e.\ $t$ (at this stage $\varrho_t\in L^\infty(0,T;H^3)$ only). Next, $\varrho_{tt}=-\bU_{v_t}\cdot\nabla\varrho-\bU_v\cdot\nabla\varrho_t-w_{v_t}\rs'$. By \eqref{A:prod} and \eqref{A:wbd}, $\nrm{v_t\cdot\nh\varrho}_{H^1}\le C\nrm{v_t}_{H^2}\nrm{\varrho}_{H^2}$, $\nrm{w_{v_t}\pz\varrho}_{H^1}\le C\nrm{\pz\varrho}_{W^{1,\infty}}\nrm{w_{v_t}}_{H^1}\le C\nrm{\varrho}_{H^4}\nrm{v_t}_{H^2}$, $\nrm{\bU_v\cdot\nabla\varrho_t}_{H^1}\le C\nrm{\bU_v}_{W^{1,\infty}}\nrm{\varrho_t}_{H^2}\le C\nrm{v}_{H^4}\nrm{\varrho_t}_{H^3}$ and $\nrm{w_{v_t}\rs'}_{H^1}\le C\nrm{v_t}_{H^2}$. Thus $\varrho_{tt}\in L^\infty(0,T;H^1)$, and with $\nrm{v}_{H^4}+\nrm{v_t}_{H^2}\le R$ for a.e.\ $t$ this gives the second bound in \eqref{A:transport}, as an essential supremum.
\end{proof}

\begin{remark}
In three dimensions the regularization step cannot be replaced by an estimate of the flow of $\bU_v$ in $W^{3,\infty}$: for $v\in H^5(\Omega)$ one only has $\nh^3w_v=-\int_0^z\nh^3\divh v$ with $\nh^4v\in H^1(\Omega)$, and $H^1(\Tl)\not\subset L^\infty(\Tl)$, so $\nh^3w_v$ need not be bounded. The commutator estimate \eqref{A:comm} needs only $\nrm{\nabla\bU_v}_{L^\infty}\lesssim\nrm{v}_{H^4}$ and $\nrm{w_v}_{H^4}\lesssim\nrm{v}_{H^5}$.
\end{remark}

\subsection{The linear velocity problem}

For given $\rho$, $c$, $v$ and a forcing $\eta$ consider
\begin{equation}\label{A:L}
\rho\big(u_t+\bU_v\cdot\nabla u\big)+f\rho u^\perp+cu-\mu\Lap u+\nh P=\eta,\qquad u(t)\in V_\sigma,\qquad u(0)=u_0,
\end{equation}
where $P=P(t,x_h)$ is the surface pressure (the Lagrange multiplier of the lid constraint). Equivalently, $u(t)\in V_\sigma$ and
\begin{equation*}
\big(\rho(u_t+\bU_v\cdot\nabla u)+f\rho u^\perp+cu-\eta,\phi\big)+\mu(\nabla u,\nabla\phi)=0\qquad\text{for all }\phi\in V_\sigma .
\end{equation*}
The Coriolis term is of order zero. It drops out when \eqref{A:L} is tested with $u$, and, after differentiation in time, its leading part $f\rho\pt^ju^\perp$ drops out when tested with $\pt^ju$.

\begin{lemma}\label{lem:A2}
Let $T\in(0,1]$, $\rho\in W^{1,\infty}(\Omega\times(0,T))$ with $\rho_*\le\rho\le\rho^*$, $c\in L^\infty(\Omega\times(0,T))$ with $c_t\in L^1(0,T;L^2)$, $v\in L^\infty(0,T;H^3)$ with values in $V_\sigma$ and $v_t\in L^1(0,T;H^2)$, $\eta\in W^{1,1}(0,T;L^2)$, and $u_0\in H^2\cap V_\sigma$. Then \eqref{A:L} has a unique solution with
\begin{equation*}
u\in L^\infty(0,T;H^2\cap V_\sigma),\qquad u_t\in L^\infty(0,T;L^2)\cap L^2(0,T;V_\sigma),\qquad \nh P\in L^\infty(0,T;L^2).
\end{equation*}
The solution is unique even among weak solutions $u\in L^\infty(0,T;L^2)\cap L^2(0,T;V_\sigma)$ for which $t\mapsto(\rho u(t),\phi)$ is absolutely continuous for every $\phi\in V_\sigma$, the weak formulation holds with $(\rho u_t,\phi)$ read as $\ddt(\rho u,\phi)-(\rho_tu,\phi)$, and $(\rho u)(0)=\rho(0)u_0$.
\end{lemma}

\begin{proof}
Note that $\eta\in C([0,T];L^2)$, $c\in C([0,T];L^2)$ with $\nrm{c(0)}_{L^\infty}\le\nrm{c}_{L^\infty}$, and $v\in C([0,T];H^2)$, so that $\bU_v(0)$ is defined and $\nrm{\bU_v}_{L^\infty}\le C\nrm{v}_{L^\infty(0,T;H^3)}$ on $[0,T]$ by \eqref{A:wbd}.

\emph{Galerkin approximation.} Let $u^m(t)=\sum_{k\le m}c_k^m(t)e_k$ solve
\begin{equation*}
\big(\rho(u^m_t+\bU_v\cdot\nabla u^m)+f\rho(u^m)^\perp+cu^m-\eta,e_j\big)+\mu(\nabla u^m,\nabla e_j)=0,\qquad j\le m,
\end{equation*}
with $u^m(0)=\mathcal P_mu_0$. The mass matrix $(\rho e_j,e_k)$ is uniformly positive definite and Lipschitz in $t$; the coefficients $(\rho\bU_v\cdot\nabla e_k,e_j)$, $(f\rho e_k^\perp,e_j)$, $(ce_k,e_j)$ and $(\eta,e_j)$ lie in $W^{1,1}(0,T)\subset C([0,T])$, because $\pt\bU_v=\bU_{v_t}\in L^1(0,T;L^2)$ by \eqref{A:wbd}. Hence the linear system has a unique solution $c^m\in C^1([0,T])$ with $(c^m)'\in W^{1,1}(0,T)$, and the computations below are legitimate. Since $\mathcal P_m$ commutes with $A_h$, $\nrm{A_hu^m(0)}\le\nrm{A_hu_0}$ and $\mu\nrm{\nabla u^m(0)}^2=(A_h\mathcal P_mu_0,\mathcal P_mu_0)\le\nrm{A_hu_0}\nrm{u_0}$. Since $A_hu^m\in\mathrm{span}\{e_k\}_{k\le m}$ and $\mu(\nabla u^m,\nabla e_j)=(A_hu^m,e_j)$, the Galerkin equation gives $A_hu^m=\mathcal P_m\mathcal Q\big[\eta-\rho(u^m_t+\bU_v\cdot\nabla u^m)-f\rho(u^m)^\perp-cu^m\big]$, and therefore, by \eqref{A:Ah} and the Poincar\'e inequality $\nrm{u^m}\le C\nrm{\nabla u^m}$,
\begin{equation}\label{A:H2gal}
\nrm{u^m}_{H^2}\le C\big(\nrm{\eta}_{L^2}+\rho^*\nrm{u^m_t}_{L^2}+\big(\rho^*\nrm{\bU_v}_{L^\infty}+|f|\rho^*+\nrm{c}_{L^\infty}\big)\nrm{\nabla u^m}_{L^2}\big).
\end{equation}

\emph{Uniform estimates.} Testing with $u^m_t$ gives
\begin{equation*}
\rho_*\nrm{u^m_t}^2+\frac\mu2\ddt\nrm{\nabla u^m}^2\le\Big(\nrm\eta+C\big(\rho^*\nrm{\bU_v}_{L^\infty}+|f|\rho^*+\nrm{c}_{L^\infty}\big)\nrm{\nabla u^m}\Big)\nrm{u^m_t}.
\end{equation*}
Next we differentiate the Galerkin system in $t$ and test with $u^m_t$. Since $(f\rho(u^m_t)^\perp,u^m_t)=0$ and $(\rho u^m_{tt},u^m_t)=\frac12\ddt\nrm{\sqrt\rho u^m_t}^2-\frac12(\rho_tu^m_t,u^m_t)$, we obtain
\begin{equation*}
\begin{aligned}
\frac12\ddt\nrm{\sqrt\rho u^m_t}^2+\mu\nrm{\nabla u^m_t}^2&=(\eta_t,u^m_t)-\frac12(\rho_tu^m_t,u^m_t)-\big(f\rho_t(u^m)^\perp+c_tu^m+cu^m_t,u^m_t\big)\\
&\quad-\big(\rho_t\bU_v\cdot\nabla u^m+\rho\bU_{v_t}\cdot\nabla u^m+\rho\bU_v\cdot\nabla u^m_t,u^m_t\big).
\end{aligned}
\end{equation*}
We bound the terms on the right. First, $|(\rho\bU_v\cdot\nabla u^m_t,u^m_t)|\le\frac\mu4\nrm{\nabla u^m_t}^2+C(\rho^*)^2\nrm{\bU_v}_{L^\infty}^2\nrm{u^m_t}^2$, and $|(\rho_t\bU_v\cdot\nabla u^m+f\rho_t(u^m)^\perp,u^m_t)|\le C\nrm{\rho_t}_{L^\infty}(\nrm{\bU_v}_{L^\infty}+|f|)\nrm{\nabla u^m}\nrm{u^m_t}$. For the term with $\bU_{v_t}=(v_t,w_{v_t})$ we use $H^2\subset L^\infty$ for $v_t$ and \eqref{A:prod}, \eqref{A:wbd} for $w_{v_t}$:
\begin{equation*}
|(\rho\bU_{v_t}\cdot\nabla u^m,u^m_t)|\le\rho^*\big(\nrm{v_t}_{L^\infty}\nrm{\nabla u^m}+C\nrm{w_{v_t}}_{L^\infty_zL^4_h}\nrm{u^m}_{H^2}\big)\nrm{u^m_t}\le C\rho^*\nrm{v_t}_{H^2}\nrm{u^m}_{H^2}\nrm{u^m_t}.
\end{equation*}
Finally $|(c_tu^m,u^m_t)|\le\nrm{c_t}_{L^2}\nrm{u^m}_{L^\infty}\nrm{u^m_t}\le C\nrm{c_t}_{L^2}\nrm{u^m}_{H^2}\nrm{u^m_t}$, $|(cu^m_t,u^m_t)|+\frac12|(\rho_tu^m_t,u^m_t)|\le(\nrm{c}_{L^\infty}+\nrm{\rho_t}_{L^\infty})\nrm{u^m_t}^2$ and $|(\eta_t,u^m_t)|\le\frac12\nrm{\eta_t}(1+\nrm{u^m_t}^2)$. We insert \eqref{A:H2gal} for $\nrm{u^m}_{H^2}$. With $Y_m\eqdef\nrm{\sqrt\rho u^m_t}^2+\mu\nrm{\nabla u^m}^2$, the two inequalities combine into
\begin{equation*}
Y_m'+\mu\nrm{\nabla u^m_t}^2\le\beta(t)(1+Y_m),
\end{equation*}
where $\beta\in L^1(0,T)$ is independent of $m$ and depends only on $\rho_*$, $\rho^*$, $\nrm{\rho_t}_{L^\infty}$, $\nrm{c}_{L^\infty}$, $\nrm{v}_{L^\infty(H^3)}$, $\nrm{\eta}_{L^\infty(L^2)}$, and linearly on $\nrm{\eta_t}_{L^2}+\nrm{v_t}_{H^2}+\nrm{c_t}_{L^2}$. The initial value is bounded independently of $m$: taking $\phi=u^m_t(0)$ in the Galerkin system at $t=0$ and using $(A_hu^m(0),\phi)=-(\mu\Lap u^m(0),\phi)$ for $\phi\in H_\sigma$,
\begin{equation*}
\nrm{u^m_t(0)}\le\rho_*^{-1}\Big(\nrm{\eta(0)}+\nrm{A_hu_0}+C\big(\rho^*\nrm{\bU_v(0)}_{L^\infty}+|f|\rho^*+\nrm{c}_{L^\infty}\big)\nrm{\nabla\mathcal P_mu_0}\Big).
\end{equation*}
Gronwall's inequality then bounds $u^m_t$ in $L^\infty(0,T;L^2)\cap L^2(0,T;H^1)$ and $\nabla u^m$ in $L^\infty(0,T;L^2)$, and \eqref{A:H2gal} bounds $u^m$ in $L^\infty(0,T;H^2)$, uniformly in $m$.

\emph{Passage to the limit and uniqueness.} A subsequence converges weakly-$*$ to a solution of the weak formulation with the stated regularity, first for $\phi=e_j$ and then, by completeness of $\{e_k\}$ in $V_\sigma$, for all $\phi\in V_\sigma$; $u(0)=u_0$ since $u^m(0)=\mathcal P_mu_0\to u_0$. Since $\eta-\rho(u_t+\bU_v\cdot\nabla u)-f\rho u^\perp-cu\in L^\infty(0,T;L^2)$ and $u\in L^\infty(0,T;\cD(A_h))$, the weak formulation says $A_hu=\mathcal Q[\eta-\rho(u_t+\bU_v\cdot\nabla u)-f\rho u^\perp-cu]$. Hence \eqref{A:L} holds a.e.\ with
\begin{equation*}
\nh P\eqdef(I-\mathcal Q)\big[\eta-\rho(u_t+\bU_v\cdot\nabla u)-f\rho u^\perp-cu+\mu\Lap u\big],
\end{equation*}
which is a horizontal gradient of a function of $x_h$ by the characterization of $I-\mathcal Q$ and $A_hu=-\mu\Lap u+\nh p_u$. For uniqueness in the larger class, the difference $\delta$ of two solutions satisfies
\begin{equation*}
\frac12\ddt\nrm{\sqrt\rho\delta}^2+\mu\nrm{\nabla\delta}^2=\frac12(\rho_t\delta,\delta)-(\rho\bU_v\cdot\nabla\delta,\delta)-(c\delta,\delta)\le\frac\mu2\nrm{\nabla\delta}^2+C\nrm\delta^2 ,
\end{equation*}
where the Coriolis term has dropped out because $\delta^\perp\cdot\delta=0$. Since $\rho\delta$ does not satisfy the lid constraint, the Lions--Magenes lemma \cite[Lemma III.1.2]{Temam2001} is not applied verbatim; instead we use its proof by mollification in time. With a time mollifier $J_\epsilon$ we test the weak formulation with $J_\epsilon J_\epsilon\delta$, which takes values in $V_\sigma$ because $V_\sigma$ does not depend on $t$; since $\rho\in W^{1,\infty}$, the Friedrichs commutator $J_\epsilon(\rho\delta)-\rho J_\epsilon\delta$ tends to zero in $L^2(0,T;H^1)$, and the identity follows as $\epsilon\to0$. Since $\delta(0)=0$, Gronwall's inequality gives $\delta=0$.
\end{proof}

\begin{lemma}\label{lem:A3}
Let $T\in(0,1]$. Let $\rho=\rs+\varrho$ with $\rho_*\le\rho\le\rho^*$, $\varrho\in C([0,T];H^4)$, $\varrho(0)=\vin$,
\begin{equation*}
\sup_{[0,T]}\nrm{\varrho}_{H^4}\le K_0,\qquad \esssup_{[0,T]}\big(\nrm{\varrho_t}_{H^3}+\nrm{\varrho_{tt}}_{H^1}\big)\le K_1,
\end{equation*}
and let $v$ take values in $V_\sigma$ with $v(0)=\uin$ and $\nrm{v}_{Y_T'}\le R$. Let $h\eqdef g\int_0^z\nh\varrho\dd z'$. Then the problem
\begin{equation}\label{A:Lmain}
\rho\big(u_t+\bU_v\cdot\nabla u\big)+f\rho u^\perp-\mu\Lap u+\nh P=h,\qquad u(t)\in V_\sigma,\qquad u(0)=\uin,
\end{equation}
has a unique solution with $\nrm{u}_{Y_T}<\infty$. Moreover $u\in C([0,T];H^4)$, $u_t\in C([0,T];H^2)$, $u_{tt}\in C([0,T];L^2)$, $P\in C([0,T];H^3(\Tl))$, $u_t(0)=u_1$, and
\begin{equation}\label{A:main}
\nrm{u}_{Y_T}^2\le C_0+C_1T^{\frac12}\big(1+\nrm{u}_{Y_T}^2\big).
\end{equation}
Here $C_0$ depends only on $\rho_*,\rho^*,K_0,M_0,\nrm{v_t(0)}_{H^2}$ and $\nrm{\varrho_t(0)}_{H^2}$, and $C_1$ depends on these and on $R,K_1$.
\end{lemma}

\begin{proof}
Note that $\varrho_t\in C([0,T];H^2)$ by interpolation, that $\nrm{\nabla\rho}_{L^\infty}\le C(1+K_0)$ and $\nrm{\rho_t}_{L^\infty}\le CK_1$, so $\rho\in W^{1,\infty}(\Omega\times(0,T))$, and that $\nrm{h}_{H^3}\le C\nrm{\varrho}_{H^4}$, $\nrm{h_t}_{H^2}\le C\nrm{\varrho_t}_{H^3}$, $\nrm{h_{tt}}_{L^2}\le C\nrm{\varrho_{tt}}_{H^1}$. Whenever \Cref{lem:stokes} is applied below at a fixed time to a field in $V_\sigma$ satisfying a Stokes system weakly, we use that weak solutions in $V_\sigma$ are unique (Lax--Milgram), so that the field coincides with the solution provided by \Cref{lem:stokes}. By \eqref{A:wbd} and \eqref{A:prod},
\begin{equation}\label{A:conv}
\nrm{\bU_v\cdot\nabla u}_{H^k}\le C\nrm{v}_{H^{\max\{k+1,3\}}}\nrm{u}_{H^{k+1}},\qquad k=0,1,2,3;
\end{equation}
for $k=0$ use $\nrm{\bU_v}_{L^\infty}\le C\nrm{v}_{H^3}$, and for $k\ge1$ use \eqref{A:prod} for $v\cdot\nh u$ and for $w_v\pz u$ with $\nrm{w_v}_{H^k}\le C\nrm{v}_{H^{k+1}}$ ($k=2,3$), respectively $\nrm{w_v\pz u}_{H^1}\le C\nrm{w_v}_{H^2}\nrm{\pz u}_{H^1}$.

\emph{Step 1 (first application of \Cref{lem:A2}).} \Cref{lem:A2} applies with $c=0$, $\eta=h\in W^{1,\infty}(0,T;L^2)$ and $u_0=\uin\in H^4\cap V_\sigma$, since $v\in L^\infty(0,T;H^4)$ and $v_t\in L^\infty(0,T;H^2)$. It gives $u\in L^\infty(0,T;H^2)$ and $u_t\in L^\infty(0,T;L^2)\cap L^2(0,T;V_\sigma)$. Applying \Cref{lem:stokes} at a.e.\ $t$ to
\begin{equation*}
-\mu\Lap u+\nh P=\zeta\eqdef h-\rho u_t-\rho\bU_v\cdot\nabla u-f\rho u^\perp,
\end{equation*}
with $\nrm{\zeta}_{H^1}\le C(1+K_0)\big(K_0+\nrm{u_t}_{H^1}+\nrm{v}_{H^3}\nrm{u}_{H^2}+\nrm{u}_{H^1}\big)$ by \eqref{A:conv}, we also get $u\in L^2(0,T;H^3)$.

\emph{Step 2 (time derivative).} Let $U$ be the solution given by \Cref{lem:A2} of
\begin{equation*}
\begin{gathered}
\rho(U_t+\bU_v\cdot\nabla U)+f\rho U^\perp+\rho_tU-\mu\Lap U+\nh\Pi=\eta_1,\qquad U(0)=u_1,\\
\eta_1\eqdef h_t-\rho_t\bU_v\cdot\nabla u-\rho\bU_{v_t}\cdot\nabla u-f\rho_tu^\perp .
\end{gathered}
\end{equation*}
The hypotheses hold: $u_1\in H^2\cap V_\sigma$, $c=\rho_t\in L^\infty$, $c_t=\varrho_{tt}\in L^\infty(0,T;H^1)$, and $\eta_1\in W^{1,1}(0,T;L^2)$. For the last claim, $\eta_1\in L^\infty(0,T;L^2)$ by \eqref{A:conv}, \eqref{A:wbd} and \eqref{A:prod} (e.g.\ $\nrm{w_{v_t}\pz u}\le C\nrm{w_{v_t}}_{L^\infty_zL^4_h}\nrm{u}_{H^2}\le C\nrm{v_t}_{H^2}\nrm{u}_{H^2}$), and
\begin{equation*}
\begin{aligned}
\pt\eta_1&=h_{tt}-\rho_{tt}\bU_v\cdot\nabla u-2\rho_t\bU_{v_t}\cdot\nabla u-\rho_t\bU_v\cdot\nabla u_t\\
&\quad-\rho\bU_{v_{tt}}\cdot\nabla u-\rho\bU_{v_t}\cdot\nabla u_t-f\rho_{tt}u^\perp-f\rho_tu_t^\perp .
\end{aligned}
\end{equation*}
All terms lie in $L^2(0,T;L^2)$, using $H^1\subset L^6$ for $\varrho_{tt}$, $v_{tt}$ and $\nabla u\in L^\infty(0,T;L^3)$, except possibly
\begin{equation*}
\begin{aligned}
&\nrm{w_{v_{tt}}\pz u}_{L^2}\le\nrm{w_{v_{tt}}}_{L^2_hL^\infty_z}\nrm{\pz u}_{L^\infty_hL^2_z}\le C\nrm{v_{tt}}_{H^1}\nrm{u}_{H^3},\\
&\nrm{w_{v_t}\pz u_t}_{L^2}\le\nrm{w_{v_t}}_{L^\infty}\nrm{\pz u_t}_{L^2}\le C\nrm{v_t}_{H^3}\nrm{u_t}_{H^1}
\end{aligned}
\end{equation*}
(\Cref{lem:aniso}(e) and \eqref{A:wbd}), which are products of two functions in $L^2(0,T)$, hence in $L^1(0,T)$.

We show $U=u_t$. By the weak form of \eqref{A:Lmain}, for $\phi\in V_\sigma$ and a.e.\ $t$,
\begin{equation*}
(\rho u_t,\phi)=r_\phi(t)\eqdef\big(h-\rho\bU_v\cdot\nabla u-f\rho u^\perp,\phi\big)-\mu(\nabla u,\nabla\phi).
\end{equation*}
Since $u\in H^1(0,T;H^1)\cap L^\infty(0,T;H^2)$, $\rho\in W^{1,\infty}$, $h\in W^{1,\infty}(0,T;L^2)$ and $\bU_v\in W^{1,\infty}(0,T;L^2)\cap L^\infty((0,T)\times\Omega)$, $r_\phi\in W^{1,1}(0,T)$, with
\begin{equation*}
r_\phi'=\big(\eta_1-\rho\bU_v\cdot\nabla u_t-f\rho u_t^\perp,\phi\big)-\mu(\nabla u_t,\nabla\phi).
\end{equation*}
Hence $t\mapsto(\rho u_t,\phi)$ is absolutely continuous and $u_t\in L^\infty(0,T;L^2)\cap L^2(0,T;V_\sigma)$ satisfies the weak form of the equation for $U$, in the sense of the uniqueness class of \Cref{lem:A2}. Its initial value is $r_\phi(0)=(h(0)-\rin\bU_{\mathrm{in}}\cdot\nabla\uin-f\rin\uin^\perp,\phi)-\mu(\nabla\uin,\nabla\phi)=(F_0,\phi)=(\rin u_1,\phi)$, by \eqref{A:u1} and $(\nh P_0,\phi)=0$. By uniqueness, $U=u_t$. Consequently $u_t\in L^\infty(0,T;H^2)$, $u_{tt}\in L^\infty(0,T;L^2)\cap L^2(0,T;V_\sigma)$, $u_t\in C([0,T];L^2)$ and $u_t(0)=u_1$.

\emph{Step 3 (elliptic regularity).} We apply \Cref{lem:stokes} at a.e.\ $t$ to
\begin{equation*}
-\mu\Lap u+\nh P=\zeta,\qquad -\mu\Lap u_t+\nh P_t=\zeta_t=h_t-\rho_tu_t-\rho u_{tt}-(\rho\bU_v\cdot\nabla u)_t-f(\rho u^\perp)_t ,
\end{equation*}
in the following order. (a) Since $u_t\in L^\infty(0,T;H^2)$, the bound for $\nrm{\zeta}_{H^1}$ in Step 1 gives $u\in L^\infty(0,T;H^3)$. (b) By \eqref{A:prod} with $k=2$ and \eqref{A:conv},
\begin{equation}\label{A:zeta}
\nrm{\zeta}_{H^k}\le C(1+K_0)\big(K_0+\nrm{u_t}_{H^k}+\nrm{v}_{H^{k+1}}\nrm{u}_{H^{k+1}}+\nrm{u}_{H^k}\big),\qquad k=2,3,
\end{equation}
so $u\in L^\infty(0,T;H^4)$. (c) The terms of $\zeta_t$ satisfy
\begin{equation*}
\begin{aligned}
&\nrm{h_t}_{H^1}+\nrm{\rho_tu_t}_{H^1}+\nrm{f(\rho u^\perp)_t}_{H^1}\le C(K_0,K_1)\big(1+\nrm{u_t}_{H^1}+\nrm{u}_{H^1}\big),\\
&\nrm{\rho u_{tt}}_{H^1}\le C(1+K_0)\nrm{u_{tt}}_{H^1},\\
&\nrm{(\rho\bU_v\cdot\nabla u)_t}_{H^1}\le C(K_0,K_1)\big(\nrm{v}_{H^3}\nrm{u}_{H^2}+\nrm{v_t}_{H^2}\nrm{u}_{H^4}+\nrm{v}_{H^4}\nrm{u_t}_{H^2}\big);
\end{aligned}
\end{equation*}
here we used $\nrm{\nabla\varrho_t}_{L^\infty}\le C\nrm{\varrho_t}_{H^3}$, \eqref{A:conv} with $k=1$, $\nrm{v_t\cdot\nh u}_{H^1}\le C\nrm{v_t}_{H^2}\nrm{u}_{H^2}$, $\nrm{w_{v_t}\pz u}_{H^1}\le C\nrm{w_{v_t}}_{H^1}\nrm{\pz u}_{W^{1,\infty}}\le C\nrm{v_t}_{H^2}\nrm{u}_{H^4}$ and $\nrm{\bU_v\cdot\nabla u_t}_{H^1}\le C\nrm{\bU_v}_{W^{1,\infty}}\nrm{u_t}_{H^2}\le C\nrm{v}_{H^4}\nrm{u_t}_{H^2}$. Since $u_{tt}\in L^2(0,T;H^1)$, we get $u_t\in L^2(0,T;H^3)$. (d) By \eqref{A:zeta} with $k=3$, $u\in L^2(0,T;H^5)$. Hence $\nrm{u}_{Y_T}<\infty$.

Time continuity follows from the Lions--Magenes lemma for $u$ ($u\in L^2H^5$, $u_t\in L^2H^3$) and for $u_t$ ($u_t\in L^2H^3$, $u_{tt}\in L^2H^1$). For $u_{tt}$, note that after (a)--(b), $\pt\eta_1\in L^2(0,T;L^2)$, since now $\nrm{w_{v_{tt}}\pz u}\le\nrm{w_{v_{tt}}}\nrm{\pz u}_{L^\infty}\le C\nrm{v_{tt}}_{H^1}\nrm{u}_{H^3}$ and $\nrm{w_{v_t}\pz u_t}\le\nrm{w_{v_t}}_{L^\infty_zL^4_h}\nrm{\pz u_t}_{L^2_zL^4_h}\le C\nrm{v_t}_{H^2}\nrm{u_t}_{H^2}$ by \eqref{A:wbd} and \Cref{lem:aniso}(f). Differentiating the weak form of the equation for $U=u_t$ once more in $t$, we see that $t\mapsto(\rho u_{tt},\phi)$ is absolutely continuous for $\phi\in V_\sigma$, with derivative bounded by $a(t)\nrm{\phi}_{H^1}$, $a\in L^2(0,T)$. As in the uniqueness part of \Cref{lem:A2} (cf.\ \Cref{rem:justify}(b)), this justifies the energy identity of Step 4(i) below. For the continuity of $u_{tt}$ we write the equation for $U=u_t$ as $\rho u_{tt}+\nh\Pi=X$ with
\begin{equation*}
X\eqdef\eta_1-\rho\bU_v\cdot\nabla u_t-f\rho u_t^\perp-\rho_tu_t+\mu\Lap u_t .
\end{equation*}
After (a)--(c), $u_t\in C([0,T];H^2)$; moreover $\eta_1\in W^{1,1}(0,T;L^2)\subset C([0,T];L^2)$, $\rho\in C(\ol\Omega\times[0,T])$, $\rho_t\in C([0,T];H^2)$ and $\bU_v\in C([0,T];L^\infty)$. Hence every term of the equation except $\rho u_{tt}$ and $\nh\Pi$ lies in $C([0,T];L^2)$, i.e.\ $X\in C([0,T];L^2)$. Since $u_{tt}(t)\in V_\sigma$ for a.e.\ $t$, $\Pi(t)$ solves $\divh(a\nh\Pi)=\divh\int_0^1\rho^{-1}X\dd z$ on $\Tl$, with $a\eqdef\int_0^1\rho^{-1}\dd z\in C([0,T];L^\infty(\Tl))$ and $(\rho^*)^{-1}\le a\le\rho_*^{-1}$. For solutions of $\divh(a_i\nh\Pi_i)=\divh F_i$, $i=1,2$, the Lax--Milgram lemma gives $\nrm{\nh\Pi_1-\nh\Pi_2}_{L^2(\Tl)}\le C(\rho_*,\rho^*)\big(\nrm{F_1-F_2}_{L^2(\Tl)}+\nrm{a_1-a_2}_{L^\infty(\Tl)}\nrm{\nh\Pi_2}_{L^2(\Tl)}\big)$. Therefore $\nh\Pi\in C([0,T];L^2)$, and $u_{tt}=\rho^{-1}(X-\nh\Pi)\in C([0,T];L^2)$. Finally, $\zeta\in C([0,T];H^2)$, because $\varrho\in C([0,T];H^4)$, $u_t\in C([0,T];H^2)$, $u,v\in C([0,T];H^3)$ and \eqref{A:conv} is bilinear; \Cref{lem:stokes} with $m=2$ gives $\nh P\in C([0,T];H^2(\Tl))$, hence $P\in C([0,T];H^3(\Tl))$ after normalization to zero mean.

\emph{Step 4 (the estimate \eqref{A:main}).} Write $N\eqdef\nrm{u}_{Y_T}$. Quantities bounded in terms of $\rho_*,\rho^*,K_0,M_0,\nrm{v_t(0)}_{H^2}$ and $\nrm{\varrho_t(0)}_{H^2}$ are called \emph{data-type} and denoted by $\Lambda$; by \eqref{A:u1bd}, $\nrm{u_1}_{H^2}\le\Lambda$. Since $u=\uin+\int_0^tu_t$, $v=\uin+\int_0^tv_t$, $u_t=u_1+\int_0^tu_{tt}$, $v_t=v_t(0)+\int_0^tv_{tt}$ and $\varrho_t=\varrho_t(0)+\int_0^t\varrho_{tt}$, we have for $t\le T\le1$:
\begin{equation}\label{A:pointwise}
\begin{aligned}
&\nrm{u(t)}_{H^3}\le M_0+\sqrt TN,\quad \nrm{u_t(t)}_{H^1}\le\Lambda+\sqrt TN,\\
&\nrm{v(t)}_{H^3}\le M_0+\sqrt TR,\quad \nrm{v_t(t)}_{H^1}\le\Lambda+\sqrt TR,\\
&\nrm{\varrho_t(t)}_{H^1}\le\Lambda+TK_1,\qquad \nrm{\varrho_t(t)}_{H^2}\le\Lambda+C\sqrt TK_1 .
\end{aligned}
\end{equation}
The last bound follows by interpolation: $\nrm{\varrho_t(t)-\varrho_t(0)}_{H^2}\le C\nrm{\varrho_t(t)-\varrho_t(0)}_{H^1}^{1/2}\nrm{\varrho_t(t)-\varrho_t(0)}_{H^3}^{1/2}\le C(TK_1)^{1/2}(2K_1)^{1/2}$. In addition, $\nrm{u}_{H^4}+\nrm{u_t}_{H^2}+\nrm{u_{tt}}\le CN$ and $\nrm{v}_{H^4}+\nrm{v_t}_{H^2}+\nrm{v_{tt}}\le CR$ on $[0,T]$.

(i) \emph{Energy for $u_{tt}$.} Differentiating the equation for $U=u_t$ in $t$ and testing with $u_{tt}\in V_\sigma$ (the surface pressure drops out; the identity is justified in Step 3), we obtain, since $(f\rho u_{tt}^\perp,u_{tt})=0$,
\begin{equation*}
\begin{aligned}
\frac12\ddt\nrm{\sqrt\rho u_{tt}}^2+\mu\nrm{\nabla u_{tt}}^2&=-g\int_\Omega\Big(\int_0^z\varrho_{tt}\dd z'\Big)\divh u_{tt}\dx-\frac32\int_\Omega\rho_t|u_{tt}|^2\dx-\int_\Omega\rho_{tt}u_t\cdot u_{tt}\dx\\
&\quad-\int_\Omega\big(\rho\bU_v\cdot\nabla u_{tt}+2(\rho\bU_v)_t\cdot\nabla u_t+(\rho\bU_v)_{tt}\cdot\nabla u\big)\cdot u_{tt}\dx\\
&\quad-f\int_\Omega\big(\rho_{tt}u^\perp+2\rho_tu_t^\perp\big)\cdot u_{tt}\dx ,
\end{aligned}
\end{equation*}
where $(\rho\bU_v)_t=\rho_t\bU_v+\rho\bU_{v_t}$ and $(\rho\bU_v)_{tt}=\rho_{tt}\bU_v+2\rho_t\bU_{v_t}+\rho\bU_{v_{tt}}$. The first term is $\le\frac\mu4\nrm{\nabla u_{tt}}^2+CK_1^2$, and the convective term $\rho\bU_v\cdot\nabla u_{tt}$ contributes at most $\frac\mu4\nrm{\nabla u_{tt}}^2+C(\rho^*)^2R^2\nrm{u_{tt}}^2$. All other terms are bounded by $C(K_0,K_1,R)(1+\nrm{v_{tt}}_{H^1})N^2$. For instance, by $H^2\subset L^\infty$, \eqref{A:wbd} and \eqref{A:prod},
\begin{equation*}
\begin{aligned}
&|(\rho_{tt}u_t,u_{tt})|+|f|\,|(\rho_{tt}u^\perp,u_{tt})|\le C\nrm{\varrho_{tt}}\big(\nrm{u_t}_{L^\infty}+\nrm{u}_{L^\infty}\big)\nrm{u_{tt}}\le CK_1N^2,\\
&|(\rho w_{v_t}\pz u_t,u_{tt})|\le C\rho^*\nrm{w_{v_t}}_{L^\infty_zL^4_h}\nrm{u_t}_{H^2}\nrm{u_{tt}}\le C\nrm{v_t}_{H^2}N^2,\\
&|(\rho_{tt}\bU_v\cdot\nabla u,u_{tt})|\le\nrm{\varrho_{tt}}\nrm{\bU_v}_{L^\infty}\nrm{\nabla u}_{L^\infty}\nrm{u_{tt}}\le CK_1RN^2,\\
&|(\rho\bU_{v_{tt}}\cdot\nabla u,u_{tt})|\le\rho^*\big(\nrm{v_{tt}}+\nrm{w_{v_{tt}}}\big)\nrm{\nabla u}_{L^\infty}\nrm{u_{tt}}\le C\nrm{v_{tt}}_{H^1}N^2,
\end{aligned}
\end{equation*}
and the remaining terms ($\rho_t|u_{tt}|^2$, $\rho_t\bU_v\cdot\nabla u_t$, $\rho v_t\cdot\nh u_t$, $\rho_t\bU_{v_t}\cdot\nabla u$, $\rho_tu_t^\perp$) are estimated in the same way with $\nrm{\rho_t}_{L^\infty}\le CK_1$, $\nrm{v_t}_{L^\infty}\le CR$, $\nrm{w_{v_t}}\le\nrm{\divh v_t}$ and $\nrm{\pz u}_{L^\infty}\le CN$. Integrating over $(0,t)$, using $\int_0^T\nrm{v_{tt}}_{H^1}\le\sqrt TR$, $T\le\sqrt T$ and the Poincar\'e inequality $\nrm{u_{tt}}_{H^1}\le C\nrm{\nabla u_{tt}}$, we obtain
\begin{equation}\label{A:Utt}
\sup_{[0,T]}\nrm{u_{tt}}^2+\int_0^T\nrm{u_{tt}}_{H^1}^2\le\Lambda\big(1+\nrm{u_{tt}(0)}^2\big)+C_1T^{\frac12}(1+N^2),\qquad \nrm{u_{tt}(0)}\le\Lambda .
\end{equation}
Here $\nrm{u_{tt}(0)}$ is data-type: testing the equation for $U=u_t$ with $u_{tt}(t)\in V_\sigma$ gives, for a.e.\ $t$, $\rho_*\nrm{u_{tt}}\le\nrm{\eta_1-\rho\bU_v\cdot\nabla u_t-\rho_tu_t-f\rho u_t^\perp+\mu\Lap u_t}$, and both sides are continuous in $t$; at $t=0$, with $u_t(0)=u_1$, $\bU_v(0)=\bU_{\mathrm{in}}$ and $\rho(0)=\rin$,
\begin{equation*}
\begin{aligned}
\nrm{u_{tt}(0)}&\le C(\rho_*,\rho^*)\Big(\nrm{\varrho_t(0)}_{H^1}+\nrm{\varrho_t(0)}_{H^2}\big(M_0^2+\nrm{u_1}+M_0\big)\\
&\qquad+\nrm{v_t(0)}_{H^2}M_0+M_0\nrm{u_1}_{H^1}+\nrm{u_1}_{H^2}\Big),
\end{aligned}
\end{equation*}
where we used $\nrm{\bU_{v_t}(0)\cdot\nabla\uin}\le\nrm{v_t(0)}_{L^\infty}\nrm{\nabla\uin}+\nrm{w_{v_t}(0)}\nrm{\pz\uin}_{L^\infty}\le C\nrm{v_t(0)}_{H^2}M_0$.

(ii) \emph{$\sup\nrm{u_t}_{H^2}$.} By \Cref{lem:stokes}, $\nrm{u_t}_{H^2}\le C\big(\nrm{h_t}+\nrm{\rho_tu_t}+\rho^*\nrm{u_{tt}}+\nrm{(\rho\bU_v\cdot\nabla u)_t}+|f|\nrm{(\rho u^\perp)_t}\big)$. By \eqref{A:pointwise}, \eqref{A:wbd} and \eqref{A:prod}, for example
\begin{equation*}
\nrm{\rho\bU_{v_t}\cdot\nabla u}\le\rho^*\big(\nrm{v_t}_{L^4}\nrm{\nabla u}_{L^4}+\nrm{w_{v_t}}\nrm{\pz u}_{L^\infty}\big)\le C\nrm{v_t}_{H^1}\nrm{u}_{H^3}\le C(\Lambda+\sqrt TR)(M_0+\sqrt TN),
\end{equation*}
each term other than $\rho^*\nrm{u_{tt}}$ is a product of factors of the form $\Lambda+\sqrt T(N+R+K_1)$, exactly one of which contains $N$ (the problem is linear in $u$). Together with \eqref{A:Utt}, this gives $\sup\nrm{u_t}_{H^2}^2\le\Lambda+C_1T^{\frac12}(1+N^2)$.

(iii) \emph{$\sup\nrm{u}_{H^4}$.} By \Cref{lem:stokes} and \eqref{A:zeta} with $k=2$, $\nrm{u}_{H^4}\le C(K_0)\big(K_0+\nrm{u_t}_{H^2}+\nrm{v}_{H^3}\nrm{u}_{H^3}+\nrm{u}_{H^2}\big)$. By (ii) and \eqref{A:pointwise}, $\sup\nrm{u}_{H^4}^2\le\Lambda+C_1T^{\frac12}(1+N^2)$.

(iv) \emph{$\int\nrm{u_t}_{H^3}^2$.} By \Cref{lem:stokes}, $\nrm{u_t}_{H^3}\le C\nrm{\zeta_t}_{H^1}$. By the bounds in Step 3(c), the term $\rho u_{tt}$ contributes $C(K_0)\int_0^T\nrm{u_{tt}}_{H^1}^2$, which is controlled by \eqref{A:Utt}, and the other terms are bounded pointwise by $C(K_0,K_1,R)(1+N)$, so their contribution is $\le C_1T(1+N^2)$.

(v) \emph{$\int\nrm{u}_{H^5}^2$.} By \Cref{lem:stokes} and \eqref{A:zeta} with $k=3$, $\nrm{u}_{H^5}\le C(K_0)\big(K_0+\nrm{u_t}_{H^3}+\nrm{v}_{H^4}\nrm{u}_{H^4}+\nrm{u}_{H^3}\big)$; we integrate in time and use (iv).

Adding (i)--(v) gives \eqref{A:main}.
\end{proof}

\subsection{Difference estimate}

\begin{lemma}\label{lem:A4}
Let $v,v'\in\mathcal B_{R,T}$, let $\rho,\rho'$ be given by \Cref{lem:A1} with velocities $\bU_v,\bU_{v'}$, let $u,u'$ be the solutions of \eqref{A:Lmain} given by \Cref{lem:A3} with $(\rho,\bU_v)$ and $(\rho',\bU_{v'})$, and suppose $\nrm{u'}_{Y_T}\le R$. Then
\begin{equation*}
\sup_{[0,T]}\nrm{\rho-\rho'}_{L^2}+\nrm{u-u'}_{W_T}\le C(M_0,\rho_*,R)\,T^{\frac12}\nrm{v-v'}_{W_T}.
\end{equation*}
\end{lemma}

\begin{proof}
Write $\delta\rho=\delta\varrho=\rho-\rho'$, $\delta v=v-v'$, $\delta\bU=\bU_v-\bU_{v'}=\bU_{\delta v}$, $\delta u=u-u'$ and $\delta P=P-P'$. By \eqref{A:transport}, $\sup\nrm{\varrho'}_{H^4}\le C(M_0,R)$, and $\rho,\rho'\in C^1([0,T];L^2)$ by \Cref{lem:A1}. The difference solves $\delta\rho_t+\bU_v\cdot\nabla\delta\rho=-\delta\bU\cdot\nabla\rho'$, $\delta\rho(0)=0$. We test with $\delta\rho$; the transport term drops out because $\bU_v$ is divergence-free and tangent to $\pa\Omega$. Using $\nrm{\bU_{\delta v}}_{L^2}\le\nrm{\delta v}+\nrm{\divh\delta v}$ (by \eqref{A:wbd}) and $\nrm{\delta v}\le C\nrm{\nabla\delta v}$, we obtain
\begin{equation*}
\ddt\nrm{\delta\rho}\le\nrm{\delta\bU\cdot\nabla\rho'}\le\nrm{\delta\bU}\nrm{\nabla\rho'}_{L^\infty}\le C(1+\nrm{\varrho'}_{H^4})\nrm{\nabla\delta v},
\end{equation*}
so $\sup\nrm{\delta\rho}\le C(M_0,R)\sqrt T\nrm{\delta v}_{W_T}$. Next, using $\rho u^\perp-\rho'u'^\perp=\rho\,\delta u^\perp+\delta\rho\,u'^\perp$,
\begin{equation*}
\rho(\delta u_t+\bU_v\cdot\nabla\delta u)+f\rho\,\delta u^\perp-\mu\Lap\delta u+\nh\delta P=g\int_0^z\nh\delta\varrho\dd z'-\delta\rho\,\big(u'_t+\bU_v\cdot\nabla u'+fu'^\perp\big)-\rho'\delta\bU\cdot\nabla u',
\end{equation*}
with $\delta u(0)=0$. We test with $\delta u\in V_\sigma$; this is legitimate since $\delta u\in H^1(0,T;H^1_0)$. Since $\rho_t+\divv(\rho\bU_v)=0$, the first term on the left gives $\frac12\ddt\nrm{\sqrt\rho\delta u}^2$; the Coriolis term drops out because $\delta u^\perp\cdot\delta u=0$, and the pressure term vanishes. On the right we use $|g\int(\int_0^z\nh\delta\varrho)\cdot\delta u|=|g\int(\int_0^z\delta\varrho)\divh\delta u|\le g\nrm{\delta\varrho}\nrm{\nabla\delta u}$; $\nrm{u'_t+\bU_v\cdot\nabla u'+fu'^\perp}_{L^\infty}\le C(\nrm{u'_t}_{H^2}+\nrm{v}_{H^3}\nrm{u'}_{H^3}+\nrm{u'}_{H^2})\le C(1+R)^2$ by $H^2\subset L^\infty$ and \eqref{A:wbd}; and $|(\rho'\delta\bU\cdot\nabla u',\delta u)|\le\rho^*\nrm{\delta\bU}\nrm{\nabla u'}_{L^\infty}\nrm{\delta u}\le C\nrm{u'}_{H^3}\nrm{\nabla\delta v}\nrm{\delta u}$. Therefore
\begin{equation*}
\ddt\nrm{\sqrt\rho\delta u}^2+\mu\nrm{\nabla\delta u}^2\le C(R)\big(\nrm{\delta\rho}^2+\nrm{\delta u}^2+\nrm{\nabla\delta v}\nrm{\delta u}\big),
\end{equation*}
and Gronwall's inequality, the bound for $\delta\rho$ and $\int_0^T\nrm{\nabla\delta v}\nrm{\delta u}\le\sqrt T\nrm{\delta v}_{W_T}\sup\nrm{\delta u}$ give the claim.
\end{proof}

\subsection{Proof of \texorpdfstring{\Cref{prop:local}}{the local well-posedness proposition}}

Let $K_0\eqdef2M_0+1$. Let $\Lambda_*$ be the maximum of the constants $C_0$ of \Cref{lem:A3} with $K_0$ as above, $\nrm{v_t(0)}_{H^2}\le\nrm{u_1}_{H^2}$ and $\nrm{\varrho_t(0)}_{H^2}\le\nrm{\varrho_1}_{H^3}$; by \eqref{A:u1bd} and $\rho^*\le\nrm{\rs}_{L^\infty}+CM_0$, $\Lambda_*=\Lambda_*(M_0,\rho_*)$. Set $R^2\eqdef\max\{2\Lambda_*+1,\,M_0^2\}$, and let $K_1\eqdef C(1+R)^2(1+K_0)$ be the bound in \eqref{A:transport}. Choose $T_0=T_0(M_0,\rho_*)\in(0,1]$ so small that the following hold:
\begin{enumerate}[label=(\alph*)]
\item $(M_0+C\sqrt{T_0}R)e^{C\sqrt{T_0}R}\le K_0$;
\item $C_1T_0^{1/2}\le\frac12$, where $C_1$ is the constant of \Cref{lem:A3}; it depends on $\rho_*,\rho^*,K_0,M_0,R,K_1$ and on $\nrm{v_t(0)}_{H^2}$, $\nrm{\varrho_t(0)}_{H^2}$, which below are always $\nrm{u_1}_{H^2}$ or $0$, respectively $\nrm{\varrho_1}_{H^2}$ or $0$, and hence are bounded by $C(M_0,\rho_*)$ by \eqref{A:u1bd};
\item the constant of \Cref{lem:A4} satisfies $C(M_0,\rho_*,R)T_0^{1/2}\le\frac12$.
\end{enumerate}
Write $T=T_0$.

\emph{Starting point.} Let $u^0$ be the solution of \eqref{A:Lmain} with the time-independent coefficients $\rho\equiv\rin$ and $v\equiv\uin$ (so $\bU_v\equiv\bU_{\mathrm{in}}$); here $\varrho\equiv\vin$, $K_1$ may be taken $0$, $\nrm{v}_{Y_T'}=\nrm{\uin}_{H^4}\le R$, $v_t(0)=0$ and $\varrho_t(0)=0$. By \Cref{lem:A3}, $u^0(0)=\uin$ and $u^0_t(0)=u_1$, and by (b), \eqref{A:main} gives $\nrm{u^0}_{Y_T}^2\le2\Lambda_*+1\le R^2$, so $u^0\in\mathcal B_{R,T}$.

\emph{Iteration.} Given $u^n\in\mathcal B_{R,T}$, let $\rho^{n+1}=\rs+\varrho^{n+1}$ be given by \Cref{lem:A1} with velocity $\bU_{u^n}$. By (a), $\sup\nrm{\varrho^{n+1}}_{H^4}\le K_0$, and by \eqref{A:transport} the time derivatives are bounded by $K_1$; moreover $\varrho^{n+1}\in C([0,T];H^4)$, $\rho_*\le\rho^{n+1}\le\rho^*$ and $\varrho^{n+1}_t(0)=\varrho_1$. Let $u^{n+1}$ solve \eqref{A:Lmain} with $(\rho^{n+1},\bU_{u^n})$; \Cref{lem:A3} applies with $v=u^n$, since $\nrm{u^n}_{Y_T'}\le R$ and $\nrm{u^n_t(0)}_{H^2}=\nrm{u_1}_{H^2}$. Then $u^{n+1}(0)=\uin$, $u^{n+1}_t(0)=u_1$, and $\nrm{u^{n+1}}_{Y_T}\le R$ by \eqref{A:main} and (b). Hence $u^{n+1}\in\mathcal B_{R,T}$.

\emph{Convergence.} By \Cref{lem:A4} (with $v=u^n$, $v'=u^{n-1}$) and (c), $\nrm{u^{n+1}-u^n}_{W_T}\le\frac12\nrm{u^n-u^{n-1}}_{W_T}$ and $\sup\nrm{\rho^{n+1}-\rho^n}_{L^2}\le\frac12\nrm{u^n-u^{n-1}}_{W_T}$. Thus $u^n\to u$ in $C([0,T];L^2)\cap L^2(0,T;H^1)$ and $\rho^n\to\rho$ in $C([0,T];L^2)$. By the uniform bounds and the lower semicontinuity of norms under weak-$*$ convergence, $\nrm{u}_{Y_T}\le R$ and $\sup\nrm{\varrho}_{H^4}\le K_0$, where $\varrho=\rho-\rs$; moreover $u(t)\in V_\sigma$, $u(0)=\uin$, and $u_t(0)=u_1$ because $u^n_t\to u_t$ in $C([0,T];H^1)$ by the Aubin--Simon lemma. Passing to the limit in the transport equations and in the weak form of \eqref{A:Lmain}, we find that $\rho$ solves $\rho_t+\bU_u\cdot\nabla\rho=0$, $\rho(0)=\rin$, and that $u$ is a weak solution of \eqref{A:Lmain} with coefficients $(\rho,\bU_u)$. By the uniqueness in \Cref{lem:A1}, $\rho$ is the solution given there with $v=u$; hence $\varrho\in C([0,T];H^4)$, $\rho_*\le\rho\le\rho^*$, and \eqref{A:transport} holds. By the uniqueness in \Cref{lem:A2} (applied in \Cref{lem:A3}), $u$ is the solution given by \Cref{lem:A3} with $(\rho,\bU_u)$; hence $u\in C([0,T];H^4)\cap L^2(0,T;H^5)$, $u_t\in C([0,T];H^2)\cap L^2(0,T;H^3)$, $u_{tt}\in C([0,T];L^2)\cap L^2(0,T;H^1)$, and \eqref{A:Lmain} holds with a surface pressure $P\in C([0,T];H^3(\Tl))$. With $w=w_u$ and $q\eqdef P-g\int_0^z\varrho\dd z'$, we have $\nh q=\nh P-g\int_0^z\nh\varrho\dd z'$ and $\pz q=-g\varrho$, so $(\varrho,u,q)$ solves \eqref{PE} with \eqref{eq:ID} and \eqref{eq:w}; the transport equation for $\rho$ is \eqref{PE}$_1$ since $\bU_u\cdot\nabla\rs=w\rs'$. Since $u\in C([0,T];H^4)$, we have $w_u\in C([0,T];H^3)$ by \eqref{A:wbd}, and $\varrho_t=-\bU_u\cdot\nabla\varrho-w_u\rs'\in C([0,T];H^3)$ because $H^3(\Omega)$ is an algebra. Moreover $\nabla q=(\nh P-g\int_0^z\nh\varrho\dd z',-g\varrho)\in C([0,T];H^2)$ (cf.\ \Cref{rem:gradq}). The no-slip condition and the lid constraint hold since $u(t)\in V_\sigma$; in particular $P$ is the surface pressure of \Cref{lem:propagation}(i), and the clamped conditions $\varrho=\pz\varrho=0$ on $\pa\Omega$ are propagated by \Cref{lem:propagation}(ii) (the solution is classical: $\varrho\in C([0,T];C^2(\ol\Omega))$, $\varrho_t\in C([0,T];H^3)\subset C([0,T];C^1(\ol\Omega))$ and $u\in C([0,T];C^2(\ol\Omega))$). Thus $(\varrho,u,q)$ solves \eqref{PE} with \eqref{eq:ID}--\eqref{eq:lid} in the class \eqref{eq:localclass}.

\emph{Uniqueness.} Let $(\varrho,u,q)$ and $(\varrho',u',q')$ be two solutions of \eqref{PE} with \eqref{eq:ID}--\eqref{eq:lid} in the class \eqref{eq:localclass} on $[0,T']$. Then $u(t),u'(t)\in V_\sigma$ by \eqref{eq:BC} and \eqref{eq:lid}, $w=w_u$ and $w'=w_{u'}$ by \eqref{eq:w}, $\nrm{u}_{Y_{T'}}+\nrm{u'}_{Y_{T'}}\le R'$ for some $R'$, and $\sup\nrm{\varrho}_{H^4}+\sup\nrm{\varrho'}_{H^4}<\infty$. By $\pz q=-g\varrho$, $q=P-g\int_0^z\varrho\dd z'$ with a surface pressure $P=P(t,x_h)$, which we normalize to have zero mean, and similarly for $q'$. The density $\rho=\rs+\varrho$ solves the transport problem of \Cref{lem:A1} with $v=u$; by the flow representation in the proof of \Cref{lem:A1}, $\rho(t)=\rin\circ X(0,t,\cdot)$, so $\inf\rho=\inf\rin=\rho_*>0$ and $\sup\rho=\rho^*$, and likewise for $\rho'$. Moreover $u$ solves \eqref{A:Lmain} with $(\rho,\bU_u)$, and similarly for the primed solution. The proof of \Cref{lem:A4}, applied with $v=u$ and $v'=u'$ on an interval $[t_0,t_0+\tau]$ on which $\rho-\rho'$ and $u-u'$ vanish at $t_0$, gives $\nrm{u-u'}_{W}\le C\tau^{1/2}\nrm{u-u'}_{W}$, where $\nrm{\cdot}_W$ is the norm $\nrm{\cdot}_{W_T}$ of \eqref{A:Ynorm} taken over $[t_0,t_0+\tau]$, with $C$ depending only on $R'$ and the bounds on $\varrho,\varrho'$. Hence $u=u'$ on a short interval of fixed length, and, by repetition, on $[0,T']$; then $\rho=\rho'$ by \Cref{lem:A1}. Finally, $\nh P$ is determined by $(\rho,u)$ through \eqref{A:Lmain}, so $\nh P=\nh P'$, hence $P=P'$ by the normalization and $\nabla q=\nabla q'$. Thus $(\varrho,v,\nabla q)$ is unique, and $q$ is unique up to a function of $t$. The uniqueness argument does not need the solutions to lie in $\mathcal B_{R,T}$ with the particular $R$ fixed above.

\emph{Dependence of $T_0$.} By construction $T_0$ depends only on $M_0$ and $\rho_*$ (and on $l,\mu,g,f,\rs$); the upper bound $\rho^*\le\nrm{\rs}_{L^\infty}+CM_0$ is controlled by $M_0$. This proves \Cref{prop:local}, with the additional information that the existence time is bounded below in terms of an upper bound for $\nrm{\vin}_{H^4}+\nrm{\uin}_{H^4}$ and a positive lower bound for $\rin$. \qed

\section*{Acknowledgement}
Q. Wang was supported by the Natural Science Foundation of Sichuan Province (No.~2025ZNSFSC0072) and the National Natural Science Foundation of China (No. 12301131).

\end{document}